\documentclass[a4paper, 10pt, oneside, reqno]{amsart}

\usepackage[utf8]{inputenc}
\usepackage{amsfonts,amssymb,amsthm,amsmath}
\usepackage{mathrsfs}
\usepackage{geometry}
\usepackage{graphicx,graphics}
\usepackage{color}
\usepackage{tikz}
\usepackage{tikz-cd}
\usepackage{comment}
\usepackage{subfiles}
\usepackage{stmaryrd}
\usepackage{hyperref}
\usepackage{todonotes}
\usepackage{cleveref}

\hypersetup{colorlinks = true, linkbordercolor = {white}, linkcolor = {blue}, citecolor = {blue}}

\theoremstyle{plain}
\newtheorem{thm}{Theorem}[subsection]
\newtheorem{lem}[thm]{Lemma}
\newtheorem{prop}[thm]{Proposition}
\newtheorem{corol}[thm]{Corollary}
\newtheorem{conj}[thm]{Conjecture}
\theoremstyle{definition}
\newtheorem{defi}[thm]{Definition}

\newtheorem{defilem}[thm]{Lemma-Definition}
\newtheorem{construction}[thm]{Construction}
\theoremstyle{remark}
\newtheorem{rque}[thm]{Remark}
\newtheorem{rques}[thm]{Remarks}
\newtheorem{ex}[thm]{Example}
\newtheorem{nota}[thm]{Notation}

\newcommand{\F}{\mathbb{F}}
\newcommand{\Fp}{\F_p}
\newcommand{\Fq}{\F_q}
\newcommand{\Fl}{\F_{\ell}}

\newcommand{\Fqb}{\overline{\F}_q}
\newcommand{\Flb}{\overline{\F}_{\ell}}

\newcommand{\Z}{\mathbb{Z}}
\newcommand{\Zl}{\Z_{\ell}}
\newcommand{\Zlb}{\overline{\Z}_{\ell}}

\newcommand{\Q}{\mathbb{Q}}
\newcommand{\Ql}{\Q_{\ell}}
\newcommand{\Qlb}{\overline{\Q}_{\ell}}

\newcommand{\RR}{\mathrm{R}}

\newcommand{\Ind}{\mathrm{Ind}}

\newcommand{\Gm}{\mathbb{G}_m}
\newcommand{\Ga}{\mathbb{G}_a}
\newcommand{\Spec}{\mathrm{Spec}}

\newcommand{\A}{\mathbb{A}}

\newcommand{\Gbf}{\mathbf{G}}
\newcommand{\Bbf}{\mathbf{B}}
\newcommand{\Ubf}{\mathbf{U}}
\newcommand{\Tbf}{\mathbf{T}}
\newcommand{\Wbf}{\mathbf{W}}
\newcommand{\Nbf}{\mathbf{N}}
\newcommand{\BwB}{\Bbf w\Bbf}
\newcommand{\Pbf}{\mathbf{P}}
\newcommand{\Vbf}{\mathbf{V}}
\newcommand{\Lbf}{\mathbf{L}}
\newcommand{\Zbf}{\mathbf{Z}}
\newcommand{\Xibf}{\mathbf{\Xi}}
\newcommand{\Gammabf}{\mathbf{\Gamma}}
\newcommand{\Hbf}{\mathbf{H}}
\newcommand{\Sbf}{\mathbf{S}}
\newcommand{\Mbf}{\mathbf{M}}
\newcommand{\Qbf}{\mathbf{Q}}
\newcommand{\Rbf}{\mathbf{R}}
\newcommand{\Xbf}{\mathbf{X}}
\newcommand{\UP}{\Ubf_{\Pbf}}
\newcommand{\UQ}{\Ubf_{\Qbf}}
\newcommand{\UR}{\Ubf_{\Rbf}}
\newcommand{\LP}{\Lbf_{\Pbf}}
\newcommand{\LQ}{\Lbf_{\Qbf}}
\newcommand{\LR}{\Lbf_{\Rbf}}

\newcommand{\pfrak}{\mathfrak{p}}
\newcommand{\qfrak}{\mathfrak{q}}

\newcommand{\ofrak}{\mathfrak{o}}

\newcommand{\Hom}{\mathrm{Hom}}
\newcommand{\D}{\mathbb{D}}
\newcommand{\id}{\mathrm{id}}

\newcommand{\Tr}{\mathrm{Tr}}

\newcommand{\ev}{\mathrm{ev}}

\newcommand{\End}{\mathrm{End}}

\newcommand{\ind}{\mathrm{ind}}
\newcommand{\Cur}{\mathrm{Cur}}

\newcommand{\Frob}{\mathrm{F}}

\newcommand{\Gal}{\mathrm{Gal}}
\newcommand{\Coh}{\mathrm{Coh}}

\newcommand{\Whitt}{\mathrm{Whitt}}
\newcommand{\Fun}{\mathrm{Fun}}

\newcommand{\Ad}{\mathrm{Ad}}
\newcommand{\Av}{\mathrm{Av}}

\newcommand{\opp}{\mathrm{op}}

\newcommand{\Spr}{\mathrm{Spr}}

\newcommand{\pr}{\mathrm{pr}}
\newcommand{\RGamma}{\mathrm{R\Gamma}}

\newcommand{\DD}{\mathrm{D}}
\newcommand{\pt}{\mathrm{pt}}
\newcommand{\spec}{\mathrm{spec}}

\newcommand{\unip}{\mathrm{unip}}
\newcommand{\Perv}{\mathrm{Perv}}
\newcommand{\Stab}{\mathrm{Stab}}

\newcommand{\V}{\mathbb{V}}

\newcommand{\can}{\mathrm{can}}

\newcommand{\Ext}{\mathrm{Ext}}

\newcommand{\Rep}{\mathrm{Rep}}

\newcommand{\der}{\mathrm{der}}

\newcommand{\HH}{\mathbb{H}}

\newcommand{\rev}{\mathrm{rev}}
\newcommand{\Irr}{\mathrm{Irr}}
\newcommand{\coh}{\mathrm{coh}}
\newcommand{\qcoh}{\mathrm{qcoh}}

\newcommand{\Ccal}{\mathcal{C}}
\newcommand{\ab}{\mathrm{ab}}

\newcommand{\Ocal}{\mathcal{O}}

\newcommand{\Exc}{\mathrm{Exc}}

\newcommand{\tot}{\mathrm{tot}}

\newcommand{\ad}{\mathrm{ad}}

\newcommand{\tors}{\mathrm{tors}}

\newcommand{\Alg}{\mathrm{Alg}}

\newcommand{\Tcal}{\mathcal{T}}

\newcommand{\loc}{\mathrm{loc}}

\newcommand{\Isog}{\mathrm{Isog}}
\newcommand{\Aut}{\mathrm{Aut}}
\newcommand{\Lcal}{\mathcal{L}}
\newcommand{\Zcal}{\mathcal{Z}}

\newcommand{\Mcal}{\mathcal{M}}
\newcommand{\BSBim}{\mathrm{BSBim}}
\newcommand{\SBim}{\mathrm{SBim}}

\newcommand{\tr}{\mathrm{tr}}
\newcommand{\ch}{\mathrm{ch}}

\newcommand{\sss}{\mathrm{ss}}
\newcommand{\Xcal}{\mathcal{X}}
\newcommand{\Ycal}{\mathcal{Y}}
\newcommand{\Morita}{\mathbf{Morita}}
\newcommand{\hc}{\mathrm{hc}}
\newcommand{\Ch}{\mathrm{Ch}}
\newcommand{\tilt}{\mathrm{tilt}}
\newcommand{\BS}{\mathrm{BS}}
\newcommand{\AS}{\mathrm{AS}}
\newcommand{\Pic}{\mathrm{Pic}}
\newcommand{\geo}{\mathrm{geo}}
\newcommand{\rat}{\mathrm{rat}}
\newcommand{\Bim}{\mathrm{Bim}}
\newcommand{\Loc}{\mathrm{Loc}}
\newcommand{\pos}{\mathrm{pos}}
\newcommand{\Int}{\mathrm{Int}}
\newcommand{\Arr}{\mathrm{Arr}}
\newcommand{\Ybf}{\mathbf{Y}}
\newcommand{\EE}{\mathbb{E}}
\newcommand{\simpconn}{\mathrm{sc}}
\newcommand{\tor}{\mathrm{tor}}
\newcommand{\Hbffrak}{\boldsymbol{\mathcal{H}}}
\newcommand{\KK}{\Ind\mathrm{K}^b}
\newcommand{\MM}{\mathbb{M}}
\newcommand{\TT}{\mathbb{T}}
\newcommand{\rigid}{\mathrm{rig}}
\newcommand{\Pin}{\mathrm{Pin}}

\newcommand{\incl}{\mathrm{incl}}
\newcommand{\mult}{\mathrm{mult}}
\newcommand{\Acal}{\mathcal{A}}
\newcommand{\Ncal}{\mathcal{N}}
\newcommand{\Dcal}{\mathcal{D}}
\newcommand{\Out}{\mathrm{Out}}
\newcommand{\unit}{\mathrm{unit}}
\newcommand{\counit}{\mathrm{counit}}
\newcommand{\cycl}{\mathrm{cycl}}
\newcommand{\colondef}{\mathpunct{:}}

\title{Jordan decomposition, categorical traces and Deligne--Lusztig theory}
\author{Arnaud Eteve}

\begin{document}

\maketitle

\textbf{Abstract.} This paper is a continuation of the author's previous work. We use the categorical trace formalism to construct the categorical Jordan decomposition for representations of finite groups of Lie type. As a second application, we study the endomorphism algebra of the Gelfand--Graev representation and recover a result of Li and Shotton--Li. 

\tableofcontents

\section{Introduction}

\subsection{Overview}

Let $\Gbf$ be a reductive group over $\Fqb$ and let $\Frob : \Gbf \to \Gbf$ be a Frobenius endomorphism coming from some $\Fq$-structure. The goal of this paper is to begin drawing consequences, for the representation theory of $\Gbf^{\Frob}$, of the main theorem of \cite{EteveFreeMonodromic} about the categorical trace of Frobenius on the Hecke category of $\Gbf$. Let $\ell \neq p$ be a prime and let $\Lambda \in \{\Flb, \Zlb, \Qlb\}$ be a coefficient ring. The Hecke category $\HH^{\Gbf}_{\Lambda}$ is a $\Lambda$-linear category which geometrizes the Yokonuma-Hecke algebra \cite{Yokonuma}. Let $\Bbf = \Tbf\Ubf$ be a Borel pair of $\Gbf$ that is stable under $\Frob$. The Hecke category is defined as a certain category of étale sheaves of $\Lambda$-modules on $\Ubf \backslash \Gbf/\Ubf$. It comes equipped with an endomorphism 
$$\Frob^* : \HH^{\Gbf}_{\Lambda} \to \HH^{\Gbf}_{\Lambda}.$$
The categorical trace of a monoidal endomorphism is a categorification of the classical notion of Hochschild homology. The categorical trace $\Tr(\Frob^*, \HH^{\Gbf}_{\Lambda})$ is a category equipped with a functor
$$\tr : \HH^{\Gbf}_{\Lambda} \to \Tr(\Frob^*, \HH^{\Gbf}_{\Lambda})$$ 
together with isomorphisms $\tr(A * B) = \tr(\Frob^*(B) * A)$ for all objects $A,B \in \HH^{\Gbf}_{\Lambda}$ (and equipped with appropriate compatibilities between them). Such a functor $\tr$ is called an $\Frob^*$-trace functor and the categorical trace is the universal category receiving a trace functor. To make sense of all the higher structures involved, we work with $\infty$-categories in the sense of \cite{HigherTopos}. Categorical traces in this context have been studied in \cite{BenZviNadlerCharacterTheory}, \cite{HoyoisScherotzkeSibilla}, \cite{ToyModel}, and the idea of considering the trace of Frobenius was motivated by the Langlands program, see \cite{HowToInventShtukas} and \cite{Zhu}. The main theorem of \cite{EteveFreeMonodromic} can be stated as follows. 
\begin{thm}[\cite{EteveFreeMonodromic}]
There is a natural equivalence of categories 
$$\Tr(\Frob^*, \HH^{\Gbf}_{\Lambda}) \cong \DD(\pt/\Gbf^{\Frob}, \Lambda).$$
\end{thm}
The category $\DD(\pt/\Gbf^{\Frob}, \Lambda)$ is the category of étale $\Lambda$-sheaves on the finite stack $\pt/\Gbf^{\Frob}$ which by étale descent is canonically equivalent to $\DD(\Rep_{\Lambda} \Gbf^{\Frob})$, the full derived category of representations of $\Gbf^{\Frob}$ on $\Lambda$-modules. Furthermore, the canonical functor 
$$\HH^{\Gbf}_{\Lambda} \to \Tr(\Frob^*, \HH^{\Gbf}_{\Lambda})$$
is also described in \emph{loc. cit.} in terms of the geometry of $\Gbf$ and we show how to recover the usual constructions of Deligne--Lusztig theory. 

More recently, Zhu \cite{ZhuCategoricalTame} and Gaitsgory--Rozenblyum--Varshavsky \cite{DennisCategoricalTraceDLTheory} have also independently established categorical trace statements of the same shape and explained how to recast many Deligne--Lusztig type constructions within this abstract framework. The formalism we use is slightly different from theirs: we use free monodromic sheaves instead of cofree monodromic sheaves for the foundations of the theory. 

This paper has four main sections. 
\begin{enumerate}
\item In Section \ref{sec:HeckeCat}, we review the construction and the key structures of the free monodromic Hecke categories as previously studied by \cite{BezrukavnikovYun}, \cite{BezrukavnikovRicheSoergelTheory}, \cite{Gouttard}, \cite{EteveFreeMonodromic}. The main goal of this section is to prove the endoscopy result, Theorem \ref{thm:EndoscopyHeckeCat}, which improves on \cite{Gouttard} and \cite{LusztigYun}.
\item In Section \ref{sec:parabolic-character-sheaves}, we review the formalism of parabolic character sheaves of Lusztig \cite{LusztigParabolicCS1}. The main use of the constructions of this section is to provide the tools to explain the compatibility of the categorical trace formalism with Deligne--Lusztig induction. 
\item In Section \ref{sectionJordan}, we construct a full Jordan decomposition, refining the ones of \cite{LusztigYun}, \cite{BonnafeRouquier} and \cite{BonnafeDatRouquier}. This construction is obtained by taking the trace of Frobenius on the equivalence of Theorem \ref{thm:EndoscopyHeckeCat}. 
\item In Section \ref{sectionGelfandGraev}, we apply the formalism to the study of the Gelfand--Graev representation of $\Gbf^{\Frob}$ when $\Gbf$ has connected center and $\ell$ is good for $\Gbf$. We give new geometric proofs of the main theorems of \cite{Dudas} and \cite{ShottonLi}. 
\end{enumerate}

We now review in more detail each of these sections. 

\subsection{Hecke categories}

The goal of Section \ref{sec:HeckeCat} is to review the construction and the key properties of the free monodromic Hecke categories as constructed in \cite{BezrukavnikovYun}, \cite{BezrukavnikovRicheSoergelTheory}, \cite{Gouttard} and \cite{EteveFreeMonodromic} and also to extend the endoscopy results of \cite{LusztigYun} and \cite{Gouttard}. 

\begin{defi}
The (free monodromic) Hecke category is the category 
$$\HH^{\Gbf}_{\Lambda, \tot} = \bigoplus_{s, s' \in \Ch(\Tbf)} \DD(\Ubf \backslash \Gbf/\Ubf)_{[s,s']},$$
where $\DD(\Ubf \backslash \Gbf/\Ubf)_{[s,s']}$ denotes the category of $[s,s']$-free monodromic sheaves on the stack $\Ubf \backslash \Gbf/\Ubf$ and $\Ch(\Tbf)$ denotes the set of characters of $\pi_1(\Tbf,1) \to \Lambda^{\times}$ of finite order prime to $p$ (if $\Lambda = \Qlb$) or $p\ell$ (if $\Lambda = \Flb, \Zlb)$. We denote by $\HH_{[s,s']}^{\Gbf}$ the summand indexed by $[s,s']$. 
\end{defi} 
We review some of the construction of the free monodromic sheaves in Section \ref{sec:SemisimplePoints}. We fix a pinning of $(\Gbf,\Bbf,\Tbf,X)$ and we denote by $(\Gbf^*,\Bbf^*,\Tbf^*,X^*)$ the pinned reductive group that is Langlands dual to $\Gbf$. In Section \ref{sec:SemisimplePoints}, we also identify the set $\Ch(\Tbf)$ with a subset of $\Tbf^*$. 

Let $s \in \Tbf^*$, we consider the centralizer $\Zbf(s)$ of $s$ in $\Gbf^*$ and $\Zbf^{\circ}(s)$ the connected centralizer of $s$. Following \cite{LusztigYun}, we construct in Section \ref{sec:EndoscopicGroup} the endoscopic group of $s$, which we denote by $\Hbf_s$. 
This group satisfies 
\begin{enumerate}
\item the Langlands dual group of $\Hbf_s^{\circ}$ is $\Zbf^{\circ}(s)$,
\item there is a canonical isomorphism $\pi_0(\Zbf(s)) = \pi_0(\Hbf_s)$. 
\end{enumerate}
The group $\Hbf_s$ is defined canonically and $\Hbf^{\circ}_s$ comes equipped with a (relative) pinning $(\Hbf_s^{\circ}, \Bbf_s,\Tbf, X_s)$, see Definition \ref{def:RelativePinning}. Similarly to the Hecke category of $\Gbf$, we define
$$\HH^{\Hbf}_{[s,s]} = \DD(\Ubf_s \backslash \Hbf_s /\Ubf_s)_{[1,1]}$$ 
which is the category of free monodromic unipotent sheaves on $\Ubf_s \backslash \Hbf_s /\Ubf_s$. Note that this stack is potentially disconnected. Given two $\Wbf$-conjugate elements $s,s' \in \Ch(\Tbf)$, we construct a scheme $\Hbf_{[s,s']}$ which is naturally an $(\Hbf_s,\Hbf_{s'})$-bitorsor and we define 
$$\HH^{\Hbf}_{[s,s']} = \DD(\Ubf_s \backslash \Hbf_{[s,s']}/\Ubf_{s'})_{[1,1]}.$$
This category is naturally an $(\HH^{\Hbf}_{[s,s]}$-$\HH^{\Hbf}_{[s',s']})$ bimodule. 

Lusztig--Yun \cite{LusztigYun} produced a very efficient way to organize the data of all the categories $(\HH^{\Gbf}_{[s,s']})$ and $(\HH^{\Hbf}_{[s,s']})$. There is a groupoid $\Xibf$ (viewed as a discrete stack over $\Fqb$) with $\mathrm{Ob}(\Xibf) = \Ch(\Tbf)$ and morphisms given by 
$$\Hom_{\Xibf}(s,s') = \pi_0(\Hbf_{[s,s']}).$$
The set of connected components of $\Xibf$ is the set of $\Wbf$-orbits on $\Ch(\Tbf)$ and for such an orbit $\ofrak \subset \Ch(\Tbf)$, we denote by $\Xibf_{\ofrak}$ the corresponding connected component. We review the construction of this groupoid in \Cref{sec:blocks}. 

Given a $\Wbf$-orbit $\ofrak$, the collection of bitorsors $(\Hbf_{[s,s']})$ is then organized in a gerbe 
$$\Hbffrak_{\ofrak} \to \Xibf_{\ofrak}$$
which is banded by $\Hbf_{s}^{\circ}$ (for some choice $s \in \ofrak$), see \Cref{sec:EndoscopicGroup}. 

Let us also denote by $\Morita$ the $(\infty, 1)$-category whose objects are monoidal presentable $\Lambda$-categories, whose morphisms are bimodules, and whose composition is given by relative tensor products. 

The collection of categories $(\HH^{\Gbf}_{[s,s']})$ and $(\HH^{\Hbf}_{[s,s']})$ then canonically defines functors 
$$\HH_{\ofrak}^{\Gbf}, \HH_{\ofrak}^{\Hbf} : \Xibf_{\ofrak}^{\opp} \to \Morita,$$
we refer to \Cref{sec:block-and-morita}.

\begin{thm}[Endoscopy for Hecke categories]\label{thm:EndoscopyHeckeCatIntro}
If $\Lambda = \Flb, \Zlb$, assume that $\ell$ is not a torsion prime for $\Gbf^*$ (see Definition \ref{def:torsion-prime}). We fix a $\Wbf$-orbit $\ofrak \subset \Ch(\Tbf)$. 
\begin{enumerate}
\item There is an isomorphism of functors $\Xibf^{\opp}_{\ofrak} \to \Morita$
$$\EE_{\ofrak} : \HH_{\ofrak}^{\Gbf} \cong \HH_{\ofrak}^{\Hbf}.$$
\item The equivalence $\EE_{\ofrak}$ depends only on a pinning of $\Gbf$ and no other extra data. 
\end{enumerate}
\end{thm}

If $\Lambda = \Qlb$, the first part of this Theorem was proved in \cite{LusztigYun} for `equivariant' Hecke categories instead of the free monodromic ones. Although the setting is slightly different, the same proof applies here. The second part of the Theorem is \Cref{corol:dependency-on-pinning}. Unraveling the data contained in such an isomorphism, we have for all $s,s' \in \ofrak$, a monoidal equivalence 
$$\EE_s : \HH^{\Gbf}_{[s,s]} \cong \HH^{\Hbf}_{[s,s]}$$
and an isomorphism of bimodules 
$$\EE_{[s,s']} : \HH^{\Gbf}_{[s,s']} \cong \HH^{\Hbf}_{[s,s']}.$$

The core idea here is to relate both sides appearing in the equivalences of Theorem \ref{thm:EndoscopyHeckeCatIntro} to suitable categories of Soergel bimodules. In Section \ref{sec:SoergelBim}, we construct for every pair of elements $(s,s') \in \Ch(\Tbf)$ a category $\SBim_{[s,s']}$ of `Soergel bimodules' in the same way that \cite{Gouttard} and \cite{LusztigYun} construct one. As before, we organize these into a functor 
$$\SBim : \Xibf^{\opp} \to \Morita.$$

\begin{thm}[Endoscopy for Hecke categories, Theorem \ref{thm:EndoscopyHeckeCat}] 
If $\Lambda = \Flb, \Zlb$, assume that $\ell$ is not a torsion prime for $\Gbf^*$ (see Definition \ref{def:torsion-prime}). 
There is an isomorphism of functors $\Xibf_{\ofrak}^{\opp} \to \Morita$,
$$\HH_{\ofrak}^{\Gbf} \cong \SBim \cong \HH_{\ofrak}^{\Hbf}.$$
\end{thm}

The proof of this theorem combines standard ideas, most of them due to Soergel \cite{Soergel} and \cite{BezrukavnikovYun}; we review its structure in Section \ref{sec:SoergelRealizationForG}. The equivalences on the $\Gbf$-side are due to \cite{Gouttard}, while the construction on the endoscopic side is new. The main difficulty stems from the disconnectedness of $\Hbf_s$ and, consequently, from the absence of \emph{a priori} canonical trivializations of the Whittaker model on the endoscopic side. We show in Section \ref{sec:TrivializationWhittModelOfH} that, under our hypothesis on $\ell$, such a trivialization can be canonically constructed. In recent work, Sandvik \cite{EndoscopyHeckeCatSandvik} and Dhillon--Li--Yun--Zhu \cite{DLZY} have produced endoscopic equivalences for affine Hecke categories, addressing similar difficulties caused by disconnected endoscopic groups.

\subsection{Jordan decomposition}

\subsubsection{Deligne--Lusztig theory, parabolic character sheaves and traces}

Before stating the Jordan decomposition Theorem, we explain how to recast Deligne--Lusztig theory using the trace formalism. Recall that we have fixed $\Gbf, \Bbf, \Tbf$ and $\Frob : \Gbf \to \Gbf$ that stabilizes $\Bbf, \Tbf$. In this situation Lang's theorem guarantees that there is a canonical isomorphism of stacks:
$$\frac{\Gbf}{\Ad_{\Frob}\Gbf} \cong \pt/\Gbf^{\Frob},$$ 
where $\Ad_{\Frob}$ denotes $\Frob$-twisted conjugation, given by $\Ad_{\Frob}(g)(x) = gx\Frob(g^{-1})$. 
Consequently there is an equivalence of categories 
$$\DD(\frac{\Gbf}{\Ad_{\Frob}\Gbf}, \Lambda) \cong \DD(\Rep_{\Lambda}\Gbf^{\Frob})$$
where, as before, the LHS denotes the category of étale sheaves of $\Lambda$-modules on the corresponding quotient. 

In \cite{EteveFreeMonodromic} and \cite{EteveTilting}, we have studied the stack $\frac{\Ubf \backslash \Gbf /\Ubf}{\Ad_{\Frob}(\Tbf)}$ and the horocycle correspondence 
$$\frac{\Ubf \backslash \Gbf /\Ubf}{\Ad_{\Frob}(\Tbf)} \xleftarrow{\alpha} \frac{\Gbf}{\Ad_{\Frob}\Bbf} \xrightarrow{\beta} \frac{\Gbf}{\Ad_{\Frob}\Gbf}.$$
In \emph{loc. cit.}, we observe that the Bruhat stratification of $\Gbf$ defines a stratification of $\frac{\Ubf \backslash \Gbf /\Ubf}{\Ad_{\Frob}(\Tbf)}$ whose strata are isomorphic to $\pt/\Tbf^{w\Frob} \rtimes (\Ubf \cap {^w}\Ubf)$ for varying elements $w \in \Wbf$. Let us denote by $j_w$ the corresponding inclusion.  Furthermore, the functor 
\begin{align*}
\DD(\Rep_{\Lambda}(\Tbf^{w\Frob})) \cong \DD(\pt/\Tbf^{w\Frob} \rtimes (\Ubf \cap {^w}\Ubf), \Lambda) &\to \DD(\frac{\Gbf}{\Ad_{\Frob}\Gbf}, \Lambda) \\
A &\mapsto \beta_!\alpha^*j_{w,!}A,
\end{align*}
is isomorphic to the usual Deligne--Lusztig induction functor, see \Cref{prop:ComparisonTraceAndDeligneLusztig}. 

The theory of parabolic character sheaves of Lusztig \cite{LusztigParabolicCS1} studies the parabolic version of the stack $\frac{\Ubf \backslash \Gbf /\Ubf}{\Ad_{\Frob}(\Tbf)}$. Let $\Pbf = \Lbf\UP$ be a standard parabolic with standard Levi. The corresponding parabolic horocycle stack is 
$$\frac{\UP \backslash \Gbf /\Frob(\UP)}{\Ad_{\Frob}\Lbf}.$$
The main result of \cite{LusztigParabolicCS1} proves that this stack admits a useful stratification, which we record in the following Proposition. 

\begin{prop}[\Cref{corol:StructureOfStrata}, \Cref{lem:ConnectednessParabCS}]
The stack $\frac{\UP \backslash \Gbf /\Frob(\UP)}{\Ad_{\Frob}\Lbf}$ has a stratification, which refines the one coming from $(\Pbf, \Frob(\Pbf))$-double cosets, and whose strata are isomorphic to stacks of the form 
$$\pt/\Mbf^{w\Frob} \rtimes \Vbf$$
where $\Mbf \subset \Lbf$ is a standard Levi subgroup, $w \in \Wbf$ and $\Vbf$ is a connected unipotent group. 
\end{prop} 

\begin{rque} 
The construction of the stratification is due to Lusztig. As far as we know, the description of the strata is new, specifically the statement about the connectedness of the unipotent groups that appear. 
\end{rque} 

Generalizing the horocycle correspondence, there is a parabolic version, given by 
$$\frac{\UP \backslash \Gbf /\Frob(\UP)}{\Ad_{\Frob}(\Lbf)} \xleftarrow{\alpha_{\Pbf}} \frac{\Gbf}{\Ad_{\Frob}\Pbf} \xrightarrow{\beta_{\Pbf}} \frac{\Gbf}{\Ad_{\Frob}\Gbf}.$$

 Let $v \in \Wbf$ and $\dot{v}$ be a lift of $v$ such that $\Ad(\dot{v}) \circ \Frob(\Lbf) = \Lbf$. Then the corresponding double coset is a stratum of Lusztig's stratification. Let us denote by $j_v$ the inclusion of the corresponding stratum. The corresponding stratum is isomorphic to 
$$\pt/(\Lbf^{\dot{v}\Frob} \rtimes (\UP \cap {^v}\Frob(\UP)))$$
and the functor 
\begin{align*}
\DD(\Rep_{\Lambda}(\Lbf^{\dot{v}\Frob})) \cong \DD(\pt/\Lbf^{\dot{v}\Frob} \rtimes (\UP \cap {^v}\Frob(\UP)), \Lambda) &\to \DD(\frac{\Gbf}{\Ad_{\Frob}\Gbf}, \Lambda) \\
A &\mapsto \beta_{\Pbf,!}\alpha^*_{\Pbf}j_{v,!}A,
\end{align*}
is isomorphic to the usual parabolic Deligne--Lusztig induction, see \Cref{prop:ComparisonTraceAndDeligneLusztig}.

The link with the categorical trace formalism is then given as follows. Let us denote by $\HH^{\Lbf}_{\tot} \subset \HH^{\Gbf}_{\tot}$ the Hecke category for the Levi $\Lbf$, see \Cref{sec:CompactWithParabInduction}.
\begin{thm}[\Cref{thm:calculation-traces}]
There are canonical equivalences 
\begin{enumerate}
\item $\Tr(\Frob^*,\HH^{\Gbf}_{\tot}) \cong \DD(\Rep_{\Lambda}\Gbf^{\Frob})$, 
\item $\Tr(\HH^{\Gbf}_{\tot, \Frob},\HH^{\Lbf}_{\tot}) \cong \DD(\frac{\UP \backslash \Gbf /\Frob(\UP)}{\Ad_{\Frob}(\Lbf)})$, where the LHS is the trace of $\HH^{\Gbf}_{\tot, \Frob}$ viewed as an $\HH^{\Lbf}$-bimodule (see \Cref{sec:recollection-traces}). 
\end{enumerate}
Furthermore, the inclusion $\HH^{\Lbf}_{\tot} \subset \HH^{\Gbf}_{\tot}$ induces a morphism on traces 
$$\Tr(\HH^{\Gbf}_{\tot, \Frob},\HH^{\Lbf}_{\tot}) \to \Tr(\Frob^*,\HH^{\Gbf}_{\tot})$$
which under the above equivalences is isomorphic to $\beta_{\Pbf,!}\alpha_{\Pbf}^*$.
\end{thm} 

\begin{rque} 
The above theorem is an $\Frob$-twisted version of the corresponding Theorem of \cite{LiNadlerYun}. We prove it in \Cref{sec:traces-calculations}. This Theorem can also be deduced from the results of \cite{ZhuCategoricalTame} and \cite{DennisCategoricalTraceDLTheory}, using a different formalism for Hecke categories.
\end{rque}

\subsubsection{Statements}

As a first application of the trace formalism, we construct a Jordan decomposition for the representation theory of $\Gbf^{\Frob}$. Recall that we denote by $\Gbf^*$ the group that is Langlands dual to $\Gbf$ and let us denote by $\Frob^* : \Gbf^* \to \Gbf^*$ the isogeny dual to $\Frob$. As is standard in Deligne--Lusztig theory, let us introduce the following equivalence relations on $(\Gbf^*_{\mathrm{ss}})^{\Frob^*}$, the set of semisimple elements in $(\Gbf^*)^{\Frob^*}$. 

\begin{defi}[\Cref{def:geometric-and-rational-conj}]
Let $s,t \in (\Gbf^*_{\mathrm{ss}})^{\Frob^*}$.
\begin{enumerate}
\item If $\Lambda = \Qlb$, we say that $s$ and $t$ are geometrically conjugate if $s$ and $t$ are conjugate in $\Gbf^*$, and we write $s \sim_{\Qlb, \geo} t$. We say that $s$ and $t$ are rationally conjugate if $s$ and $t$ are $(\Gbf^*)^{\Frob^*}$-conjugate, and we write $s \sim_{\Qlb, \rat} t$.
\item If $\Lambda = \Zlb, \Flb$, we say that $s$ and $t$ are geometrically conjugate if the prime-to-$\ell$ parts of $s$ and $t$ are $\Gbf^*$-conjugate, and we write $s \sim_{\Lambda, \geo} t$. We say that $s$ and $t$ are rationally conjugate if the prime-to-$\ell$ parts of $s$ and $t$ are $(\Gbf^*)^{\Frob^*}$-conjugate, and we write $s \sim_{\Lambda, \rat} t$. 
\end{enumerate}
\end{defi}

\begin{rque}[see \Cref{sec:the-group-gg}]
There is a bijection 
$$\Ch(\Tbf)/\Wbf \cong \Gbf^*_{\sss, \Lambda}/\sim_{\Lambda, \geo}$$
between the set of $\Wbf$-orbits in $\Ch(\Tbf)$ and the set of $\Gbf^*$-conjugacy classes of semisimple elements of order invertible in $\Lambda$. We shall identify those two sets using this bijection. Furthermore, the set of $\Gbf^*$-orbits that intersect $(\Gbf^*)^{\Frob^*}$ are exactly the $\Wbf$-orbits $\ofrak \subset \Ch(\Tbf)$ that are stable under $\Frob$, i.e. $\Frob(\ofrak) = \ofrak$. For such an orbit, we shall usually write $\ofrak'_{\rat} \subset \ofrak$ for a $(\Gbf^*)^{\Frob^*}$-stable orbit inside. 
\end{rque}

Let us fix a $\Frob$-stable $\Wbf$-orbit $\ofrak \subset \Ch(\Tbf)$. Recall that we have a morphism of stacks 
$$\Hbffrak_{\ofrak} \to \Xibf_{\ofrak}$$
and since $\ofrak$ is $\Frob$-stable, they both come equipped with a morphism $\Frob : \Hbffrak_{\ofrak}\to \Hbffrak_{\ofrak}$. We shall denote by $\Lcal_{\Frob}\Hbffrak_{\ofrak}$ the $\Frob$-loop stack, this is the pullback 
\[\begin{tikzcd}
	{\Lcal_{\Frob}\Hbffrak_{\ofrak}} & {\Hbffrak_{\ofrak}} \\
	{\Hbffrak_{\ofrak}} & {\Hbffrak_{\ofrak} \times \Hbffrak_{\ofrak}.}
	\arrow[from=1-1, to=1-2]
	\arrow[from=1-1, to=2-1]
	\arrow["{\id \times \Frob}", from=1-2, to=2-2]
	\arrow["\Delta"', from=2-1, to=2-2]
\end{tikzcd}\]
In particular, there is a well-defined morphism of stacks 
$$\Lcal_{\Frob}\Hbffrak_{\ofrak} \to \Lcal_{\Frob}\Xibf_{\ofrak}$$
which equips the category $\DD(\Lcal_{\Frob}\Hbffrak_{\ofrak}, \Lambda)$ with an action of $\DD(\Lcal_{\Frob}\Xibf_{\ofrak}, \Lambda)$.

After choosing some $s \in \ofrak$, this loop stack is isomorphic to 
$$\Lcal_{\Frob}\Hbffrak_{\ofrak} \cong \frac{\Hbf_{[s,\Frob(s)]}}{\Ad_{\Frob}\Hbf_s}.$$
Note that the stack $\frac{\Hbf_{[s,\Frob(s)]}}{\Ad_{\Frob}\Hbf_s}$ need not be connected. In general, its connected components are in bijection with the set of $(\Gbf^*)^{\Frob^*}$-orbits in $\ofrak$ (see \Cref{lem:RationalConjClasses}). Given some element $\dot{\beta} \in \Hbf_{[s, \Frob(s)]}$, there is an isomorphism 
$$\frac{\Hbf_{[s,\Frob(s)]}}{\Ad_{\Frob}\Hbf_s} = \bigsqcup_{[\gamma] \in H^1(\Ad(\dot{\beta})\Frob, \pi_0(\Hbf_{[s,s]}))} \pt/\Hbf_{[s,s]}^{\gamma\Ad(\dot{\beta})\Frob}.$$
Here, for each cohomology class $[\gamma]$, one chooses a cocycle representative $\gamma$. Changing that representative conjugates the corresponding twisted Frobenius, so the displayed classifying stack is well-defined up to canonical equivalence.

\begin{thm}[Jordan decomposition]\label{thm:EndoscopyRepIntro}
If $\Lambda = \Zlb, \Flb$, we assume that $\ell$ is not a torsion prime for $\Gbf^*$ (see Definition \ref{def:torsion-prime}). 
\begin{enumerate}
\item There is a direct sum decomposition into geometric series
$$\DD(\Rep_{\Lambda} \Gbf^{\Frob}) = \bigoplus_{\ofrak \in (\Gbf^*)^{\Frob^*}/\sim_{\Lambda, \geo}} \DD^{\ofrak}(\Rep_{\Lambda}\Gbf^{\Frob}).$$
Furthermore, the category $\DD^{\ofrak}(\Rep_{\Lambda}\Gbf^{\Frob})$ is generated by the collection of objects of the form $\tr_{\Frob}(\Delta_{w,s})$ for $w \in \Wbf$ and $s \in \ofrak$ such that $w\Frob(s)=s$, where $\Delta_{w,s} \in \HH^{\Gbf}_{\tot}$ is the corresponding standard object, see \Cref{cor:TraceVsDeligneLusztig}. 
\item Each geometric series, associated with a $\Frob$-stable $\Wbf$-orbit $\ofrak$, is equipped with a canonical action of $\DD(\Lcal_{\Frob}\Xibf_{\ofrak}, \Lambda)$. This yields in particular a direct sum decomposition in rational series 
$$\DD^{\ofrak}(\Rep_{\Lambda}\Gbf^{\Frob}) = \bigoplus_{\ofrak'_{\rat} \subset \ofrak}\DD^{\ofrak'_{\rat}}(\Rep_{\Lambda}\Gbf^{\Frob}).$$
Each of the direct summands $\DD^{\ofrak'_{\rat}}(\Rep_{\Lambda}\Gbf^{\Frob})$ is generated by the collection of objects of the form $\tr_{\Frob}(\Delta_{w,\Frob(s)})$ for $(s,w) \in \Lcal_{\Frob}\Xibf_{\ofrak}$ in the connected component corresponding to $\ofrak_{\rat}' \subset \ofrak$. Here $s$ is the base point of the loop and $w\Frob(s)=s$, so that $\Delta_{w,\Frob(s)}$ belongs to the block $\HH_{[s,\Frob(s)]}$.
\item Let $\ofrak \subset \Ch(\Tbf)$ be a $\Frob$-stable $\Wbf$-orbit; there exists a $\DD(\Lcal_{\Frob}\Xibf_{\ofrak})$-linear equivalence, which we call the endoscopic equivalence, 
$$\EE_{\ofrak}^{\Frob} : \DD^{\ofrak}(\Rep_{\Lambda}\Gbf^{\Frob}) \cong \DD^{\unip}(\Lcal_{\Frob}\Hbffrak_{\ofrak}, \Lambda)$$
where the RHS is the category of sheaves that, for some choice of $s$ and $\dot{\beta}$ as above, give unipotent representations of $\Hbf_{[s,s]}^{\gamma\Ad(\dot{\beta})\Frob}$ for all $\gamma$, see \Cref{sec:unipotent-representations-endoscopic-side} for the definition. 
\item Given a rational orbit $\ofrak_{\rat}' \subset \ofrak$, the above equivalence $\EE_{\ofrak}^{\Frob}$ yields equivalences
$$\DD^{\ofrak_{\rat}'}(\Rep_{\Lambda}\Gbf^{\Frob}) \cong \DD^{\unip}(\Lcal_{\Frob}\Hbffrak_{\ofrak}^{\ofrak_{\rat}'}, \Lambda)$$
where $\Lcal_{\Frob}\Hbffrak_{\ofrak}^{\ofrak_{\rat}'}$ denotes the connected component associated with $\ofrak'_{\rat}$. 
\item Given a standard parabolic $\Pbf = \Lbf\UP$, there is a decomposition in geometric series 
$$\DD(\frac{\UP \backslash \Gbf/\Frob(\UP)}{\Ad_{\Frob}\Lbf}, \Lambda) = \bigoplus_{\ofrak \in (\Gbf^*)^{\Frob^*}/\sim_{\Lambda, \geo}} \DD^{\ofrak}(\frac{\UP \backslash \Gbf/\Frob(\UP)}{\Ad_{\Frob}\Lbf}, \Lambda)$$
that glues together the geometric series of the various strata. There exists an endoscopic equivalence 
$$\EE_{\ofrak}^{\Frob, \Lbf} : \DD^{\ofrak}(\frac{\UP \backslash \Gbf/\Frob(\UP)}{\Ad_{\Frob}\Lbf}, \Lambda) \cong \DD^{\unip}_{(\UP, \UP)}(\Lcal_{\Frob}^{\Lbf}\Hbffrak_{\ofrak}, \Lambda)$$
where the RHS is defined in \Cref{defi:unipotent-parabolic-cs-endoscopic-side}. 
\item Let $v \in \Wbf$ and $\dot{v}$ be an element such that $\Ad(\dot{v})\Frob(\Lbf) = \Lbf$ and let $\ofrak^{\Lbf} \subset \ofrak$ be a $\Ad(v)\Frob$-stable $\Wbf^{\Lbf}$-orbit, where $\Wbf^{\Lbf}$ is the Weyl group of $\Lbf$. The following diagram is canonically commutative 
\[\begin{tikzcd}
	{\DD^{\ofrak^{\Lbf}}(\Rep_{\Lambda}\Lbf^{\dot{v}\Frob}) } & {\DD^{\unip}(\Lcal_{v\Frob}\Hbffrak_{\ofrak^{\Lbf}}^{\Lbf}, \Lambda)} \\
	{\DD^{\ofrak}(\Rep_{\Lambda}\Gbf^{\Frob})} & {\DD^{\unip}(\Lcal_{\Frob}\Hbffrak_{\ofrak}, \Lambda)}
	\arrow["{\EE_{\ofrak^{\Lbf}}^{v\Frob}}", from=1-1, to=1-2]
	\arrow["{\RR_{\dot{v}, \Pbf}}"', from=1-1, to=2-1]
	\arrow["{\RR_{\dot{v}, \Pbf}^{\Hbf}}", from=1-2, to=2-2]
	\arrow["{\EE_{\ofrak}^{\Frob}}"', from=2-1, to=2-2]
\end{tikzcd}\]
where the vertical functor on the RHS is defined in \Cref{defi:DL-induction-endoscopic-side}.
\item The equivalence $\EE_{\ofrak}^{\Frob}$ sends the (compactly supported) cohomology of a Deligne--Lusztig variety to the (compactly supported) cohomology of a Deligne--Lusztig variety.
\item The equivalence $\EE_{\ofrak}^{\Frob}$ is $t$-exact.
\item The equivalence $\EE_{\ofrak}^{\Frob}$ is compatible with inversion of $\ell$, see \Cref{sec:compatibility-with-inversion-of-ell} for what this means. 
\item The equivalence $\EE_{\ofrak}^{\Frob}$ is compatible with isogenies, see \Cref{sec:compatibility-with-isogenies} for what this means. 
\item The equivalence $\EE_{\ofrak}^{\Frob}$ is compatible with automorphisms of $(\Gbf, \Frob)$, see \Cref{thm:compat-automorphisms} for what this means. 
\end{enumerate}
\end{thm}

Not all parts of this statement are new, but we do give a proof of all points of this theorem in Section \ref{sectionJordan} using the formalism of categorical traces. Let us now indicate which parts of this theorem are new and which parts were already known. 
\begin{enumerate}
\item Points $(i)$ and $(ii)$ are the well-known partitions in geometric and rational series. If $\Lambda = \Qlb$, they were constructed in \cite{DeligneLusztig}, while for $\Lambda = \Flb, \Zlb$ this was done in \cite{BroueMichel}. 
\item If $\Lambda = \Qlb$, then point $(iii)$ was proved in \cite{LusztigYun}, building on \cite{LusztigCatCenter1} and \cite{LusztigCatCenter2}. However, they do not address the compatibility with Deligne--Lusztig induction.
\item For a general $\Lambda$, one of the first results is \cite{BroueTorus} when $\Zbf(s)$ is a maximal torus. If $\Zbf(s)$ is a Levi of $\Gbf^*$, then points $(iii)$--$(xi)$ are proved in \cite{BonnafeRouquier}, \cite{BonnafeDatRouquier} and \cite{Ruhstorfer} without restrictions on $\ell$, and the resulting equivalence is splendid in the sense of Rickard \cite{Rickard}.	
\end{enumerate}

\begin{rque}
In \cite{BonnafeDatRouquier} and \cite{BonnafeRouquier}, the functor realizing the equivalence $(iii)$ is given by Deligne--Lusztig induction. More precisely, the equivalence $(iii)$ is given by a bimodule; in those works, this bimodule is the top cohomology group of a parabolic Deligne--Lusztig variety. In contrast, our functor is obtained using abstract methods, which do not allow us to keep track of the bimodule realizing the equivalence. However, the compatibility with Deligne--Lusztig induction implies that our functor is the same as the one constructed by \cite{BonnafeRouquier} and \cite{BonnafeDatRouquier}, see \Cref{thm:comparison-with-BDR}.
\end{rque}

Let us now indicate where each point of Theorem \ref{thm:EndoscopyRepIntro} is proved. 
\begin{enumerate}
\item Point $(i)$ is proved in Section \ref{sec:GeometricSeries}.
\item Point $(ii)$ is proved in Section \ref{sec:RationalSeries}.
\item Points $(iii), (iv)$, $(v)$ and $(vi)$ are proved in Sections \ref{sec:EndoscopyReps} and \ref{sec:Variant}. 
\item Points $(vii)$, $(viii)$ and $(ix)$ are proved in Section \ref{sec:PropertiesEndoscopic}.
\item Point $(x)$ is proved in \Cref{sec:compatibility-with-isogenies}.
\item Point $(xi)$ is proved in \Cref{sec:compat-with-automorphisms}.
\end{enumerate}

The proof of Theorem \ref{thm:EndoscopyRepIntro} is obtained by taking the categorical trace of Theorem \ref{thm:EndoscopyHeckeCatIntro}. Using the trace formalism, we obtain additional structure on $\Rep_{\Lambda}\Gbf^{\Frob}$, which appears relevant to the modular representation theory of $\Gbf^{\Frob}$. 

\begin{thm}[Theorem \ref{thm:LinearityForReps}]\label{thm:LinearityIntro}
There is a canonical $\Ocal((\Tbf^{\vee}\sslash \Wbf)^{\Frob^{\vee}})$-linear structure on $\Rep_{\Lambda} \Gbf^{\Frob}$, where $\Tbf^{\vee}$ is the torus dual to $\Tbf$ defined over $\Lambda$ and $\Frob^{\vee}$ is the morphism dual to the restriction of $\Frob$ to $\Tbf$. 
\end{thm}

\begin{rque}
The decomposition in geometric series is a formal consequence of this structure, as the connected components of the scheme $(\Tbf^{\vee}\sslash \Wbf)^{\Frob^{\vee}}$ are in canonical bijection with the $\Frob^*$-stable conjugacy classes of semisimple elements in $\Gbf^*$.
\end{rque}

\begin{rque}
If $\Lambda = \Qlb$, the algebra $\Ocal((\Tbf^{\vee}\sslash \Wbf)^{\Frob^{\vee}})$ is isomorphic to $\prod \Qlb$ and this structure yields nothing more than the decomposition in geometric series.
\end{rque}

\begin{rque}
The structure of Theorem \ref{thm:LinearityIntro} is a map 
$$\Ocal((\Tbf^{\vee}\sslash \Wbf)^{\Frob^{\vee}}) \to Z(\Lambda[\Gbf^{\Frob}]).$$ 
We understand this map as an analogue for finite groups of Lie type of the map $\Exc \to \mathfrak{Z}(\Gbf(E))$ defined by \cite{FarguesScholze} and \cite{GenestierLafforgue} from the excursion algebra to the Bernstein center of a $p$-adic group defined over a non-archimedean local field $E$. 
\end{rque}

\subsubsection{Applications}

Let us conclude this section with two (related) applications of Jordan decomposition: 
\begin{enumerate}
\item the unitriangularity of the decomposition matrix of $\Gbf^{\Frob}$, which was originally conjectured by Geck \cite{GeckConjecture}, 
\item the classification of modular irreducible representations of $\Gbf^{\Frob}$. 
\end{enumerate}

\begin{defi}[\Cref{defi:JordanIsoClasses}]

Let $\mathbb{J}(\Gbf^{\Frob}, \Flb)$ be the set of equivalence classes of pairs $(s, \rho)$ where $s \in (\Gbf^*)^{\Frob^*}$ is a semisimple element of order prime to $\ell$ and $\rho \in \Irr_{\Flb}^{\unip}(\Hbf_s^{\Frob})$ up to rational conjugation of $s$. 
\end{defi} 

\begin{thm}[\Cref{thm:classificationOfIrrep}]
Assume that $\ell$ is not a torsion prime for $\Gbf^*$, 
\begin{enumerate}
\item there is a bijection depending only on the chosen pinning of $\Gbf$ 
$$\Irr_{\Flb}(\Gbf^{\Frob}) \cong \mathbb{J}(\Gbf^{\Frob}, \Flb).$$
\item Assume further that $p$ and $\ell$ are good for $\Gbf$ and that $\Gbf$ has connected center. Then $\Irr_{\Flb}^{\unip}(\Hbf_s^{\Frob}) \cong \Irr_{\Qlb}^{\unip}(\Hbf_s^{\Frob})$.
\item Assume that $p$ and $\ell$ are good for $\Gbf$. Then Geck's conjecture holds, i.e., the decomposition matrix of $\Gbf^{\Frob}$ is unitriangular.
\end{enumerate}
\end{thm} 

Point $(ii)$ of the Theorem is the main Theorem of \cite{DudasBrunatTaylor}. 

\begin{conj}
\Cref{thm:classificationOfIrrep} holds for all $p$, all $\ell \neq p$, and all $\Gbf$. 
\end{conj}

\subsection{Endomorphisms of the Gelfand--Graev representation}

We give a second application of the trace formalism to the structure of the endomorphism algebra of the Gelfand--Graev representation. In this section, we assume that $\Gbf$ has connected center, that $\ell$ is good for $\Gbf$ (see Definition \ref{def:torsion-prime}), and that $\Lambda = \Zlb$. Let $\Ubf^-$ be the unipotent radical of the Borel opposite to $\Bbf$. Since $\Bbf$ is stable under $\Frob$, so is $\Ubf^-$. We let 
$$\psi : (\Ubf^-)^{\Frob} \to \Lambda^{\times}$$
be a generic character of $(\Ubf^-)^{\Frob}$. The Gelfand--Graev representation is the representation
$$\Gamma_{\psi} = \ind_{(\Ubf^-)^{\Frob}}^{\Gbf^{\Frob}} \psi.$$

The goal of Section \ref{sectionGelfandGraev} is to use the formalism of Hecke categories, horocycle correspondences and categorical traces to give geometric proofs of the following two theorems. 

\begin{thm}[\cite{Dudas}]\label{thmDudasintro}\footnote{Zhu \cite{ZhuCategoricalTame} gave a proof of this result in a similar way to our proof.}
Let $w \in \Wbf$ and let ${^*}R_w : \DD(\Lambda[\Gbf^{\Frob}]) \to \DD(\Lambda[\Tbf^{w\Frob}])$ be the Deligne--Lusztig restriction functor, see 
Proposition \ref{prop:ComparisonTraceAndDeligneLusztig}. Then there exists a $\Tbf^{w\Frob}$-linear isomorphism 
$${^*}R_w(\Gamma_{\psi}) = \Lambda[\Tbf^{w\Frob}][-\ell(w)].$$
\end{thm}

By functoriality, this theorem provides a map
$$\Cur_w : \End(\Gamma_{\psi}) \to \Lambda[\Tbf^{w\Frob}]$$
called the $w$-Curtis map. 
After choosing a trivialization of the roots of unity in $\Fqb^{\times}$, we can produce an isomorphism 
$$\Lambda[\Tbf^{w\Frob}] = \Ocal((\Tbf^{\vee})^{w\Frob^{\vee}}).$$
Here $\Tbf^{\vee}$ denotes the torus dual to $\Tbf$ defined over $\Lambda$, $\Frob^{\vee}$ denotes the morphism dual to $\Frob$, and $(\Tbf^{\vee})^{w\Frob^{\vee}}$ denotes the subscheme of $w\Frob^{\vee}$-fixed points; see Section \ref{sec:NotationGeometry}. 
There is a natural map of schemes over $\Lambda$
$$(\Tbf^{\vee})^{w\Frob^{\vee}} \to (\Tbf^{\vee}\sslash \Wbf)^{\Frob^{\vee}},$$
where $\Tbf^{\vee}\sslash \Wbf$ is the GIT quotient of $\Tbf^{\vee}$ by $\Wbf$. On rings of functions, this defines the $w$-spectral Curtis map 
$$\Cur_w^{\spec} : \Ocal((\Tbf^{\vee}\sslash \Wbf)^{\Frob^{\vee}}) \to \Ocal((\Tbf^{\vee})^{w\Frob^{\vee}}).$$ We set 
$$\Cur^{\spec} = \bigoplus_w \Cur_w^{\spec} : \Ocal((\Tbf^{\vee}\sslash \Wbf)^{\Frob^{\vee}}) \to \bigoplus_w \Ocal((\Tbf^{\vee})^{w\Frob^{\vee}}).$$

\begin{thm}[\cite{Li}, \cite{ShottonLi}]\label{thmShottonLiIntro}
There is a unique isomorphism $\End(\Gamma_{\psi}) \cong \Ocal((\Tbf^{\vee}\sslash \Wbf)^{\Frob^{\vee}})$ such that the following diagram commutes 
\[\begin{tikzcd}
	{\End(\Gamma_{\psi})} & {\Ocal((\Tbf^{\vee}\sslash \Wbf)^{\Frob^{\vee}})} \\
	{\bigoplus_{w\in \Wbf}\Lambda[\Tbf^{w\Frob}]} & {\bigoplus_{w\in \Wbf}\Ocal((\Tbf^{\vee})^{w\Frob^{\vee}}).}
	\arrow[from=1-1, to=1-2]
	\arrow["\Cur"', from=1-1, to=2-1]
	\arrow["{\Cur^{\spec}}", from=1-2, to=2-2]
	\arrow[Rightarrow, no head, from=2-1, to=2-2]
\end{tikzcd}\]
\end{thm}

We provide proofs of Theorems \ref{thmDudasintro} and \ref{thmShottonLiIntro} in Sections \ref{sec:ProofDudas} and \ref{sec:ProofShottonLi} respectively.  

\begin{rque}
It should be noted that our construction provides explicit maps in Theorem \ref{thmShottonLiIntro} which are automatically defined over $\Lambda$. By contrast, the existing proof first establishes the statement over $\Qlb$ and then examines the denominators that appear in specific bases to descend it to $\Zlb$. 
\end{rque}


\subsection{Notations}

\subsubsection{Geometry}\label{sec:NotationGeometry}
We fix $k = \Fqb$, the algebraic closure of a finite field. We will consider algebraic stacks of finite type over $k$. 

For a stack $X$ together with an endomorphism $f : X \to X$, we denote by $X^f$ the stack of $f$-fixed points. By definition, this is the stack fitting into the following Cartesian diagram 
\[\begin{tikzcd}
	{X^f} & X \\
	X & {X \times X}
	\arrow[from=1-1, to=1-2]
	\arrow[from=1-1, to=2-1]
	\arrow["{\id \times f}", from=1-2, to=2-2]
	\arrow["\Delta"', from=2-1, to=2-2]
\end{tikzcd}\]
where $\Delta$ is the diagonal of $X$. Given a group $\Gamma$ acting on an affine scheme $X = \Spec R$, we denote by $X\sslash \Gamma = \Spec R^{\Gamma}$ the GIT quotient of $X$ by $\Gamma$. 

\subsubsection{$\infty$-Categories}

We will work with $\infty$-categories as constructed in \cite{HigherTopos} and \cite{HigherAlgebra}. All categories are considered to be $\infty$-categories. We denote by $\Pr$ the category of presentable categories and for a ring $\Lambda$, we denote by $\Pr_{\Lambda}$ the category of $\Lambda$-linear presentable categories. For a scheme $X$ over $\Lambda$, we denote by $\DD_{\qcoh}(X)$ the derived category of quasi-coherent sheaves on $X$, by $\DD_{\coh}(X)$ the full subcategory of $\DD_{\qcoh}(X)$ of cohomologically bounded complexes with coherent cohomology and by $\Coh(X)$ the abelian category of coherent sheaves on $X$.

\subsubsection{Reductive groups}\label{sec:NotationsReductive}

We fix once and for all a connected reductive group $\Gbf$ over $k$, together with a Borel pair $\Bbf = \Tbf\Ubf \subset \Gbf$. We denote by $\Wbf$ the associated Weyl group. We fix a set of lifts $\dot{w} \in \Nbf(\Tbf)$ of the elements of $\Wbf$. 

We denote by $\Gbf_{\der}, \Gbf_{\simpconn}, \Gbf_{\ad}$ the derived subgroup, the simply connected cover of $\Gbf_{\der}$, and the adjoint quotient of $\Gbf$, respectively. We denote by $\Tbf_{\der}, \Tbf_{\simpconn}$, and $\Tbf_{\ad}$ the (inverse) images of $\Tbf$ in the respective groups. 

\subsubsection{Roots and pinning}
We introduce the following notation for root systems. With the chosen Borel pair $\Bbf = \Tbf\Ubf$ of $\Gbf$, we denote by 
$$(X^*, X_*, \Phi, \Phi^{\vee}, \Delta, \Delta^{\vee})$$ 
the based root datum of $\Gbf$, where $\Phi \subset X^* = X^*(\Tbf)$ (resp. $\Phi^{\vee} \subset X_* = X_*(\Tbf)$) denotes the set of roots (resp. coroots), $\Phi^{+}$ the positive roots (resp. $\Phi^{\vee,+}$ the positive coroots) and $\Delta$ the simple roots (resp. $\Delta^{\vee}$ the simple coroots). We fix a pinning $(x_{\alpha})_{\alpha \in \Delta}$ of $\Gbf$ where 
$$x_{\alpha} : \Ga \cong \Gbf_{\alpha}.$$ 

We shall denote by $\Z\Phi$ (resp. $\Z\Phi^{\vee}$) the lattices generated by $\Phi$ in $X^*(\Tbf)$ (resp. by $\Phi^{\vee}$ in $X_*(\Tbf)$). 

\subsubsection{Extended pinning}\label{sec:ExtendedPinning}
On top of the pinning, we fix a trivialization of all the root groups 
$$x_{\alpha} : \Ga \cong \Gbf_{\alpha}.$$
We do not require that these satisfy Chevalley-type relations, but we will assume that they are Frobenius-equivariant. We call this datum an extended pinning of $\Gbf$. Note that after choosing an ordering of the simple roots, using the Chevalley relations, one can produce such a datum out of the already chosen pinning of $\Gbf$. It is not clear that an extended pinning produced using the Chevalley relations can be made to be $\Frob$-equivariant.\footnote{The extended pinning will be used to make certain calculations explicit, but the resulting Theorem \ref{thm:EndoscopyRepIntro} does not depend on the extended part of the pinning.}

\subsubsection{Fundamental group and bad primes}

We recall that, for a reductive group, the Borovoi fundamental group $\pi_1(\Gbf)$ is defined to be 
$$\pi_1(\Gbf) = X_*(\Tbf)/X_*(\Tbf_{\simpconn}).$$
Note that this presentation depends on the choice of $\Tbf$; however, the group $\pi_1(\Gbf)$ depends only on $\Gbf$; see \cite{Borovoi}. 

\begin{defi}[\protect{\cite[Chapter E, Section I.4]{SeminarAlgebraicGroups}}]\label{def:torsion-prime}
A prime $p$ is a bad prime (resp. a torsion prime) if $p$ divides 
$$|(X^*(\Tbf)/\Z\Phi_1)_{\tor}| \ (\text{resp.} |(X_*(\Tbf)/\Z\Phi_1^{\vee})_{\tor}|)$$
for some subroot datum $(X^*, X_*, \Phi_1, \Phi_1^{\vee})$ of $(X^*, X_*, \Phi, \Phi^{\vee})$. 
\end{defi} 

\begin{rques}
\begin{enumerate}
\item As pointed out in \cite[Chapter E, Section I.4]{SeminarAlgebraicGroups}, for a subroot datum $(X^*, X_*, \Phi_1, \Phi_1^{\vee})$, only $\Phi_1$ is required to be integrally closed in $\Phi$.
\item If $\Gbf$ is semisimple and simply connected, then bad/torsion primes depend only on the types of the simple quotient of $\Gbf$, the primes in question are in $\{2,3,5\}$ as indicated in the tables in \emph{loc. cit.}.
\item In general, for $\Gbf$ an arbitrary reductive group, a prime is a torsion prime if it is either a torsion prime for $\Gbf_{\simpconn}$ or if it divides $|\pi_1(\Gbf)_{\tor}|$. 
\end{enumerate}
\end{rques}

\subsubsection{Frobenius and finite groups}

A Frobenius endomorphism $\Frob : \Gbf \to \Gbf$ is a finite purely inseparable isogeny such that a power of $\Frob$ is a $q^{\delta}$-Frobenius endomorphism of $\Gbf$ for some integer $\delta \geq 1$. We denote by $\Gbf^{\Frob}$ the finite discrete group of its fixed points. We assume that all the data we have chosen thus far is $\Frob$-equivariant, i.e., $\Frob(\Bbf, \Tbf) = (\Bbf, \Tbf)$ and the trivialization of all root groups makes the following diagram commutative 
\[\begin{tikzcd}
	\Ga & \Ga \\
	{\Gbf_{\alpha}} & {\Gbf_{\Frob(\alpha)}.}
	\arrow["{t \mapsto t^q}", from=1-1, to=1-2]
	\arrow["{x_{\alpha}}"', from=1-1, to=2-1]
	\arrow["{x_{\Frob(\alpha)}}", from=1-2, to=2-2]
	\arrow["\Frob"', from=2-1, to=2-2]
\end{tikzcd}\]

\subsubsection{Sheaves}

We fix a ring $\Lambda \in \{\Flb, \Zlb, \Qlb\}$ and we denote by $\Lambda_0 \in \{\Fl,\Zl, \Ql\} $ the prime subring of $\Lambda$. For a stack $X$ over $k$, we denote by $\DD(X, \Lambda)$ the category of ind-constructible sheaves of $\Lambda$-modules on $X$. If $\Lambda$ is clear from the context, we will write $\DD(X)$. 

\subsection{Acknowledgments}

We thank Jean-François Dat and Olivier Dudas for their continuous support during the preparation of this paper. Section \ref{sectionGelfandGraev} is an improvement on the corresponding section of the thesis of the author \cite{EteveThesis} which was directed by Jean-François Dat and of which Olivier Dudas was a referee. The author was supported by the Max Planck Institute for Mathematics during the preparation of this paper. The author thanks Tom Gannon, Jens Eberhardt, Simon Riche, Daniel Juteau, Dennis Gaitsgory, Jack Shotton, Cédric Bonnafé, Gunter Malle, Gerhard Hiss, Britta Späth, Marc Cabanes, Jay Taylor and Zhicheng Feng for comments and discussions during the preparation of this paper.

\section{Preliminaries}

\subsection{Characters of tori and semisimple points}\label{sec:SemisimplePoints}

Let $\Tbf$ be a torus over $k = \Fqb$. We recall some notation from \cite{EteveFreeMonodromic}. We let $\pi_1^t(\Tbf)$ be the tame fundamental group of $\Tbf$ at the point $1 \in \Tbf$, i.e., the largest prime-to-$p$ quotient of $\pi_1(\Tbf, 1)$. There is a canonical isomorphism 
$$\pi_1^t(\Tbf) = X_*(\Tbf) \otimes_{\Z} \hat{\Z}^{(p)}(1).$$
We denote by $\pi_1^{\ell}(\Tbf)$ the largest pro-$\ell$ quotient of $\pi_1(\Tbf)$, so that $\pi_1^{\ell}(\Tbf) = X_*(\Tbf) \otimes_{\Z} \Zl(1)$. 

We introduce the following ring $\RR_{\Tbf, \Lambda}$.
\begin{enumerate}
\item $\RR_{\Tbf, \Lambda} := \Lambda_0\llbracket \pi_1^{\ell}(\Tbf) \rrbracket \otimes_{\Lambda_0} \Lambda$ if $\Lambda = \Zlb, \Flb$, 
\item $\RR_{\Tbf, \Qlb} := \left(\RR_{\Tbf, \Zl}\otimes_{\Zl}\Ql\right)^{\wedge}\otimes_{\Ql}\Qlb$, where $\wedge$ refers to the completion along the augmentation ideal of $\RR_{\Tbf, \Zl}\otimes_{\Zl}\Ql$.
\end{enumerate}
We will write $\RR_{\Tbf}$ if $\Lambda$ is clear from the context. 

There is a canonical map 
$$\pi_1(\Tbf) \to \RR_{\Tbf}^{\times}.$$
Corresponding to this map, there is a locally constant sheaf $L_{\Tbf}$ of free rank-one $\RR_{\Tbf}$-modules on $\Tbf$, which is called the free monodromic unipotent local system on $\Tbf$. 

We denote by 
$$\Ch_{\Lambda}(\Tbf) := \{\chi : \pi_1^t(\Tbf) \to \Lambda^{\times} \ \text{of finite order prime to the residue characteristic of} \ \Lambda \}$$
i.e., the set of characters of $\pi_1^t(\Tbf)$ of order prime to $p$ if $\Lambda = \Qlb$ and of order prime to $p\ell$ if $\Lambda = \Zlb, \Flb$. We will drop the index $\Lambda$ if it is clear from the context.

We denote by $\Tbf^{\vee}$ (resp. $\Tbf^{\vee}_{\Lambda_0}$) the torus dual to $\Tbf$ defined over $\Lambda$ (resp. over $\Lambda_0$). 
\begin{enumerate}
\item There is a canonical isomorphism $\RR_{\Tbf} = \Ocal(\Tbf^{\vee, \wedge}_{\Lambda_0, 1}) \otimes_{\Lambda_0} \Lambda$, where $(-)^{\wedge}_1$ denotes the completion at the point $1 \in \Tbf^{\vee}_{\Lambda_0}$. 
\item After choosing a topological generator of $\pi_1^t(\Gm)$, we obtain a bijection 
\begin{equation}\label{eq:CharactersVsPoints}
\Ch(\Tbf) = \Tbf^{\vee, \tors}(\Lambda)
\end{equation}
where the right-hand side is the set of $\Lambda$-points of $\Tbf^{\vee}$ that are torsion of order prime to $p$ if $\Lambda = \Qlb$ and prime to $p\ell$ if $\Lambda = \Flb, \Zlb$. 
\end{enumerate}

Let $\chi : \pi_1^t(\Tbf) \to \Lambda^{\times}$ be an element of $\Ch(\Tbf)$ and let $s \in \Tbf^{\vee}$ be the point corresponding to $\chi$. We denote by $\Lcal_{s}$ the Kummer sheaf over $\Tbf$ corresponding to $\chi$. 

For a stack $X$ with an action of $\Tbf$, we denote by 
$$\DD(X/(\Tbf, L_{\Tbf} \otimes_{\Lambda} \Lcal_s), \RR_{\Tbf})$$ 
the category of free $s$-monodromic sheaves on $X$ as constructed in \cite{EteveFreeMonodromic}, \cite{Gouttard} and \cite[Appendix A]{BezrukavnikovYun}. 

\subsection{The group $\Gbf^*$}\label{sec:the-group-gg}

Recall that we have a pinned reductive group
$$(\Gbf, \Bbf, \Tbf, (x_{\alpha})_{\alpha \in \Delta}).$$

Corresponding to the pinned reductive group $\Gbf$, there is a unique pinned reductive group $\Gbf^*$ over $k$, equipped with a Borel pair $\Bbf^* = \Tbf^*\Ubf^*$ whose root datum is 
$$(X_*, X^*, \Phi^{\vee}, \Phi, \Delta^{\vee}, \Delta).$$
The torus $\Tbf^*$ is the torus dual to $\Tbf$ defined over $k$ and we have $X^*(\Tbf^*) = X_*(\Tbf)$. 
The choice of the topological generator of $\pi_1^t(\Gm)$ yields a trivialization of $k^{\times} \cong (\Q/\Z)^{(p)}$ and a bijection
$$\Ch(\Tbf)/\Wbf \cong \Gbf^*_{\mathrm{ss}, \Lambda}/\Gbf^*$$
where the right-hand side is the set of conjugacy classes of semisimple points of $\Gbf^*$ of order prime to $p$ if $\Lambda = \Qlb$ and of order prime to $p\ell$ if $\Lambda = \Zlb, \Flb$. 

\subsection{The root system attached to a semisimple point}

Let $s \in \Tbf^*$ be a semisimple point. We denote by $\Zbf = \Zbf_{\Gbf^*}(s)$ the centralizer of $s$ and by $\Zbf^{\circ}$ the neutral connected component of $\Zbf$. Proceeding as in \cite[Section 2.4]{LusztigYun}, we discuss the structure of the root system of $\Zbf^{\circ}$. Let $\Phi^{\vee}_{s}$ be the set of roots $\alpha^{\vee}$ of $\Tbf^*$ such that 
$$\alpha^{\vee}(s) = 1.$$
Note that this condition is equivalent to the condition that $\alpha^{\vee,*}\Lcal_s \cong \Lambda$. 
We let $\Phi_s$ be the set of roots that are dual to the roots in $\Phi^{\vee}_s$. Lemmas \ref{lem:RootDatumCentralizer} and \ref{lem:Centralizers} are standard. We refer to \cite[Section 9.1]{LusztigYun}.

\begin{lem}\label{lem:RootDatumCentralizer}
The tuple $(X^*(\Tbf), X_*(\Tbf), \Phi_s, \Phi^{\vee}_s)$ is the root datum dual to the root datum of $(\Zbf^{\circ}, \Tbf^*)$. 
\end{lem}

The Borel $\Bbf^* \subset \Gbf^*$ determines a Borel $\Bbf_s^* = \Bbf^* \cap \Zbf^{\circ}$ of $\Zbf^{\circ}$ and therefore a notion of positive roots for the root system $(X^*(\Tbf), X_*(\Tbf), \Phi_s, \Phi^{\vee}_s)$. Note that a root of $\Phi_s$ is positive if and only if it is positive in $\Phi$. We denote by $\Delta_s \subset \Phi_s$ and $\Delta_{s}^{\vee} \subset \Phi_s^{\vee}$ the associated sets of simple roots and coroots.

We denote by 
\begin{enumerate}
\item $\Wbf_s \colondef= \Stab_{\Wbf}(s)$, 
\item $\Wbf_s^{\circ} \subset \Wbf_s$ the subgroup generated by the reflections $s_{\alpha} \in \Wbf$ such that $\alpha \in \Delta_s$. 
\end{enumerate}

\begin{lem}\label{lem:Centralizers}
\begin{enumerate}
\item The subgroup $\Wbf_s^{\circ} \subset \Wbf_s$ is normal and the quotient $\Gammabf_s \colon = \Wbf_s/\Wbf_s^{\circ}$ is abelian. 
\item The group $\Wbf_s^{\circ}$ is the Weyl group of the pair $(\Zbf^{\circ}, \Tbf^*)$. There is a canonical isomorphism $\pi_0(\Zbf) = \Gammabf_s$.
\end{enumerate}
\end{lem}

\begin{lem}[\cite{Steinberg}]\label{lem:ControlCentralizer}
Let $\Gbf^*_{\mathrm{sc}}$ be the simply connected cover of the derived subgroup $\Gbf^*_{\der}$ of $\Gbf^*$. Let $s \in \Tbf^*$. Then there is a canonical embedding $\Gammabf_s \subset \pi_1(\Gbf^*_{\der})$, where
$$\pi_1(\Gbf^*_{\der})=X_*(\Tbf^*\cap\Gbf^*_{\der})/\mathbf Z\Phi^{\vee}
\cong X^*\!\left(\ker(\Gbf^*_{\mathrm{sc}}\to\Gbf^*_{\der})\right).$$
In particular, $|\Gammabf_s|$ divides $|\pi_1(\Gbf^*_{\der})|$. 
\end{lem}

\subsection{Blocks}\label{sec:blocks}

Let $s,s' \in \Tbf^*$. Following \cite{LusztigYun}, we denote by 
$${_s}\Wbf_{s'} := \{w \in \Wbf, ws' = s\} \subset \Wbf.$$
The set ${_s}\Wbf_{s'}$ is a $(\Wbf_s, \Wbf_{s'})$-bi-torsor. We denote by 
$${_s}\underline{\Wbf}_{s'} := \Wbf_s^{\circ} \backslash {_s}\Wbf_{s'} \cong {_s}\Wbf_{s'}/\Wbf_{s'}^{\circ},$$
and call its elements blocks. This set is naturally a $(\Gammabf_s, \Gammabf_{s'})$-bi-torsor, where as before $\Gammabf_s$ (resp. $\Gammabf_{s'})$ denotes the quotient $\Wbf_s/\Wbf_s^{\circ}$ (resp. $\Wbf_{s'}/\Wbf_{s'}^{\circ}$). Given a block $\beta \in {_s}\underline{\Wbf}_{s'}$, we can consider $\beta$ as a subset of $\Wbf$. We equip $\beta$ with the induced Bruhat order. 

\begin{lem}[\protect{\cite[Lemma 4.8]{LusztigYun}}]\label{lem:ProductOfBlocks}
\begin{enumerate}
\item Let $\beta \in {_s}\underline{\Wbf}_{s'}$ be a block. Then $\beta$ has a unique minimal element $w_{\beta}$ and a unique maximal element $w^{\beta}$. 
\item Let $\beta \in {_s}\underline{\Wbf}_{s'}$ and $\gamma \in {_{s'}}\underline{\Wbf}_{s^{\prime\prime}}$ be two blocks. Then there exists a unique block denoted by $\beta \gamma$ such that for all $w \in \beta$ and $w' \in \gamma$ we have $ww' \in \beta \gamma$. 
\item Let $\beta \in {_s}\underline{\Wbf}_{s'}$ and $\gamma \in {_{s'}}\underline{\Wbf}_{s^{\prime\prime}}$ be two blocks. There are equalities
$$w_{\beta}w_{\gamma} = w_{\beta\gamma}, w_{\beta}w^{\gamma} = w^{\beta\gamma}.$$
\end{enumerate}
\end{lem}

Following \cite[Section 4.5]{LusztigYun}, we define $\Xibf$ to be the following groupoid:
\begin{enumerate}
\item its objects are the elements of $\Ch(\Tbf)$,
\item given $s,s'\in \Ch(\Tbf)$, the set of maps from $s$ to $s'$ is
$$\Hom_{\Xibf}(s,s') \colondef = {_s}\underline{\Wbf}_{s'},$$ 
\item the composition is given by the product of blocks of Lemma \ref{lem:ProductOfBlocks}. 
\end{enumerate}
It follows from \cite{LusztigYun} that the composition is associative and that the block containing $1$ acts by the identity; hence $\Xibf$ is a well-defined category. It is clear that the connected components of $\Xibf$ are in bijection with the $\Wbf$-orbits of $\Ch(\Tbf)$. Given such an orbit $\mathfrak{o}$, we denote by $\Xibf_{\mathfrak{o}}$ the corresponding connected component of $\Xibf$. We will also consider $\Xibf$ as a locally finite stack over $\Fqb$. 

We let $\widetilde{\Xibf} \to \Xibf$ be the groupoid object in schemes defined by 
\begin{enumerate}
\item its objects are the elements of $\Ch(\Tbf)$, 
\item the scheme of morphisms above $\beta \in {_s}\underline{\Wbf}_{s'}$ is the scheme of all lifts $\dot{w}_{\beta} \in \Nbf(\Tbf)$  of the minimal element $w_{\beta}$, we denote this scheme by $\widetilde{\Xibf}^{\beta}_{[s,s']}$,
\item the composition is induced by the multiplication in $\Gbf$. 
\end{enumerate}
It follows from \emph{loc. cit.} that $\widetilde{\Xibf}$ is a well-defined groupoid object in schemes and that the map $\widetilde{\Xibf} \to \Xibf$ is a $\Tbf$-gerbe. Note that on a given component $\widetilde{\Xibf}^{\beta}_{[s,s']}$ the left and right actions of $\Tbf$ differ by the action of $w_{\beta}$. 

\subsection{Endoscopic groups}\label{sec:EndoscopicGroup}

\begin{defi}[\protect{\cite[Section 10.2]{LusztigYun}, Relative Pinning}]\label{def:RelativePinning}
Let $(\Gbf,\Bbf,\Tbf)$ and $(\Hbf,\Bbf_{\Hbf}, \Tbf)$ be two reductive groups with Borel pairs and the same maximal torus such that the root datum $(X^*, \Phi_{\Hbf}, \Phi^{\vee}_{\Hbf})$ is a subroot datum of $(X^*, \Phi, \Phi^{\vee})$ satisfying $\Phi^+_{\Hbf} \subset \Phi^+$. A pinning of $\Hbf$ relative to $\Gbf$ consists, for every simple root $\alpha \in \Delta_{\Hbf}$ of $\Hbf$, of an isomorphism 
$$\iota_{\alpha} : \Hbf_{\alpha} \cong \Gbf_{\alpha}$$
of root groups. 
\end{defi}

\begin{rque}
Given $(\Gbf,\Bbf,\Tbf)$, the group $(\Hbf,\Bbf_{\Hbf}, \Tbf)$ is canonically determined by a pinning relative to $\Gbf$. 
\end{rque}

Let $s \in \Ch(\Tbf)$ and let $(\Hbf_s^{\circ}, \Bbf_s, \Tbf_s, \{\iota_{\alpha}\}_{\alpha \in \Delta_s})$ be the group dual to $\Zbf_s^{\circ}$ and relatively pinned to $\Gbf$. This determines $\Hbf_s^{\circ}$ canonically and $\Hbf_s^{\circ}$ is usually called the endoscopic group (with respect to $s$). 

Let $\mathfrak{o} \subset \Ch(\Tbf)$ be a $\Wbf$-orbit. 
Following \cite[Section 10.3]{LusztigYun}, we consider for all blocks $\beta \in {_s}\underline{\Wbf}_{s'}$ for $s,s'\in \mathfrak{o}$,
$$\Hbf_{[s,s']}^{\beta} \colondef = \widetilde{\Xibf}_{[s,s']}^{\beta} \times^{\Tbf} \Hbf_{s'}^{\circ} \cong  \Hbf_s^{\circ} \times^{\Tbf} \widetilde{\Xibf}_{[s,s']}^{\beta}.$$
The scheme $\Hbf_{[s,s']}^{\beta}$ is naturally a $(\Hbf_s^{\circ}, \Hbf_{s'}^{\circ})$-bi-torsor. Furthermore, given $\beta \in {_s}\underline{\Wbf}_{s'}$ and $\gamma \in {_{s'}}\underline{\Wbf}_{s^{\prime\prime}}$, two composable blocks, there is a well-defined multiplication
$$\Hbf_{[s,s']}^{\beta} \times \Hbf_{[s',s^{\prime\prime}]}^{\gamma} \to \Hbf_{[s,s^{\prime\prime}]}^{\beta\gamma}.$$
The collection of all $\Hbf_{[s,s']}^{\beta}$ assembles into a groupoid $\Hbffrak_{\mathfrak{o}}$ in schemes equipped with a morphism of groupoids
$$\Hbffrak_{\mathfrak{o}} \to \Xibf_{\mathfrak{o}}$$
whose fiber at a block $\beta \in {_s}\underline{\Wbf}_{s'}$ is $\Hbf_{[s,s']}^{\beta}$.
Let $s \in \mathfrak{o}$. We set
$$\Hbf_{[s,s]} = \bigsqcup_{\gamma \in \Gammabf_s} \Hbf^{\gamma}_{[s,s]}.$$
The groupoid structure of $\Hbffrak_{\mathfrak{o}}$ restricts to a group structure on $\Hbf_{[s,s]}$. The group of components of $\Hbf_{[s,s]}$ is canonically identified with $\Gammabf_s$. 

\begin{rque}\label{rque:PinningvsRelativePinning}
Recall that we have a trivialization of all root groups of $\Gbf$ (see \ref{sec:NotationsReductive}), composing the relative pinning with these trivializations yields a compatible system of pinnings of $\Hbf_s^{\circ}$ for all $s$ simultaneously. 
\end{rque} 

\subsection{Endoscopic groups and parabolic subgroups}\label{sec:endoscopic-parabolic-groups}

Let $\Lbf$ be a standard Levi subgroup of $\Gbf$, and let $\Pbf=\Lbf\Ubf_{\Pbf}$ be its associated standard parabolic subgroup. Let $\Lbf^{*} \subset \Pbf^* \subset \Gbf^*$ be the corresponding dual standard Levi and parabolic subgroup. Let $s \in \Ch(\Tbf)$ and consider the following tower of groups 
$$\Zbf_{\Lbf^*}(s) \subset \Zbf_{\Pbf^*}(s) \subset \Zbf_{\Gbf^*}(s).$$

\begin{lem}\label{lem:ParabEndoscopicGroups}
\begin{enumerate}
\item The group $\Zbf_{\Pbf^*}^{\circ}(s)$ is a standard parabolic subgroup of $\Zbf^{\circ}_{\Gbf^*}(s)$ with standard Levi $\Zbf^{\circ}_{\Lbf^*}(s)$ and unipotent radical $\Zbf_{\Ubf_{\Pbf^*}}(s)$. 
\item $\pi_0(\Zbf_{\Pbf^*}(s)) = \pi_0(\Zbf_{\Lbf^*}(s))$, and $\Zbf_{\Pbf^*}(s) = \Zbf_{\Lbf^*}(s)\Zbf_{\Ubf_{\Pbf^*}}(s)$, with the second factor being connected.
\item Let $\Wbf_s^{\Lbf}$ be the Weyl group of $\Zbf_{\Lbf^*}(s)$ and $\Gammabf_s^{\Lbf} = \pi_0(\Zbf_{\Pbf^*}(s))$, under the canonical splittings 
$$\Wbf_s^{\Lbf} = \Wbf_{s}^{\Lbf, \circ} \rtimes \Gammabf_s^{\Lbf} \subset \Wbf_{s} = \Wbf_s^{\circ}\rtimes \Gammabf_s$$
the subgroup $\Gammabf_s^{\Lbf}$ maps injectively to $\Gammabf_s$; we identify it with its image. 
\end{enumerate}
\end{lem}

\begin{proof}
The first two points are well known; let us prove the third. Let $\beta \in \Gammabf_s^{\Lbf}$ and let $w_{\beta}^{\Lbf} \in \Wbf_{s}^{\Lbf}$ be the corresponding minimal element. Let $\beta' \in \Gammabf_s$ be the image of the block $\beta$, and let $w_{\beta'}$ be the corresponding minimal element of $\Wbf_s$. Since $w_{\beta}^{\Lbf}$ is minimal, it sends all positive roots of $\Zbf_{\Lbf^*}(s)$ to positive roots. Moreover, since $\Zbf_{\Lbf^*}(s)$ normalizes the unipotent radical of $\Zbf_{\Pbf^*}(s)$, $w_{\beta}^{\Lbf}$ sends all positive roots in the unipotent radical of $\Zbf_{\Pbf^*}(s)$ to positive roots. Thus $w_{\beta}^{\Lbf}$ sends all positive roots of $\Zbf^{\circ}_{\Gbf^*}(s)$ to positive roots, so it has length $0$ in $\Wbf_s$ and we have $w_{\beta}^{\Lbf} = w_{\beta'}$. 
\end{proof}

Let $s \in \Ch(\Tbf)$, dualizing the above sequence we have a canonical standard Levi and parabolic subgroup  
$$\Hbf^{\Lbf, \circ}_s \subset \Hbf_s^{\Pbf,\circ} \subset \Hbf_s^{\Gbf,\circ}$$
and extensions 
$$\Hbf^{\Lbf}_{[s,s]} \subset \Hbf^{\Pbf}_{[s,s]} \subset \Hbf^{\Gbf}_{[s,s]}.$$ 
We write
$$\Ubf_s^{\Hbf,\Pbf}:=\mathrm{R}_u(\Hbf_s^{\Pbf,\circ})$$
for the unipotent radical of the endoscopic parabolic. It is normalized by $\Hbf^{\Pbf}_{[s,s]}$, and Frobenius carries it to $\Ubf_{\Frob(s)}^{\Hbf,\Frob(\Pbf)}$. This notation distinguishes the endoscopic unipotent radical from the subgroup $\Ubf_{\Pbf}\subset\Gbf$.

Let $s,s'$ be two elements that are $\Wbf^{\Lbf}$-conjugate. If $\beta\in {_s}\underline{\Wbf}^{\Lbf}_{s'}$ is a Levi block, we denote its image in ${_s}\underline{\Wbf}_{s'}$ by the same letter. Then there are naturally defined bi-torsors 
$$\Hbf^{\Lbf}_{[s,s']} \subset \Hbf^{\Pbf}_{[s,s']} \subset \Hbf^{\Gbf}_{[s,s']}.$$

Let $\ofrak \subset \Ch(\Tbf)$ be the $\Wbf$-orbit containing $s$ and let $\ofrak^{\Lbf} \subset \ofrak$ be the corresponding $\Wbf^{\Lbf}$-orbit. Set $\ofrak^{\Gbf}=\ofrak$. Putting together all the bitorsors $\Hbf^{?}_{[s,s']}$, we get a commutative diagram of groupoids 

\begin{equation}\label{diag:functoriality-for-endoscopic-groups}
\begin{tikzcd}
	{\Hbffrak^{\Lbf}_{\ofrak^{\Lbf}}} & {\Hbffrak^{\Pbf}_{\ofrak^{\Lbf}}} & {\Hbffrak^{\Gbf}_{\ofrak^{\Gbf}}} \\
	{\Xibf_{\ofrak^{\Lbf}}^{\Lbf}} & {\Xibf_{\ofrak^{\Lbf}}^{\Lbf}} & {\Xibf_{\ofrak^{\Gbf}}^{\Gbf}}
	\arrow[from=1-1, to=2-1]
	\arrow["{\alpha_{\Pbf}^{\Hbf}}"', from=1-2, to=1-1]
	\arrow["{\beta_{\Pbf}^{\Hbf}}", from=1-2, to=1-3]
	\arrow[from=1-2, to=2-2]
	\arrow[from=1-3, to=2-3]
	\arrow[equals, from=2-2, to=2-1]
	\arrow[from=2-2, to=2-3]
\end{tikzcd}
\end{equation}

\subsection{Endoscopic groups and isogenies}\label{sec:PrelimEndoscopicGroupsAndIsogeny}
Let $f : \Gbf \to \Gbf'$ be a central isogeny and denote by $f^* : \Gbf'^* \to \Gbf^*$ the dual isogeny. Let $\Zbf_f = \ker(f^*)$ be the kernel of $f^*$. In this subsection, we shall denote objects defined with respect to $\Gbf'$ with a $'$ (for instance, we write $\Bbf', \Tbf', \Wbf', \dots$). 

\subsubsection{Orbits} 
Let $\ofrak \subset \Ch(\Tbf)$ be a $\Wbf$-orbit. Let $\ofrak'\subset \Ch(\Tbf')$ be a $\Wbf$-orbit such that $f^*(\ofrak') = \ofrak$. For $z \in \Zbf_f$, set 
$$z\ofrak' \colondef = \{zs' \in \Ch(\Tbf'), s'\in \ofrak'\}$$ 
is a $\Wbf$-orbit in $\Ch(\Tbf')$. This operation defines a transitive action of $\Zbf_f$ on $\{\ofrak' \subset \Ch(\Tbf'), f^*(\ofrak') = \ofrak\}$.

\subsubsection{Stabilizer and transporters}
Let $s'_1,s'_2 \in \Ch(\Tbf')$ with image $s_1,s_2 \in \Ch(\Tbf)$. The natural isomorphism $\Wbf' \to \Wbf$ yields an injective map 
\begin{equation}\label{eq:IsogenyTransporter}
{_{s'_1}}\Wbf_{s'_2} \to {_{s_1}}\Wbf_{s_2}.
\end{equation}

\begin{lem}\label{lem:iso-connected-stabilizer}
The map (\ref{eq:IsogenyTransporter}), for $s_1' = s'_2$, induces an isomorphism 
$\Wbf_{s'_1}^{'\circ} = \Wbf_{s_1}^{\circ}.$
\end{lem}

\begin{proof}
The map $\Gbf'^* \to \Gbf^*$ induces an isogeny 
$\Zbf_{\Gbf'^*}^{\circ}(s_1') \to \Zbf_{\Gbf^*}^{\circ}(s_1)$
and therefore
$\Wbf_{s'_1}^{'\circ} = \Wbf_{s_1}^{\circ}.$
\end{proof}

It follows from \Cref{lem:iso-connected-stabilizer} that we have a well-defined injective map 
$${_{s'_1}}\underline{\Wbf}'_{s'_2} \to {_{s_1}}\underline{\Wbf}_{s_2}.$$
Recall that we denote by $\Gammabf_s = \pi_0( \Zbf_{\Gbf^*}(s)) = {_{s}}\underline{\Wbf}_{s}$. We therefore have an injective map $\Gammabf'_{s'_1} \to \Gammabf_{s_1}.$
We denote by 
$$\Gammabf_{f, s'_1} =  \Gammabf_{s_1}/\Gammabf'_{s'_1}.$$
The choice of $s'_1$ yields an embedding 
$$i_{s'_1} : \Gammabf_{f,s'_1} \to \Zbf_f$$
given by sending $\gamma \mapsto [\dot{\gamma}, s'_1]$ where $\dot{\gamma}$ is any lift to $\Gbf'^*$ of $\gamma$.
More generally, given $s_1,s_2,s'_1,s'_2$ as above, there is a well-defined map 
$$i_{s'_1,s'_2} : {_{s_1}}\underline{\Wbf}_{s_2} \to \Zbf_f$$
given by $i_{s'_1,s'_2}(\beta) = \dot{\beta}s'_2\dot{\beta}^{-1}s_1^{'-1}$, where $\dot{\beta}$ is any lift of $\beta$ to $\Gbf^{'*}$. We note that this map sends ${_{s'_1}}\underline{\Wbf}'_{s'_2}$ to $1 \in \Zbf_f$. 

\begin{lem}\label{lem:SectionIsogenyTransporter}
With $s_1,s_2$ and $s'_1,s'_2$ chosen as above, the canonical map
$$\bigsqcup_{z \in i_{s'_1, s'_2}({_{s_1}}\underline{\Wbf}_{s_2})} {_{s'_1z}}\underline{\Wbf}'_{s'_2} \to {_{s_1}}\underline{\Wbf}_{s_2}$$
is an isomorphism.
\end{lem}

\begin{proof}
This follows by unwinding the definition. We show that 
$$\bigsqcup_{z \in i_{s'_1, s'_2}({_{s_1}}{\Wbf}_{s_2})} {_{s'_1z}}{\Wbf}'_{s'_2} = {_{s_1}}{\Wbf}_{s_2}.$$
Let $w \in {_{s_1}}{\Wbf}_{s_2}$ and let $\dot{w}$ be any lift of $w$ to $\Gbf^{'*}$. Setting $z = \dot{w}s_2'\dot{w}^{-1}s_1^{'-1}$, we have $w \in {_{s'_1z}}{\Wbf}'_{s'_2}$. This establishes surjectivity. Conversely, if $w_1' \in {_{s'_1z_1}}{\Wbf}'_{s'_2}$ and $w_2' \in {_{s'_1z_2}}{\Wbf}'_{s'_2}$ are equal in ${_{s_1}}{\Wbf}_{s_2}$, then $z_1 = z_2$ and thus $w'_1 = w'_2$, as ${_{s'_1z}}{\Wbf}'_{s'_2} \to {_{s_1}}{\Wbf}_{s_2}$ is injective. 
The equality is stable under the left and right connected stabilizers. Quotienting by $\Wbf_{s_1}^{\circ}=\Wbf_{s'_1}^{\prime\circ}$ and $\Wbf_{s_2}^{\circ}=\Wbf_{s'_2}^{\prime\circ}$ therefore gives the asserted equality of underlined block sets.
\end{proof}

\begin{lem}\label{lem:IsogenyPreservesBlock}
The map ${_{s'_1}}{\Wbf}'_{s'_2} \subset {_{s_1}}{\Wbf}_{s_2}$ sends the minimal length elements onto minimal length elements. 
\end{lem}

\begin{proof}
This is clear since both sets are $\Wbf_{s'_1}^{\circ} = \Wbf_{s_1}^{\circ}$-stable subsets of $\Wbf$ and their inclusion in $\Wbf$ preserves the length.
\end{proof}

\subsubsection{Endoscopic groups}

Let $s' \in \Ch(\Tbf^{'})$ have image $s \in \Ch(\Tbf)$. Then the isogeny
$$\Zbf_{\Gbf^{'*}}^{\circ}(s') \to \Zbf_{\Gbf^*}^{\circ}(s)$$
yields a dual isogeny of (relatively) pinned reductive groups
$$\Hbf_{s}^{\circ} \to \Hbf^{'\circ}_{s'}.$$

More generally, if $\beta' \in \Gammabf'_{s'} \subset \Gammabf_s$, by Lemma \ref{lem:IsogenyPreservesBlock}, we have a well-defined map 
\begin{equation}\label{eq:maps:endoscpic-groups-isogeny}
\Hbf^{\beta}_{[s,s]} \to \Hbf^{' \beta'}_{[s',s']}.
\end{equation}
Note that the maps (\ref{eq:maps:endoscpic-groups-isogeny}) \emph{do not} assemble to a map of groups
$$\Hbf_{[s,s]} \to \Hbf'_{[s',s']}.$$

For blocks, we have the following maps. Let $\beta' \in {_{s'_1}}\underline{\Wbf}'_{s'_2}$ have image $\beta \in {_{s_1}}\underline{\Wbf}_{s_2}$. Then we have a well-defined map 
\begin{equation}\label{eq:ProjectionIsogenyEndoscopicSide}
\Hbf^{\beta}_{[s_1,s_2]} \to \Hbf^{'\beta'}_{[s'_1, s'_2]}.
\end{equation}
This map is $(\Hbf^{\circ}_{s_1}, \Hbf^{\circ}_{s_2}) \to (\Hbf^{'\circ}_{s'_1}, \Hbf^{'\circ}_{s'_2})$-linear.

\subsubsection{A bimodule on the endoscopic side}\label{sec:BimoduleEndoscopicSide}

Let $s' \in \Ch(\Tbf^{'})$ with image $s \in \Ch(\Tbf)$ be as above. We define
$$\Hbf'_{[s,s']} = \bigsqcup_{z \in \Zbf_f} \Hbf'_{[zs', s']}$$
and we now equip this scheme with a structure of a $(\Hbf_{[s,s]}, \Hbf'_{[s',s']})$-bimodule. The right action of $\Hbf'_{[s',s']}$ is given by right translation on each of the summands. We note that 
\begin{align*}
 \Hbf'_{[s,s']} &= \bigsqcup_{z} \Hbf'_{[zs', s']} \\
 &= \bigsqcup_z \bigsqcup_{\beta' \in {_{zs'}}\underline{\Wbf'}_{s'}} \Hbf^{',\beta'}_{[zs', s']}, 
\end{align*}
and under the bijection of Lemma \ref{lem:SectionIsogenyTransporter}, we will simply write this as 
$$\Hbf'_{[s,s']} = \bigsqcup_{\beta \in {_{s}}\underline{\Wbf}_{s}} \Hbf^{',\beta}_{[z(\beta)s', s']},$$
where $z(\beta)\in\Zbf_f$ is the unique element determined by that bijection.
Let $\beta_1, \beta_2 \in {_{s}}\underline{\Wbf}_{s}$. We define a map
$$a_{\beta_1,\beta_2} : \Hbf_{[s,s]}^{\beta_1} \times \Hbf_{[s,s']}^{', \beta_2} \to \Hbf_{[s,s']}^{',\beta_1\beta_2}.$$
Let $z_1 \in \Zbf_f$ and $\beta_2'$ be such that $\Hbf_{[s,s']}^{',\beta_2} = \Hbf_{[z_1s',s']}^{',\beta_2'}$. By Lemma \ref{lem:SectionIsogenyTransporter}, there exists a unique $z_2 \in \Zbf_f$ and a unique $\beta_1'$ such that $\beta_1'$ maps to $\beta_1$ under the inclusion 
$${_{z_2s'}}\underline{\Wbf'}_{z_1s'} \hookrightarrow {_s}\underline{\Wbf}_s.$$
We have a well-defined map 
\begin{equation}\label{eq:DefineActionPart1}
\Hbf^{\beta_1}_{[s,s]} \to \Hbf^{',\beta_1'}_{[z_2s',z_1s']}.
\end{equation}
We set $a_{\beta_1,\beta_2}$ to be the following composition 
$$\Hbf_{[s,s]}^{\beta_1} \times \Hbf_{[s,s']}^{', \beta_2} \to \Hbf_{[z_2s',z_1s']}^{',\beta'_1} \times \Hbf_{[z_1s',s']}^{',\beta_2'} \to \Hbf_{[z_2s',s']}^{',\beta_1'\beta_2'},$$
where the first map is (\ref{eq:DefineActionPart1}) and the second map is the multiplication in the groupoid $\Hbf'$. 

\begin{lem}\label{lem:ConstructionBimodEndoscopicSide}
The collection of maps $(a_{\beta_1,\beta_2})_{\beta_1,\beta_2}$ defines a left action of $\Hbf_{[s,s]}$ on $\Hbf'_{[s,s']}$. 
\end{lem}

\begin{proof}
The fact that $a_{1, \beta} = \id_{\Hbf'_{[z_1s',s']}}$ follows from the fact that the projection (\ref{eq:ProjectionIsogenyEndoscopicSide}) is $\Hbf_{[s,s]}^{\circ} \to \Hbf^{',\circ}_{[s',s']}$ linear. It remains to prove the commutativity of the following diagram for all $\beta_1,\beta_2,\beta_3 \in {_s}\underline{\Wbf}_{s}$ 
\[\begin{tikzcd}
	{\Hbf_{[s,s]}^{\beta_1} \times \Hbf_{[s,s]}^{\beta_2} \times\Hbf_{[s,s']}^{',\beta_3}} & {\Hbf_{[s,s]}^{\beta_1}\times \Hbf_{[s,s']}^{',\beta_2\beta_3}} \\
	{\Hbf_{[s,s]}^{\beta_1\beta_2} \times \Hbf_{[s,s']}^{',\beta_3}} & {\Hbf_{[s,s']}^{',\beta_1\beta_2\beta_3}.}
	\arrow["{\id \times a_{\beta_2, \beta_3}}", from=1-1, to=1-2]
	\arrow["{m \times \id_{\Hbf_{[s,s']}^{',\beta_3}}}"', from=1-1, to=2-1]
	\arrow["{a_{\beta_1, \beta_2\beta_3}}", from=1-2, to=2-2]
	\arrow["{a_{\beta_1\beta_2, \beta_3}}"', from=2-1, to=2-2]
\end{tikzcd}\]
Let $(z_3, \beta_3')$ be as before, chosen such that $\Hbf_{[s,s']}^{',\beta_3} = \Hbf_{[z_3s',s']}^{',\beta_3'}$. Let $z_2$ be such that $\beta_2$ lies in ${_{z_2s'}}\underline{\Wbf'}_{z_3s'} \subset {_s}\underline{\Wbf}_{s}$ (where the embedding is given by Lemma \ref{lem:SectionIsogenyTransporter} for the pair $(s_1',s_2') = (s', z_3s')$). Similarly let $z_1$ be such that $\beta_1 \in {_{z_1s'}}\underline{\Wbf'}_{z_2s'} \subset {_s}\underline{\Wbf}_{s}$, where the embedding is again given by Lemma \ref{lem:SectionIsogenyTransporter} for the pair $(s_1',s_2') = (s', z_2s')$, and let $z_{12}$ be such that $\beta_1\beta_2 \in {_{z_{12}s'}}\underline{\Wbf'}_{z_3s'}$ (for the pair $(s_1',s_2') = (s', z_3s')$). By construction, we have $z_{12} = z_1$. This implies that we have a commutative diagram 
\[\begin{tikzcd}
	{\Hbf_{[s,s]}^{\beta_1} \times \Hbf_{[s,s]}^{\beta_2}} & {\Hbf_{[z_1s', z_2s']}^{',\beta_1'} \times \Hbf_{[z_2s',z_3s']}^{',\beta_2'}} \\
	{\Hbf_{[s,s]}^{\beta_1\beta_2}} & {\Hbf_{[z_{12}s',z_3s']}^{',\beta_1'\beta_2'}}
	\arrow[from=1-1, to=1-2]
	\arrow[from=1-1, to=2-1]
	\arrow[from=1-2, to=2-2]
	\arrow[from=2-1, to=2-2]
\end{tikzcd}\]
where the vertical maps are given by multiplication and the horizontal maps are of the form (\ref{eq:ProjectionIsogenyEndoscopicSide}). The lemma follows. 
\end{proof}

\subsection{Rationality}

Recall that we have equipped $\Gbf$ with a purely inseparable isogeny $\Frob : \Gbf \to \Gbf$. We denote by $\Frob^* : \Gbf^* \to \Gbf^*$ the dual isogeny which is also purely inseparable. 

Let $s,t \in \Gbf^*$. We denote by
$${_s}\Gbf^*_{t} \colondef = \{g \in \Gbf^*, \Ad(g)(t) = s\}.$$
It is clear that this is a $(\Zbf(s), \Zbf(t))$-bi-torsor, where, as before, $\Zbf(s)$ is the centralizer of $s$ in $\Gbf^*$. We denote 
$${_s} \underline{\Gbf}^*_t \colondef= \Zbf^{\circ}(s) \backslash {_s}\Gbf^*_{t} \cong {_s}\Gbf^*_{t} / \Zbf^{\circ}(t).$$
Note that if both $s,t$ lie in $\Tbf^*$, then there is a canonical isomorphism of $(\Gammabf_s, \Gammabf_t)$-bi-torsors
$${_s} \underline{\Gbf}^*_t \cong {_s}\underline{\Wbf}_t.$$

\begin{defi}\label{def:geometric-and-rational-conj}
We say that a semisimple point $s \in \Gbf^*$ is rational if $s \in (\Gbf^*)^{\Frob^*}$. We say that the $\Gbf^*$-conjugacy class of $s$ is rational if $\Frob^*(\Ad(\Gbf^*)(s)) = \Ad(\Gbf^*)(s)$.
\end{defi}

Let $\Lambda \in \{\Flb, \Zlb, \Qlb\}$ be a coefficient ring. We define two equivalence relations on $(\Gbf^*_{\mathrm{ss}})^{\Frob^*}$, the set of rational semisimple points of $(\Gbf^*)^{\Frob^*}$.
\begin{enumerate}
\item if $\Lambda = \Qlb$, 
\begin{enumerate}
\item $s \sim_{\Lambda, \geo} t$ if and only if $s$ and $t$ are conjugate in $\Gbf^*$,
\item $s \sim_{\Lambda, \rat} t$ if and only if $s$ and $t$ are conjugate in $(\Gbf^*)^{\Frob^*}$,
\end{enumerate}
\item if $\Lambda = \Flb, \Zlb$, 
\begin{enumerate}
\item $s \sim_{\Lambda, \geo} t$ if and only if the prime-to-$\ell$ parts of $s$ and $t$ are conjugate in $\Gbf^*$,
\item $s \sim_{\Lambda, \rat} t$ if and only if the prime-to-$\ell$ parts of $s$ and $t$ are conjugate in $(\Gbf^*)^{\Frob^*}$.
\end{enumerate}
\end{enumerate}
Clearly if $s \sim_{\Lambda,\rat} t$ then $s \sim_{\Lambda, \geo} t$ and if $s \sim_{\Zlb, ?} t$ then $s \sim_{\Qlb, ?} t$ for $? \in \{\geo, \rat\}$. 

\begin{lem}[\protect{\cite[Proposition 3.21]{DigneMichelBook}}]\label{lem:RationalConjClassFaciles}
Let $s \in (\Gbf^*_{\mathrm{ss}})^{\Frob^*}$, then there is a canonical bijection 
$$\{(\Gbf^*)^{\Frob^*}\text{-orbits in} \ \Ad(\Gbf^*)(s)\} \cong H^1(\Frob^*, \Gammabf_s).$$
\end{lem}

\begin{lem}\label{lem:RationalConjClasses}
Let $s \in \Tbf^*$ be such that its $\Gbf^*$-conjugacy class is rational. Then there is a canonical bijection
$$\{(\Gbf^*)^{\Frob^*}\text{-orbits in} \ \Ad(\Gbf^*)(s)\} \cong |\frac{{_s}\underline{\Wbf}_{\Frob(s)}}{\Ad_{\Frob}\Gammabf_s}|$$
where the right-hand side is the total space of the groupoid quotient of ${_s}\underline{\Wbf}_{\Frob(s)}$ by the action of $\Gammabf_s$ given by $\Ad_{\Frob}(\gamma)(\beta) = \gamma \beta \Frob(\gamma^{-1})$. 
\end{lem}

\begin{proof}
The proof is essentially the same as the proof of Lemma \ref{lem:RationalConjClassFaciles}. Let us indicate the main steps. We first note that there is a canonical bijection 
$$|\frac{{_s}\underline{\Wbf}_{\Frob(s)}}{\Ad_{\Frob}\Gammabf_s}| = |\frac{{_s}\Gbf^*_{\Frob^*(s)}}{\Ad_{\Frob^*}(\Zbf(s))}|.$$
Let $t \in (\Gbf^*)^{\Frob^*}$ be $\Gbf^*$-conjugate to $s$, then there exists $g \in \Gbf^*$ such that $gsg^{-1} = t$. It follows that $\Frob^*(gsg^{-1}) = gsg^{-1}$ and that the element $g^{-1}\Frob^*(g)$ belongs to ${_s}\Gbf^*_{\Frob^*(s)}$. Arguing as in \cite[Proposition 3.21]{DigneMichelBook}, it is clear that the corresponding point of $|\frac{{_s}\Gbf^*_{\Frob^*(s)}}{\Ad_{\Frob^*}(\Zbf(s))}|$ is independent of the choice of $g$. 

Conversely, let $g \in {_s}\Gbf^*_{\Frob^*(s)}$ represent an element in $|\frac{{_s}\Gbf^*_{\Frob^*(s)}}{\Ad_{\Frob^*}(\Zbf(s))}|$. There exists $h \in \Gbf^*$ such that $h^{-1}\Frob^*(h) = g$. Let $t = hsh^{-1}$. We have 
$$t^{-1}\Frob^*(t) =  hs^{-1}h^{-1}\Frob^*(hsh^{-1}) = hs^{-1}g\Frob^*(s)\Frob^*(h^{-1}) = hg\Frob^*(h^{-1}) = 1$$
where the third equality comes from the fact that $g\Frob^*(s) = sg$. It follows that $t \in (\Gbf^*)^{\Frob^*}$. The same argument as in \emph{loc. cit.} shows that the $(\Gbf^*)^{\Frob^*}$-conjugacy class of $t$ depends only on the point in $|\frac{{_s}\Gbf^*_{\Frob^*(s)}}{\Ad_{\Frob^*}(\Zbf(s))}|$. This gives an inverse to the previous map. 
\end{proof}

\subsection{The Morita Category}\label{sec:MoritaCategory}

Recall that we denote by $\Pr_{\Lambda}$ the $\infty$-category of presentable $(\infty,1)$-categories that are $\DD(\Lambda)$-linear. We denote by 
$$\Morita$$ 
the category whose objects are monoidal categories in $\Pr_{\Lambda}$, i.e., objects in $\Alg(\Pr_{\Lambda})$, and whose morphism space $\Hom_{\Morita}(A,B)$ is the groupoid $(A,B)-\Bim^{\cong}$ of $(A,B)$-bimodules, with only isomorphisms retained (so that $\Morita$ is an $(\infty,1)$-category and not an $(\infty,2)$-category). The composition in $\Morita$ is given by relative tensor product, namely 
\begin{align*}
\Hom_{\Morita}(A,B) \times \Hom_{\Morita}(B,C) &\to \Hom_{\Morita}(A,C) \\
(M,N) &\mapsto M \otimes_B N.
\end{align*}
This category is constructed in \cite[Section 4.4.3]{HigherAlgebra}; see in particular Remark 4.4.3.11 of \emph{loc.\ cit.}

\section{Free monodromic Hecke categories}\label{sec:HeckeCat}

\subsection{The Hecke category $\HH^{\Gbf}$}\label{sec:Hecke-cat-generalities}

We start by recalling the definition and construction of the free monodromic Hecke categories. 
\begin{defi}
We denote by $\HH_{\tot}$ (or $\HH^{\Gbf}_{\tot,\Lambda}$ if we want to put some emphasis on $\Gbf$ and/or $\Lambda$) the category
$$\HH_{\tot} \colondef = \bigoplus_{s,s' \in \Ch(\Tbf)} \DD(\Ubf \backslash \Gbf/\Ubf, \RR_{\Tbf \times \Tbf})_{[s,s']},$$
where $\HH_{[s,s']} = \DD(\Ubf \backslash \Gbf/\Ubf, \RR_{\Tbf \times \Tbf})_{[s,s']}$ is the category of free monodromic $(\Lcal_{s}, \Lcal_{s'})$ sheaves on $\Ubf \backslash \Gbf/\Ubf$. We refer to \cite{Gouttard} and \cite[Section 4.1]{EteveFreeMonodromic} for a construction. Let $\ofrak \subset \Ch(\Tbf)$ be a $\Wbf$-orbit. We set
$$\HH_{\ofrak, \tot} \colondef = \oplus_{s,s' \in \ofrak} \HH_{[s,s']},$$
This is a direct summand of $\HH_{\tot}$.
\end{defi}

\begin{lem}[\protect{\cite[Theorem 1.2.1]{EteveFreeMonodromic}}]
\begin{enumerate}
\item The category $\HH_{\tot}$ is compactly generated and an object is compact if and only if it is constructible for the Bruhat stratification. 
\item The category $\HH_{[s,s']}$ is nonzero only if $s$ and $s'$ are $\Wbf$-conjugate.
\item The category $\HH_{\tot}$ is equipped with a well-defined perverse $t$-structure. 
\end{enumerate}
\end{lem}

Consider the following convolution diagram 
\[\begin{tikzcd}
	& {\Ubf \backslash \Gbf \times^{\Ubf} \Gbf/\Ubf} & {\Ubf \backslash \Gbf/\Ubf} \\
	{\Ubf \backslash \Gbf/\Ubf} && {\Ubf \backslash \Gbf/\Ubf}
	\arrow["m", from=1-2, to=1-3]
	\arrow["{\pr_1}"{description}, from=1-2, to=2-1]
	\arrow["{\pr_2}"{description}, from=1-2, to=2-3]
\end{tikzcd}\]
where $\pr_i$ is induced by the $i$th projection for $i = 1,2$, and $m$ is induced by multiplication. 
The category $\HH$ is then equipped with the monoidal structure 
$$A * B = m_!(A \widehat{\boxtimes}_{\Lambda} B)[\dim \Tbf].$$
We refer to \cite[Section 4.2]{EteveFreeMonodromic} for the notations. We also recall that for $A \in \HH_{[s_1,s_2]}$ and $B \in \HH_{[s_3,s_4]}$ we have
\begin{enumerate}
\item $A * B = 0$ if $s_2 \neq s_3$, 
\item $A * B \in \HH_{[s_1, s_4]}$. 
\end{enumerate}

Let $n \in \Nbf(\Tbf)$ be a lift of some element $w \in \Wbf$. This element yields a splitting 
$$\Ubf \backslash \Bbf n\Bbf/\Ubf = \Tbf.n/(\Ubf \cap {^n}\Ubf),$$
and therefore a map $\nu_n : \Ubf \backslash \Bbf n\Bbf/\Ubf \to \Tbf$.  Let us also denote by $j_w : \Ubf \backslash \Bbf n\Bbf/\Ubf \to \Ubf \backslash \Gbf/\Ubf$ the inclusion of the stratum. 
We denote by 
\begin{enumerate}
\item $\Delta_{n,s} = j_{w,!}\nu_n^*(L_{\Tbf} \otimes_{\Lambda} \Lcal_s)[\ell(w) + \dim \Tbf]$, the standard sheaves,
\item $\nabla_{n,s} = j_{w,*}\nu_n^*(L_{\Tbf} \otimes_{\Lambda} \Lcal_s)[\ell(w) + \dim \Tbf]$, the costandard sheaves,
\end{enumerate}
where $L_{\Tbf}$ denotes the free monodromic unipotent sheaf, as recalled in Section \ref{sec:SemisimplePoints}. It is clear from the definition that if $n,n'$ are two lifts of some element $w \in \Wbf$, the sheaves $\Delta_{n,s}$ and $\Delta_{n',s}$ are (noncanonically) isomorphic.

\begin{nota}
Recall that we have fixed a set of lifts $\dot{w} \in \Nbf(\Tbf)$ of the elements of $\Wbf$. We will denote by 
$$\Delta_{w,s} \colondef = \Delta_{\dot{w},s}, \nabla_{w,s} \colondef = \nabla_{\dot{w},s}.$$
We also denote by 
$$\Delta_{w} \colondef = \oplus_{s \in \Ch(\Tbf)} \Delta_{w,s}, \nabla_w \colondef = \oplus_{s \in \Ch(\Tbf)} \nabla_{w,s}.$$ 
\end{nota}

\begin{defi}\label{def:TiltingSheaves}
Let $A \in \HH^{\heartsuit}_{\tot}$ be a perverse sheaf.
\begin{enumerate}
\item A $\Delta$-filtration on $A$ is a filtration whose graded pieces are isomorphic to sheaves of the form $\Delta_{n,s}$.
\item A $\nabla$-filtration on $A$ is a filtration whose graded pieces are isomorphic to sheaves of the form $\nabla_{n,s}$.
\item The sheaf $A$ is tilting if it has both a $\Delta$ and a $\nabla$-filtration.  
\end{enumerate}
We denote by $\HH_{\tilt}$ the full subcategory of $\HH^{\heartsuit}_{\tot}$ of tilting objects.
\end{defi}

\begin{lem}[\protect{\cite[7.7]{BezrukavnikovRicheSoergelTheory}, \cite[8.4.1]{Gouttard}}]
Let $w,w' \in \Wbf$ and $n,n' \in \Nbf(\Tbf)$ be lifts of $w$ and $w'$, respectively. Let $s \in \Ch(\Tbf)$. 
\begin{enumerate}
\item if $\ell(w) + \ell(w') = \ell(ww')$ then there are canonical isomorphisms
$$\Delta_{n,n's} * \Delta_{n',s} \cong \Delta_{nn',s}, \nabla_{n,n's} * \nabla_{n',s} \cong \nabla_{nn',s}.$$
\item There is an isomorphism 
$$\Delta_{n^{-1}, ns} * \nabla_{n,s} \cong \Delta_{1,s} \cong \nabla_{n^{-1},ns} * \Delta_{n,s}.$$
\end{enumerate}
\end{lem}

\begin{lem}\label{lem:VanishingBetweenStandardAndCostandard}
For all $n,n' \in \Nbf(\Tbf)$, lifts of $w$ and $w'$, respectively, and $s,s' \in \Ch(\Tbf)$, we have 
\begin{enumerate}
\item $\Hom(\Delta_{n,s}, \nabla_{n',s'}) = 0$ if $w \neq w'$ or $s \neq s'$,
\item $\Hom(\Delta_{n,s}, \nabla_{n',s'}) \cong \RR_{\Tbf}$ if $w = w'$ and $s = s'$.
\end{enumerate}
\end{lem}

\begin{proof}
Clear from the definitions. 
\end{proof}

\begin{thm}[Structure of tilting sheaves]\label{thm:StructureCategoryviaTilting}
\begin{enumerate}
\item Let $T,T' \in \HH_{\tilt}$. Then $\Hom(T,T')$ is concentrated in degree $0$. 
\item The convolution of two tilting sheaves is tilting.
\item A sheaf is tilting if and only if its $!$ and $*$-restrictions to all strata are free monodromic perverse sheaves. 
\item For all $(w,s)$ there exists an indecomposable tilting sheaf $T_{w,s}$, unique up to isomorphism, such that $T_{w,s}$ contains $\Delta_{w,s}$ in a $\Delta$-filtration with multiplicity one and $T_{w,s}$ is supported on the closure of $\Ubf \backslash \BwB/\Ubf$. 
\item The category $\HH_{\tilt}$ is Krull-Schmidt, Karoubian and generates $\HH_{\tot}$ under colimits. 
\item There exists a unique weight structure on $\HH_{\tot}$ whose heart is $\HH_{\tilt}$ and the weight complex functor $\HH_{\tot} \to \KK(\HH_{\tilt})$ is an equivalence. 
\end{enumerate}
\end{thm}

\begin{proof}
\begin{enumerate}
\item Immediate from Lemma \ref{lem:VanishingBetweenStandardAndCostandard}.
\item This is \cite[Proposition 4.3.4]{BezrukavnikovYun} and \cite[Lemma 9.7.5]{Gouttard}.
\item This is a classical statement about tilting perverse sheaves, see \cite{TiltingExercises} and \cite[Section 5.4]{BezrukavnikovRicheSoergelTheory} in the free monodromic context. 
\item This is \cite[Lemma A.7.3]{BezrukavnikovYun}, \cite[Proposition 5.11]{BezrukavnikovRicheSoergelTheory} and \cite[Proposition 9.5.2]{Gouttard}.
\item The Krull-Schmidt property is \cite[Corollary 5.4]{BezrukavnikovRicheSoergelTheory}. The generation follows from point $(iv)$ and the Karoubian property follows from $(iii)$. 
\item Since the tilting objects generate the category $\HH_{\tot}$ under colimits and the mapping space between two tilting sheaves is concentrated in degree $0$, we can argue as in \cite[Corollary 6.22]{EberhardtEteve}.
\end{enumerate}
The reader should note that the references \cite{BezrukavnikovYun}, \cite{BezrukavnikovRicheSoergelTheory}, and \cite{Gouttard} consider the case where $\Lambda$ is a field ($\Qlb$ for the first reference and $\Flb$ for the last two). The proofs of the $\Zlb$ case are essentially the same as in the $\Flb$ case; details can be found in the thesis of the author \cite[Section 3.B]{EteveThesis}.
\end{proof}

\begin{rque}
Since the category $\HH_{\tilt}$ is contained in the perverse heart of $\HH_{\tot}$, it follows that the realization functor $\Ind\mathrm{K}^b(\HH_{\tilt}) \to \DD(\HH^{\heartsuit}_{\tot}) \to \HH_{\tot}$ is an equivalence whose inverse is the weight complex functor.
\end{rque}

\subsection{Blocks and Morita equivalences}\label{sec:block-and-morita}

\begin{lem}\label{lem:ForgetSupport}
Let $s,s' \in \Ch(\Tbf)$ be $\Wbf$-conjugate, and let $\beta \in {_s}\underline{\Wbf}_{s'}$ be a block. Then the natural map 
$$\Delta_{\dot{w}_{\beta},s'} \to \nabla_{\dot{w}_{\beta}, s'}$$
is an isomorphism. 
\end{lem}

\begin{proof}
This is \cite[Proposition 5.2]{LusztigYun} and \cite[Corollary 8.5.5]{Gouttard}.
\end{proof}

\begin{rque} 
It follows that there are isomorphisms $T_{w_{\beta}, s'} \cong \Delta_{\dot{w}_{\beta},s'} \cong \nabla_{\dot{w}_{\beta}, s'}$. 
\end{rque} 

\begin{lem}
The category $\HH_{[s,s']}$ splits into a direct sum 
$$\HH_{[s,s']} = \bigoplus_{\beta \in {_s}\underline{\Wbf}_{s'}} \HH^{\beta}_{[s,s']}$$
where $\HH^{\beta}_{[s,s']}$ is the full subcategory of $\HH_{[s,s']}$ generated by all the objects $\Delta_{w, s'}$ where $w \in \beta$. 
\end{lem}

\begin{proof}
This is \cite[Proposition 4.11]{LusztigYun} and \cite[Lemma 8.5.6]{Gouttard}.
\end{proof}

\begin{defi}
We denote by $\HH^{\circ}_{[s,s]}$ the category corresponding to $\beta = 1 \in {_s}\underline{\Wbf}_s = \Gammabf_s$. 
\end{defi}

\begin{lem}\label{lem:ConvolRespectsBlocks}
The convolution of $\HH_{\tot}$ respects blocks: if $\beta \in {_s}\underline{\Wbf}_{s'}$ and $\gamma \in {_{s'}}\underline{\Wbf}_{s^{\prime\prime}}$ are composable blocks, then for $A \in \HH^{\beta}_{[s,s']}$ and $B \in \HH^{\gamma}_{[s',s^{\prime\prime}]}$, we have 
$$A * B \in \HH^{\beta\gamma}_{[s,s^{\prime\prime}]}.$$ 
\end{lem}

\begin{proof}
This is \cite[Proposition 4.13]{LusztigYun} and \cite[Lemma 8.5.1]{Gouttard}.
\end{proof}

\begin{rque}
It follows from Lemma \ref{lem:ConvolRespectsBlocks} that the category $\HH^{\circ}_{[s,s]}$ is monoidal and that for all $\beta \in {_s}\underline{\Wbf}_{s'}$ and $\gamma \in {_s'}\underline{\Wbf}_{s}$ the category $\HH^{\beta}_{[s,s']}$ is a left module over $\HH^{\circ}_{[s,s]}$ and $\HH^{\gamma}_{[s',s]}$ is a right module over $\HH^{\circ}_{[s,s]}$. 
\end{rque}

\begin{lem}\label{lem:ModuleOfRankOne}
The category $\HH^{\beta}_{[s,s']}$ is a free left-module of rank one over $\HH^{\circ}_{[s,s]}$.
\end{lem}

\begin{proof}
Consider the functor $\HH^{\circ}_{[s,s]} \to \HH^{\beta}_{[s,s']}$ given by 
$$A \mapsto A * \Delta_{w_{\beta}, s'}.$$
By \cite[Proposition 5.2]{LusztigYun} and \cite[Lemma 8.5.2]{Gouttard} this functor is an equivalence. It is clearly linear over $\HH^{\circ}_{[s,s]}$ (which acts by left convolution). 
\end{proof}

\begin{lem}\label{lem:MoritaFunctorG}
There is a well-defined functor 
$$\HH : \Xibf^{\opp} \to \Morita$$ 
which sends the object $s \in \Ch(\Tbf)$ to the monoidal category $\HH^{\circ}_{[s,s]}$ and the morphism $\beta \in {_s}\underline{\Wbf}_{s'}$ to the $\HH^{\circ}_{[s,s]}$-$\HH^{\circ}_{[s',s']}$-bimodule $\HH^{\beta}_{[s,s']}$. 
\end{lem}

\begin{proof}
There is a priori a lax functor provided by convolution 
$$\HH^{\beta}_{[s,s']} \otimes_{\HH^{\circ}_{[s',s']}} \HH^{\gamma}_{[s',s^{\prime\prime}]} \to \HH^{\beta\gamma}_{[s,s^{\prime\prime}]}.$$ 
By Lemma \ref{lem:ModuleOfRankOne} this morphism is an isomorphism which yields the desired morphism $\Xibf^{\opp} \to \Morita$. 
\end{proof}

\begin{lem}\label{lem:MoritaEquivalence1}
Let $\mathfrak{o}$ be a $\Wbf$-orbit and let $s \in \mathfrak{o}$. Consider the $\HH_{\mathfrak{o},{\tot}}$-$\HH_{[s,s]}$-bimodule $\MM_s = \oplus_{s' \in \mathfrak{o}} \HH_{[s',s]}$. This bimodule defines a Morita equivalence between $\HH_{\mathfrak{o},{\tot}}$ and $\HH_{[s,s]}$. 
\end{lem}

\begin{proof}
Denote by ${_s}\MM = \oplus_{s' \in \mathfrak{o}} \HH_{[s,s']}$. The convolution in $\HH$ yields functors
$${_s}\MM \otimes \MM_s \to \HH_{[s,s]}$$ 
and 
$$\MM_s \otimes {_s}\MM \to \HH_{\mathfrak{o},{\tot}}.$$
They factor over $\HH_{\mathfrak{o},{\tot}}$ and $\HH_{[s,s]}$ respectively and induce functors
$${_s}\MM \otimes_{\HH_{\mathfrak{o},{\tot}}} \MM_s \to \HH_{[s,s]}$$ 
and 
$$\MM_s \otimes_{\HH_{[s,s]}} {_s}\MM \to \HH_{\mathfrak{o},{\tot}}.$$
By \cite[4.6.2.1]{HigherAlgebra}, it is enough to show that these functors are equivalences. 

For the first one, consider the map 
$$\HH_{[s,s]} \to {_s}\MM \otimes_{\HH_{\mathfrak{o}}} \MM_s, A \mapsto A \otimes \Delta_{1,s}.$$
Clearly, the composition $\HH_{[s,s]} \to  {_s}\MM \otimes_{\HH_{\mathfrak{o}}} \MM_s \to \HH_{[s,s]}$ is the identity, hence the first map is fully faithful. Moreover for all $A \otimes B \in {_s}\MM \otimes_{\HH_{\mathfrak{o}}} \MM_s$, we have $A \otimes B = (A * B) \otimes \Delta_{1,s}$. It follows that the first map in the above composition is essentially surjective and therefore an isomorphism. 

For the second map, note that 
$$\MM_s \otimes_{\HH_{[s,s]}} {_s}\MM = \bigoplus_{s',s^{\prime\prime}} \HH_{[s',s]} \otimes_{\HH_{[s,s]}} \HH_{[s,s^{\prime\prime}]}.$$
Arguing as in Lemma \ref{lem:ModuleOfRankOne}, the category $\HH_{[s',s]}$ is free of rank one over $\HH_{[s,s]}$ hence $\HH_{[s',s]} \otimes_{\HH_{[s,s]}} \HH_{[s,s^{\prime\prime}]} = \HH_{[s',s^{\prime\prime}]}$ and this second map is an isomorphism. 
\end{proof}

\subsection{The Hecke category $\HH^{\Hbf}$}

We fix a $\Wbf$-orbit $\mathfrak{o} \subset \Ch(\Tbf)$. Recall that we have a groupoid $\Hbffrak_{\mathfrak{o}} \to \Xibf_{\mathfrak{o}}$ constructed in Section \ref{sec:EndoscopicGroup}. For $s \in \mathfrak{o}$, the fiber at $s$ is the group $\Hbf^{\circ}_s$ and the fiber over a block $\beta \in {_s}\underline{\Wbf}_{s'}$ yields a $\Hbf^{\circ}_s$-$\Hbf^{\circ}_{s'}$ bitorsor $\Hbf^{\beta}_{[s,s']}$. Moreover all the groups $\Hbf^{\circ}_{s}$ come equipped with a pinning (relative to $\Gbf$) and with a Borel pair $\Bbf_s = \Ubf_s\Tbf_s$. 

For an element $s \in \mathfrak{o}$, we denote by 
$$\HH^{\Hbf,\circ}_{[s,s]} = \DD(\Ubf_s \backslash \Hbf^{\circ}_{s} / \Ubf_s, \RR_{\Tbf_s \times \Tbf_s})_{\unip, \unip},$$
and for a block $\beta \in {_s}\underline{\Wbf}_{s'}$, we define the category 
$$\HH^{\Hbf, \beta}_{[s,s']} = \DD(\Ubf_s\backslash \Hbf^{\beta}_{[s,s']}/\Ubf_{s'}, \RR_{\Tbf_s \times \Tbf_{s'}})_{\unip, \unip}$$
of free monodromic unipotent sheaves on $\Ubf_s\backslash \Hbf^{\beta}_{[s,s']}/\Ubf_{s'}$. 

Let $\beta \in {_s}\underline{\Wbf}_{s'}$ and $\gamma \in {_{s'}}\underline{\Wbf}_{s^{\prime\prime}}$ be composable blocks, and consider the diagram 
\[\begin{tikzcd}
	& {\Ubf_s \backslash \Hbf^{\beta}_{[s,s']} \times^{\Ubf_{s'}}\Hbf_{[s',s^{\prime\prime}]}^{\gamma} /\Ubf_{s^{\prime\prime}}} & {\Ubf_s \backslash \Hbf^{\beta\gamma}_{[s,s^{\prime\prime}]}/\Ubf_{s^{\prime\prime}}} \\
	{\Ubf_s \backslash \Hbf^{\beta}_{[s,s']}/\Ubf_{s'}} && {\Ubf_{s'} \backslash \Hbf_{[s',s^{\prime\prime}]}^{\gamma} /\Ubf_{s^{\prime\prime}}}
	\arrow["m", from=1-2, to=1-3]
	\arrow["{\pr_1}"{description}, from=1-2, to=2-1]
	\arrow["{\pr_2}"{description}, from=1-2, to=2-3]
\end{tikzcd}\]
We define the convolution product as follows:
$$\HH^{\Hbf,\beta} \times \HH^{\Hbf,\gamma} \to \HH^{\Hbf, \beta\gamma}, (A,B) \mapsto m_!(A \widehat{\boxtimes} B)[\dim \Tbf].$$ 
Consider the category 
$$\HH^{\Hbf}_{\tot} = \oplus_{\beta} \HH^{\Hbf, \beta}.$$
We can equip this category with a monoidal structure such that for $A \in \HH^{\beta}$ and $B \in \HH^{\gamma}$ the convolution $A * B$ is zero if $\beta$ and $\gamma$ are not composable and given by the formula above if they are. 
Similarly, we can fix an element $s \in \mathfrak{o}$ and consider the category 
$$\HH^{\Hbf}_{[s,s]} = \oplus_{\beta \in \Gammabf_s} \HH^{\Hbf, \beta}_{[s,s]}.$$

As in the case of $\Gbf$, we define the standard, costandard, and tilting objects. First recall that 
$$\Hbf^{\beta}_{[s,s']} = \widetilde{\Xibf}^{\beta}_{[s,s']} \times^{\Tbf_{s'}} \Hbf^{\circ}_{s'} = \Hbf^{\circ}_s \times^{\Tbf_{s}} \widetilde{\Xibf}_{[s,s']}^{\beta}.$$
We denote by $\Nbf^{\beta}_{[s,s']}$ the following subscheme of $\Hbf^{\beta}_{[s,s']}$:
$$\Nbf^{\beta}_{[s,s']} = \widetilde{\Xibf}^{\beta}_{[s,s']} \times^{\Tbf_{s'}} \Nbf^{\circ}_{s'} = \Nbf^{\circ}_s \times^{\Tbf_{s}} \widetilde{\Xibf}_{[s,s']}^{\beta}$$
where $\Nbf^{\circ}_s \subset \Hbf^{\circ}_s$ is the normalizer of $\Tbf_{s}$. Note that we have a canonical isomorphism of sets
$$\Nbf^{\beta}_{[s,s']}/\Tbf_{s'} = \Tbf_{s} \backslash \Nbf^{\beta}_{[s,s']} = \beta.$$

Let $w \in \beta$ and $n \in \Nbf^{\beta}_{[s,s']}$. We obtain a splitting
$$\Ubf_s \backslash \Bbf_sn\Bbf_{s'} /\Ubf_{s'} = n\Tbf_{s'}/(\Ubf_{s'} \cap {^n}\Ubf_s)$$
where ${^n}\Ubf_s \subset \Hbf^{\circ}_{s'}$. In particular, we obtain a map $\nu_n : \Ubf_s \backslash \Bbf_sn\Bbf_{s'} /\Ubf_{s'} \to \Tbf_{s'}$.
We denote by $j_w : \Ubf_s \backslash \Bbf_sn\Bbf_{s'} /\Ubf_{s'} \subset \Ubf_s \backslash \Hbf^{\beta}_{[s,s']} /\Ubf_{s'}$ the inclusion of the stratum and by
\begin{enumerate}
\item $\Delta_{n,s'} = j_{w,!}\nu_n^*L_{\Tbf_{s'}}[\ell_{\beta}(w) + \dim \Tbf]$, the standard sheaves, 
\item $\nabla_{n,s'} = j_{w,*}\nu_n^*L_{\Tbf_{s'}}[\ell_{\beta}(w) + \dim \Tbf]$, the costandard sheaves.
\end{enumerate}
As before, we get a notion of tilting sheaves, and it is clear that the analogue of Theorem \ref{thm:StructureCategoryviaTilting} holds for $\HH^{\Hbf, \beta}_{[s,s']}$ in place of $\HH^{\Gbf}$. 

\begin{lem}
Let $w \in \beta$ and $w' \in \gamma$ where $\beta \in {_s}\underline{\Wbf}_{s'}$ and $\gamma \in {_{s'}}\underline{\Wbf}_{s^{\prime\prime}}$ are two composable blocks. Let $n,n'$ be lifts of $w \in \Nbf^{\beta}_{[s,s']}$ and $w'\in \Nbf^{\gamma}_{[s',s^{\prime\prime}]}$. 
\begin{enumerate}
\item If $\ell_{\beta}(w) + \ell_{\gamma}(w') = \ell_{\beta\gamma}(ww')$, then there are canonical isomorphisms 
$$\Delta_{n,s'} * \Delta_{n',s^{\prime\prime}} \cong \Delta_{nn',s^{\prime\prime}}, \nabla_{n,s'} * \nabla_{n',s^{\prime\prime}} \cong \nabla_{nn',s^{\prime\prime}}.$$
\item There are isomorphisms 
$$\Delta_{n,n^{-1}s} * \nabla_{n^{-1},s} \cong \Delta_{1,s} \cong \nabla_{n,n^{-1}s} * \Delta_{n^{-1},s}.$$
\end{enumerate}
\end{lem}

Let $\beta \in {_s}\underline{\Wbf}_{s'}$ be a block. Since $\Hbf^{\beta}_{[s,s']}$ is a $\Hbf^{\circ}_{s}$-$\Hbf^{\circ}_{s'}$-bitorsor, the following lemma is immediate (compare with Lemmas \ref{lem:ModuleOfRankOne} and \ref{lem:MoritaFunctorG}).
\begin{lem}
The category $\HH^{\Hbf,\beta}_{[s,s']}$ is a free $\HH^{\Hbf, \circ}_{[s,s]}$-module of rank $1$. Moreover, there is a well-defined functor 
$$\HH^{\Hbf} : \Xibf_{\mathfrak{o}}^{\opp} \to \Morita$$
such that the image of $s \in \mathfrak{o}$ is $\HH^{\Hbf, \circ}_{[s,s]}$ and the image of $\beta \in {_s}\underline{\Wbf}_{s'}$ is $\HH^{\Hbf, \beta}_{[s,s']}$. 
\end{lem}

\subsection{Soergel bimodules}\label{sec:SoergelBim}

In this section, we assume that $\ell$ is not a torsion prime for $\Gbf^*$ (see Definition \ref{def:torsion-prime}). This hypothesis allows us to use the Pittie-Steinberg Theorem even if $\Gbf$ does not have connected center. 
\begin{thm}[\protect{\cite{PittieSteinberg}, \cite[Proposition 10.2.3]{Gouttard}}]\label{thm:PittieSteinbergCompleted}
Let $s \in \Ch(\Tbf)$,
\begin{enumerate}
\item the ring $\RR_{\Tbf}$ is free of finite rank $|\Wbf_s^{\circ}|$ over $\RR_{\Tbf}^{\Wbf_{s}^{\circ}}$, 
\item there exists a basis $(e_w)_{w \in \Wbf_{s}^{\circ}}$ of $\RR_{\Tbf}$ such that 
$$\det(v(e_w))_{v,w \in \Wbf_s^{\circ}}$$
is nonzero in the fraction field of $\RR_{\Tbf}^{\Wbf_{s}^{\circ}}$. 
\end{enumerate}
\end{thm}

We define the schemes over which Soergel bimodules will be defined.
\begin{enumerate}
\item For $s \in \Ch(\Tbf)$, we denote by
$$\Ccal(\Tbf)_s = \Ccal_s \colondef= s \times \Spec(\RR_{\Tbf}).$$
\item We denote by 
$$\Ccal(\Tbf) = \bigsqcup_{s \in \Ch(\Tbf)} \Ccal_s.$$
\end{enumerate}
Let us specify the meaning of $\Ccal_s$. To $s \in \Ch(\Tbf)$  corresponds a $\Lambda$-point of $\Tbf^{\vee}$ which is defined over a finite extension $\Lambda_1$ of $\Lambda_0$ (recall that $\Lambda_0$ denotes the `prime subring' of $\Lambda$). The scheme $\Ccal_s$ is then canonically identified with $\Tbf_{\Lambda_1, s}^{\vee, \wedge} \times_{\Lambda_1} \Lambda$, where $\Tbf_{\Lambda_1, s}^{\vee, \wedge}$ is the completion of $\Tbf_{\Lambda_1}^{\vee}$ at the point $s$. It follows that there is a canonical map $\Ccal(\Tbf) \to \Tbf^{\vee}$. 
Up to the extensions $\Lambda \times_{\Lambda_1} -$ we will refer to $\Ccal(\Tbf)$ as the formal completion of $\Tbf^{\vee}$ at all the torsion points of $\Tbf^{\vee}$. Let $\mathfrak{o} \subset \Ch(\Tbf)$ be a $\Wbf$-orbit. We denote by $\Ccal(\Tbf)_{\mathfrak{o}}$ the union of all connected components of $\Ccal(\Tbf)$ indexed by elements $s \in \mathfrak{o}$. 

For $s,s' \in \Ch(\Tbf)$ and $\beta \in {_s}\underline{\Wbf}_{s'}$, we denote by 
$$\Ccal^{\beta}_{[s,s']} \subset \Ccal_s \times \Ccal_{s'}$$ 
the reduced closed subscheme of $\Ccal_s \times \Ccal_{s'}$ which is the union of all the graphs of the elements $w \in \beta$. Note that an element $w \in \beta$ satisfies $w(s') = s$ hence defines a morphism $\Ccal_{s'} \to \Ccal_s$. Given $s \in \Ch(\Tbf)$, we also denote by 
$$\Ccal^{\circ}_{s}$$
the scheme $\Ccal^{\beta}_{[s,s']}$ for $s = s'$ and $\beta = 1$. Note that all the schemes $\Ccal_{[s,s']}^{\beta}$ for $s, s'$ in a single orbit $\mathfrak{o}$ are all (noncanonically) isomorphic.

\begin{lem}
We have 
$$\Ccal^{\circ}_{s} = s \times \Spec(\RR_{\Tbf}) \times_{\Spec(\RR_{\Tbf}^{\Wbf_s^{\circ}})} \Spec(\RR_{\Tbf}) \times s.$$ 
\end{lem}

\begin{proof}
Both schemes are closed subschemes of $\Ccal_s \times \Ccal_s$ and they have the same points. It is enough to show that the right-hand side is reduced. This follows from Theorem \ref{thm:PittieSteinbergCompleted}. 
\end{proof}

Let $\beta$ and $\gamma$ be composable blocks in $\mathfrak{o}$. We define a convolution product
$$\DD_{\qcoh}(\Ccal_{[s,s']}^{\beta}) \times \DD_{\qcoh}(\Ccal_{[s',s^{\prime\prime}]}^{\gamma}) \to \DD_{\qcoh}(\Ccal_{[s,s^{\prime\prime}]}^{\beta\gamma})$$ 
as follows. Consider the following diagram
\[\begin{tikzcd}
	& {\Ccal_{s} \times \Ccal_{s'} \times \Ccal_{s^{\prime\prime}}} & {\Ccal_{s} \times \Ccal_{s^{\prime\prime}}} \\
	& {\Ccal_{[s,s']}^{\beta} \times_{\Ccal_{s'}} \Ccal_{[s',s^{\prime\prime}]}^{\gamma}} & {\Ccal_{[s,s^{\prime\prime}]}^{\beta\gamma}} \\
	{\Ccal_{[s,s']}^{\beta}} && {\Ccal_{[s',s^{\prime\prime}]}^{\gamma}}
	\arrow["{\pr_{13}}"{description}, from=1-2, to=1-3]
	\arrow[from=2-2, to=1-2]
	\arrow["{\pr_{13}}"{description}, from=2-2, to=2-3]
	\arrow["{\pr_1}"{description}, from=2-2, to=3-1]
	\arrow["{\pr_2}"{description}, from=2-2, to=3-3]
	\arrow[hook, from=2-3, to=1-3]
\end{tikzcd}\]
where $\pr_1, \pr_2$ and $\pr_{13}$ are the obvious projections. The convolution is defined as 
$$(A,B) \mapsto \pr_{13,*}(\pr_1^*A \otimes \pr_2^*B).$$
Noting that objects of $\DD_{\qcoh}(\Ccal_{[s,s']}^{\beta})$ are $\Ocal(\Ccal_s)$-$\Ocal(\Ccal_{s'})$-bimodules, this convolution is nothing else than the convolution of bimodules.

The following statements are standard.
\begin{enumerate}
\item The convolution product defined above equips $\DD_{\qcoh}(\Ccal_{[s,s]}^{\circ})$ with a monoidal structure.
\item The category $\DD_{\qcoh}(\Ccal_{[s,s']}^{\beta})$ is a bimodule over $\DD_{\qcoh}(\Ccal_{[s,s]}^{\circ})$-$\DD_{\qcoh}(\Ccal_{[s',s']}^{\circ})$ that is free of rank one as a left and right module respectively. 
\end{enumerate}

As before, we can organize these categories into a functor 
$$\DD_{\qcoh}(\Ccal \times_{\Ccal\sslash W} \Ccal) : \Xibf^{\opp} \to \Morita$$ 
such that the image of $s$ is $\DD_{\qcoh}(\Ccal_{[s,s]}^{\circ})$ and the image of $\beta$ is $\DD_{\qcoh}(\Ccal_{[s,s']}^{\beta})$. 

\begin{defi}
Let $s \in \mathfrak{o}$ and let $w$ be a simple reflection of $\Wbf$ (which we do not write as $s$ to avoid a conflict of notation). We define the Bott-Samelson bimodule 
$$\BS(w,s) \colondef= \RR_{\Tbf} \otimes_{\RR_{\Tbf}^w} \RR_{\Tbf}$$ 
if $ws = s$, and 
$$\BS(w,s) \colondef= \RR_{\Tbf}$$
if $ws \neq s$. In both cases, this object is understood as an object of $\Coh(\Ccal_{[ws,s]}^{\beta})$, where $\beta$ is the block such that $w \in \beta$. 

Let $\underline{w} = (w_1, \dots, w_n)$ be a sequence of simple reflections. We define the Bott-Samelson bimodule 
$$\BS(\underline{w},s) \colondef= \BS(w_1, (w_2\dots w_n)(s)) * \dots * \BS(w_{n - 1}, w_ns) *\BS(w_n, s).$$
This is an object of $\Coh(\Ccal_{[ws,s]}^{\beta})$ where $\beta \in {_{ws}}\underline{\Wbf}_s$ is the block containing $w = w_1\dots w_n$. 
\end{defi}

\begin{lem}
For all sequences $(w_1, \dots, w_n)$ of simple reflections, the Bott-Samelson bimodule $\BS(\underline{w},s)$ is concentrated in degree $0$. 
\end{lem}

\begin{proof}
Writing that $\Ccal_{[s,s']}^{\beta}$ is a subscheme of $\Ccal_s \times \Ccal_{s'}$, we check by induction on the length of $\underline{w}$ that the $\BS(\underline{w},s)$ is free as a $\Ocal(\Ccal_s)$-module and as a $\Ocal(\Ccal_{s'})$-module separately. If $\ell(\underline{w}) = 1$ then the statement is clear by definition. Assume the statement is true for sequences of length $n-1$, then $\BS(\underline{w},s) = \BS((w_1, \dots, w_{n-1}), w_ns) \otimes_{\Ocal({\Ccal_{w_ns}})} \BS(w_n,s)$ which is clearly free as a left (resp. right) module over $\Ocal({\Ccal_{ws}})$ (resp. $\Ocal({\Ccal_s})$). 
\end{proof}

\begin{defi}
Let $\beta \in {_s}\underline{\Wbf}_{s'}$ be a block in $\mathfrak{o}$. We define the category of Soergel bimodules $\SBim^{\beta}_{[s,s']}$ as the full subcategory of $\Coh(\Ccal_{[s,s']}^{\beta})$ generated by all Bott-Samelson bimodules $\BS(\underline{w}, s')$ where $ws' = s$ and $w \in \beta$ under isomorphisms, direct sums, and direct summands. 

If $\beta = 1$, then we denote this category by $\SBim^{\circ}_{[s,s]}$. In other words, the category of Soergel bimodules is the Karoubian completion of the category of Bott-Samelson bimodules.  
\end{defi} 

By construction the class of Bott-Samelson bimodules is closed under convolution, it follows that the same property holds for the category of Soergel bimodules. We therefore have the following structures 
\begin{enumerate}
\item $\SBim^{\circ}_{[s,s]}$ is a monoidal additive category, 
\item the category $\SBim^{\beta}_{[s,s']}$ is a $\SBim^{\circ}_{[s,s]}$-$\SBim^{\circ}_{[s',s']}$-bimodule. 
\end{enumerate}
We consider the homotopy categories $\KK(\SBim^{\beta}_{[s,s']})$ of the categories of Soergel bimodules. They are equipped with the appropriate monoidal/module structures as above.

\begin{lem}
The category $\KK(\SBim^{\beta}_{[s,s']})$ is free of rank one over $\KK(\SBim^{\circ}_{[s,s]})$. 
\end{lem}

\begin{proof}
Let $w_{\beta}$ be the minimal element of $\beta$ and write $w_{\beta} = w_1\dots w_n$ as a product of simple reflections. For all $i = 1, \dots,n$ we have $w_i \dots w_n s' \neq w_{i+1} \dots w_n s'$, where the empty product $w_{n+1}\cdots w_n$ is the identity; hence $\BS((w_1, \dots, w_n), s')$ is the structure sheaf of the graph of $w_{\beta}$ in $\Ccal_s \times \Ccal_{s'}$, hence $- * \BS((w_1, \dots, w_n), s') : \SBim^{\circ}_{[s,s]} \to \SBim^{\beta}_{[s,s']}$ is an exact equivalence.
\end{proof}

As before, we can organize all these categories into a single functor
$$\SBim : \Xibf^{\opp} \to \Morita$$
which sends an element $s$ to $\KK(\SBim^{\circ}_{[s,s]})$ and $\beta$ to $\KK(\SBim^{\beta}_{[s,s']})$. 

Let $\mathfrak{o} \subset \Ch(\Tbf)$ be a $\Wbf$-orbit. We denote by 
$$\Ccal_{\mathfrak{o}} \colondef= \bigsqcup_{\beta} \Ccal_{[s,s']}^{\beta}$$
the disjoint union of all schemes $\Ccal_{[s,s']}^{\beta}$ where $\beta \in {_s}\underline{\Wbf}_{s'}$ and $s,s' \in \mathfrak{o}$. We can equip the category $\DD_{\qcoh}(\Ccal_{\mathfrak{o}})$ with the following monoidal structure. Let $A \in \DD_{\qcoh}(\Ccal^{\beta}_{[s_1,s_2]})$ and $B \in \DD_{\qcoh}(\Ccal^{\gamma}_{[s_3,s_4]})$. Then we set 
\begin{enumerate}
\item $A * B = 0$ if $s_2 \neq s_3$, i.e., $\beta$ and $\gamma$ are not composable, 
\item $A * B$ is the convolution defined above if $s_2 = s_3$. 
\end{enumerate}

\begin{rque}\label{rque:ConnectedCenterSpaceOfComponents}
Assume that $\Gbf$ has connected center, in this case all the sets ${_s}\underline{\Wbf}_{s'}$ are reduced to points by Lemma \ref{lem:ControlCentralizer}. Consider the scheme $\Ccal(\Tbf) \times_{\Ccal(\Tbf)\sslash \Wbf} \Ccal(\Tbf)$, its connected components are indexed by pairs $(s,s')$ such that there exists $w \in \Wbf$ with $ws' = s$. The connected component attached to such a pair $(s,s')$ is canonically identified with $\Ccal_{[s,s']}^{\beta}$ where $\beta$ is the unique element of ${_s}\underline{\Wbf}_{s'}$. 
\end{rque}

\subsection{Soergel realization for $\HH^{\Gbf}$}\label{sec:SoergelRealizationForG}

\subsubsection{Whittaker model}\label{sec:WhittakerModelForG}
Recall that we have fixed a pinning of the group $\Gbf$. We denote by $\Ubf^-$ the unipotent of $\Bbf^-$, the Borel subgroup opposite to $\Bbf$. The pinning provides a canonical generic character
$$\phi : \Ubf^- \to \Ubf^-/[\Ubf^-,\Ubf^-] = \prod_{\Delta} \Ga \xrightarrow{\Sigma} \Ga$$
where the second morphism comes from the identification of the simple root groups with $\Ga$ through the pinning. We fix a primitive $p$-th root of $1$ in $\Lambda$ and an Artin-Schreier sheaf $\AS$ on $\Ga$. Finally, we denote by 
$$\Lcal_{\psi} \colondef= \phi^*\AS$$
the pullback of the Artin-Schreier sheaf to $\Ubf^-$.

We denote by $\Whitt(\HH^{\Gbf}_{\tot})$ the category 
$$\Whitt(\HH^{\Gbf}_{\tot}) \colondef= \bigoplus_{s \in \Ch(\Tbf)} \DD((\Ubf^-,\Lcal_{\psi}) \backslash \Gbf/\Ubf, \RR_{\Tbf})_{[-,s]}$$ 
of free monodromic sheaves (with respect to the right action of $\Tbf$) that are left $(\Ubf^-, \Lcal_{\psi})$-equivariant. 

\begin{lem}[\cite{BBMWhitttaker}]\label{lem:WhittakerSheavesSupport}
Any sheaf in $\Whitt(\HH^{\Gbf}_{\tot})$ is supported on the open cell $\Ubf^-\Bbf$.
\end{lem}

Consider the inclusion $i : \Tbf \hookrightarrow \Gbf$. It follows from Lemma \ref{lem:WhittakerSheavesSupport} that the functor $i^*$ induces an equivalence 
$$i^* : \Whitt(\HH^{\Gbf}_{\tot}) \cong \bigoplus_s \DD(\Tbf, \RR_{\Tbf})_{s} \cong \DD_{\qcoh}(\Ccal(\Tbf))$$
where the second equivalence follows from taking the fiber at $1 \in \Tbf$. 

\begin{defi}[Averaging functors]\label{def:AveragingFunctors}
There is a pair of functors 
$$\Av_{\Ubf} : \Whitt(\HH^{\Gbf}_{\tot}) \leftrightarrows \HH^{\Gbf}_{\tot} : \Av_{\psi}$$
given by the following formulas
$$\Av_{\Ubf} = a_!(\Lambda \boxtimes -)[\dim \Ubf], \Av_{\psi} = \bar{a}_!(\Lcal_{\psi} \boxtimes -)[\dim \Ubf]$$ 
where $a : \Ubf \times \Gbf/\Ubf \to \Gbf/\Ubf$ and $\bar{a} : \Ubf^- \times \Gbf/\Ubf \to \Gbf/\Ubf$ are the action maps. 
\end{defi}

\begin{lem}[\protect{\cite[Theorem 1.5]{BBMWhitttaker}}]\label{lem:AveragingPsiIsRightAdjoints}
The functor $\Av_{\psi}$ is perverse $t$-exact and right adjoint to $\Av_{\Ubf}$. 
\end{lem}

\subsubsection{The big tilting sheaf}\label{subsec:BigTilting}

Recall that we denote by $\Delta_1 = \oplus_{s} \Delta_{1,s}$ the unit in $\HH^{\Gbf}_{\tot}$. We denote by 
$$\mathbb{T} \colondef= \Av_{\Ubf} \Av_{\psi} \Delta_1 = \bigoplus_s \mathbb{T}_s.$$

\begin{lem}\label{lem:MultiplicityOne}
The object $\mathbb{T}$ is a direct sum of tilting objects, more precisely we have
$$\mathbb{T} = \bigoplus_{s, \beta} T_{w^{\beta},s}$$
where $s$ ranges through all elements of $\Ch(\Tbf)$ and $\beta$ through all blocks. (Recall that $w^{\beta}$ denotes the maximal element of the block $\beta$ and $T_{w^{\beta},s}$ the unique indecomposable tilting sheaf supported on the closure of $\Bbf w^{\beta}\Bbf$ and containing $\Delta_{w^{\beta},s}$ with multiplicity one.) 
\end{lem}

\begin{proof}
This is \cite[Lemma 4.4.11]{BezrukavnikovYun}, \cite[Lemma 10.1]{BezrukavnikovRicheSoergelTheory} and \cite[Proposition 12.9.3]{Gouttard}.
\end{proof}

\begin{rque}
The objects $T_{w^{\beta},s}$ provided by Theorem \ref{thm:StructureCategoryviaTilting} are not canonical. However the object $\mathbb{T}$ depends only on the pinning of $\Gbf$ and on no other choice. Taking the various direct summands of Lemma \ref{lem:MultiplicityOne}, we get a canonical construction of these objects. 
\end{rque}

\begin{lem}
There is a canonical isomorphism of functors $\HH^{\Gbf}_{\tot} \to \HH^{\Gbf}_{\tot}$
$$ \Av_{\Ubf} \Av_{\psi} \cong \mathbb{T} * -.$$
\end{lem}

\begin{proof}
This is \cite[Lemma 4.4.11]{BezrukavnikovYun} and \cite[Proposition 12.9.3]{Gouttard}. 
\end{proof}

\begin{corol}
The object $\mathbb{T}$ has a canonical coalgebra structure. 
\end{corol}

\begin{proof}
Since $\Av_{\Ubf}$ is left adjoint to $\Av_{\psi}$, it follows that $\Av_{\Ubf} \Av_{\psi}$ is a comonad acting on $\HH^{\Gbf}$. The comultiplication of $\mathbb{T}$ is given by the unit of the adjunction $\id \to \Av_{\psi}\Av_{\Ubf}$
$$\mathbb{T} \cong \Av_{\Ubf} \Av_{\psi}\Delta_1 \to \Av_{\Ubf} \Av_{\psi}\Av_{\Ubf} \Av_{\psi}\Delta_1 \cong \mathbb{T} * \mathbb{T}.$$ 
The counit of $\mathbb{T}$ 
$$\mathbb{T} \cong \Av_{\Ubf}\Av_{\psi}\Delta_1 \to \Delta_1$$ 
comes from the adjunction morphism $\Av_{\Ubf}\Av_{\psi} \to \id$. 
\end{proof}

\subsubsection{Linearity}\label{subsect:Linearity}

\begin{lem}
The category $\HH^{\Gbf}_{\tot}$ is $\Ccal(\Tbf) \times \Ccal(\Tbf)$-linear. This structure induces a $\Ccal_{[s,s']}^{\beta}$-linear structure on the direct summand $\HH^{\beta}_{[s,s']}$.
\end{lem}

\begin{proof}
The two $\Ccal(\Tbf)$-linear structures come from the left and right monodromy actions of Verdier, we refer to \cite[Section 2.5]{EteveFreeMonodromic} for a discussion. The second statement follows from the formality of $\HH^{\Gbf}_{\tot}$ by Theorem \ref{thm:StructureCategoryviaTilting} and the fact that for standard sheaves $\End(\Delta_{w,s}) = \RR_{\Tbf}$ with the bimodule structures given by the graph of $w$. 
\end{proof}

\begin{lem}\label{lem:LinearityAlgebra}
The category $\HH^{\Gbf}_{\tot} \in \Pr_{\Lambda}$ naturally lifts to a monoidal category in $\Pr_{\Ccal(\Tbf)\sslash \Wbf}$.
\end{lem}

\begin{proof}
By Theorem \ref{thm:StructureCategoryviaTilting}, the category $\HH_{\tot}$ is generated under shifts, cones, retracts, and colimits by tilting objects. It is therefore enough to show that the convolution functor $\HH_{\tilt} \times \HH_{\tilt} \to \HH_{\tilt}$ is $\Ccal(\Tbf)\sslash \Wbf$-linear. Since tilting sheaves are perverse, this can be checked in the abelian category $\HH^{\heartsuit}_{\tot}$. It follows from the fact that the left and right $\Ccal(\Tbf)$-linear structures on $\HH_{\tot}$ agree after restriction to $\Ccal(\Tbf)\sslash \Wbf$.
\end{proof}

\begin{rque}
Note that we have a natural map $\Ccal(\Tbf)\sslash \Wbf \to \Tbf^{\vee} \sslash \Wbf$, hence we can also consider the category $\HH^{\Gbf}_{\tot}$ as a monoidal category in $\Pr_{\Tbf^{\vee} \sslash \Wbf}$.
\end{rque}

\subsubsection{Endomorphismensatz}

\begin{thm}[Endomorphismensatz]\label{thm:Endomorphismensatz}
Assume that $\ell$ is not a torsion prime for $\Gbf^*$ (see Definition \ref{def:torsion-prime}). Then the canonical map 
$$\Ocal(\Ccal_{[s,s']}^{\beta}) \to \End(T_{[s,s']}^{\beta})$$ 
is an isomorphism. 
\end{thm}

\begin{proof}
This is \cite[Proposition 4.7.3]{BezrukavnikovYun}, \cite[Theorem 9.1]{BezrukavnikovRicheSoergelTheory} and \cite[Theorem 10.4.1]{Gouttard}. Again, we note that in the references we have provided, the authors work in the setting where $\Lambda$ is a field. Since the map in Theorem \ref{thm:Endomorphismensatz} is compatible with reduction mod $\ell$, the $\Zlb$-statement follows from the $\Flb$-statement by Nakayama's lemma. 
\end{proof}

\begin{rque}\label{rque:WWAction}
In the proof of Theorem \ref{thm:Endomorphismensatz}, a key point is to construct a $\Wbf \times \Wbf$-action on $\End(T_{\mathfrak{o}})$, where $T_{\mathfrak{o}} = \oplus_{s,s' \in \mathfrak{o}, \beta} T_{w^{\beta},s'}$. This action is constructed as follows. Set $\Delta_{w, \mathfrak{o}} = \oplus_{s \in \mathfrak{o}} \Delta_{w,s}$. Then there are isomorphisms
$$\Delta_{w, \mathfrak{o}} * T_{\mathfrak{o}} * \Delta_{v, \mathfrak{o}} \cong T_{\mathfrak{o}}$$
for all $v,w \in \Wbf$, by \cite[Proposition 7.10]{BezrukavnikovRicheSoergelTheory} and \cite[Proposition 9.7.6]{Gouttard}. For all $w,v \in \Wbf$ this defines a morphism 
$$\End(T_{\mathfrak{o}}) \to \End(\Delta_{v,\mathfrak{o}} * T_{\mathfrak{o}} * \Delta_{w,\mathfrak{o}}) \cong \End(T_{\mathfrak{o}})$$
where the first map is induced by the functor $A \mapsto \Delta_{v,\mathfrak{o}} * A * \Delta_{w,\mathfrak{o}}$ and the second by the above isomorphism. It is shown in \cite[Lemma 9.6]{BezrukavnikovRicheSoergelTheory} and \cite[Lemma 10.4.5]{Gouttard} that this morphism does not depend on any  previous choices and defines a canonical action of $\Wbf \times \Wbf$ on $\End(T_{\mathfrak{o}})$. 
\end{rque}

\begin{rque}
Assume that $\Gbf$ has connected center, see Remark \ref{rque:ConnectedCenterSpaceOfComponents}, then we can reformulate Theorem \ref{thm:Endomorphismensatz} as 
\begin{equation}
\End(\mathbb{T}) \cong \Ocal(\Ccal(\Tbf) \times_{\Ccal(\Tbf)\sslash \Wbf} \Ccal(\Tbf)).
\end{equation}
\end{rque}

\subsubsection{Rigidified minimal standard sheaves}\label{sec:Rigidified-minimal-standard}

\begin{defi}[\protect{\cite[Section 5.11]{LusztigYunV3}}]\label{def:min-standard}
Let $\beta \in {_s}\underline{\Wbf}_{s'}$ be a block. A rigidified minimal standard sheaf for $\beta$ is a pair $(\Delta^{\beta}, \alpha)$ where
\begin{enumerate}
\item $\Delta^{\beta} \in \HH^{\beta}_{[s,s']}$ is isomorphic to $\Delta_{w_{\beta},s'}$ 
\item $\alpha$ is an isomorphism $\RR_{\Tbf} \cong 1^*\Av_{\psi}\Delta^{\beta}[\dim \Tbf]$. 
\end{enumerate}
\end{defi}

\begin{lem}\label{lem:min-standard-existence}
For all blocks $\beta \in {_s}\underline{\Wbf}_{s'}$, there exists, up to unique isomorphism, a unique rigidified minimal standard sheaf $(\Delta^{\beta}, \alpha)$. 
\end{lem}

\begin{proof}
By \cite[Lemma 12.9.2]{Gouttard}, we have that $1^*\Av_{\psi}(\Delta_{w_{\beta},s'}) \cong \RR_{\Tbf}$. Hence, there exists a rigidified minimal standard sheaf. 
Let $(\Delta^{\beta}, \alpha)$ and $(\Delta'^{\beta}, \alpha')$ be two rigidified minimal standard sheaves. Let $\lambda : \Delta^{\beta} \cong \Delta'^{\beta}$ be an isomorphism. Applying the functor $1^*\Av_{\psi}(-)[\dim \Tbf]$ to the isomorphism $\lambda$, we get 
$$\RR_{\Tbf} \xleftarrow{\alpha} 1^*\Av_{\psi}(\Delta^{\beta})[\dim \Tbf] \cong 1^*\Av_{\psi}(\Delta'^{\beta})[\dim \Tbf] \xrightarrow{\alpha'} \RR_{\Tbf}.$$
The composition yields an element $\mu \in \RR_{\Tbf}^{\times}$. After rescaling by $\frac{1}{\mu}$, we get a commutative diagram 
\[\begin{tikzcd}
	{1^*\Av_{\psi}(\Delta^{\beta})[\dim \Tbf]} && {1^*\Av_{\psi}(\Delta'^{\beta})[\dim \Tbf]} \\
	& {\RR_{\Tbf}}
	\arrow["{\frac{1}{\mu}\lambda}", Rightarrow, no head, from=1-1, to=1-3]
	\arrow["\alpha"', from=1-1, to=2-2]
	\arrow["{\alpha'}", from=1-3, to=2-2]
\end{tikzcd}\]
and thus an isomorphism of pairs $(\Delta^{\beta}, \alpha) \cong (\Delta'^{\beta}, \alpha')$. Finally, since $\End(\Delta^{\beta}) = \RR_{\Tbf}$, this isomorphism is unique. 
\end{proof}

Consider the constant functor 
$$C : \Xibf^{\opp} \to \Morita$$ 
where $C(s) = \DD(\Lambda)$ and $C(\beta) = \DD(\Lambda)$ (seen as a $\DD(\Lambda)$-bimodule in an obvious way). 

\begin{lem}\label{lem:RigidifiedMinStandard}
Let $\beta \in {_s}\underline{\Wbf}_{s'}$ and $\gamma \in {_{s'}}\underline{\Wbf}_{s^{\prime\prime}}$ be two composable blocks. Then there is a canonical isomorphism 
$$\Delta^{\beta} * \Delta^{\gamma} \cong \Delta^{\beta\gamma}.$$
Furthermore, there is a well-defined natural transformation of functors $\Xibf^{\opp} \to \Morita$ 
$$\Delta^{\min} : C \Rightarrow \HH^{\Gbf}$$
such that $\Delta^{\min}(\beta) \cong \Delta^{\beta}$. 
\end{lem}

\begin{proof}
Denote by $\delta_{\psi, s} \in \DD((\Ubf^-,\Lcal_{\psi}) \backslash \Gbf /\Ubf, \RR_{\Tbf})_{[-,s]}$ the unique perverse sheaf such that $1^*\delta_{\psi,s}[\dim \Tbf] = \RR_{\Tbf}$. Let $(\Delta^{\beta}, \alpha)$ be the unique minimal rigidified standard sheaf corresponding to $\beta$. The data of the isomorphism $\alpha$ is equivalent to the data of an isomorphism $\Av_{\psi}\Delta^{\beta} = \delta_{\psi,s'}$. 

Clearly the sheaf $\Delta^{\beta} * \Delta^{\gamma}$ is isomorphic to $\Delta^{\beta\gamma}$. To construct this isomorphism canonically, we have to trivialize $\Av_{\psi}(\Delta^{\beta} * \Delta^{\gamma})$. Arguing as in \cite[Lemma 5.12]{LusztigYunV3}, we have a sequence of isomorphisms
\begin{align*}
\Av_{\psi}(\Delta^{\beta} * \Delta^{\gamma}) &\cong \Av_{\psi}(\Delta^{\beta}) * \Delta^{\gamma} \\
&\xrightarrow{\alpha} \delta_{\psi, s'} * \Delta^{\gamma} \\
&\cong \Av_{\psi}(\Delta^{\gamma}) \\
&\xrightarrow{\alpha'} \delta_{\psi,s^{\prime\prime}}.
\end{align*}
This provides us with a trivialization of $\Av_{\psi}(\Delta^{\beta} * \Delta^{\gamma})$ and thus there is a canonical isomorphism of minimal rigidified standard sheaves $\Delta^{\beta} * \Delta^{\gamma} = \Delta^{\beta\gamma}$. 

Finally, $\Delta^{\min}$ is well-defined if the isomorphisms $\Delta^{\beta} * \Delta^{\gamma} = \Delta^{\beta\gamma}$ respect the associativity constraints of the category $\HH$. This is clear by the uniqueness of the isomorphism of rigidified minimal standard sheaves. 
\end{proof}

\begin{rque}\label{remark:rigidified-standard}
Definition \ref{def:min-standard} and Lemmas \ref{lem:min-standard-existence} and \ref{lem:RigidifiedMinStandard} do not require that we only deal with minimal elements in the various blocks (though they play a special role). In particular, the same proof yields for all pairs $(w, s)$ where $w\in \Wbf$ and $s \in \Ch(\Tbf)$ a pair $(\Delta_{w,s}^{\rigid}, \alpha_{w,s})$ where $\Delta_{w,s}^{\rigid}$ is a standard sheaf isomorphic to $\Delta_{w,s}$ and $\alpha$ is an isomorphism $\alpha: 1^*\Av_{\psi}\Delta_{w,s}^{\rigid} \cong \RR_{\Tbf}$. Moreover, if $\underline{w} = (w_1, \dots, w_n)$ is a reduced expression for $w$, set $s_i=w_{i+1}\cdots w_n(s)$, with the empty product equal to the identity. There is then a canonical isomorphism 
$$\Delta_{w_1,s_1}^{\rigid} * \dots * \Delta_{w_n,s_n}^{\rigid} \cong \Delta_{w,s}^{\rigid}.$$
\end{rque}

\begin{rque}\label{remark:rigidified-costandard}
As in Remark \ref{remark:rigidified-standard}, we can define rigidified costandard sheaves $\nabla_{w,s}^{\rigid}$ for all $w \in \Wbf$ and $s \in \Ch(\Tbf)$. We note that as there is a canonical isomorphism $\Av_{\psi}(\nabla_{w,w^{-1}s}^{\rigid} * \Delta_{w^{-1},s}^{\rigid}) \cong \delta_{\psi,s}$, there is a canonical isomorphism $\nabla_{w,w^{-1}s}^{\rigid} * \Delta_{w^{-1},s}^{\rigid} \cong \Delta_{1,s}$ (and similarly if we swap the roles of $\Delta$ and $\nabla$). 
\end{rque} 

\subsubsection{The functor $\V$}

We assume that $\ell$ is not a torsion prime for $\Gbf^*$ (see Definition \ref{def:torsion-prime}).

\begin{defi}[Bott-Samelson tilting sheaves]\label{def:Bott-Samelson-Tilting}
Let $w$ be a simple reflection and $s \in \Ch(\Tbf)$. We denote by $i_w : \Gbf_w \subset \Gbf$ the inclusion of the minimal parabolic containing $w$. 
\begin{enumerate}
\item if $ws = s$, we let $T_{w,s} = i_{w,*}i_w^*\mathbb{T}_s$,
\item if $ws \neq s$, then the object $i_w^*\mathbb{T}_s \in \HH^{\Gbf_w}_{[-,s]}$ splits canonically as the direct sum $\Delta_{1,s} \oplus \Delta_{w,s}$. We set $T_{w,s} = i_{w,*}\Delta_{w,s}$.
\end{enumerate}
For a sequence of simple reflections $\underline{w} = (w_1, \dots, w_n)$ we denote by 
$$T_{\underline{w},s} \colondef= T_{w_1, w_2\dots w_n(s)} * \dots * T_{w_n, s}$$
and we call those objects Bott-Samelson tilting sheaves.
\end{defi} 

\begin{rque}
The objects $T_{\underline{w},s}$ depend only on the pinning of $\Gbf$ and not on extra choices. 
\end{rque}

\begin{rque}\label{remark:canon-short-exact-seq}
Let $w$ be a simple reflection of $\Wbf$, and let $s$ satisfy $ws = s$. Note that $T_{w,s}$ is the maximal tilting sheaf (with right monodromy $s$) in the Hecke category $\HH^{\Gbf_w}_{[s,s]}$ for the minimal Levi subgroup $\Gbf_w$. In particular, $T_{w,s}$ represents, on $\HH^{\Gbf_w}_{[s,s]}$, the functor $\Av_{\psi}$. It follows that there is a unique map 
$$T_{w,s} \to \nabla_{w,s}^{\rigid}$$ 
such that the map $\alpha : \delta_{\psi,s} \cong \Av_{\psi}\nabla_{w,s}^{\rigid}$ is adjoint to $T_{w,s} \cong \Av_{U}\delta_{w,s} \to \nabla_{w,s}^{\rigid}$.
\end{rque}

\begin{lem}
The Bott-Samelson tilting sheaves are tilting sheaves and the category $\HH_{\tilt}$ is generated under direct sums, isomorphisms and direct summands by the Bott-Samelson tilting sheaves. 
\end{lem}

\begin{proof}
Since tilting sheaves are stable under convolution by Theorem \ref{thm:StructureCategoryviaTilting}, it follows that all Bott-Samelson tilting sheaves are tilting. For the rest of the statement it is enough to check that $T_{w, s}$ is a direct summand of a Bott-Samelson tilting sheaf. Let $w= w_1\dots w_n$ be a minimal expression for $w$. Then $T_{\underline{w}, s}$ is supported on the closure of $\BwB$ and contains $\Delta_{w,s}$ with multiplicity one by \cite[Lemma 9.9.1]{Gouttard} hence contains $T_{w,s}$ with multiplicity one. 
\end{proof}

\begin{defi}
The $\V^{\beta}_{[s,s']}$ functor is defined as 
$$\V^{\beta}_{[s,s']} \colondef= \Hom(T^{\beta}_{[s,s']}, -) : \HH^{\beta}_{[s,s']} \to \DD_{\qcoh}(\Ccal^{\beta}_{[s,s']}).$$
Taking a direct sum over all blocks $\beta$ in a single $\Wbf$-orbit $\mathfrak{o}$ yields a functor 
$$\V_{\mathfrak{o}} \colondef= \Hom(T_{\mathfrak{o}}, -) : \HH_{\mathfrak{o}} \to \DD_{\qcoh}(\Ccal_{\mathfrak{o}}).$$
\end{defi}

The coalgebra structure of $\mathbb{T}$ (and thus of $T_{\mathfrak{o}}$) yields a canonical map for all $A,B \in \HH_{\mathfrak{o}}$
$$\V_{\mathfrak{o}}(A) \otimes_{\Lambda} \V_{\mathfrak{o}}(B) \to \V_{\mathfrak{o}}(A * B)$$
obtained as the following composition 
\begin{align*}
\V_{\mathfrak{o}}(A) \otimes_{\Lambda} \V_{\mathfrak{o}}(B)  &\cong \Hom(T_{\mathfrak{o}}, A) \otimes_{\Lambda} \Hom(T_{\mathfrak{o}}, B) \\
&\to \Hom(T_{\mathfrak{o}} * T_{\mathfrak{o}}, A*B) \\
&\to \Hom(T_{\mathfrak{o}}, A*B).
\end{align*}
Here, the second line comes from the bifunctoriality of $-*-$, and the third line comes from the comultiplication of $T_{\mathfrak{o}}$. 

\begin{lem}\label{lem:MonoidalityG}
For all $A \in \HH^{\beta}_{[s,s']}$ and $B \in \HH^{\gamma}_{[s',s^{\prime\prime}]}$ the map 
$$\V_{\mathfrak{o}}(A) \otimes_{\Lambda} \V_{\mathfrak{o}}(B) \to \V_{\mathfrak{o}}(A * B)$$ 
canonically factors through 
$$\V_{\mathfrak{o}}(A) \otimes_{\Lambda} \V_{\mathfrak{o}}(B) \to \V_{\mathfrak{o}}(A) \otimes_{\Ocal(\Ccal_{s'})} \V_{\mathfrak{o}}(B) \to \V_{\mathfrak{o}}(A * B).$$
Furthermore, the map $\V_{\mathfrak{o}}(A) \otimes_{\Ocal(\Ccal_{s'})} \V_{\mathfrak{o}}(B) \to \V_{\mathfrak{o}}(A * B)$ is an isomorphism. In particular the functor $\V_{\mathfrak{o}}$ is equipped with a canonical monoidal structure.  
\end{lem}

\begin{proof}
This \cite[Proposition 4.6.4]{BezrukavnikovYun}, \cite[Proposition 11.5]{BezrukavnikovRicheSoergelTheory} and \cite[Proposition 10.12.1]{Gouttard}. 
\end{proof}

It follows from Theorem \ref{thm:StructureCategoryviaTilting} that if $T \in \HH^{\beta}_{[s,s'], \tilt}$, the object $\V^{\beta}_{[s,s']}(T)$ is a coherent sheaf concentrated in degree $0$. 

\begin{lem}\label{lem:VRank1}
Let $w$ be a simple reflection, and let $\beta \in {_s}\underline{\Wbf}_{s'}$ be the block containing it. Then there is a canonical isomorphism 
$$\V^{\beta}_{[s,s']}(T_{w,s'}) \cong \BS(w,s')$$
as objects of $\Coh(\Ccal_{[s,s']}^{\beta})$. 
\end{lem}

\begin{proof}
This is \cite[Lemma 10.4.9]{Gouttard}. 
\end{proof}

\begin{lem}\label{lem:VBottSamelson}
Let $\underline{w}$ be a sequence of simple reflections. Then there is a canonical isomorphism 
$$\V^{\beta}_{[s,s']}(T_{\underline{w},s'}) \cong \BS(\underline{w}, s').$$ 
\end{lem}

\begin{proof}
This follows from the monoidality of the functor $\V$ and the rank one case \ref{lem:VRank1}. 
\end{proof}

\begin{thm}\label{thm:SoergelRealizationOfABlock}
For all $\beta \in {_s}\underline{\Wbf}_{s'}$ the functor $\V$ induces an equivalence 
$$\V^{\beta}_{[s,s']} : \HH^{\beta}_{[s,s'], \tilt} \to \SBim^{\beta}_{[s,s']}.$$
\end{thm}

\begin{proof}
The essential surjectivity follows from the fact that both categories are generated under direct sums, direct summands and isomorphisms by Bott-Samelson objects by Lemma \ref{lem:VBottSamelson}. The fully faithfulness is \cite[Proposition 9.10.2]{Gouttard}. 
\end{proof}

\begin{thm}\label{thm:SoergelRealizationForG}
The functor $\V$ defines a natural isomorphism of functors 
$$\HH^{\Gbf} \Rightarrow \SBim : \Xibf^{\opp} \to \Morita.$$
\end{thm}

\begin{proof}
The monoidality of the functor $\V$ yields a natural transformation of functors $\HH^{\Gbf} \Rightarrow \SBim$ the fact that it is an isomorphism is Theorem \ref{thm:SoergelRealizationOfABlock}. 
\end{proof}

\subsection{Soergel realization for $\HH^{\Hbf}$}

We assume that $\ell$ is not a torsion prime for $\Gbf^*$ (see Definition \ref{def:torsion-prime}). We fix a $\Wbf$-orbit $\mathfrak{o} \subset \Ch(\Tbf)$. In this section, we establish the analogue of Theorem \ref{thm:SoergelRealizationForG} for the functor $\HH^{\Hbf}$. However, there is an extra difficulty coming from the fact that the $(\Hbf^{\circ}_{[s,s]}, \Hbf^{\circ}_{[s',s']})$-bitorsors $\Hbf^{\beta}_{[s,s']}$ are not canonically trivialized.

\subsubsection{The Whittaker model}

The pinning of $\Hbf^{\circ}_{[s,s]}$ relative to $\Gbf$, together with the trivialization of the root groups of $\Gbf$, determines a character (see Remark \ref{rque:PinningvsRelativePinning})
$$\phi_s : \Ubf^-_s \to \A^1$$
as in Section \ref{sec:WhittakerModelForG}. 
Hence we can consider $\Lcal_{\psi} = \phi_s^*\AS$ which is a rank one character sheaf on $\Ubf_s^-$. We then define, as in the $\HH^{\Gbf}$ case, the following Whittaker model 
$$\Whitt(\HH^{\Hbf, \beta}_{[s,s']}) \colondef= \DD((\Ubf_s^-, \Lcal_{\psi})\backslash \Hbf^{\beta}_{[s,s']}/\Ubf_{s'}, \RR_{\Tbf})_{\unip}.$$ 

Recall from the construction of $\Hbf^{\beta}_{[s,s']}$ that there is a natural inclusion $i_{[s,s']}^{\beta} : \widetilde{\Xibf}^{\beta}_{[s,s']} \to \Hbf^{\beta}_{[s,s']}$. 
\begin{lem}
The pullback $i_{[s,s']}^{\beta,*} : \Whitt(\HH^{\Hbf, \beta}_{[s,s']}) \to \DD(\widetilde{\Xibf}^{\beta}_{[s,s']}, \RR_{\Tbf})_{\unip}$ is an equivalence. 
\end{lem} 

\begin{proof}
This is a special case of Lemma \ref{lem:WhittakerSheavesSupport}.
\end{proof}

\begin{rque}
It is clear that the category $\DD(\widetilde{\Xibf}^{\beta}_{[s,s']}, \RR_{\Tbf})_{\unip}$ is equivalent to $\DD(\RR_{\Tbf})$, with the identification given by pullback at any point. However, it is not yet clear that this identification yields an identification of functors $\Xibf_{\mathfrak{o}} \to \Morita$ between the functor $\beta \mapsto \DD(\RR_{\Tbf})$ and $\beta \mapsto \DD(\widetilde{\Xibf}^{\beta}_{[s,s']}, \RR_{\Tbf})_{\unip}$. We should note that in \cite{LusztigYun}, the authors work with categories of sheaves that are equivariant under $\Tbf$ (and not free monodromic). This setting permits a canonical object in their analog of the category $\DD(\widetilde{\Xibf}^{\beta}_{[s,s']}, \RR_{\Tbf})_{\unip}$, namely, the constant sheaf. 
\end{rque}

\subsubsection{Trivializing minimal standard objects of $\HH^{\Hbf}$}\label{sec:TrivializationWhittakerModel}\label{sec:TrivializationWhittModelOfH}

Let us introduce the following notation:
\begin{enumerate}
\item Let $\DD : \Xibf_{\mathfrak{o}}^{\opp} \to \Morita$ be the functor $\beta \mapsto \DD(\widetilde{\Xibf}^{\beta}_{[s,s']}, \RR_{\Tbf})_{\unip}$, 
\item Let $C : \Xibf_{\mathfrak{o}}^{\opp} \to \Morita$ be the functor $\beta \mapsto \DD(\RR_{\Tbf})$, the constant diagram. 
\end{enumerate}
We aim to produce a canonical isomorphism $C \Rightarrow \DD$. We consider the category $\mathcal{N}$ of such natural isomorphisms $N : C \Rightarrow \DD$ which are subject to the following conditions.
\begin{enumerate}
\item The map $N^{\beta}$ is $\RR_{\Tbf}$-linear.
\item For a block $\beta$, the map $N^{\beta} : \DD(\RR_{\Tbf}) = C^{\beta} \to \DD^{\beta}$ is a $t$-exact equivalence.
\item Let $I$ be the augmentation ideal of $\RR_{\Tbf}$. The category $\DD^{\beta} \otimes_{\RR_{\Tbf}} \RR_{\Tbf}/I$ is canonically equivalent to $\DD(\widetilde{\Xibf}^{\beta}_{[s,s']}/T, \Lambda) = \DD(\Lambda)$, where the equivalence is provided by the constant sheaf on $\widetilde{\Xibf}^{\beta}_{[s,s']}$. The map $N^{\beta} : \DD(\RR_{\Tbf}) \cong C^{\beta}$ induces the identity on $\DD(\RR_{\Tbf})\otimes_{\RR_{\Tbf}} \RR_{\Tbf}/I = \DD(\Lambda) = C^{\beta} \otimes_{\RR_{\Tbf}} \RR_{\Tbf}/I$. 
\end{enumerate}
Let us unfold the definition of the category $\mathcal{N}$. By linearity, the data of $N^{\beta}$ is equivalent to the data of an object $y_{\beta} \in \DD^{\beta}$ for all blocks $\beta$. The second condition states that $y_{\beta} \in \DD^{\beta, \heartsuit}$, and the last one states that there is a canonical isomorphism $y_{\beta}\otimes_{\RR_{\Tbf}} \RR_{\Tbf}/I = \Lambda_{\widetilde{\Xibf}^{\beta}_{[s,s']}}$. Furthermore, the collection of objects $(y_{\beta})$ is equipped with a collection of coherent isomorphisms 
\begin{equation}\label{eq:NaturalityIntoMorita}
y_{\beta} * y_{\beta'} = y_{\beta\beta'}
\end{equation}
for composable blocks $\beta$ and $\beta'$. 

\begin{lem}\label{lem:ExistenceTrivialization}
The category $\mathcal{N}$ is not empty.
\end{lem}

\begin{proof}
For each block $\beta$, choose a geometric point of the $\Tbf_s$-$\Tbf_{s'}$-torsor $\widetilde{\Xibf}^{\beta}_{[s,s']}$. It identifies $\DD^{\beta}$ with $\DD(\RR_{\Tbf})$; let $y_{\beta}$ be the object corresponding to the free rank-one module $\RR_{\Tbf}$. Then $y_{\beta}$ lies in the heart and its reduction modulo $I$ is the constant sheaf. Thus a collection satisfying $(ii)$ and $(iii)$ exists. Let $\beta, \beta'$ be two composable blocks, the two objects $y_{\beta} * y_{\beta'}$ and $y_{\beta\beta'}$ both yield an equivalence of $\DD^{\beta\beta', \heartsuit} = \DD(\RR_{\Tbf})^{\heartsuit}$. They differ by a free $\RR_{\Tbf}$-module of rank one $\Lcal_{\beta,\beta'}$, that is, there is a canonical isomorphism 
$$y_{\beta} * y_{\beta'} = \Lcal_{\beta, \beta'} \otimes_{\RR_{\Tbf}} y_{\beta\beta'}.$$
Furthermore, the $\RR_{\Tbf}/I = \Lambda$-module $\Lcal_{\beta, \beta'}/I$ is canonically trivialized as both $(y_{\beta} * y_{\beta'})/I$ and $y_{\beta\beta'}/I$ are trivialized.
The assignment $(\beta, \beta') \mapsto \Lcal_{\beta, \beta'}$ defines a $2$-cocycle on $\Xibf_{\mathfrak{o}}$ valued in the groupoid $\Pic(\RR_{\Tbf})$; that is, for all triples of composable blocks $(\beta_1, \beta_2, \beta_3)$, there is a canonical isomorphism 
\begin{equation}\label{eq:twoCocycle}
\Lcal_{\beta_2, \beta_3} \otimes_{\RR_{\Tbf}} \Lcal_{\beta_1, \beta_2\beta_3}^{-1} \otimes_{\RR_{\Tbf}} \Lcal_{\beta_1\beta_2, \beta_3} \otimes_{\RR_{\Tbf}} \Lcal_{\beta_1,\beta_2} \cong \RR_{\Tbf}.
\end{equation}
For all composable pairs $(\beta_1, \beta_2)$, we choose a generator $f_{\beta_1, \beta_2} \in \Lcal_{\beta_1, \beta_2}$ as an $\RR_{\Tbf}$-module which reduces to $1 \in \Lcal_{\beta_1, \beta_2}/I$. The choice of this element yields an isomorphism $\Lcal_{\beta_1, \beta_2} \cong \RR_{\Tbf}$. For all triples of composable morphisms $(\beta_1, \beta_2, \beta_3)$, we define
$$h_{\beta_1,\beta_2,\beta_3} = f_{\beta_2, \beta_3}f_{\beta_1, \beta_2\beta_3}^{-1}f_{\beta_1\beta_2, \beta_3}f_{\beta_1, \beta_2}^{-1} \in \RR_{\Tbf}^{\times}.$$
It is an immediate verification that the assignment $(\beta_1,\beta_2,\beta_3) \mapsto h_{\beta_1,\beta_2,\beta_3}$ defines a $3$-cocycle on $\Xibf_{\mathfrak{o}}$ valued in $\RR_{\Tbf}^{\times}$. Furthermore, since all $f_{\beta_i,\beta_j}$ reduce to $1$ mod $I$, the cocycle $h$ is valued in $1 + I\RR_{\Tbf} \subset \RR_{\Tbf}^{\times}$. 
By Lemma \ref{lem:TrivializationCohomology}, this cocycle is cohomologically trivial. Hence, there exists a $2$-cochain $g$ valued in $(1 + I\RR_{\Tbf})$ whose differential is $h$. Consider the elements $f'_{\beta_1,\beta_2} = g_{\beta_1,\beta_2}f_{\beta_1,\beta_2} \in \Lcal_{\beta_1,\beta_2}$. This family provides a trivialization of the $2$-cocycle $\Lcal_{\beta_1,\beta_2}$, and we can therefore make coherent choices of isomorphisms 
$$y_{\beta_1} * y_{\beta_2} = y_{\beta_1\beta_2}.$$
Hence there exists an isomorphism $C \Rightarrow \DD$. 
\end{proof}

\begin{lem}\label{lem:Connectedness}
The category $\mathcal{N}$ is a $1$-truncated connected groupoid.
\end{lem}

\begin{proof}
The data of an object in $\mathcal{N}$ is determined by the data of the collection of objects $y_{\beta} \in \DD^{\beta, \heartsuit}$. Since these all lie in an abelian category, there are no nontrivial higher morphisms; this proves the $1$-truncated statement. It is also clear that this category is a groupoid. Let us prove that it is connected. Let $(y_{\beta})$ and $(z_{\beta})$ be two natural isomorphisms $C \Rightarrow \DD$. They differ by an $\RR_{\Tbf}$-line bundle $\Lcal_{\beta}$, i.e., $z_{\beta} = \Lcal_{\beta} \otimes_{\RR_{\Tbf}} y_{\beta}$. Arguing as in the proof of Lemma \ref{lem:ExistenceTrivialization}, the map $\beta \mapsto \Lcal_{\beta}$ defines a $1$-cocycle on $\Xibf_{\mathfrak{o}}$ valued in $\Pic(\RR_{\Tbf})$. As before, this cocycle is controlled by a $2$-cocycle valued in $(1 + I\RR_{\Tbf})$ which is cohomologically trivial by Lemma \ref{lem:TrivializationCohomology}. Therefore, the $1$-cocycle $\Lcal_{\beta}$ can be trivialized, and the two objects $(y_{\beta})$ and $(z_{\beta})$ are isomorphic.
\end{proof}

\begin{lem}\label{lem:NoNontrivialIso}
There are no nontrivial morphisms in $\mathcal{N}$.
\end{lem}

\begin{proof}
Since $\mathcal{N}$ is a connected groupoid, it is enough to show that there are no nontrivial automorphisms.
Let $(y_{\beta})$ be an object in $\mathcal{N}$. Let $(\alpha_{\beta}) : (y_{\beta}) \to (y_{\beta})$ be an automorphism of $(y_{\beta})$. Since $y_{\beta}$ corresponds to an $\RR_{\Tbf}$-module of rank one $\End(y_{\beta}) = \RR_{\Tbf}$ and so $\alpha_{\beta} \in 1 + I\RR_{\Tbf}$ as before. The relation $y_{\beta}*y_{\beta'} = y_{\beta\beta'}$ shows that $\beta \mapsto \alpha_{\beta}$ is a morphism of groupoids which is trivial by Lemma \ref{lem:TrivializationCohomology}. Hence $\alpha_{\beta} = \id$ for all $\beta$. 
\end{proof}

\begin{lem}\label{lem:TrivializationCohomology}
For all $i > 0$, the cohomology group $H^i(\Xibf_{\mathfrak{o}}, 1 + I\RR_{\Tbf})$ vanishes. 
\end{lem}

\begin{proof}
Since the groupoid $\Xibf_{\mathfrak{o}}$ is equivalent to the classifying stack $\pt/\Gammabf_s$ for any choice of $s \in \mathfrak{o}$, we get that 
$$H^i(\Xibf_{\mathfrak{o}}, 1 + I\RR_{\Tbf}) = H^i(\Gammabf_s, 1 + I\RR_{\Tbf})$$
where the right-hand side is group cohomology. We equip the group $(1 + I\RR_{\Tbf})$ with the decreasing filtration $(1 + I^n\RR_{\Tbf})_n$. It is clear that $(1 + I\RR_{\Tbf})$ is complete with respect to the topology defined by this filtration. Furthermore, the graded pieces are isomorphic to $\Lambda$-modules of the form $\Lambda^N$ for some nonnegative integer $N$. 
If $\Lambda = \Qlb$, then $H^i(\Gammabf_s, 1 + I\RR_{\Tbf}) = 0$ for all primes $\ell$. If $\Lambda = \Zlb, \Flb$, then our assumption on $\ell$ implies that $\ell$ does not divide $|\Gammabf_s|$ by Lemma \ref{lem:ControlCentralizer} and thus $H^i(\Gammabf_s, 1 + I\RR_{\Tbf}) = 0$. 
\end{proof}

The combination of Lemmas \ref{lem:ExistenceTrivialization}, \ref{lem:NoNontrivialIso}, and \ref{lem:Connectedness} shows that there exists a canonical trivialization 
\begin{equation}\label{eq:Delta-on-the-endoscopic-side}
\Delta : C \Rightarrow \DD.
\end{equation}

\subsubsection{The functor $\V$}\label{sec:functor-V-endoscopic-side}

Now that we have constructed a canonical trivialization, we can extend the theory of the $\V$-functor to $\HH^{\Hbf}$. We consider the equivalence (\ref{eq:Delta-on-the-endoscopic-side})
$$\Delta : C \cong \DD$$
of Section \ref{sec:TrivializationWhittakerModel}. For all blocks $\beta \in {_s}\underline{\Wbf}_{s'}$, we then have an object $R^{\beta} \in \DD(\widetilde{\Xibf}^{\beta}_{[s,s']}, \RR_{\Tbf})_{\unip}$ such that under the equivalence $\DD^{\beta} \cong C^{\beta}$ the object $R^{\beta}$ corresponds to $\RR_{\Tbf}$. The inclusion $i_{[s,s']}^{\beta} :\widetilde{\Xibf}_{[s,s']}^{\beta} \to \Hbf_{[s,s']}^{\beta}$ yields a morphism of functors $\Xibf_{\mathfrak{o}} \to \Morita$ 
$$i_{[s,s'],*}^{\beta} : \DD(\widetilde{\Xibf}_{[s,s']}^{\beta}, \RR_{\Tbf})_{\unip} \to \HH^{\Hbf, \beta}_{[s,s']}.$$
We set $\Delta^{\Hbf,\beta} =  i_{[s,s'],*}^{\beta}R^{\beta}$. The isomorphism $C \cong \DD$ shows that 
$$\Delta^{\Hbf,\beta} * \Delta^{\Hbf,\gamma} \cong \Delta^{\Hbf, \beta \gamma}$$
for two composable blocks $\beta$ and $\gamma$. 

We then define the big tilting sheaf $\mathbb{T}^{\Hbf, \beta} \in \HH^{\Hbf, \beta}_{[s,s']}$ to be 
$$\mathbb{T}^{\Hbf, \beta} \cong \Av_{\Ubf} \Av_{\psi} \Delta^{\Hbf,\beta}.$$
Note that  $\Av_{\psi} \Delta^{\Hbf,\beta} \in \Whitt(\HH^{\beta}_{[s,s']})$ corresponds to $R^{\beta}$ under the equivalence $\Whitt(\HH^{\beta}_{[s,s']}) = \DD(\widetilde{\Xibf}_{[s,s']}^{\beta}, \RR_{\Tbf})_{\unip}.$

\begin{thm}[Endomorphismensatz] 
There is a natural isomorphism 
$$\End(\mathbb{T}^{\Hbf, \beta}) \cong \Ocal(\Ccal^{\beta}).$$
\end{thm}

\begin{proof}
The map 
$$\Ocal(\Ccal^{\beta}) \to \End(\mathbb{T}^{\Hbf, \beta})$$ 
is simply given by the monodromy. We can now choose a point of the $(\Hbf^{\circ}_{[s,s]},\Hbf^{\circ}_{[s',s']})$-bitorsor $\Hbf^{\beta}_{[s,s']}$. Left and right translation identify it with the corresponding group torsors and transport both torus and monodromy actions. Under these identifications, the statement that this map is an isomorphism is an instance of Theorem \ref{thm:Endomorphismensatz}. 
\end{proof}

\begin{thm}
For all tilting sheaves $T \in \HH^{\Hbf, \beta}_{[s,s']}$, the object
$$\Hom(\mathbb{T}^{\Hbf, \beta}, T) \in \DD_{\qcoh}(\Ccal^{\beta}_{[s,s']})$$ 
is concentrated in degree $0$ and is a Soergel bimodule. 
\end{thm}

\begin{proof}
After trivializing $\Hbf^{\beta}_{[s,s']}$, this is an instance of Theorem \ref{thm:SoergelRealizationOfABlock}. 
\end{proof}

\begin{defi}
We define the functor 
$$\V^{\Hbf, \beta} : \HH^{\Hbf, \beta}_{[s,s']} \to \KK(\SBim^{\beta}_{[s,s']})$$
as the functor obtained by passing to homotopy categories the functors
$$\V^{\Hbf, \beta} = \Hom(\mathbb{T}^{\Hbf, \beta}, -).$$
\end{defi} 

Consider the category 
$$\HH^{\Hbf}_{\mathfrak{o}, \tot} = \oplus_{\beta} \HH^{\Hbf,\beta}_{[s,s']}$$
where the sum ranges over all blocks $\beta \in {_s}\underline{\Wbf}_{s'}$ where $s,s' \in \mathfrak{o}$. We equip this category with a convolution structure. Let $A  \in \HH^{\Hbf,\beta}_{[s_1,s_2]}$ and $B \in \HH^{\Hbf,\gamma}_{[s_3,s_4]}$. Then we set 
\begin{enumerate}
\item $A * B = 0$ if $s_2 \neq s_3$,
\item $A * B \in \HH^{\Hbf, \beta\gamma}_{[s_1, s_4]}$ is the convolution defined previously. 
\end{enumerate}
Similarly we consider the category 
$$\Whitt(\HH^{\Hbf}_{\mathfrak{o}, \tot}) := \oplus_{\beta} \Whitt(\HH^{\Hbf,\beta}_{[s,s']}).$$
As before, there is a pair of adjoint functors 
$$\Av_{\Ubf} : \Whitt(\HH^{\Hbf}_{\mathfrak{o}, \tot}) \leftrightarrows \HH^{\Hbf}_{\mathfrak{o}, \tot} : \Av_{\psi}.$$
Set $\Delta^{\Hbf} = \oplus_{\beta} \Delta^{\Hbf, \beta}$. Then we have 
$$\mathbb{T}^{\Hbf} = \oplus_{\beta} \mathbb{T}^{\Hbf, \beta} = \Av_{\Ubf} \Av_{\psi} \Delta^{\Hbf}.$$
Arguing as in the case of $\HH^{\Gbf}$, the object $\mathbb{T}^{\Hbf}$ has a canonical coalgebra structure in $\HH^{\Hbf}_{\mathfrak{o}, \tot}$ and thus the functor $$\V^{\Hbf} : \HH^{\Hbf}_{\mathfrak{o}, \tot} \to \DD_{\qcoh}(\Ccal_{\mathfrak{o}})$$
is equipped with a canonical lax-monoidal structure. 

\begin{lem}\label{lem:MonoidalityHSide}
The functor 
$$\V^{\Hbf} : \HH^{\Hbf}_{\mathfrak{o}, \tot} \to \DD_{\qcoh}(\Ccal_{\mathfrak{o}})$$
is strong monoidal. 
\end{lem}

\begin{proof}
This is a particular case of Lemma \ref{lem:MonoidalityG}. 
\end{proof}

\begin{thm}\label{thm:SoergelRealizationHSide}
The functor $\V^{\Hbf, \beta} : \HH^{\Hbf, \beta}_{[s,s']} \to \KK(\SBim^{\beta}_{[s,s']})$ is an equivalence. Moreover, the collection of all $\V^{\Hbf, \beta}$ yields an isomorphism of functors 
\begin{align*}
\Xibf_{\mathfrak{o}}^{\opp} &\to \Morita \\
\HH^{\Hbf} &\cong \SBim.
\end{align*}
\end{thm}

\begin{proof}
The monoidality of the functor $\V^{\Hbf}$ of Lemma \ref{lem:MonoidalityHSide} yields a natural transformation 
$$\HH^{\Hbf} \Rightarrow \SBim.$$
The fact that the functor $\V^{\Hbf, \beta}$ is an equivalence can be checked after trivializing $\Hbf^{\beta}_{[s,s']}$ and using Theorem \ref{thm:SoergelRealizationOfABlock}.
\end{proof}

\subsection{Endoscopy}

\begin{thm}\label{thm:EndoscopyHeckeCat}
Let $\mathfrak{o} \subset \Ch(\Tbf)$ be a $\Wbf$-orbit and assume that $\ell$ is not a torsion prime for $\Gbf^*$ (see Definition \ref{def:torsion-prime}). Then there is a canonical isomorphism, depending only on the pinning of $\Gbf$, of functors $\Xibf_{\mathfrak{o}}^{\opp} \to \Morita$
\begin{equation}\label{eq:EndoscopyForHeckeCat}
\EE : \HH^{\Gbf} \cong \HH^{\Hbf}.
\end{equation}
For all blocks $\beta$, the equivalence $\EE^{\beta} : \HH^{\Gbf,\beta}_{[s,s']} \cong \HH^{\Hbf,\beta}_{[s,s']}$ is $t$-exact and preserves standard, costandard, and tilting sheaves (we shall also write $\EE^{\circ}_s, \EE_{[s,s]}, \EE^{\beta}_{[s,s']}$, \dots for the corresponding endoscopic equivalences between categories with the same decorations).
\end{thm}

\begin{proof}
Theorems \ref{thm:SoergelRealizationForG} and \ref{thm:SoergelRealizationHSide} show that both functors are isomorphic to $\SBim$. Moreover this equivalence preserve tilting sheaves by construction. By monoidality since the costandard sheaves are the duals of the standard sheaves, it is enough to check the statement for the standard sheaves. For the standard sheaves, again by monoidality we can reduce to checking that $\Delta_{w,s}^{\Gbf}$ is mapped onto $\Delta_{w,s}^{\Hbf}$ under this equivalence for $w$ a simple reflection. If $ws \neq s$ then $\Delta_{w,s}^{\Gbf}$ and $\Delta_{w,s}^{\Hbf}$ are both tilting sheaves and the statement is clear by Lemma \ref{lem:ForgetSupport}. If $ws = s$, then there are canonical short exact sequences of perverse sheave, see \cite[Appendix C]{BezrukavnikovYun}
$$0 \to \Delta_{w,s}^{\Gbf} \to T_{w,s}^{\Gbf} \to \Delta_{1,s}^{\Gbf} \to 0$$
and 
$$0 \to \Delta_{w,s}^{\Hbf} \to T_{w,s}^{\Hbf} \to \Delta_{1,s}^{\Hbf} \to 0.$$
As both $T_{w,s}^{\Gbf}$ and $\Delta_{1,s}^{\Gbf}$ (resp. $T_{w,s}^{\Hbf}$ and $\Delta_{1,s}^{\Hbf}$) are tilting, we deduce that the image of $\Delta_{w,s}^{\Gbf}$ under the equivalence is $\Delta_{w,s}^{\Hbf}$. 
\end{proof}

\begin{rque}
Since both functors are valued in $\Morita$, this also implies that there are monoidal equivalences 
$$\HH^{\Gbf, \circ}_{[s,s]} \cong \HH^{\Hbf, \circ}_{[s,s]}$$
for all elements $s \in \mathfrak{o}$. 
\end{rque}

\begin{rque}
Under the equivalence of Theorem \ref{thm:EndoscopyHeckeCat}, the rigidified minimal standard sheaf $\Delta^{\beta}$ corresponds to the sheaf $\Delta^{\Hbf, \beta}$. 
\end{rque} 

\begin{rque}
Theorem \ref{thm:EndoscopyHeckeCat} establishes the endoscopy for free monodromic Hecke categories and generalizes the results of \cite{Gouttard}. 
\end{rque}

\subsection{Properties of the endoscopic equivalence}

In this section, we assume that $\ell$ is not a torsion prime for $\Gbf^*$ (see Definition \ref{def:torsion-prime}).

\subsubsection{Inversion of $\ell$}

In this section, we assume that $\Gbf$ has connected center for simplicity.
Let $s \in \Ch_{\Zlb}(\Tbf) \subset \Tbf^*$. It is an element of $\Tbf^*$ of order prime to $\ell$. Let $s_{\ell} \in \Tbf^*$ be an element of order a power $\ell$. Then we have 
$$\Zbf_{\Gbf^*}(ss_{\ell}) = \Zbf_{\Zbf_{\Gbf^*}(s)}(s_{\ell}).$$
Let $\Hbf_s$ and $\Hbf_{ss_{\ell}}$ be the endoscopic groups corresponding to $s$ and $ss_{\ell}$ respectively. 

In \cite[Section 2.4]{EteveFreeMonodromic}, we introduced a functor 
$$\HH^{\Gbf}_{[s,s], \Zlb} \xrightarrow{\Loc_{s_{\ell},s_{\ell}}} \HH^{\Gbf}_{[ss_{\ell}, ss_{\ell}], \Qlb}$$ 
which is essentially obtained by inverting $\ell$ and completing at a certain ideal of $\RR_{\Tbf, \Zlb}[\frac{1}{\ell}]$. It follows from \emph{loc. cit.} that this functor sends free monodromic standard (resp. costandard, resp. tilting) sheaves to standard (resp. costandard, resp. tilting) sheaves and is equipped with a canonical monoidal structure. 

Let us denote by 
$$\Hbf_{[s,s], [s_{\ell}, s_{\ell}], \Qlb}^{\Hbf_s} = \DD(\Ubf_s\backslash \Hbf_s /\Ubf_s, \RR_{\Tbf} \otimes \RR_{\Tbf})_{[s_{\ell}, s_{\ell}]},$$
the corresponding Hecke category over $\Qlb$. 

\begin{lem}\label{lem:CompatInvertEll}
There is a commutative diagram 
\[\begin{tikzcd}
	{\HH^{\Gbf}_{[s,s], \Zlb}} & {\HH^{\Gbf}_{[ss_{\ell}, ss_{\ell}], \Qlb}} \\
	{\HH^{\Hbf_s}_{[s,s], \Zlb}} & {\HH^{\Hbf_s}_{[s,s],[s_{\ell}, s_{\ell}], \Qlb}} & {\HH^{\Hbf_{ss_{\ell}}}_{[ss_{\ell}, ss_{\ell}], \Qlb}}
	\arrow["{{\Loc_{s_{\ell}, s_{\ell}}}}", from=1-1, to=1-2]
	\arrow["{\EE_{\Zlb,s}}"', from=1-1, to=2-1]
	\arrow["{\EE_{\Qlb,ss_{\ell}}}", from=1-2, to=2-3]
	\arrow["{{\Loc_{s_{\ell}, s_{\ell}}}}"', from=2-1, to=2-2]
	\arrow[from=2-2, to=2-3]
\end{tikzcd}\]
where the left vertical map and the two rightmost maps are the endoscopic equivalences for $(\Gbf,s)$ over $\Zlb$, $(\Gbf,ss_{\ell})$ over $\Qlb$ and $(\Hbf_s, s_{\ell})$ over $\Qlb$ respectively.
\end{lem}

\begin{proof}
First note that the functor $\Loc_{[s_{\ell}, s_{\ell}]}$ preserves the maximal tilting object as it commutes with $\Av_{\Ubf}$ and $\Av_{\psi}$. We therefore have a commutative diagram 
\[\begin{tikzcd}
	{\HH^{\Gbf}_{[s,s], \Zlb, \tilt}} && {\HH^{\Gbf}_{[ss_{\ell},ss_{\ell}], \Qlb, \tilt}} & \\
	{\SBim_{[s,s],\Zlb}} && {\SBim_{[ss_{\ell},ss_{\ell}], \Qlb}} & {\HH^{\Hbf_{ss_{\ell}}}_{[ss_{\ell}, ss_{\ell}], \Qlb, \tilt}} \\
	{\HH^{\Hbf_s}_{[s,s], \Zlb, \tilt}} && {\HH^{\Hbf_s}_{[s,s],[s_{\ell},s_{\ell}],\Qlb,\tilt}}
	\arrow["{{\Loc_{[s_{\ell}, s_{\ell}]}}}", from=1-1, to=1-3]
	\arrow["{\V^{\Gbf}_{[s,s]}}"', from=1-1, to=2-1]
	\arrow["{\V^{\Gbf}_{[ss_{\ell},ss_{\ell}]}}", from=1-3, to=2-3]
	\arrow["{{\otimes_{\Ocal(\Ccal_{[s,s]}^{\Zlb})}\Ocal(\Ccal_{[ss_{\ell},ss_{\ell}]}^{\Qlb})}}", from=2-1, to=2-3]
	\arrow["{\V^{\Hbf_{ss_{\ell}}}_{[ss_{\ell},ss_{\ell}]}}", from=2-4, to=2-3]
	\arrow["{\V^{\Hbf}_{[s,s]}}", from=3-1, to=2-1]
	\arrow["{{\Loc_{[s_{\ell}, s_{\ell}]}}}"', from=3-1, to=3-3]
	\arrow["{\V^{\Hbf_s}_{[s,s],[s_{\ell},s_{\ell}]}}"', from=3-3, to=2-3]
\end{tikzcd}\]
where the equality signs are given by the functor $\V$. This proves the lemma.
\end{proof}

\subsubsection{Compatibility with parabolic induction}\label{sec:CompactWithParabInduction}

Let $\Lbf \subset \Pbf \subset \Gbf$ be a standard Levi with standard parabolic $\Pbf$. Let us denote by $i : \Pbf \subset \Gbf$ the inclusion and by $r : \Pbf \to \Lbf$ the projection. The parabolic induction functor 
$$i_*r^* : \HH^{\Lbf}_{\tot} \to \HH^{\Gbf}_{\tot}$$
is fully faithful. 
Dually there is a standard parabolic and Levi $\Lbf^* \subset \Pbf^* \subset \Gbf^*$. 
Let $s \in \Ch_{\Lambda}(\Tbf) \subset \Tbf^*$ and consider the inclusion 
$$\Zbf_{\Lbf^*}(s) \subset \Zbf_{\Pbf^*}(s) \subset \Zbf_{\Gbf^*}(s)$$
of centralizers. The group $\Zbf_{\Lbf^*}^{\circ}(s)$ is a Levi subgroup of $\Zbf_{\Gbf^*}^{\circ}(s)$ and we have an induced morphism of component groups $\Gammabf_s^{\Lbf} \to \Gammabf_s^{\Gbf}$, see Lemma \ref{lem:ParabEndoscopicGroups}.
Let $\mathfrak{o}^{\Lbf} \subset \mathfrak{o}^{\Gbf}$ be the $\Wbf^{\Lbf}$ and $\Wbf^{\Gbf}$ orbits containing $s$. There is a canonical morphism of groupoids 
$$\Xibf_{\mathfrak{o}^{\Lbf}}^{\Lbf} \to \Xibf^{\Gbf}_{\mathfrak{o}^{\Gbf}}.$$
This induces a morphism of groupoids
$$\widetilde{\Xibf}_{\mathfrak{o}^{\Lbf}}^{\Lbf} \to \widetilde{\Xibf}^{\Gbf}_{\mathfrak{o}^{\Gbf}}.$$
Furthermore, upon passing to endoscopic groups, we obtain a Levi subgroup $\Hbf_s^{\circ, \Lbf}$ of $\Hbf_s^{\circ, \Gbf}$ which is pinned relatively to $\Lbf$ using the pinning of $\Hbf_s^{\circ,\Gbf}$ relative to $\Gbf$. This morphism induces a correspondence of groupoids
$$\Hbffrak_{\mathfrak{o}^{\Lbf}}^{\Lbf} \xleftarrow{r^{\Hbf}} \Hbffrak_{\mathfrak{o}^{\Pbf}}^{\Pbf} \xrightarrow{i^{\Hbf}} \Hbffrak_{\mathfrak{o}^{\Gbf}}^{\Gbf}$$
and therefore a fully faithful functor 
$$i_*^{\Hbf}r^{*,\Hbf} : \HH^{\Hbf^{\Lbf}}_{\mathfrak{o}, \tot} \to \HH^{\Hbf^{\Gbf}}_{\mathfrak{o}, \tot}.$$

\begin{lem}\label{lem:CompactEndoHeckeCatParabInduction}
There is a commutative diagram 
\[\begin{tikzcd}
	{\HH^{\Lbf}_{\mathfrak{o}^{\Lbf},\tot}} & {\HH^{\Gbf}_{\mathfrak{o}^{\Gbf}, \tot}} \\
	{\HH^{\Hbf^{\Lbf}}_{\mathfrak{o}^{\Lbf}, \tot}} & {\HH^{\Hbf^{\Gbf}}_{\mathfrak{o}^{\Gbf}, \tot}}
	\arrow["{{i_*r^*}}", from=1-1, to=1-2]
	\arrow["{\EE_{\ofrak^{\Lbf}}^{\Lbf}}"', from=1-1, to=2-1]
	\arrow["{\EE_{\ofrak}^{\Gbf}}", from=1-2, to=2-2]
	\arrow["{{i_*^{\Hbf}r^{*,\Hbf}}}"', from=2-1, to=2-2]
\end{tikzcd}\]
where the vertical maps are the endoscopic equivalences for $\Lbf$ and $\Gbf$ of Theorem \ref{thm:EndoscopyHeckeCat} respectively. 
\end{lem}
 
\begin{proof}
Since the morphisms $i_*$ are fully faithful, it is enough to check that the endoscopic equivalence maps $\HH^{\Lbf}_{\mathfrak{o}^{\Lbf}, \tot}$ onto $\HH^{\Hbf^{\Lbf}}_{\mathfrak{o}^{\Lbf},\tot}$. But $\HH^{\Lbf}_{\mathfrak{o}^{\Lbf},\tot}$ is the full subcategory of $\HH^{\Gbf}_{\mathfrak{o}^{\Gbf},\tot}$ spanned by objects $\Delta_{w,s}$ where $s \in \mathfrak{o}^{\Lbf}$ and $w \in \Wbf^{\Lbf}$, and similarly for $\HH^{\Hbf_{\Lbf}}_{\mathfrak{o}_{\Lbf},\tot}$. The statement now follows from the fact that the endoscopic equivalence preserves standard objects. 
\end{proof}

\subsubsection{A special case: $\Hbf = \Gbf$}

Let $s \in \Ch(\Tbf)$ and assume that $s \in \Zbf(\Gbf^*)$. It immediately implies that 
\begin{enumerate}
\item the $\Wbf$-orbit containing $s$ is $\{s\}$, 
\item the endoscopic group $\Hbf_{[s,s]} = \Gbf$ (compatibly with the relative pinning).
\end{enumerate}

\begin{lem}[\protect{\cite[Lemma 2.3]{LusztigYun}}]
The character sheaf $\Lcal_s$ on $\Tbf$ extends uniquely to a rank-one character sheaf $\Lcal_s^{\Gbf}$ on $\Gbf$. 
\end{lem}

\begin{rque}
It follows from \cite{LusztigYun} that the sheaf $\Lcal_s^{\Gbf}$ is canonically trivialized when restricted to $\Ubf$ and $\Ubf^-$. In particular, as a sheaf, it descends to a (uniquely defined) sheaf on $\Ubf \backslash \Gbf/\Ubf$.
\end{rque} 

There is then a well-defined monoidal functor
$$\EE^{\mathrm{naive}} : \HH^{\Gbf}_{[s,s]} \xrightarrow{\Lcal_s^{\Gbf, -1}} \HH^{\Gbf}_{[1,1]} = \HH^{\Hbf}_{[s,s]}.$$

\begin{lem}\label{lem:endoscopy-trivial-case}
There is a canonical isomorphism of monoidal functors $\HH^{\Gbf}_{[s,s]} \to \HH^{\Hbf}_{[s,s]}$,
$$\EE^{\mathrm{naive}} \cong \EE_{[s,s]}.$$
\end{lem} 

\begin{proof}
It is enough to prove that there is an isomorphism of coalgebras 
$$\TT_s \cong \TT_1 \otimes_{\Lambda} \Lcal_s^{\Gbf}.$$
As both objects $\TT_s$ and $\TT_1$ corepresent $\Av_{\psi}$, it is enough to construct a canonical isomorphism for all $A \in \HH_{[1,1]}^{\Gbf}$
$$\Av_{\psi}(A \otimes_{\Lambda} \Lcal_s^{\Gbf}) = \Av_{\psi}(A) \otimes_{\Lambda} \Lcal_s^{\Gbf}.$$

Let us introduce the following diagram 
\[\begin{tikzcd}
	{\Ubf^- \times \Gbf/\Ubf} & {\Gbf \times \Gbf/\Ubf} \\
	{\Gbf/\Ubf}
	\arrow["{i \times \id}", from=1-1, to=1-2]
	\arrow["a"', from=1-1, to=2-1]
	\arrow["m", from=1-2, to=2-1]
\end{tikzcd}\]\\
where $i : \Ubf^- \to \Gbf$ is the inclusion, $a$ is the action map, and $m$ is the multiplication in $\Gbf$.

We then have $\Av_{\psi}(A) = a_!(\Lcal_{\psi} \boxtimes A)[\dim \Ubf]$. Finally we can directly compute 
\begin{align*}
\Av_{\psi}(A) \otimes_{\Lambda} \Lcal_s^{\Gbf} &= a_!(\Lcal_{\psi} \boxtimes A)[\dim \Ubf] \otimes_{\Lambda} \Lcal_s^{\Gbf} \\
&= m_!((i_!\Lcal_{\psi} \boxtimes A) \otimes_{\Lambda} m^*\Lcal_s^{\Gbf} )[\dim \Ubf] \\
&= m_!((i_!\Lcal_{\psi} \boxtimes A) \otimes_{\Lambda} (\Lcal_s^{\Gbf} \boxtimes \Lcal_s^{\Gbf} ))[\dim \Ubf] \\
&= m_!((i_!\Lcal_{\psi} \otimes \Lcal_s^{\Gbf}) \boxtimes (A \otimes_{\Lambda} \Lcal_s^{\Gbf}))[\dim \Ubf] \\
&= m_!((i_!\Lcal_{\psi}) \boxtimes (A \otimes_{\Lambda} \Lcal_s^{\Gbf}))[\dim \Ubf] \\
&= \Av_{\psi}(A \otimes_{\Lambda} \Lcal_s^{\Gbf}).
\end{align*}
The maps are as follows 
\begin{enumerate}
\item the first and last ones are by definition, 
\item the second is the projection formula, 
\item the third comes from the isomorphism $m^*\Lcal_s^{\Gbf} = \Lcal_{s}^{\Gbf} \boxtimes \Lcal_s^{\Gbf}$ coming from its definition as a character sheaf on $\Gbf$, 
\item the fourth one is again the projection formula, 
\item the fifth one comes from the fact that $i^*\Lcal_s^{\Gbf}$ is trivial.
\end{enumerate}
\end{proof}
\subsubsection{Independence on the extended pinning}\label{sec:indep-extended-pinning}

The equivalence $\EE$ constructed in Theorem \ref{thm:EndoscopyHeckeCat} depends on the datum of an extended pinning of $\Gbf$; see Section \ref{sec:ExtendedPinning}. We now remove this dependence and prove that we only need a pinning of $\Gbf$. Let us choose a $\Wbf$-orbit $\ofrak \subset \Ch(\Tbf)$. 

Let $s \in \ofrak$. There is a unique maximal Levi $\Lbf_s^* \subset \Zbf_{\Gbf^*}(s)$ of $\Gbf^*$, namely, the Levi of $\Gbf^*$ whose root datum is spanned by simple roots of $(\Gbf^*, \Tbf^*)$ that are contained in the root datum of $(\Zbf_{\Gbf^*}^{\circ}(s), \Tbf^*)$. We shall denote by $\Lbf_s \subset \Gbf$ the corresponding Levi of $\Gbf$. Clearly, $s \in \Zbf(\Lbf_s^*)$, and by Lemma \ref{lem:endoscopy-trivial-case}, the endoscopic equivalence is given by 
$$\HH^{\Lbf_s}_{[s,s]} \xrightarrow{\otimes \Lcal_s^{\Lbf_s,-1}} \HH^{\Hbf^{\Lbf_s}}_{[s,s]}.$$
Let us denote this functor by $\EE_s^{\Lbf}$. Note that this functor does not depend on any choices. 

Let $s,t \in \ofrak$ and let $\sigma \in {_s}\Wbf_t$ be a simple reflection of $\Wbf$ and denote by $\beta$ its image in ${_s}\underline{\Wbf}_t$. There are canonically defined objects 
\begin{enumerate}
\item $\Delta_t^{\beta, \Gbf} \in \HH^{\Gbf, \beta}_{[s,t]}$, given by rigidified minimal standard sheaves, see Section \ref{sec:Rigidified-minimal-standard}, 
\item $\Delta_t^{\beta, \Hbf} \in \HH^{\Hbf, \beta}_{[s,t]}$, constructed in Section \ref{sec:functor-V-endoscopic-side}.  
\end{enumerate}
Furthermore there is a canonical isomorphism 
$$\EE(\Delta_t^{\beta, \Gbf}) \cong \Delta_t^{\beta, \Hbf}.$$

Consider the groupoid $\mathfrak{E}$ of all natural transformations of functors $\Xibf_{\ofrak}^{\opp} \to \Morita$
$$E : \HH_{\ofrak}^{\Gbf} \Rightarrow \HH_{\ofrak}^{\Hbf},$$
such that for all $\beta \in {_s}\underline{\Wbf}_t$, the functor $E^{\beta} : \HH^{\Gbf, \beta}_{[s,t]} \to \HH^{\Hbf, \beta}_{[s,t]}$ is $t$-exact, sends standard/costandard/tilting sheaves to standard/costandard/tilting sheaves, and is equipped with isomorphisms
\begin{enumerate}
\item $E_{|\HH_{[s,s]}^{\Lbf_s}} \cong \EE_s$,
\item $E(\Delta_t^{\beta, \Gbf}) \cong \Delta_t^{\beta, \Hbf}$, compatible with the isomorphisms $\Delta^{\beta_1, \Gbf} * \Delta^{\beta_2, \Gbf} \cong \Delta^{\beta_1\beta_2, \Gbf}$ for all composable blocks $\beta_1,\beta_2$. 
\end{enumerate}

\begin{thm}\label{thm:uniqueness-of-endoscopic-equiv}
The groupoid $\mathfrak{E}$ is contractible and nonempty.
\end{thm}

\begin{proof}
The fact that this groupoid is nonempty follows from Theorem \ref{thm:EndoscopyHeckeCat}, which provides a functor $\EE \in \mathfrak{E}$. 

Let $E \in \mathfrak{E}$ and define $A := \EE^{-1} \circ E$. This is an automorphism of $\HH_{\ofrak}^{\Gbf}$. It is enough to prove that there is a unique isomorphism $A \cong \id_{\HH_{\ofrak}^{\Gbf}}$ compatible with the above data. For all pairs $(\sigma, s)$ where $\sigma$ is a simple reflection of $\Wbf$ and $s \in \ofrak$, the functor $A$ is equipped with an isomorphism $a_{\sigma,s} : A(T_{\sigma,s}) \cong T_{\sigma,s}$. As $A$ yields a monoidal self-equivalence of $\HH_{\ofrak}^{\Gbf}$, for all Bott-Samelson tilting sheaves $T_{\underline{w}, s}$ (see Definition \ref{def:Bott-Samelson-Tilting}), we have a uniquely determined isomorphism 
$$a_{\underline{w},s} : A(T_{\underline{w},s}) \cong T_{\underline{w},s}.$$

Our strategy is to prove that the collection of maps $(a_{\underline{w},s})_{\underline{w},s}$ are natural, that is, for all pairs $(\underline{w},s), (\underline{w}',s')$ the composition 
\begin{equation}\label{eq:composition-natural}
\Hom(T_{\underline{w},s}, T_{\underline{w}',s'}) \xrightarrow{A} \Hom(A(T_{\underline{w},s}), A(T_{\underline{w}',s'})) \xrightarrow{a} \Hom(T_{\underline{w},s}, T_{\underline{w}',s'})
\end{equation}
is the identity, where the first map is induced by the functor $A$ and the second by the isomorphisms $a_{\underline{w}, s}$ and $a_{\underline{w}', s'}$. Assuming this, the maps $(a_{\underline{w},s})$ define a natural isomorphism $A \Rightarrow \id$ on the subcategory of Bott-Samelson objects. Since the Hecke category is obtained by taking $\Ind\mathrm{K}^b(\BSBim)$, this natural transformation extends to a natural isomorphism on all of $\HH_{\ofrak}^{\Gbf}$. 

Let $w \in \Wbf$ and $s \in \ofrak$. We claim that there is a canonical isomorphism $b_{w,s} : A(\Delta_{w,s}^{\rigid}) \cong \Delta_{w,s}^{\rigid}$. If $w$ is a simple reflection, this is part of the defining data of $A$. In general, write $w = \sigma v$ where $\ell(v) + 1 = \ell(w)$ and $\sigma$ is a simple reflection. Then we have 
\begin{align*}
\Av_{\psi}(A(\Delta_{w,s}^{\rigid})) &\cong 1^*\Av_{\psi}(A(\Delta_{\sigma,vs}^{\rigid} * \Delta_{v,s}^{\rigid})) \\
&\cong \Av_{\psi}(A(\Delta_{\sigma,vs}^{\rigid}) * A(\Delta_{v,s}^{\rigid})) \\
&\cong \Av_{\psi}(\Delta_{\sigma,vs}^{\rigid}) * A(\Delta_{v,s}^{\rigid}) \\
&\cong \Av_{\psi}A(\Delta_{v,s}^{\rigid}),
\end{align*}
where the first isomorphism is by the multiplicativity of rigidified standard sheaves, the second by the monoidality of $A$ and the rest of the calculation proceeds by induction as in the proof of Lemma \ref{lem:RigidifiedMinStandard}. The isomorphisms $(b_{w,s})$ are compatible with the convolution of rigidified standard sheaves. 

We now prove the naturality of the collection of maps $(a_{\underline{w},s})$. For $(\underline{w},s)$ and $(\underline{w}',s')$, using that $T_{\sigma,s}^{\vee} \cong T_{\sigma, \sigma s}$ for a simple reflection $\sigma$, we can assume that $\underline{w} = 1$. Therefore, it is enough to show that for all $(\underline{w}', s')$, the composition (\ref{eq:composition-natural}) $\Hom(\Delta_{1,s}, T_{\underline{w}'},s') \to 
\Hom(\Delta_{1,s}, T_{\underline{w}',s'})$ is the identity. We show the more general statement that for all $w \in \Wbf$, the composition 
$$\Hom(\Delta_{w,s}^{\rigid}, T_{\underline{w}',s'}) \xrightarrow{A} \Hom(A(\Delta_{w,s}^{\rigid}), A(T_{\underline{w}',s'})) \rightarrow  \Hom(\Delta_{w,s}^{\rigid}, T_{\underline{w}',s'})$$
is the identity, where the second map is induced by the isomorphisms $a_{\underline{w}',s'}$ and $b_{w,s}$. We can proceed by induction on the length of $\underline{w}'$. Let us write $T_{\underline{w}',s'} = T_n = T_{n-1} * T_{\sigma_n, s'}$ and use the short exact sequence 
$$0 \to \Delta_{1,s'} \to T_{\sigma_n,s'} \to \nabla_{\sigma_n,s'}^{\rigid} \to 0,$$
we then have a commutative diagram 
\[\begin{tikzcd}
	{\Hom(\Delta_{v,s'}^{\rigid}, T_{n-1}) } & {\Hom(A(\Delta_{v,s'}^{\rigid}), A(T_{n-1})) } & {\Hom(\Delta_{v,s'}^{\rigid}, T_{n-1}) } \\
	{\Hom(\Delta_{v,s'}^{\rigid}, T_n) } & {\Hom(A(\Delta_{v,s'}^{\rigid}), A(T_n)) } & {\Hom(\Delta_{v,s'}^{\rigid}, T_n) } \\
	{\Hom(\Delta_{v,s'}^{\rigid}, T_{n-1} * \nabla_{\sigma_n,s'}^{\rigid}) } & {\Hom(A(\Delta_{v,s'}^{\rigid}), A(T_{n-1}* \nabla_{\sigma_n,s'}^{\rigid})) } & {\Hom(\Delta_{v,s'}^{\rigid}, T_{n-1} * \nabla_{\sigma_n,s'}^{\rigid})}
	\arrow["A", from=1-1, to=1-2]
	\arrow[from=1-1, to=2-1]
	\arrow["a", from=1-2, to=1-3]
	\arrow[from=1-2, to=2-2]
	\arrow[from=1-3, to=2-3]
	\arrow["A", from=2-1, to=2-2]
	\arrow[from=2-1, to=3-1]
	\arrow["a", from=2-2, to=2-3]
	\arrow[from=2-2, to=3-2]
	\arrow[from=2-3, to=3-3]
	\arrow["A"', from=3-1, to=3-2]
	\arrow["a"', from=3-2, to=3-3]
\end{tikzcd}\]\\
Using Remark \ref{remark:rigidified-costandard}, it follows that the third row of the diagram is isomorphic to
$$0 \to \Hom(\Delta_{v\sigma_n, s_n}^{\rigid}, T_{n-1}) \to \Hom(A\Delta_{v\sigma_n, s_n}^{\rigid}), A(T_{n-1})) \to \Hom(\Delta_{v\sigma_n, s_n}^{\rigid}, T_{n-1}) \to 0,$$
the theorem then follows by induction. 
\end{proof}

\begin{corol}\label{corol:dependency-on-pinning}
The equivalence $\EE$ depends only on the pinning of $\Gbf$ and not on the extended pinning chosen in Section \ref{sec:ExtendedPinning}.
\end{corol}

\begin{rque}
While the equivalence $\EE$ depends only on the pinning, the Soergel description of Theorem \ref{thm:SoergelRealizationHSide} depends on the extended pinning of $\Gbf$. 
\end{rque} 

\subsubsection{Changing the pinning of $\Gbf$}\label{sec:changing-the-pinning}

The endoscopic equivalence $\EE$ depends on the pinning of $\Gbf$ that we have chosen. Let $\Pin$ be the variety of all pinnings of $(\Gbf, \Ubf^-)$. Let $\Vbf = \Ubf^-_{\ab} = \prod_{\alpha \in \Delta} \Gbf_{-\alpha}$ be the abelianization of $\Ubf^-$; this is a vector group. 
Let $\Vbf^{\vee}$ be the dual of $\Vbf$ and $\Vbf^{\vee,*} = \prod_{\alpha\in \Delta}\Gbf_{-\alpha}^{\vee, *} \subset \Vbf^{\vee}$ and $\Gbf^{\vee,*}_{-\alpha} = \Gbf_{-\alpha}^{\vee} - 0$. The group $\Tbf_{\ad}$ acts on both $\Pin$ and $\Vbf^{\vee,*}$.

\begin{lem}
Both $\Pin$ and $\Vbf^{\vee,*}$ are isomorphic as $\Tbf_{\ad}$-torsors. 
\end{lem}

\begin{proof}
Any element $\phi' \in \Vbf^{\vee, *}$ is a map $\Vbf \to \Ga$ which yields a nonzero map 
$$\Gbf_{-\alpha} \hookrightarrow \Vbf \to \Ga$$ and therefore a pinning of $(\Gbf, \Ubf^-)$. This defines a map $\Vbf^{\vee,*} \to \Pin$ which is the desired isomorphism. 
\end{proof}

Let $t \in \Tbf_{\ad}$, we denote by $\phi_t = \Ad(t)(\phi)$ and by $\Lcal_{\psi_t} = \phi_t^*\AS$ the corresponding Artin-Schreier sheaf. For a pinning $\phi' \in \Pin$, we shall write $(\Delta_{\phi'}^{\Gbf, \beta})$ the collection of rigidified minimal standard sheaves defined using the pinning $\phi'$ in place of $\phi$. Similarly, we will write $\EE_{t}$ for the endoscopic equivalence induced by $\phi_t$.

\begin{lem}\label{lem:calculation-conj-on-min-std}
Let $\beta \in {_{\beta s}}\underline{\Wbf}_s$ be a block, and let $t \in \Tbf$. There is a canonical $\RR_{\Tbf \times \Tbf}$-linear isomorphism 
$$\Delta_{\phi_t}^{\Gbf, \beta} \cong (L_{\Tbf \times \Tbf}(t, t^{-1}) \otimes_{\Lambda} \Lcal_{s}(w_{\beta}(t)/t)) \otimes_{\RR_{\Tbf \times \Tbf}} \Delta_{\phi}^{\Gbf, \beta},$$
where $L_{\Tbf \times \Tbf}(t, t^{-1})$ (resp. $\Lcal_{s}(w_{\beta}(t)/t)$)  is the fiber of $L_{\Tbf \times \Tbf}$ (resp. $\Lcal_s$) at the point $(t,t^{-1}) \in \Tbf \times \Tbf$ (resp. $w_{\beta}(t)/t) \in \Tbf$). 
\end{lem}

\begin{proof}
We have $\Av_{\psi} = \Av_{\psi_t} \circ \Ad(t)_*$, and hence $\Delta_{\phi_t}^{\Gbf,\beta} = \Ad(t)_*\Delta_{\phi}^{\Gbf, \beta}$. The lemma now follows from the fact that $\Delta_{\phi}^{\Gbf,\beta}$ is left $L_{\Tbf} \otimes_{\Lambda} \Lcal_{\beta s}$-equivariant and right $\Lcal_{\Tbf} \otimes_{\Lambda} \Lcal_s$-equivariant. 
\end{proof}

Let $\ofrak \subset \Ch(\Tbf)$, $s \in \ofrak$, and $\beta$ be a block, and consider the following collection of sheaves on $\Tbf$:
$$\Lcal_s^{\beta} = \Ad(w_{\beta})^*\Lcal_s \otimes_{\Lambda} \Lcal_s^{-1}.$$ 
As the elements $w_{\beta}$ are multiplicative for composable blocks (see Lemma \ref{lem:ProductOfBlocks}) there is a canonical isomorphism 
$$\Lcal_{\beta_2(s)}^{\beta_1} \otimes \Lcal_{s}^{\beta_2} \cong \Lcal_s^{\beta_1\beta_2}.$$
Moreover these isomorphisms are compatible with composition of blocks and therefore the map 
$$\otimes \Lcal(t) : \HH_{\ofrak,\tot}^{\Hbf} = \bigoplus_{s,\beta} \HH^{\Hbf, \beta}_{[\beta(s),s]} \xrightarrow{\bigoplus_{s,\beta} \otimes_{\Lambda} \Lcal_{s}^{\beta}(t)} \bigoplus_{s,\beta} \HH^{\Hbf, \beta}_{[\beta(s),s]} = \HH_{\ofrak,\tot}^{\Hbf}$$
is monoidal. 

\begin{lem}\label{lem:changing-pinning}
Let $t \in \Tbf$ and let $\ofrak \subset \Ch(\Tbf)$. There is a canonical isomorphism of functors $ \Lcal(t) \otimes \Ad(t)_*\circ \EE \cong \EE \circ \Ad(t)_*$.
\end{lem}

\begin{proof}
The strategy is to apply Theorem \ref{thm:uniqueness-of-endoscopic-equiv}. When restricted to $\Hbf^{\Lbf_s}_{[s,s]}$, the functor $\otimes \Lcal(t)$ is trivial as $\Hbf^{\Lbf_s}_{[s,s]} \subset \Hbf^{\Gbf, \circ}_{[s,s]}$. On the other hand, by the calculation of Lemma \ref{lem:calculation-conj-on-min-std}, we have $\Ad(t)_*\Delta^{\Hbf, \beta} \cong L_{\Tbf \times \Tbf}(t,t^{-1}) \otimes_{\RR_{\Tbf \times \Tbf}} \Delta^{\Hbf, \beta}$ and $\Ad(t)_*\Delta^{\Gbf, \beta} \cong \Lcal^{\beta}(t) \otimes_{\Lambda} L_{\Tbf \times \Tbf}(t,t^{-1}) \otimes_{\RR_{\Tbf \times \Tbf}} \Delta^{\Gbf, \beta}$. It follows that there are canonical isomorphisms, compatible with convolution
$$\Lcal(t) \otimes \Ad(t)_*\circ \EE(\Delta^{\Gbf, \beta}) \cong \EE \circ \Ad(t)_*(\Delta^{\Gbf, \beta}).$$
\end{proof}

Let $s \in \Ch(\Tbf)$, and let $\beta \in {_{s'}}\underline{\Wbf}_s$ be a block. Consider the sheaf $\Lcal_{s}^{\beta} \in \DD(\Tbf, \Lambda)$. It is immediate from its definition that it descends to a sheaf on $\Tbf_{\ad}$; let $\Lcal_{s,\ad}^{\beta} \in \DD(\Tbf_{\ad}, \Lambda)$ denote the corresponding sheaf. Since, on $\Tbf$, it is a rank-one monodromic sheaf, the same is true on $\Tbf_{\ad}$. Consequently, taking its fiber at $1$, we have a map
$$\pi_1(\Tbf_{\ad}) = X_*(\Tbf_{\ad}) \otimes_{\Z} \hat{\Z}^{(p)}(1) \to \Lambda^{\times}.$$
Considering only the case where $s' = s$ and using the topological generator of $\hat{\Z}^{(p)}(1)$, we get a bilinear map 
$$\ev^{\geo} : \Gammabf_s \times X_*(\Tbf_{\ad}) \to \Lambda^{\times}.$$
On the other hand, recall that, by Lemma \ref{lem:ControlCentralizer}, the group $\Gammabf_s$ is naturally a subgroup of $\Zbf(\Gbf^*_{\simpconn})$. Since $X_*(\Tbf_{\ad}) = X^*(\Tbf^*_{\simpconn})$, we get a pairing
$$\ev : \Gammabf_s \times X_*(\Tbf_{\ad}) \to \Fqb^{\times} \to \Lambda^{\times}$$
where the second map is the embedding $\Fqb^{\times} \subset \Lambda^{\times}$ induced by the topological generator of $\hat{\Z}^{(p)}(1)$. 

\begin{lem}
We have $\ev = \ev^{\geo}$.
\end{lem}

\begin{proof}
Let $\beta \in \Gammabf_s$. Since on $\Tbf$, we have $\Lcal_s^{\beta} = \Ad(w_{\beta})^*\Lcal_s \otimes_{\Lambda} \Lcal_s^{-1}$, it follows that the induced map $\ev^{\geo}(\beta, -) : X_*(\Tbf) \to X_*(\Tbf_{\ad}) \to \Lambda^{\times}$ is given by $\lambda \mapsto \Ad(w_{\beta})(\lambda)(s)/\lambda(s)$. Let $z_{\beta} \in \Zbf(\Gbf^*_{\simpconn})$ be the image of $\beta$. We have that $z_{\beta} = [\dot{n}_{w_{\beta}}, \dot{s}] \in \Zbf(\Gbf^*_{\simpconn})$ where $\dot{s}$ is any lift of $s \in \Tbf^*$ to $\widetilde{\Gbf}^*$ any central extension of $\Gbf^*$ whose derived subgroup is $\Gbf^*_{\simpconn}$ and $\dot{n}_{w_{\beta}}$ is any lift in $\widetilde{\Gbf}^*$ of $w_{\beta}$. It follows that 
$$\ev(\beta, -) : X^*(\Tbf^*_{\simpconn}) \to \Lambda^{\times}$$
is given by $\lambda \mapsto \lambda([\dot{n}_{w_{\beta}}, \dot{s}])$, and when restricted to $X^*(\Tbf^*)$, this is the map $\lambda \mapsto \Ad(w_{\beta})\lambda(s)/\lambda(s)$. The lemma follows.
\end{proof}

Recall that we have a pairing 
$$\ev^{\can} : \pi_0(\Zbf(\Gbf)) \times \pi_1(\Gbf^*_{\der}) \to \Lambda^{\times}$$
after embedding $\Fqb^{\times}$ in $\Lambda^{\times}$. 
By Lemma \ref{lem:ControlCentralizer}, we have a canonical map $\Gammabf_s \to \pi_1(\Gbf_{\der}^*)$. 
Consider the short exact sequence $1\to\Zbf(\Gbf)\to\Tbf\to\Tbf_{\ad}\to1$. The resulting exact sequence of character groups and the functor $\Hom(-,\Z)$ give a boundary map
$$X_*(\Tbf_{\ad})\longrightarrow \Ext^1\!\bigl(X^*(\Zbf(\Gbf)),\Z\bigr)
\cong \Ext^1\!\bigl(X^*(\pi_0\Zbf(\Gbf)),\Z\bigr).$$
For the finite character group $A=X^*(\pi_0\Zbf(\Gbf))$, the canonical identification $\Ext^1(A,\Z)=\Hom(A,\Q/\Z)$, followed by the fixed prime-to-$p$ Tate generator, identifies this group with the geometric points of $\pi_0(\Zbf(\Gbf))$. We therefore obtain the required map
$$X_*(\Tbf_{\ad})\longrightarrow \pi_0(\Zbf(\Gbf)).$$

\begin{lem}
The following diagram is commutative 
\[\begin{tikzcd}
	{\pi_0(\Zbf(\Gbf)) \times \pi_1(\Gbf^*_{\der}) } & {\Lambda^{\times}} \\
	{X_*(\Tbf_{\ad}) \times \Gammabf_s.}
	\arrow["{\ev^{\can}}", from=1-1, to=1-2]
	\arrow[from=2-1, to=1-1]
	\arrow["\ev"', from=2-1, to=1-2]
\end{tikzcd}\]
\end{lem}

\begin{proof}
This is a standard calculation.
\end{proof}

\subsubsection{Frobenius equivariance}

We discuss the Frobenius equivariance of the endoscopic equivalence. Let $\ofrak \subset \Ch(\Tbf)$ be a $\Wbf$-orbit which we assume to be stable under $\Frob$ (these are the only relevant orbits for the categorical trace construction of Section \ref{sectionJordan}). 

Hence, for all $s,s'\in\Ch(\Tbf)$, there are well-defined functors 
$$\HH_{[s,s']}^{\Gbf} \xrightarrow{\Frob^*} \HH^{\Gbf}_{[\Frob(s), \Frob(s')]}.$$ 
On the dual side, there is a well-defined morphism 
$$\Hbf_{[s,s']} \xrightarrow{\Frob} \Hbf_{[\Frob(s), \Frob(s')]}$$ 
which induces a well-defined morphism 
$$\HH^{\Hbf}_{[s,s']} \xrightarrow{\Frob^*} \HH^{\Hbf}_{[\Frob(s), \Frob(s')]}.$$ 
These functors are compatible with the monoidal structures and they induce endomorphisms
$$\HH^{\Gbf}_{\ofrak,\tot} \to \HH^{\Gbf}_{\ofrak, \tot}$$ 
and 
$$\HH^{\Hbf}_{\ofrak, \tot} \to \HH^{\Hbf}_{\ofrak, \tot}$$
as functors $\Xibf_{\ofrak}^{\opp} \to \Morita$.

\begin{lem}\label{lem:FrobeniusEquivariance}
The equivalence of Theorem \ref{thm:EndoscopyHeckeCat} is naturally $\Frob$-equivariant. 
\end{lem} 

\begin{proof}
The morphism $\Frob$ also induces an endomorphism $\Frob : \RR_{\Tbf} \to \RR_{\Tbf}$ and therefore also an endomorphism of $\SBim$ (again as a functor $\Xibf_{\ofrak}^{\opp} \to \Morita$). It is therefore enough to prove that both equivalences 
$$\HH^{\Gbf} \cong \SBim \cong \HH^{\Hbf}$$
are Frobenius equivariant. The functor $\Frob^* : \HH_{[s,s']}^{\Gbf} \to \HH^{\Gbf}_{[\Frob(s), \Frob(s')]}$ sends tilting objects to tilting objects and is $t$-exact. Hence, it is enough to prove that the functor $\V$ is $\Frob$-equivariant. Unwinding its construction, it is enough to prove that both functors $\Av_{\psi}$ and $\Av_{\Ubf}$ are $\Frob$-equivariant; this is clear, as both $\Ubf$ and the pinning of $\Gbf$ are Frobenius equivariant. 

On the endoscopic side, since the conditions of Section \ref{sec:TrivializationWhittakerModel} are stable under Frobenius, the formation of the objects $\Delta^{\Hbf, \beta}$ is Frobenius equivariant. Finally, since the pinnings of the groups $\Hbf^{\circ}_s$ are Frobenius equivariant, it follows that the corresponding functors $\Av_{\psi}$ and $\Av_{\Ubf}$ are also Frobenius equivariant. 
\end{proof}

\subsubsection{Compatibility with isogenies}\label{sec:compat-with-isog-hecke-cat}

Let $f : \Gbf \to \Gbf'$ be a finite central isogeny. We shall use the notations of Section \ref{sec:PrelimEndoscopicGroupsAndIsogeny}. 
\begin{lem}[\protect{\cite[Lemmas 2.3.6 and 2.3.8]{EteveFreeMonodromic}}]\label{lem:PushForwardIsogeny}
Denote again by $f : \Tbf \to \Tbf'$ the map induced by $f$.
There is an isomorphism of $\Lambda$-sheaves
$$f_*(L_{\Tbf}) \cong \bigoplus_{z \in \Ch(\ker(f))}(L_{\Tbf'} \otimes_{\Lambda} \Lcal_z).$$
\end{lem}

Let $s' \in \Ch(\Tbf')$ with image $s \in \Ch(\Tbf)$. 
\begin{lem}
There are isomorphisms 
\begin{enumerate}
\item $f_*\Delta_{w,s} \cong \oplus_{z \in \Zbf_f} \Delta_{w,zs'}$, 
\item $f_*\nabla_{w,s} \cong \oplus_{z \in \Zbf_f} \nabla_{w,zs'}$, 
\item $f_*T_{w,s} \cong \oplus_{z \in \Zbf_f} T_{w,zs'}$, 
\item $f_*\Av_{\psi} \cong \Av_{\psi} f_*$.
\end{enumerate}
\end{lem}

\begin{proof}
Since $f$ is finite, the functor $f_*$ commutes with $!$ and $*$ pullback and pushforward to strata. Together with Lemma \ref{lem:PushForwardIsogeny}, this implies $(i)$ and $(ii)$. Recall that a sheaf $T$ is tilting if and only if its $!$ and $*$ restrictions to strata are perverse. Clearly, this condition is preserved by $f_*$; this shows $(iii)$. Finally, the map $f$ yields an isomorphism $\Ubf \to \Ubf'$; this implies $(iv)$. 
\end{proof}

Let $\ofrak \subset \Ch(\Tbf)$ be a $\Wbf$-orbit, and let $\ofrak' \subset \Ch(\Tbf')$ be a $\Wbf$-orbit mapping to $\ofrak$. The map $f_*$ is monoidal, and there are monoidal maps
$$\HH^{\Gbf}_{\ofrak, \tot} \to \bigoplus_{\ofrak'\mapsto \ofrak} \HH^{\Gbf'}_{\ofrak', \tot} \xrightarrow{\pr_{\ofrak'}} \HH^{\Gbf'}_{\ofrak', \tot},$$
where the second map is the projection onto the corresponding direct summand. 

\begin{lem}\label{lem:MonoidalPushIsogenies}
Let $s' \in \Ch(\Tbf')$ have image $s \in \Ch(\Tbf)$. The induced map 
$$\HH^{\Gbf, \circ}_{[s,s]} \to \HH^{\Gbf', \circ}_{[s',s']}$$
is monoidal. 
\end{lem}

\begin{proof}
Clearly if $A \in \HH^{\Gbf, \circ}_{[s,s]}$, then $f_*A \in \bigoplus_{z}\HH^{\Gbf',\circ}_{[zs',zs']}$. Since this direct sum is a direct sum of monoidal categories, the projection  $\bigoplus_{z}\HH^{\Gbf',\circ}_{[zs',zs']} \to \HH^{\circ,\Gbf'}_{[s',s']}$ is monoidal. The lemma follows. 
\end{proof}

Let $(s_1',s_2',\beta')$ be a triple with $s_1',s_2' \in \Ch(\Tbf')$ and $\beta' \in {_{s_1'}}\underline{\Wbf'}_{s_2'}$, with image $(s_1,s_2, \beta)$ under the corresponding map for $\Tbf$. There is a well-defined map 
$$\HH^{\Gbf, \beta}_{[s_1,s_2]} \to \HH^{\Gbf', \beta'}_{[s_1',s_2']}$$
which is a map of bimodules over $(\HH^{\Gbf, \circ}_{[s_1,s_1]} \to \HH^{\Gbf', \circ}_{[s'_1,s'_1]}, \HH^{\Gbf, \circ}_{[s_2,s_2]} \to \HH^{\Gbf', \circ}_{[s'_2,s'_2]})$. 

We can summarize this in the following lemma. 
\begin{lem}
Let $\ofrak' \subset \Ch(\Tbf')$ be a $\Wbf$-orbit mapping to $\ofrak \subset \Ch(\Tbf)$. The functor $f_*$ yields a natural transformation of functors $\Xibf_{\ofrak'}^{\opp} \to \Morita$ as  
$$(\Xibf_{\ofrak'}^{\opp} \to \Xibf_{\ofrak}^{\opp} \xrightarrow{\HH^{\Gbf}} \Morita) \Rightarrow (\Xibf_{\ofrak'}^{\opp}\xrightarrow{\HH^{\Gbf'}} \Morita).$$
\end{lem}

On the endoscopic side, for all $s' \in \Ch(\Tbf')$ with image $s \in \Ch(\Tbf)$, there is a well-defined monoidal functor 
$$f^{\Hbf}_* : \HH^{\Hbf,\circ}_{[s,s]} \to \HH^{\Hbf', \circ}_{[s',s']}$$
obtained by pushforward along the map $\Hbf_{[s,s]}^{\circ} \to \Hbf_{[s',s']}^{',\circ}$ and projecting onto the $\Tbf'$ unipotent part. Similarly, for a triple $(s_1',s_2',\beta')$ as above mapping to $(s_1,s_2, \beta)$, pushforward along the map $\Hbf^{\beta}_{[s_1,s_2]} \to \Hbf^{',\beta'}_{[s_1',s_2']}$ and projection onto the category of $\Tbf'$-unipotent sheaves yields a functor 
$$\HH^{\Hbf,\beta}_{[s_1,s_2]} \to \HH^{\Hbf',\beta'}_{[s_1',s_2']}.$$
Since the map $\Hbf^{\beta}_{[s_1,s_2]} \to \Hbf^{',\beta'}_{[s_1',s_2']}$ is linear over $(\Hbf_{[s_1,s_1]}^{\circ} \to \Hbf_{[s'_1,s'_1]}^{',\circ})$ and $(\Hbf_{[s_2,s_2]}^{\circ} \to \Hbf_{[s'_2,s'_2]}^{',\circ})$, the corresponding functor is bilinear over $\HH^{\Hbf,\circ}_{[s_i,s_i]} \to \HH^{\Hbf', \circ}_{[s'_i,s'_i]}$ for $i = 1,2$. The previous paragraph can be summed up in terms of the following lemma.
\begin{lem}
Let  $\ofrak' \subset \Ch(\Tbf')$ be a $\Wbf$-orbit mapping to $\ofrak \subset \Ch(\Tbf)$. Pushforward along $f^{\Hbf}$ yields a natural transformation of functors $\Xibf_{\ofrak'}^{\opp} \to \Morita$ as  
$$(\Xibf_{\ofrak'}^{\opp} \to \Xibf_{\ofrak}^{\opp} \xrightarrow{\HH^{\Hbf}} \Morita) \Rightarrow (\Xibf_{\ofrak'}^{\opp}\xrightarrow{\HH^{\Hbf'}} \Morita).$$
\end{lem}

\begin{prop}
Let $\ofrak' \subset \Ch(\Tbf')$ be a $\Wbf$-orbit mapping to $\ofrak \subset \Ch(\Tbf)$. There is a canonically commuting diagram of functors $\Xibf_{\ofrak'}^{\opp} \to \Morita$
\[\begin{tikzcd}
	{\HH^{\Gbf}_{|\Xibf_{\ofrak'}^{\opp}}} & {\HH^{\Hbf}_{|\Xibf_{\ofrak'}^{\opp}}} \\
	{\HH^{\Gbf'}} & {\HH^{\Hbf'}} \\
	\arrow["\EE", from=1-1, to=1-2]
	\arrow["{{f_*}}", from=1-1, to=2-1]
	\arrow["{{f^{\Hbf}_*}}", from=1-2, to=2-2]
	\arrow["\EE"', from=2-1, to=2-2]
\end{tikzcd}\]
where the horizontal equalities are given by the endoscopic equivalence for $\Gbf$ and $\Gbf'$ and $\HH^{\Gbf}_{|\Xibf_{\ofrak'}^{\opp}}$ denotes the composition $(\Xibf_{\ofrak'}^{\opp} \to \Xibf_{\ofrak}^{\opp} \xrightarrow{\HH^{\Gbf}} \Morita)$ (and similarly for $\HH^{\Hbf}_{|\Xibf_{\ofrak'}^{\opp}}$).
\end{prop}

\begin{proof}
Let us sketch the argument. The map $\Tbf \to \Tbf'$ induces a map $\RR_{\Tbf} \to \RR_{\Tbf'}$ and the target ring is a finite free module over the source. For $(s_1',s_2',\beta')$ mapping to $(s_1, s_2,\beta)$, the functor 
\begin{align*}
\SBim_{[s_1,s_2]}^{\beta} &\to \SBim_{[s_1',s_2']}^{\beta'}  \\
M \mapsto \RR_{\Tbf'} &\otimes_{\RR_{\Tbf}} M 
\end{align*}
is well-defined. This functor yields a natural transformation, for $\Wbf$-orbits $\ofrak'$ and $\ofrak$ as above, of functors $\Xibf_{\ofrak'}^{\opp} \to \Morita$
$$\SBim_{|\Xibf_{\ofrak'}}^{\Gbf} \to \SBim^{\Gbf'}$$ 
where the exponent refers to which group we are using to construct the corresponding category of Soergel bimodules. Since the map $f_*$ commutes with taking the functor $\Av_{\psi}$, we get a natural commutative diagram of functors $\Xibf_{\ofrak'}^{\opp} \to \Morita$ 
\[\begin{tikzcd}
	{\HH^{\Gbf}_{|\Xibf_{\ofrak'}}} & {\HH^{\Hbf}_{|\Xibf_{\ofrak'}}} \\
	{\KK(\SBim^{\Gbf})_{|\Xibf_{\ofrak'}}} & {\KK(\SBim^{\Gbf'})}
	\arrow[equals, from=1-1, to=1-2]
	\arrow["{f_*}", from=1-1, to=2-1]
	\arrow["{f^{\Hbf}_*}", from=1-2, to=2-2]
	\arrow[equals, from=2-1, to=2-2]
\end{tikzcd}\]
The same holds on the endoscopic side and the proposition follows. 
\end{proof}

Finally, let $s' \in \Ch(\Tbf')$ with image $s \in \Ch(\Tbf)$ and let $\ofrak'$ (resp. $\ofrak$) be the $\Wbf$-orbit containing $s'$ (resp. containing $s$). Recall that we have Morita equivalences (see Lemma \ref{lem:MoritaEquivalence1}) $\HH_{\ofrak} \sim \HH_{[s,s]}^{\Gbf}$ and $\HH_{\ofrak'} \sim \HH_{[s',s']}^{\Gbf'}$. The monoidal map $\HH_{\ofrak}^{\Gbf} \to \HH_{\ofrak'}^{\Gbf'}$ equips the target with a structure of a $\HH_{\ofrak}^{\Gbf} \otimes \HH_{\ofrak'}^{\Gbf'}$-module. We now describe the corresponding $(\HH_{[s,s]}^{\Gbf}, \HH_{[s',s']}^{\Gbf'})$-bimodule, we invite the reader to compare with the construction of the bimodule $\Hbf_{[s,s']}'$ (on the endoscopic side) in Section \ref{sec:BimoduleEndoscopicSide}. We define 
$$\HH_{[s,s']}^{\Gbf'} = \bigoplus_{z \in \Zbf_f} \HH_{[s'z,s']}^{\Gbf'}.$$
Proceeding as in Section \ref{sec:BimoduleEndoscopicSide}, we equip this category with the structure of a $(\HH_{[s,s]}^{\Gbf}, \HH_{[s',s']}^{\Gbf'})$-bimodule. The right module structure is induced by the action of $\HH_{[s',s']}^{\Gbf'}$ by right translation on each of the summands. The left action is defined as follows: for every $\beta_1 \in {_s}\underline{\Wbf}_{s}$, $z \in \Zbf_f$, and $\beta_2' \in {_{s'z}}\underline{\Wbf'}_{s'}$, let $(z_1, \beta'_1)$ be the unique pair given by the bijection of Lemma \ref{lem:SectionIsogenyTransporter}, where $z_1 \in \Zbf_f$ and $\beta_1' \in {_{z_1s'}}\underline{\Wbf'}_{zs'}$. We define 
$$a_{\beta_1, z, \beta_2'} : \HH_{[s,s]}^{\Gbf, \beta_1} \times \HH_{[zs',s']}^{\Gbf', \beta_2'} \to \HH_{[z_1s',s']}^{\Gbf',\beta_1'\beta_2'}$$ 
as the composition 
$$ \HH_{[s,s]}^{\Gbf, \beta_1} \times \HH_{[zs',s']}^{\Gbf', \beta_2'} \to \HH_{[z_1s',zs']}^{\Gbf', \beta_1'}\times \HH_{[zs',s']}^{\Gbf', \beta_2'} \to  \HH_{[z_1s',s']}^{\Gbf',\beta_1'\beta_2'}$$
where the first map is induced by the projection $\HH_{[s,s]}^{\Gbf, \beta_1} \to  \HH_{[z_1s',zs']}^{\Gbf', \beta_1'}$ on the first factor. Taking the direct sum over all $\beta_1, z, \beta_2'$ yields the desired action map. The argument of Lemma \ref{lem:ConstructionBimodEndoscopicSide} establishes that this defines an action of $\HH_{[s,s]}^{\Gbf}$ on $\HH_{[s,s']}^{\Gbf'}$. 

\begin{lem}
Under the Morita equivalences  $\HH_{\ofrak} \sim \HH_{[s,s]}^{\Gbf}$ and $\HH_{\ofrak'} \sim \HH_{[s',s']}^{\Gbf'}$, the bimodule $\HH_{\ofrak'}^{\Gbf'}$ corresponds to $\HH_{[s,s']}^{\Gbf'}$.
\end{lem}
\begin{proof}
Recall that the bimodules given the Morita equivalences are 
$${_s}\MM^{\Gbf} = \bigoplus_{s_2} \HH_{[s,s_2]}^{\Gbf}$$
and 
$$\MM^{\Gbf'}_{s'} = \bigoplus_{s_1'} \HH_{[s_1',s']}^{\Gbf'}.$$
It is therefore enough to construct an isomorphism of $(\HH_{[s,s]}^{\Gbf}, \HH_{[s',s']}^{\Gbf'})$-bimodule 
$${_s}\MM^{\Gbf} \otimes_{\HH_{\ofrak}^{\Gbf}} \HH_{\ofrak'}^{\Gbf'} \otimes_{\HH_{\ofrak'}^{\Gbf'}} \MM_{s'}^{\Gbf'} \cong \HH_{[s,s']}^{\Gbf'}.$$
Note that $\HH_{[s,s']}^{\Gbf'}$ is a direct summand of $\MM_{s'}^{\Gbf'}$.
Clearly, we have  ${_s}\MM^{\Gbf} \otimes_{\HH_{\ofrak}^{\Gbf}} \HH_{\ofrak'}^{\Gbf'} \otimes_{\HH_{\ofrak'}^{\Gbf'}} \MM_{s'}^{\Gbf'} = {_s}\MM^{\Gbf} \otimes_{\HH_{\ofrak}^{\Gbf}} \MM_{s'}^{\Gbf'}$. Consider now the maps 
$$ \HH_{[s,s']}^{\Gbf'} \to {_s}\MM^{\Gbf} \otimes_{\HH_{\ofrak}^{\Gbf}} \MM_{s'}^{\Gbf'}$$
given by $A \mapsto \Delta_{1,s} \otimes A$ and 
$${_s}\MM^{\Gbf} \otimes_{\HH_{\ofrak}^{\Gbf}} \MM_{s'}^{\Gbf'} \to \HH_{[s,s']}^{\Gbf'}$$
given by 
$${_s}\MM^{\Gbf} \otimes_{\HH_{\ofrak}^{\Gbf}} \MM_{s'}^{\Gbf'} \to \HH_{\ofrak}^{\Gbf} \otimes_{\HH_{\ofrak}^{\Gbf}} \MM_{s'}^{\Gbf'} \to \MM_{s'}^{\Gbf'}$$
where the map is induced by the inclusion ${_s}\MM^{\Gbf} \subset \HH_{\ofrak}^{\Gbf}$. By construction, the image of ${_s}\MM^{\Gbf} \otimes_{\HH_{\ofrak}^{\Gbf}} \MM_{s'}^{\Gbf'} $ is contained in $\HH_{[s,s']}^{\Gbf'}$. Arguing as in Lemma \ref{lem:MoritaEquivalence1}, we see that these two maps are inverse equivalences of each other. 
\end{proof}

On the endoscopic side, the $(\Hbf_{[s,s]}, \Hbf_{[s',s']}')$-bimodule $\Hbf_{[s,s']}'$ constructed in Section \ref{sec:BimoduleEndoscopicSide} yields a bimodule $\HH_{[s,s']}^{\Hbf'}$ over $(\HH_{[s,s]}^{\Hbf}, \HH_{[s',s']}^{\Hbf'})$.
\begin{corol}
There is a commutative diagram in $\Morita$
\[\begin{tikzcd}
	{\HH_{[s,s]}^{\Gbf}} & {\HH_{[s,s]}^{\Hbf}} \\
	{\HH_{[s',s']}^{\Gbf'}} & {\HH_{[s',s']}^{\Hbf'}.}
	\arrow["{\EE_s^{\Gbf}}", from=1-1, to=1-2]
	\arrow["{{ \HH_{[s,s']}^{\Gbf'}}}", from=2-1, to=1-1]
	\arrow["{\EE_s^{\Gbf'}}", from=2-1, to=2-2]
	\arrow["{{ \HH_{[s,s']}^{\Hbf'}}}"', from=2-2, to=1-2]
\end{tikzcd}\]
\end{corol}

\section{Parabolic character sheaves}\label{sec:parabolic-character-sheaves}

The purpose of this section is to review the theory of $\Frob$-twisted parabolic character sheaves \cite{LusztigParabolicCS1} and its relation with Deligne-Lusztig theory.

\subsection{Horocycle stacks and functors}\label{sec:Definition-horocycle}

Recall that we have fixed a reductive group $\Gbf$ with a Borel pair $\Bbf = \Tbf\Ubf$ and a Frobenius endomorphism $\Frob : \Gbf \to \Gbf$, and that we assume that the pair $(\Tbf, \Bbf)$ is $\Frob$-stable. We consider the $\Frob$-twisted conjugation action of $\Gbf$ on itself, which we denote by $\Ad_{\Frob}$. It is given by the following formula:
$$\Ad_{\Frob}(g)(x) \colondef = gx\Frob(g^{-1}).$$
If $h \in \Gbf$, we will denote by $\Ad_{h\Frob}$ the action given by 
$$\Ad_{h\Frob}(g)(x) \colondef = h\Ad_{\Frob}(g)(x)h^{-1}.$$

\subsubsection*{Horocycle stacks and correspondences}

Following Lusztig \cite{LusztigParabolicCS1} and \cite{SpringerParabolicCS}, let $\Pbf = \Lbf\Ubf_{\Pbf}$ be a standard parabolic with standard Levi $\Lbf$ and unipotent radical $\Ubf_{\Pbf}$. The $\Frob$-twisted horocycle stack of type $\Pbf$ is the stack 
$$\frac{\Ubf_{\Pbf} \backslash \Gbf/\Frob(\Ubf_{\Pbf})}{\Ad_{\Frob}\Lbf}.$$
There are several ways to present this stack as a global quotient stack; let us describe them.
\begin{enumerate}
\item Consider the action of $\Ad_{\Frob}(\Pbf)$ on $\Gbf/\Frob(\Ubf_{\Pbf})$. Then we have 
$$ \frac{\Gbf/\Frob(\Ubf_{\Pbf})}{\Ad_{\Frob}\Pbf} \cong \frac{\Ubf_{\Pbf} \backslash \Gbf/\Frob(\Ubf_{\Pbf})}{\Ad_{\Frob}\Lbf}.$$
\item Consider the subgroup $\Delta_{\Frob}(\Lbf)(\Ubf_{\Pbf} \times \Frob(\Ubf_{\Pbf})) \subset \Pbf \times \Frob(\Pbf)$ of pairs $(p_1,p_2) \in \Pbf \times \Frob(\Pbf)$ satisfying $\Frob(l_1) = l_2$, where $l_i$ is the image of $p_i$ in $\Lbf = \Pbf/\Ubf_{\Pbf}$ (resp. $\Frob(\Pbf)/\Frob(\Ubf_{\Pbf}) = \Frob(\Lbf)$). Also consider the diagonal subgroup $\Delta(\Gbf) \subset \Gbf \times \Gbf$. Then there is a natural isomorphism 
$$\Delta_{\Frob}(\Lbf)(\Ubf_{\Pbf} \times \Frob(\Ubf_{\Pbf})) \backslash (\Gbf \times \Gbf)  /\Delta(\Gbf) \cong \frac{\Ubf_{\Pbf} \backslash \Gbf/\Frob(\Ubf_{\Pbf})}{\Ad_{\Frob}\Lbf}$$
induced by the map $(g,h) \mapsto gh^{-1}$.
\end{enumerate}

Let $\Qbf = \Mbf\Ubf_{\Qbf} \subset \Pbf = \Lbf\Ubf_{\Pbf}$ be two standard parabolics with Levi factors $\Mbf$ and $\Lbf$ and unipotent radicals $\Ubf_{\Qbf}$ and $\Ubf_{\Pbf}$ respectively. We introduce the horocycle correspondence \cite{LusztigParabolicCS1}. 

\begin{equation}\label{eq:DefiHorocycle}
\frac{\Ubf_{\Pbf} \backslash \Gbf/\Frob(\Ubf_{\Pbf})}{\Ad_{\Frob}\Lbf} \xleftarrow{\alpha_{\Qbf \subset \Pbf}} \frac{\Ubf_{\Pbf} \backslash \Gbf/\Frob(\Ubf_{\Pbf})}{\Ad_{\Frob}(\Lbf \cap \Qbf)} \xrightarrow{\beta_{\Qbf \subset \Pbf}} \frac{\Ubf_{\Qbf} \backslash \Gbf/\Frob(\Ubf_{\Qbf})}{\Ad_{\Frob}\Mbf}.
\end{equation}

The maps are defined as follows.
\begin{enumerate}
\item The map $\alpha_{\Qbf \subset \Pbf}$ is induced by taking the quotients of $\Gbf/\Frob(\Ubf_{\Pbf})$ by the groups $\Qbf \subset \Pbf$, respectively. Thus, we have 
$$\frac{\Ubf_{\Pbf} \backslash \Gbf/\Frob(\Ubf_{\Pbf})}{\Ad_{\Frob}(\Lbf \cap \Qbf)} \colondef \frac{\Gbf/\Frob(\Ubf_{\Pbf})}{\Ad_{\Frob}(\Qbf)} \to  \frac{\Gbf/\Frob(\Ubf_{\Pbf})}{\Ad_{\Frob}(\Pbf)} \colondef \frac{\Ubf_{\Pbf} \backslash \Gbf/\Frob(\Ubf_{\Pbf})}{\Ad_{\Frob}\Lbf}$$
where the first isomorphism follows from the fact that $\Qbf = (\Qbf \cap \Lbf)\Ubf_{\Pbf}$.  
\item The map $\beta_{\Qbf \subset \Pbf}$ is induced by taking the quotient by $\Qbf$ of the map $\Gbf/\Frob(\Ubf_{\Pbf}) \to \Gbf/\Frob(\Ubf_{\Qbf})$. Thus, we have
$$\frac{\Ubf_{\Pbf} \backslash \Gbf/\Frob(\Ubf_{\Pbf})}{\Ad_{\Frob}(\Lbf \cap \Qbf)} \colondef \frac{\Gbf/\Frob(\Ubf_{\Pbf})}{\Ad_{\Frob}(\Qbf)} \to \frac{\Gbf/\Frob(\Ubf_{\Qbf})}{\Ad_{\Frob}(\Qbf)} \colondef \frac{\Ubf_{\Qbf} \backslash \Gbf/\Frob(\Ubf_{\Qbf})}{\Ad_{\Frob}\Mbf}.$$
\end{enumerate}

Clearly, the map $\alpha_{\Qbf \subset \Pbf}$ is smooth and proper of relative dimension $\dim \Pbf - \dim \Qbf = \dim (\Lbf/(\Lbf\cap \Qbf))$, and the map $\beta_{\Qbf \subset \Pbf}$ is smooth and a fibration by affine spaces of dimension $\dim(\Ubf_{\Qbf}/\Ubf_{\Pbf})$. 

\begin{ex}\label{ex:ExtremalCases}
Let us highlight two examples. 
\begin{enumerate}
\item If $\Pbf = \Gbf$, then the horocycle stack is 
$$\frac{\Gbf}{\Ad_{\Frob}(\Gbf)}.$$
Note that if $\Gbf$ is connected, then by Lang's theorem, there is a canonical isomorphism $\frac{\Gbf}{\Ad_{\Frob}(\Gbf)} \cong \pt/\Gbf^{\Frob}$.
\item If $\Pbf = \Bbf$, then the horocycle stack is 
$$\frac{\Ubf\backslash \Gbf/\Ubf}{\Ad_{\Frob}\Tbf}.$$ 
This stack was studied in detail in \cite{EteveTilting}. We equip this stack with its Bruhat stratification. The strata are then (noncanonically) isomorphic to 
$$\frac{\Ubf \backslash \BwB/\Ubf}{\Ad_{\Frob}\Tbf} = \pt/(\Tbf^{w\Frob} \ltimes (\Ubf \cap ^{w}\Ubf)).$$
\end{enumerate}
\end{ex} 

\subsubsection*{Basic examples and strata}

The stack $\frac{\Ubf_{\Pbf} \backslash \Gbf/\Frob(\Ubf_{\Pbf})}{\Ad_{\Frob}\Lbf}$ is naturally stratified by $(\Pbf, \Frob(\Pbf))$-double cosets which are indexed by ${^\Pbf}\Wbf^{\Frob(\Pbf)}$. Given $v \in {^\Pbf}\Wbf^{\Frob(\Pbf)}$, we shall denote by 
$$i_v^{\Pbf} : \frac{\Ubf_{\Pbf} \backslash \Pbf v\Frob(\Pbf)/\Frob(\Ubf_{\Pbf})}{\Ad_{\Frob}\Lbf} \subset \frac{\Ubf_{\Pbf} \backslash \Gbf/\Frob(\Ubf_{\Pbf})}{\Ad_{\Frob}\Lbf}$$
the corresponding inclusion. We will write it as $i_v$ if the parabolic $\Pbf$ is clear from context. 

\subsubsection*{Categories, adjunction, and transitivity}

We now introduce the category of $\Frob$-twisted parabolic character sheaves. 

\begin{defi}
The category of $\Frob$-twisted character sheaves of type $\Pbf$ is the category of sheaves on the horocycle stack of type $\Pbf$, namely 
$$\DD(\frac{\Ubf_{\Pbf} \backslash \Gbf/\Frob(\Ubf_{\Pbf})}{\Ad_{\Frob}\Lbf}).$$
\end{defi}

\begin{rque}
Usually for character sheaves \cite{CharacterSheaves1} and parabolic character sheaves \cite{LusztigParabolicCS1}, one only selects the subcategory of `monodromic objects' of all sheaves on the appropriate horocycle stack. However, because of the $\Frob$-twisting, all sheaves will satisfy this condition. This can be checked on the two examples \ref{ex:ExtremalCases}, where the horocycle stacks of type $\Gbf$ and $\Bbf$ are stratified by classifying stacks. This fact will be generalized in \Cref{corol:StructureOfStrata} using an appropriate stratification of Lusztig. 
\end{rque}

We define the horocycle (or Harish-Chandra) functors. Let $\Qbf \subset \Pbf$ be two standard parabolic subgroups. We set
\begin{enumerate}
\item $\hc^!_{\Qbf \subset \Pbf} \colondef= \beta_{\Qbf \subset \Pbf,!}\alpha_{\Qbf \subset \Pbf}^*[2\dim (\Ubf_{\Qbf}/\Ubf_{\Pbf})] $ 
\item $\ch_{\Qbf \subset \Pbf} \colondef= \alpha_{\Qbf \subset \Pbf, !} \beta_{\Qbf \subset \Pbf}^*$
\item $\hc^*_{\Qbf \subset \Pbf} \colondef=  \beta_{\Qbf \subset \Pbf,*}\alpha_{\Qbf \subset \Pbf}^!.$
\end{enumerate}
Using the smoothness and properness of the map $\alpha_{\Qbf \subset \Pbf}$ and the smoothness of the map $\beta_{\Qbf \subset \Pbf}$, it is clear that we have an adjunction triple
$$(\hc^!_{\Qbf \subset \Pbf}, \ch_{\Qbf \subset \Pbf}, \hc^*_{\Qbf \subset \Pbf}).$$

\begin{prop}[Transitivity of induction]\label{prop:Transitivity}
Let $\Pbf_1 = \Lbf_1\Ubf_{\Pbf_1} \subset \Pbf_2 = \Lbf_2\Ubf_{\Pbf_2} \subset \Pbf_3 = \Lbf_3\Ubf_{\Pbf_3}$ be three standard parabolic subgroups. Then there are canonical isomorphisms
\begin{enumerate}
\item $\hc^!_{\Pbf_1 \subset \Pbf_2}\hc^!_{\Pbf_2 \subset \Pbf_3} = \hc^!_{\Pbf_1 \subset \Pbf_3}$
\item $\ch_{\Pbf_2 \subset \Pbf_3}\ch_{\Pbf_1 \subset \Pbf_2} = \ch_{\Pbf_1 \subset \Pbf_3}$
\item $\hc^*_{\Pbf_1 \subset \Pbf_2}\hc^*_{\Pbf_2 \subset \Pbf_3} = \hc^*_{\Pbf_1 \subset \Pbf_3}.$
\end{enumerate}
\end{prop}

\begin{proof}
After passing to adjoints, it is enough to prove one of them. Consider the following diagram
\[\begin{tikzcd}
	&& {\frac{\Ubf_{\Pbf_3} \backslash \Gbf/\Frob(\Ubf_{\Pbf_3})}{\Ad_{\Frob}(\Lbf_3 \cap \Pbf_1)}} \\
	& {\frac{\Ubf_{\Pbf_2} \backslash \Gbf/\Frob(\Ubf_{\Pbf_2})}{\Ad_{\Frob}(\Lbf_2 \cap \Pbf_1)}} && {\frac{\Ubf_{\Pbf_3} \backslash \Gbf/\Frob(\Ubf_{\Pbf_3})}{\Ad_{\Frob}(\Lbf_3 \cap \Pbf_2)}} \\
	{\frac{\Ubf_{\Pbf_1} \backslash \Gbf/\Frob(\Ubf_{\Pbf_1})}{\Ad_{\Frob}\Lbf_1}} && {\frac{\Ubf_{\Pbf_2} \backslash \Gbf/\Frob(\Ubf_{\Pbf_2})}{\Ad_{\Frob}\Lbf_2}} && {\frac{\Ubf_{\Pbf_3} \backslash \Gbf/\Frob(\Ubf_{\Pbf_3})}{\Ad_{\Frob}\Lbf_3}}
	\arrow[from=1-3, to=2-2]
	\arrow[from=1-3, to=2-4]
	\arrow[from=2-2, to=3-1]
	\arrow[from=2-2, to=3-3]
	\arrow[from=2-4, to=3-3]
	\arrow[from=2-4, to=3-5]
\end{tikzcd}\]
where the bottom two correspondences are the horocycle correspondences for $(\Pbf_1 \subset \Pbf_2)$ and $(\Pbf_1 \subset \Pbf_3)$, respectively, and the outer one is the horocycle correspondence for $(\Pbf_1 \subset \Pbf_3)$. The middle square is Cartesian, as it can be rewritten as 
\[\begin{tikzcd}
	{\frac{\Gbf/\Ubf_{\Pbf_2}}{\Ad_{\Frob}\Pbf_1}} & {\frac{\Gbf/\Ubf_{\Pbf_2}}{\Ad_{\Frob}\Pbf_2}} \\
	{\frac{\Gbf/\Ubf_{\Pbf_3}}{\Ad_{\Frob}\Pbf_1}} & {\frac{\Gbf/\Ubf_{\Pbf_3}}{\Ad_{\Frob}\Pbf_2}}
	\arrow[from=1-1, to=1-2]
	\arrow[from=2-1, to=1-1]
	\arrow[from=2-1, to=2-2]
	\arrow[from=2-2, to=1-2]
\end{tikzcd}\]
where the vertical maps are induced by $\Gbf/\Ubf_{\Pbf_3} \to \Gbf/\Ubf_{\Pbf_2}$ and the horizontal ones by modding out along the inclusion $\Pbf_1 \subset \Pbf_2$. 
\end{proof}

\subsubsection*{Comparison with Deligne--Lusztig induction}

Let us also explain the link with the usual Deligne-Lusztig varieties. Let $\Pbf = \Lbf\Ubf_{\Pbf}$ be a standard parabolic, and let $v \in {^\Pbf} \Wbf^{\Frob(\Pbf)}$, with a lift $\dot{v}$, be such that $\Ad(\dot{v})\circ \Frob(\Lbf) = \Lbf$. Attached to the tuple $(\Gbf, \Pbf, \dot{v})$, there is a Deligne--Lusztig variety 
$$Y(\dot{v}, \Pbf) \colondef= \{g\Ubf_{\Pbf} \in \Gbf/\UP, g^{-1}\Frob(g) \in \Ubf_{\Pbf}\dot{v}\Frob(\Ubf_{\Pbf}) \}$$
which carries an action of $\Gbf^{\Frob} \times \Lbf^{\dot{v}\Frob}$. 

Correspondingly, there is an isomorphism of stacks induced by the choice of $\dot{v}$, 
$$\Ubf_{\Pbf} \backslash \Pbf v \Frob(\Pbf) /\Frob(\Ubf_{\Pbf}) \simeq \pt/(\Lbf^{\dot{v}\Frob} \rtimes \Ubf_{\Pbf} \cap \Ad(\dot{v})\Frob(\Ubf_{\Pbf})).$$
Let us denote by $k_{\dot{v}} \colondef \pt/(\Lbf^{\dot{v}\Frob} \rtimes \Ubf_{\Pbf} \cap \Ad(\dot{v})\Frob(\Ubf_{\Pbf})) \to \pt/\Lbf^{\dot{v}\Frob}$ the induced projection. 

\begin{lem}[\protect{\cite[Lemma 3.2.7]{EteveFreeMonodromic}}]\label{prop:ComparisonTraceAndDeligneLusztig}
There is an isomorphism of functors
$$\DD(\pt/\Lbf^{\dot{v}\Frob}, \Lambda) \to \DD(\pt/\Gbf^{\Frob}, \Lambda)$$ 
between 
$$\ch_{\Pbf \subset \Gbf} i^{\Pbf}_{v,!} k_{\dot{v}}^* \simeq \RGamma_c(Y(\dot{v}, \Pbf), \Lambda) \otimes_{\Lbf^{\dot{v}\Frob}} -.$$ 
\end{lem}
 
\subsection{Lusztig strata, Levi reduction, and generation}\label{sec:Lusztig-stratification}

\subsubsection*{Indexing sequences}
Let $I \subset \Delta$ be a set of simple roots. 
\begin{defi}\label{def:Sequences}
Let $\Tcal_I$ be the set of sequences $(I_n, w_n)$ where $I_n \subset \Delta$ and $w_n \in \Wbf$ satisfying 
\begin{enumerate}
\item $I_0 = I$ 
\item $I_{n+1} \colondef= I_n \cap \Ad(w_n)\Frob(I_n)$, 
\item $w_{n} \in {^{I_n}}\Wbf^{{\Frob(I_n)}}$,
\item $w_{n+1} \in \Wbf_{I_{n+1}}w_n\Wbf_{\Frob(I_n)}$.
\end{enumerate}
\end{defi}

\begin{rque}
To be consistent with our previous notations, we have reversed left and right in Lusztig's notations. 
\end{rque}

\begin{rque}\label{rque:Constancy}
Let $t \in \Tcal_I$ be such a sequence. By $(ii)$, we have $I_{n+1} \subset I_n$, and the sequence of $I_n$ is constant for $n \gg 0$. We denote by $I_{\infty}$ this constant value. If $I_{n+1} = I_n$, then we have $I_n = \Ad(w_n) \Frob(I_n)$, and combining $(iii)$ and $(iv)$, we see that $w_{n+1} = w_n$. The sequence $(w_n)$ is then also eventually constant, and we denote by $w_{\infty}$ this constant value. 
\end{rque}

\subsubsection*{Construction of the strata}
Let $I \subset \Delta$ and denote by $\Pbf = \Lbf\Ubf_{\Pbf}$ the standard parabolic of type $I$. 
We denote by $\Zbf_I = (\Delta_{\Frob}(\Lbf)(\Ubf_{\Pbf} \times \Frob(\Ubf_{\Pbf}))) \backslash (\Gbf \times \Gbf)$. This is the scheme that parametrizes triples $(\Qbf, \Qbf', u\Ubf_{\Qbf})$ where $\Qbf$ is a parabolic of type $I$, $\Qbf'$ is a parabolic of type $\Frob(I)$ and $u$ satisfies $\Frob(\Ad(u)\Qbf) = \Qbf'$. 

In \cite[Section 3]{LusztigParabolicCS1}, Lusztig constructs a stratification of $\Zbf_I$ indexed by $\Tcal_I$. Let us introduce the following notation:
\begin{enumerate}
\item Given two parabolic subgroups $\Qbf, \Qbf' \subset \Gbf$, we denote by $\Qbf^{\Qbf'} = (\Qbf \cap \Qbf')\Ubf_{\Qbf}$; this is again a parabolic subgroup of $\Gbf$. 
\item Given two parabolic subgroups $\Qbf$ and $\Qbf'$ of types $I$ and $J$, respectively, we denote by $\pos(\Qbf, \Qbf') \in {^I}\Wbf^J$ their relative position. 
\end{enumerate}

Let $(\Qbf, \Qbf', u\Ubf_{\Qbf})$ be a point of $\Zbf_I$. We construct a sequence of triples $(\Qbf_n, \Qbf'_n, u_n\Ubf_{\Qbf_n})$ by induction
\begin{enumerate}
\item $(\Qbf_0, \Qbf'_0, u_0\Ubf_{\Qbf}) \colondef= (\Qbf, \Qbf', u\Ubf_{\Qbf})$, 
\item $(\Qbf_{n+1}, \Qbf'_{n+1}, u_{n+1}\Ubf_{\Qbf_{n+1}})$ is defined as 
$$\Qbf_{n+1} = \Qbf_n^{\Qbf_{n}'},  \Qbf'_{n+1} = (\Qbf'_n)^{\Qbf_{n}}, u_{n+1}\Ubf_{\Qbf_{n+1}} = u_n\Ubf_{\Qbf_{n+1}}.$$
Note that $u_{n+1}\Ubf_{\Qbf_{n+1}}$ is well-defined, as $\Qbf_{n+1} \subset \Qbf_{n}$ and thus $\Ubf_{\Qbf_{n}} \subset \Ubf_{\Qbf_{n+1}}$.
\end{enumerate}

Let $(\Qbf, \Qbf', u\Ubf_{\Qbf})$ be as above and let $(I_n, w_n)$ be the sequence where $I_n \subset \Delta$ is the set of simple roots such that $\Qbf_n$ is of type $I_n$ and $w_n = \pos(\Qbf_n, \Qbf'_n)$. 
\begin{lem}[\protect{\cite[3.11]{LusztigParabolicCS1}}]\label{lem:SequenceFromParab}
Let $(\Qbf, \Qbf', u\Ubf_{\Qbf}) \in \Zbf_I$. Then the sequence $(I_n, w_n)$ constructed above defines an element of $\Tcal_I$. 
\end{lem}

Let $t \in \Tcal_I$. Let $\Zbf_{I,t}$ be the subvariety of triples $(\Qbf, \Qbf', u\Ubf_{\Qbf})$ whose associated sequence in $\Tcal_I$ by  \Cref{lem:SequenceFromParab} is $t$. For a sequence $t = (I_i, w_i)$, we denote by $t_{n}$ the sequence whose terms are $(I_{n+i}, w_{n+i})_i$. By construction, $t_{n}$ is an element of $\Tcal_{I_{n}}$. There is a natural map 
\begin{equation}\label{eq:MapTruncation}
\lambda_{t,1} : \Zbf_{I,t} \to \Zbf_{I_1,t_{1}}
\end{equation}
given by $(\Qbf, \Qbf', u\Ubf_{\Qbf}) \mapsto (\Qbf_1, \Qbf'_1, u_1\Ubf_{\Qbf_1})$. 

\begin{thm}[\protect{\cite[Lemma 3.12]{LusztigParabolicCS1}}]\label{thm:LusztigStratification}
\begin{enumerate}
\item The variety $\Zbf_I$ is stratified by the $\Zbf_{I,t}$. 
\item The map $\Zbf_{I,t} \to \Zbf_{I_1,t_{1}}$ induces a bijection on $\Gbf$-orbits and is a fibration by affine spaces. 
\end{enumerate}
\end{thm}

\subsubsection*{Stack presentations and truncation} Let $I \subset \Delta$ be a set of simple roots and let $t = (I_n, w_n) \in \Tcal_I$. Let us denote by $\Pbf_n = \Lbf_n\Ubf_{\Pbf_n}$ the standard parabolic subgroup of type $I_n$ and by $\Pbf_{\infty} = \Lbf_{\infty}\Ubf_{\Pbf_{\infty}}$ its value for $n \gg 0$. 

Given $(\Qbf, \Qbf', g\Ubf_{\Qbf}) \in \Zbf_{I}$, we denote by $(\Qbf_n, \Qbf'_{n})$ the sequence of parabolic subgroups obtained as before. We define a sequence of $\Gbf$-stable subschemes of $\Zbf_I$ 
$$\Zbf_{I,t} = \Zbf_{I, \infty} \subset \Zbf_{I,k} \subset \dots \Zbf_{I,1} \subset \Zbf_{I,0} \subset \Zbf_I.$$
The scheme $\Zbf_{I,k}$ is the subset of triples $(\Qbf, \Qbf', g\Ubf_{\Qbf})$ such that $\pos(\Qbf_k, \Qbf'_k) = w_k$. Let us note that this is well-defined: if $\pos(\Qbf_k, \Qbf'_k) = w_k$, then $\pos(\Qbf_{k-1}, \Qbf'_{k-1}) = w_{k-1}$ by combining $(iii)$ and $(iv)$ in Definition \ref{def:Sequences}. We also note that $\Zbf_{I,k}$ depends only on the first $k$ terms of the sequence $t$.

\begin{lem}\label{lem:StackyReal}
There is a canonical isomorphism of stacks
$$\Zbf_{I,k}/\Gbf \cong \frac{\Ubf_{\Pbf_0} \backslash \Pbf_k w_k \Frob(\Pbf_k)/\Frob(\Ubf_{\Pbf_0})}{\Ad_{\Frob}(\Lbf_0 \cap \Pbf_k)},$$
making the following diagram commutative and Cartesian
\[\begin{tikzcd}
	{\Zbf_{I,k}} & {\frac{\Ubf_{\Pbf_0} \backslash \Pbf_k w_k \Frob(\Pbf_k)/\Frob(\Ubf_{\Pbf_0})}{\Ad_{\Frob}(\Lbf_0 \cap \Pbf_k)}} \\
	{\Zbf_{I,k-1}} & {\frac{\Ubf_{\Pbf_0} \backslash \Pbf_{k-1} w_{k-1} \Frob(\Pbf_{k-1})/\Frob(\Ubf_{\Pbf_0})}{\Ad_{\Frob}(\Lbf_0 \cap \Pbf_{k-1})}}
	\arrow[from=1-1, to=1-2]
	\arrow[from=1-1, to=2-1]
	\arrow[from=1-2, to=2-2]
	\arrow[from=2-1, to=2-2]
\end{tikzcd}\]
where the right vertical map is induced by the inclusions $\Lbf_0 \cap \Pbf_k \subset \Lbf_0 \cap \Pbf_{k-1}$ and $\Pbf_k w_k \Frob(\Pbf_k) \subset \Pbf_{k-1} w_{k-1} \Frob(\Pbf_{k-1})$.
\end{lem}

\begin{proof}
Consider the scheme $\Zbf_{I, I_k, w_k}$ parametrizing triples $(\Qbf_k, \Qbf'_k, g\Ubf_{\Qbf})$ such that $\Qbf_k$ (resp. $\Qbf'_k$) is a parabolic subgroup of type $I_k$ (resp. $\Frob(I_k)$), $\Frob(\Ad(g)\Qbf_k) =  \Qbf'_k$ and $\Qbf$ is the unique parabolic subgroup of type $I$ containing $\Qbf_k$ and $\pos(\Qbf_k, \Qbf'_k) = w_k$. 

There is a well-defined map $\Zbf_{I,k} \to \Zbf_{I, I_k, w_k}$ given by 
$$(\Qbf, \Qbf', g\Ubf_{\Qbf}) \mapsto (\Qbf_k, \Qbf'_k, g\Ubf_{\Qbf}).$$
This map is clearly a $\Gbf$-equivariant isomorphism as $\Qbf$ (resp. $\Qbf'$) are the unique parabolic subgroups of type $I$ (resp. $\Frob(I)$) containing $\Qbf_k$ (resp. $\Qbf_k'$). 

Let $\Zbf_{I,I_k,w_k}' \subset \Zbf_{I,I_k,w_k}$ be the subscheme of triples where $\Qbf_k = \Pbf_k$. This scheme is isomorphic to $\Pbf_k w_k \Frob(\Pbf_k)/\Frob(\Ubf_{\Pbf_0})$. There is then an isomorphism of stacks
$$\Zbf_{I,I_k,w_k}/\Gbf \cong  \Zbf_{I,I_k,w_k}'/\Pbf_k  \cong \frac{\Pbf_k w_k \Frob(\Pbf_k)/\Frob(\Ubf_{\Pbf_0})}{\Ad_{\Frob}(\Pbf_k)}.$$
Since $\Pbf_k \subset \Pbf_0$, we have an equality $\Pbf_k = (\Pbf_k \cap \Lbf_0)\Ubf_{\Pbf_0}$ and therefore an isomorphism 
$$\Zbf_{I,k}/\Gbf \cong \frac{\Ubf_{\Pbf_0} \backslash \Pbf_k w_k \Frob(\Pbf_k)/\Frob(\Ubf_{\Pbf_0})}{\Ad_{\Frob}(\Lbf_0 \cap \Pbf_k)}.$$
The commutativity of the diagram is immediate. 
\end{proof}

Let us rewrite $\Zbf_{I,k}$ as $\Zbf_{I,t,k}$ to put some emphasis on the sequence $t$. Let $k \geq n$ be two integers. There is a natural map 
$$\Zbf_{I,t,k} \to \Zbf_{I_n,t_n,k}$$
given by $(\Qbf, \Qbf', g\Ubf_{\Qbf}) \mapsto (\Qbf_n, \Qbf'_n, g\Ubf_{\Qbf_n})$. Under the isomorphism of Lemma \ref{lem:StackyReal}, this map is the map 
\begin{align*}
\frac{\Ubf_{\Pbf_0} \backslash \Pbf_k w_k \Frob(\Pbf_k)/\Frob(\Ubf_{\Pbf_0})}{\Ad_{\Frob}(\Lbf_0 \cap \Pbf_k)} &\cong  \frac{\Pbf_k w_k \Frob(\Pbf_k)/\Frob(\Ubf_{\Pbf_0})}{\Ad_{\Frob}(\Pbf_k)} \\
&\longrightarrow \frac{\Pbf_k w_k \Frob(\Pbf_k)/\Frob(\Ubf_{\Pbf_n})}{\Ad_{\Frob}(\Pbf_k)} \cong \frac{\Ubf_{\Pbf_n} \backslash \Pbf_k w_k \Frob(\Pbf_k)/\Frob(\Ubf_{\Pbf_n})}{\Ad_{\Frob}(\Lbf_n \cap \Pbf_k)} 
\end{align*}
induced by the inclusion $\Ubf_{\Pbf_0} \subset \Ubf_{\Pbf_n}$. 

\begin{corol}\label{corol:LusztigAndHorocycle}
There is a canonical commutative diagram 
\[\begin{tikzcd}
	{\frac{\Ubf_{\Pbf_0}\backslash \Pbf_0 w_0 \Frob(\Pbf_0)/\Frob(\Ubf_{\Pbf_0})}{\Ad_{\Frob}\Lbf_0}} & {\frac{\Ubf_{\Pbf_0}\backslash \Pbf_1 w_1 \Frob(\Pbf_1)/\Frob(\Ubf_{\Pbf_0})}{\Ad_{\Frob}(\Lbf_0 \cap \Pbf_1)}} & {\frac{\Ubf_{\Pbf_1}\backslash \Pbf_1 w_1 \Frob(\Pbf_1)/\Frob(\Ubf_{\Pbf_1})}{\Ad_{\Frob}\Lbf_1}} \\
	{\frac{\Ubf_{\Pbf_0}\backslash \Gbf/\Frob(\Ubf_{\Pbf_0})}{\Ad_{\Frob}\Lbf_0}} & {\frac{\Ubf_{\Pbf_0}\backslash \Gbf/\Frob(\Ubf_{\Pbf_0})}{\Ad_{\Frob}(\Lbf_0\cap\Pbf_1)}} & {\frac{\Ubf_{\Pbf_1}\backslash \Gbf/\Frob(\Ubf_{\Pbf_1})}{\Ad_{\Frob}\Lbf_1}}
	\arrow["{i_{t,0}}", from=1-1, to=2-1]
	\arrow[hook', from=1-2, to=1-1]
	\arrow["{\lambda_{t,1}}", from=1-2, to=1-3]
	\arrow["{i_{t,1}}"{description}, from=1-2, to=2-1]
	\arrow["{i'_{t,1}}"{description}, from=1-2, to=2-2]
	\arrow["{i_{t_1,0}}", from=1-3, to=2-3]
	\arrow["{\alpha_{\Pbf_1\subset\Pbf_0}}", from=2-2, to=2-1]
	\arrow["{\beta_{\Pbf_1\subset\Pbf_0}}"', from=2-2, to=2-3]
\end{tikzcd}\]
where the right square is Cartesian and the top left map is an immersion. 
\end{corol}

\begin{proof}
This is essentially a reformulation of Lemma \ref{lem:StackyReal} for the left part of the diagram. The statement about the right part of the diagram follows from the above discussion. 
\end{proof}

\begin{corol}\label{corol:StructureOfStrata}
Let $t \in \Tcal_I$. Then the stack $\Zbf_{I,t}/\Gbf$ is isomorphic to $\pt/(\Lbf_{\infty}^{w_{\infty}\Frob} \ltimes \Vbf_{\infty})$
where $\Vbf_{\infty}$ is a unipotent group. 
\end{corol}

\begin{proof}
Denote the sequence $t$ by $(I_n, w_n)$. We first assume that $w = w_0 = w_{\infty}$ and $I = I_0 = I_{\infty}$. By Remark \ref{rque:Constancy}, this implies that $w\Frob(I_0) = I_0$. Consider now the variety $\Zbf_I$. A coarser stratification than the one constructed by Lusztig is the one given by relative positions of two parabolics, namely we have 
$$\Zbf_I = \bigsqcup_{v \in {^\Pbf}\Wbf^{\Frob(\Pbf)}} \Zbf_{I,v},$$
where $\Zbf_{I,v}$ is the variety of triples $(\Qbf, \Qbf', u\Ubf_{\Qbf})$ where $\pos(\Qbf, \Qbf') = v$. Since $w\Frob(I_0) = I_0$, the variety $\Zbf_{I,w}$ contains a single $\Gbf$-orbit and we have an isomorphism of stacks
$$\Zbf_{I,w}/\Gbf \cong \frac{\Ubf_{\Pbf} \backslash \Pbf w\Frob(\Pbf)/\Frob(\Ubf_{\Pbf})}{\Ad_{\Frob}\Lbf}.$$
Furthermore, a choice of lift $\dot{w} \in \Nbf(\Tbf)$ of $w$ yields an isomorphism
$$\Zbf_{I,w}/\Gbf \cong \pt/(\Lbf^{w\Frob} \ltimes \Ubf_{\Pbf} \cap \Ad(w)\Frob(\Ubf_{\Pbf})),$$
which proves the claim if $w_0 = w_{\infty}$.

Consider now an arbitrary $t\in\Tcal_I$. Iterating the truncation maps gives a $\Gbf$-equivariant morphism
$$q_t:\Zbf_{I,t}\longrightarrow \Zbf_{I_{\infty},t_{\infty}}$$
which induces a bijection on $\Gbf$-orbits and is an iterated affine-space fibration. Thus $\Zbf_{I,t}$ has one $\Gbf$-orbit. Choose a point $x\in\Zbf_{I,t}$ mapping to the point represented by a lift $\dot w_{\infty}$. Then
$$\Zbf_{I,t}/\Gbf\simeq \pt/\Stab_{\Gbf}(x).$$
At the stable term, the explicit quotient in \Cref{lem:StackyReal} gives
$$\Stab_{\Gbf}(q_t(x))\cong \Lbf_{\infty}^{w_{\infty}\Frob}\ltimes \Vbf'_{\infty},$$
where $\Vbf'_{\infty}$ is unipotent. The morphism $q_t$ induces an injective homomorphism of stabilizers; its image projects onto $\Lbf_{\infty}^{w_{\infty}\Frob}$, and its intersection with $\Vbf'_{\infty}$ is a unipotent subgroup $\Vbf_{\infty}$. Hence
$$\Stab_{\Gbf}(x)\cong \Lbf_{\infty}^{w_{\infty}\Frob}\ltimes\Vbf_{\infty},$$
which proves the claim.
\end{proof}

\begin{rque}
We will prove in Lemma \ref{lem:ConnectednessParabCS} that the group $\Vbf_{\infty}$ of Lemma \ref{corol:StructureOfStrata} is connected. 
\end{rque} 

\subsubsection{Compatibility with induction from a Levi}\label{sec:ComparisonWithLevi}

Let $\Pbf = \Lbf\Ubf_{\Pbf}$ be the standard parabolic of type $I \subset \Delta$. Let $v \in {^\Pbf}\Wbf^{\Frob(\Pbf)}$ be an element such that $v\Frob(I) = I$. We fix $\dot{v}$ a lift of $v$. The choice of $\dot{v}$ determines a map 
$$\nu : \Pbf v \Frob(\Pbf) \to \Lbf.$$
This map is defined as follows: for $x \in  \Pbf v \Frob(\Pbf)$, there is a unique $l = \nu(x)$ such that $x = ul\dot{v}u'$ for some elements $u \in \Ubf_{\Pbf}$ and $u' \in \Frob(\Ubf_{\Pbf})$. We will denote by $\Frob_v : \Lbf \to \Lbf$ the map defined by
$$ \Frob_v = \Ad(\dot{v}) \circ \Frob.$$
We will denote by $\Bbf_{\Lbf} = \Bbf \cap \Lbf$ and by $\Ubf_{\Lbf} = \Lbf \cap \Ubf$. 

This map defines natural maps 
\begin{enumerate}
\item $\frac{\Ubf_{\Pbf} \backslash \Pbf v \Frob(\Pbf)/\Frob(\Ubf_{\Pbf})}{\Ad_{\Frob}(\Lbf)} \to \frac{\Lbf}{\Ad_{\Frob_v}(\Lbf)}$, 
\item $\frac{\Ubf_{\Pbf} \backslash \Pbf v \Frob(\Pbf)/\Frob(\Ubf_{\Pbf})}{\Ad_{\Frob}(\Bbf_{\Lbf})} \to \frac{\Lbf}{\Ad_{\Frob_v}(\Bbf_{\Lbf})}$
\item $\frac{\Ubf \backslash \Pbf v \Frob(\Pbf)/\Frob(\Ubf)}{\Ad_{\Frob}(\Tbf)} \to \frac{\Ubf_{\Lbf} \backslash\Lbf/ \Ubf_{\Lbf}}{\Ad_{\Frob_v}(\Tbf)}$
\end{enumerate}

\begin{lem}[\protect{\cite[Proposition 2.5]{LusztigParabolicCS1}}]\label{lem:CombinatoricReduction}
The map 
\begin{align*}
\Wbf_{\Pbf} &\to \Wbf_{\Pbf}v\Wbf_{\Frob(\Pbf)} \\
w &\mapsto wv,
\end{align*}
is a bijection. 
\end{lem}

The map $\nu$ also induces a map 
$$\Bbf wv \Bbf \to \Bbf_{\Lbf} w \Bbf_{\Lbf}.$$

\begin{lem}\label{lem:DiagramReductionToLevi}
Let $w \in \Wbf_{\Pbf}$. The following diagram is commutative, and all maps from the back face to the front face are $\Ubf_{\Pbf} \cap \Frob_v(\Ubf_{\Pbf})$-split gerbes:
\[\begin{tikzcd}
	&& {\frac{\Ubf_{\Pbf} \backslash \Bbf wv \Bbf/\Frob(\Ubf_{\Pbf})}{\Ad_{\Frob}(\Bbf)}} & {\frac{\Ubf \backslash \Bbf wv \Bbf/\Frob(\Ubf)}{\Ad_{\Frob}(\Tbf)}} \\
	& {\frac{\Ubf_{\Pbf} \backslash \Pbf v \Frob(\Pbf)/\Frob(\Ubf_{\Pbf})}{\Ad_{\Frob}(\Lbf)}} & {\frac{\Ubf_{\Pbf} \backslash \Pbf v \Frob(\Pbf)/\Frob(\Ubf_{\Pbf})}{\Ad_{\Frob}(\Bbf_{\Lbf})}} & {\frac{\Ubf \backslash \Pbf v \Frob(\Pbf)/\Frob(\Ubf)}{\Ad_{\Frob}(\Tbf)}} \\
	& {\frac{\Bbf_{\Lbf}w\Bbf_{\Lbf}}{\Ad_{\Frob_v}(\Bbf_{\Lbf})}} & {\frac{\Ubf_{\Lbf} \backslash \Bbf_{\Lbf}w\Bbf_{\Lbf}/ \Ubf_{\Lbf}}{\Ad_{\Frob_v}(\Tbf)}} \\
	{\frac{\Lbf}{\Ad_{\Frob_v}(\Lbf)}} & {\frac{\Lbf}{\Ad_{\Frob_v}(\Bbf_{\Lbf})}} & {\frac{\Ubf_{\Lbf} \backslash\Lbf/ \Ubf_{\Lbf}}{\Ad_{\Frob_v}(\Tbf)}}
	\arrow[from=1-3, to=1-4]
	\arrow[from=1-3, to=2-2]
	\arrow[from=1-3, to=2-3]
	\arrow[from=1-3, to=3-2]
	\arrow[from=1-4, to=2-4]
	\arrow[from=1-4, to=3-3]
	\arrow[from=2-2, to=4-1]
	\arrow[from=2-3, to=2-2]
	\arrow[from=2-3, to=2-4]
	\arrow[from=2-3, to=4-2]
	\arrow[from=2-4, to=4-3]
	\arrow[from=3-2, to=4-1]
	\arrow[from=3-2, to=4-2]
	\arrow[from=3-2, to=3-3]
	\arrow[from=3-3, to=4-3]
	\arrow[from=4-2, to=4-1]
	\arrow[from=4-2, to=4-3]
\end{tikzcd}\]
In particular, the top squares are Cartesian. 
\end{lem}

\begin{proof}
Clear from the definitions.
\end{proof}

\subsubsection*{Conservativity, connectedness, and generation}

\begin{prop}[Conservativity of Harish-Chandra transform]\label{prop:Conservativity}
Let $\Qbf \subset \Pbf$ be two standard parabolic subgroups. The functor $\hc^!_{\Qbf \subset \Pbf}$ is conservative on bounded below objects in $\DD(\frac{\Ubf_{\Pbf} \backslash \Gbf/\Frob(\Ubf_{\Pbf})}{\Ad_{\Frob}\Lbf})$.
\end{prop}

\begin{proof}
The argument is the same as in \cite[Corollary 3.3.5]{EteveFreeMonodromic}, which reduces to \cite[Theorem 1.5]{BBMWhitttaker}. Consider the category $\DD(\frac{\Lbf}{\Ad(\Lbf)})$. This is naturally a monoidal category (for the $!$-convolution) which acts on $\DD(\frac{\Ubf_{\Pbf} \backslash \Gbf/\Frob(\Ubf_{\Pbf})}{\Ad_{\Frob}\Lbf})$. Let $\Spr_{\Pbf}^{\Qbf} \in \DD(\frac{\Lbf}{\Lbf})$ be the Springer sheaf obtained as the equivariant parabolic induction $\ind_{M}^L \delta_1$, where $\delta_1$ is the skyscraper sheaf at $1 \in \frac{\Mbf}{\Mbf}$. There is a canonical isomorphism of functors
$$\ch_{\Qbf \subset \Pbf}\hc^!_{\Qbf \subset \Pbf} \cong \Spr_{\Pbf}^{\Qbf} * -.$$
The proof of \cite{EteveFreeMonodromic} shows that this functor is conservative (if $\ell$ is large, then this is already proved in \cite{BBMWhitttaker}). It follows that the functor $\hc^!_{\Qbf \subset \Pbf}$ is conservative. 
\end{proof}

Let $t \in \Tcal_I$ and let $i_t :  \Zbf_{I,t}/\Gbf \to \frac{\Ubf_{\Pbf_0} \backslash \Gbf/\Frob(\Ubf_{\Pbf_0})}{\Ad_{\Frob}(\Lbf_0)}$ and $i_{t,k} : \Zbf_{I,t,k}/\Gbf \to \frac{\Ubf_{\Pbf_0} \backslash \Gbf/\Frob(\Ubf_{\Pbf_0})}{\Ad_{\Frob}(\Lbf_0)}$ be the inclusions.

\begin{lem}\label{lem:CommutationHCStrata}
\begin{enumerate}
\item Let $A \in \DD(\Zbf_{I,t,1}/\Gbf)$. The sheaf $\hc^!_{\Pbf_1 \subset \Pbf_0} i_{t,1,!}A$ is supported on $\Zbf_{I_1,t_1,0}/\Gbf$. 
\item There is a canonical isomorphism of functors 
$$\hc^!_{\Pbf_1 \subset \Pbf_0} i_{t,1,!} = i_{t_1,0,!}\lambda_{t,1,!}[2(\dim \Ubf_{\Pbf_1} - \dim \Ubf_{\Pbf_0})].$$
\item There is a canonical isomorphism of functors
$$\ch_{\Pbf_1 \subset \Pbf_0} i_{t_1,0,!} = i_{t,1,!}\lambda_{t,1}^*.$$
\end{enumerate}
\end{lem}

\begin{proof}
Using the diagram of Corollary \ref{corol:LusztigAndHorocycle}, point $(iii)$ of the lemma becomes an immediate application of proper base change. Point $(ii)$ is an immediate consequence of point $(i)$; hence, it remains to prove point $(i)$. Let us first refine the diagram of Corollary \ref{corol:LusztigAndHorocycle} as follows:
\[\begin{tikzcd}
	& {\frac{\Ubf_{\Pbf_0}\backslash \Pbf_1 w_1 \Frob(\Pbf_1)/\Frob(\Ubf_{\Pbf_0})}{\Ad_{\Frob}(\Lbf_0 \cap \Pbf_0)}} & {\frac{\Ubf_{\Pbf_1}\backslash \Pbf_1 w_1 \Frob(\Pbf_1)/\Frob(\Ubf_{\Pbf_1})}{\Ad_{\Frob}\Lbf_1}} \\
	{\frac{\Ubf_{\Pbf_0}\backslash \Pbf_0 w_0 \Frob(\Pbf_0)/\Frob(\Ubf_{\Pbf_0})}{\Ad_{\Frob}\Lbf_0}} & {\frac{\Ubf_{\Pbf_0}\backslash \Pbf_0 w_0 \Frob(\Pbf_0)/\Frob(\Ubf_{\Pbf_0})}{\Ad_{\Frob}(\Lbf_0 \cap\Pbf_1)}} & {\frac{\Ubf_{\Pbf_1}\backslash \Pbf_0 w_0 \Frob(\Pbf_0)/\Frob(\Ubf_{\Pbf_1})}{\Ad_{\Frob}\Lbf_1}} \\
	{\frac{\Ubf_{\Pbf_0}\backslash \Gbf/\Frob(\Ubf_{\Pbf_0})}{\Ad_{\Frob}\Lbf_0}} & {\frac{\Ubf_{\Pbf_0}\backslash \Gbf/\Frob(\Ubf_{\Pbf_0})}{\Ad_{\Frob}(\Lbf_0\cap\Pbf_1)}} & {\frac{\Ubf_{\Pbf_1}\backslash \Gbf/\Frob(\Ubf_{\Pbf_1})}{\Ad_{\Frob}\Lbf_1}.}
	\arrow["{\lambda_{t,1}}", from=1-2, to=1-3]
	\arrow["{i_{t,1}}"{description}, hook', from=1-2, to=2-1]
	\arrow["{i'_{t,1}}"{description}, from=1-2, to=2-2]
	\arrow["{i_{t_1,0}}", from=1-3, to=2-3]
	\arrow["{i_{t,0}}", from=2-1, to=3-1]
	\arrow["{\alpha_{\Pbf_1\subset\Pbf_0}}", from=2-2, to=2-1]
	\arrow["{\beta_{\Pbf_1\subset\Pbf_0}}"', from=2-2, to=2-3]
	\arrow[from=2-2, to=3-2]
	\arrow[from=2-3, to=3-3]
	\arrow["{\alpha_{\Pbf_1\subset\Pbf_0}}", from=3-2, to=3-1]
	\arrow["{\beta_{\Pbf_1\subset\Pbf_0}}"', from=3-2, to=3-3]
\end{tikzcd}\]
In this diagram, the top-right and bottom-left squares are Cartesian. This implies that $\hc^!_{\Pbf_1 \subset \Pbf_0} i_{t,1,!}A$ is supported on $\frac{\Ubf_{\Pbf_1}\backslash \Pbf_0 w_0 \Frob(\Pbf_0)/\Frob(\Ubf_{\Pbf_1})}{\Ad_{\Frob}\Lbf_1}$. Let $\tilde{t} = (\tilde{I}_n, \tilde{w}_n)$ be another sequence starting with $(I_0,w_0)$ and such that $\tilde{w}_1 \neq w_1$. Then the parabolic $\tilde{\Pbf}_1$ determined by this sequence is the parabolic $\Pbf_1$, as this one depends only on $(I_0, w_0)$. We then have an immersion 
$$\frac{\Ubf_{\Pbf_1}\backslash \Pbf_1 \tilde{w}_1 \Frob(\Pbf_1)/\Frob(\Ubf_{\Pbf_1})}{\Ad_{\Frob}\Lbf_1} \xrightarrow{i_{\tilde{t}_1, 0}} \frac{\Ubf_{\Pbf_1}\backslash \Pbf_0 w_0 \Frob(\Pbf_0)/\Frob(\Ubf_{\Pbf_1})}{\Ad_{\Frob}\Lbf_1}.$$
We can now compute (in the next computation we remove the indices $_{\Pbf_1 \subset \Pbf_0}$ from the symbols $\alpha$ and $\beta$) 
\begin{align*}
i_{\tilde{t}_1, 0}^*\hc^!_{\Pbf_1 \subset \Pbf_0} i_{t,1,!} &= i_{\tilde{t}_1, 0}^*\beta_!\alpha^* i_{t,1,!}[2(\dim \Ubf_{\Pbf_1} - \dim \Ubf_{\Pbf_0})]  \\
&\cong \lambda_{\tilde{t},1,!}i_{\tilde{t}, 1}^{'*}\alpha^*i_{t,1,!}[2(\dim \Ubf_{\Pbf_1} - \dim \Ubf_{\Pbf_0})] \\
&\cong \lambda_{\tilde{t},1,!}i_{\tilde{t}, 1}^*i_{t,1,!}[2(\dim \Ubf_{\Pbf_1} - \dim \Ubf_{\Pbf_0})] \\
&= 0.
\end{align*}
Let us give some details on the above calculation:
\begin{enumerate}
\item the first line is simply unfolding the definition, 
\item the second line follows from proper base change in the upper-right square of the above diagram with $\tilde{t}$ in place of $t$, 
\item the third line follows from the fact that $i_{\tilde{t}, 1}^{'*}\alpha^* = i_{\tilde{t}, 1}^*$, 
\item the last line follows from the fact that $i_{\tilde{t}, 1}^*i_{t,1,!} = 0$ as the two subvarieties $\Zbf_{I,t,1}$ and $\Zbf_{I,\tilde{t}, 1}$ of $\Zbf_{I}$ are disjoint.
\end{enumerate}
\end{proof}

\begin{lem}\label{lem:ConnectednessParabCS}
Let $t \in \Tcal_I$. Then the group $\Vbf_{\infty}$ of Lemma \ref{corol:StructureOfStrata} is connected. 
\end{lem}

\begin{proof}
It follows from the combination of Lemma \ref{lem:CommutationHCStrata} and Proposition \ref{prop:Conservativity} that the functor $\lambda_{t,1,!}$ is conservative. By induction, it follows that the functor $\lambda_{t,\infty,!}$ is conservative on the bounded-below part of the category, where $\lambda_{t, \infty}$ is the map $\Zbf_{I,t}/\Gbf \to \Zbf_{I_{\infty}, w_{\infty}}/\Gbf$. This implies that $\Vbf_{\infty}$ is connected. 
\end{proof}

\begin{rque}
This argument for the connectedness of the unipotent part of these stabilizers is somewhat indirect. It would be very interesting to have a more direct geometric argument.
\end{rque}

\begin{lem}\label{lem:ReductionToLevi}
Let $\Pbf = \Lbf\Ubf_{\Pbf}$ be the standard parabolic of type $I \subset \Delta$. Let $v \in {^\Pbf}\Wbf^{\Frob(\Pbf)}$ be such that $v\Frob(I) = I$, and let $\dot{v}$ be a lift of $v$. Let $w \in \Wbf_{\Pbf}$. Then there is an isomorphism of functors 
$$\ch_{\Bbf \subset \Pbf}^{\Gbf} i_{wv,!}^{\Bbf} = i_{v,!}^{\Pbf}\nu^*\ch^{\Lbf}_{\Bbf_{\Lbf} \subset \Lbf} i_{w,!}^{\Bbf_{\Lbf}}$$ 
where the exponents $^{\Lbf}$ and $^{\Gbf}$ indicate the ambient group in which these functors are computed and $\nu$ is the map 
$$\frac{\Ubf_{\Pbf} \backslash \Pbf v \Frob(\Pbf)/\Frob(\Ubf_{\Pbf})}{\Ad_{\Frob}(\Lbf)} \to \frac{\Lbf}{\Ad_{\Frob_v}(\Lbf)}.$$
\end{lem}

\begin{proof}
This is immediate by proper base change using the diagram of Lemma \ref{lem:DiagramReductionToLevi}.
\end{proof}

\begin{prop}\label{prop:Generation}
The essential image of the functor $\ch_{\Qbf \subset \Pbf}$ generates $\DD(\frac{\Ubf_{\Pbf} \backslash \Gbf/\Frob(\Ubf_{\Pbf})}{\Ad_{\Frob}\Lbf})$. 
\end{prop}

\begin{proof}
By transitivity, it is enough to prove that the essential image of $\ch_{\Bbf \subset \Pbf}$ generates the category $\DD(\frac{\Ubf_{\Pbf} \backslash \Gbf/\Frob(\Ubf_{\Pbf})}{\Ad_{\Frob}\Lbf})$. Let $t \in \Tcal_I$, let $(I_{\infty},w_{\infty})$ be its stable term, and let $a\in\Wbf_{\Pbf_{\infty}}$. By \Cref{lem:CombinatoricReduction}, $a\mapsto a w_{\infty}$ is a bijection from $\Wbf_{\Pbf_{\infty}}$ onto the double coset $\Wbf_{\Pbf_{\infty}}w_{\infty}\Frob(\Wbf_{\Pbf_{\infty}})$. By \Cref{lem:ReductionToLevi}, the collection of functors $\ch_{\Bbf \subset \Pbf_{\infty}}i_{a w_{\infty},!}$ identifies with Deligne--Lusztig induction on the stable Levi and therefore generates the category of sheaves supported on this stratum; see \cite[Théorème A]{BonnafeRouquier}. By transitivity (\Cref{prop:Transitivity}), the collection $\ch_{\Bbf \subset \Pbf}i_{a w_{\infty},!}$ generates the subcategory of sheaves supported on $\Zbf_{I,t}/\Gbf$. By \cite[Proposition 2.5]{LusztigParabolicCS1}, we have 
$$\Wbf = \sqcup_{t \in \Tcal_I} \Wbf_{\Pbf_{\infty}}w_{\infty}.$$
Hence, by varying the sequence $t$, we obtain that the category generated by the essential image of $\ch_{\Bbf \subset \Pbf}$ contains all sheaves of the form $i_{t,!}A$ for some sheaf on $\Zbf_{I,t}/\Gbf$. Since the collection of stacks $\Zbf_{I,t}/\Gbf$ stratifies $\frac{\Ubf_{\Pbf}\backslash \Gbf/\Frob(\Ubf_{\Pbf})}{\Ad_{\Frob}(\Lbf)}$, we get the desired statement. 
\end{proof}

\subsection{Mackey formulas and applications}

\begin{defi}[Intertwining operators]
Let $\Pbf, \Qbf \subset \Rbf$ be standard parabolic subgroups, and write their Levi decompositions as $\Pbf = \LP\UP, \Qbf = \LQ\Ubf_{\Qbf}, \Rbf = \LR\UR$. Let $z \in {^\Pbf}\Wbf_{\Rbf}^{\Qbf}$, and fix a lift $\dot{z} \in \Rbf$ of $z$. Consider the following correspondence:
\[\begin{tikzcd}
	& {\frac{\Ubf_{\Rbf} \backslash \Gbf /\Frob(\Ubf_{\Rbf})}{\Ad_{\Frob}(\LR \cap \Pbf \cap {{^{z^{-1}}}\Qbf})}} & {\frac{\Ubf_{\Rbf} \backslash \Gbf /\Frob(\Ubf_{\Rbf})}{\Ad_{\Frob}(\LR \cap {^z}\Pbf \cap \Qbf)}} \\
	{\frac{\Ubf_{\Pbf} \backslash \Gbf /\Frob(\Ubf_{\Pbf})}{\Ad_{\Frob}(\LP)}} &&& {\frac{\Ubf_{\Qbf} \backslash \Gbf /\Frob(\Ubf_{\Qbf})}{\Ad_{\Frob}(\LQ)}}
	\arrow["{\Ad_{\Frob}(\dot{z})}"', from=1-2, to=1-3]
	\arrow["{\alpha_z}", from=1-2, to=2-1]
	\arrow["{\beta_z}"', from=1-3, to=2-4]
\end{tikzcd}\]
We define the following functors $\DD(\frac{\Ubf_{\Pbf} \backslash \Gbf /\Frob(\Ubf_{\Pbf})}{\Ad_{\Frob}(\LP)}) \to \DD(\frac{\Ubf_{\Qbf} \backslash \Gbf /\Frob(\Ubf_{\Qbf})}{\Ad_{\Frob}(\LQ)})$
\begin{enumerate}
\item $\Int_{\Pbf, \Qbf,z}^{\Rbf, !} = \beta_{z,!}\Ad_{\Frob}(\dot z)_*\alpha_z^*$,
\item $\Int_{\Pbf, \Qbf,z}^{\Rbf, *} = \beta_{z,*}\Ad_{\Frob}(\dot z)_*\alpha_z^!$.
\end{enumerate}
\end{defi}

\begin{lem}[Mackey Formula]\label{lem:Mackey}
Let $\Pbf, \Qbf \subset \Rbf$ be standard parabolic subgroups. 
\begin{enumerate}
\item The functor $\hc^!_{\Qbf \subset \Rbf}  \ch_{\Pbf \subset \Rbf}$ has a filtration indexed by the poset ${^\Pbf}\Wbf_{\Rbf}^{\Qbf}$ whose graded pieces are the functors $\Int_{\Pbf, \Qbf,z}^{\Rbf, !}[2(\dim(\UQ) - \dim(\UR))]$.
\item The functor $\hc^*_{\Qbf \subset \Rbf}  \ch_{\Pbf \subset \Rbf}$ has a filtration indexed by the poset ${^\Pbf}\Wbf_{\Rbf}^{\Qbf}$ whose graded pieces are the functors $\Int_{\Pbf, \Qbf,z}^{\Rbf, *}$.
\end{enumerate}
\end{lem}

\begin{proof}
The argument is entirely similar for both $(i)$ and $(ii)$, so we will focus on the first one. We write the standard Levi decompositions of the parabolic subgroups as before: $\Pbf = \LP\Ubf_{\Pbf}, \Qbf = \LQ\Ubf_{\Qbf}, \Rbf = \LR\Ubf_{\Rbf}$. Consider the following diagram, where the middle square is defined to be Cartesian:
\[\begin{tikzcd}
	&& {\frac{\Ubf_{\Rbf} \backslash \Gbf /\Frob(\Ubf_{\Rbf})}{\Ad_{\Frob}(\LR)}} \\
	& {\frac{\Ubf_{\Rbf} \backslash \Gbf /\Frob(\Ubf_{\Rbf})}{\Ad_{\Frob}(\LR \cap \Pbf)}} && {\frac{\Ubf_{\Rbf} \backslash \Gbf /\Frob(\Ubf_{\Rbf})}{\Ad_{\Frob}(\LR \cap \Qbf)}} \\
	{\frac{\Ubf_{\Pbf} \backslash \Gbf /\Frob(\Ubf_{\Pbf})}{\Ad_{\Frob}(\LP)}} && \Xcal && {\frac{\Ubf_{\Qbf} \backslash \Gbf /\Frob(\Ubf_{\Qbf})}{\Ad_{\Frob}(\LQ)}}
	\arrow["{\alpha_{\Pbf \subset \Rbf}}"{description}, from=2-2, to=1-3]
	\arrow["{\beta_{\Pbf \subset \Rbf}}"{description}, from=2-2, to=3-1]
	\arrow["{\alpha_{\Qbf \subset \Rbf}}"{description}, from=2-4, to=1-3]
	\arrow["{\beta_{\Qbf \subset \Rbf}}"{description}, from=2-4, to=3-5]
	\arrow["\gamma"{description}, from=3-3, to=2-2]
	\arrow["\delta"{description}, from=3-3, to=2-4]
\end{tikzcd}\]

Let us describe the stack $\Xcal$. Consider $\Ycal$, the subscheme of $\Gbf/\Frob(\UR) \times \Gbf/\Frob(\UR) \times \Rbf$ of triples $(g\Frob(\UR), h\Frob(\UR), r)$ where
$$\Ad_{\Frob}(r)(g\Frob(\UR)) = h\Frob(\UR).$$
Consider the following action of $\Pbf \times \Qbf$ on $\Ycal$ given by 
$$(p,q).(g\Frob(\UR), h\Frob(\UR), r) = (\Ad_{\Frob}(p)(g\Frob(\UR)), \Ad_{\Frob}(q)(h\Frob(\UR)), prq^{-1}).$$
By definition of the fiber of stacks, we have 
$$\Xcal = \Ycal/(\Pbf \times \Qbf).$$
Consider the map 
$$\Xcal \to \Pbf \backslash \Rbf/\Qbf$$
given by $(g\Frob(\UR), h\Frob(\UR), r) \mapsto \Pbf r\Qbf$. We now equip $\Pbf \backslash \Rbf/\Qbf$ with the stratification given by relative positions:
$$\Pbf \backslash \Rbf/\Qbf = \bigsqcup_{z \in {^\Pbf}\Wbf_{\Rbf}^{\Qbf}} \Pbf z\Qbf.$$
We equip $\Xcal$ with the pullback of this stratification and given $z \in {^\Pbf}\Wbf_{\Rbf}^{\Qbf}$, we denote by $j_z : \Xcal_z \to \Xcal$ the inclusion of the corresponding stratum. Since the stratification is finite, any sheaf $A \in \DD(\Xcal)$ has two filtrations (one increasing and one decreasing) whose graded pieces are 
$j_{z,!}j_z^*A$ and $j_{z,*}j_z^!A$.

The composition $\hc^!_{\Qbf \subset \Rbf}  \ch_{\Pbf \subset \Rbf}$ is given by 
\begin{align*}
\hc^!_{\Qbf \subset \Rbf}  \ch_{\Pbf \subset \Rbf} &=  \beta_{\Qbf \subset \Rbf,!}\alpha_{\Qbf \subset \Rbf}^*\alpha_{\Pbf \subset \Rbf,!}\beta_{\Pbf \subset \Rbf}^*[2(\dim(\UQ) - \dim(\UR))] \\
& =\beta_{\Qbf \subset \Rbf,!}\delta_!\gamma^*\beta_{\Pbf \subset \Rbf}^*[2(\dim(\UQ) - \dim(\UR))].
\end{align*}
By the above discussion, this functor is filtered with graded pieces given by 
$$\beta_{\Qbf \subset \Rbf,!}\delta_!j_{z,!}j_z^*\gamma^*\beta_{\Pbf \subset \Rbf}^*[2(\dim(\UQ) - \dim(\UR))].$$
It remains to identify these pieces. We have 
$$\Xcal_z = \{g\Frob(\UR)\}/\Pbf \cap {^z}\Qbf = \frac{\Gbf}{\Rbf \cap \Pbf \cap {^z}\Qbf},$$
and the map $\delta j_z$ is identified with $\Ad_{\Frob}(\dot{z})$. 
\end{proof}

As an application, we give a short geometric proof of the following fundamental theorem of Deligne-Lusztig.

\begin{thm}[\cite{DeligneLusztig}]
Let $(w,\theta)$ be a pair where $w \in \Wbf$ and $\theta : \Tbf^{w\Frob} \to \Qlb^{\times}$ is a regular character. Then 
$$\ch_{\Bbf \subset \Gbf}(j_{w,!}\theta)[\ell(w)]$$
is an irreducible representation in degree $0$. 
\end{thm}

\begin{proof}
As $\theta$ is regular, the map $j_{w,!}\theta[\ell(w)] \to j_{w,*}\theta[\ell(w)]$ is an isomorphism by \cite[Proposition 5.2]{LusztigYun}. Since applying $\ch_{\Bbf \subset \Gbf}$ to these objects computes the cohomology of the corresponding Deligne-Lusztig varieties, arguing as in \cite{BonnafeRouquier}, we see that the object $\ch_{\Bbf \subset \Gbf}(j_{w,!}\theta)[\ell(w)]$ lies in degree $0$. To prove that it is an irreducible representation, it is enough to prove that $\End(\ch_{\Bbf \subset \Gbf}(j_{w,!}\theta))$ is reduced to $\Qlb$. Applying adjunctions and the Mackey formula (\ref{lem:Mackey}), it is enough to prove that 
$$\Hom(j_{w,!}\theta, \Int_{\Bbf, \Bbf,v}^{\Gbf}(j_{w,!}\theta))$$
is $0$ if $v \neq 1$. The object $\Int_{\Bbf, \Bbf,v}^{\Gbf}(j_{w,!}\theta)$ is monodromic for the action of $\Tbf \times \Tbf$ on $\frac{\Ubf \backslash \Gbf /\Ubf}{\Ad_{\Frob}(\Tbf)}$ given by left and right translation with generalized monodromy $(\Ad(v)(\theta), \Ad(w\Frob(v^{-1}))(\theta))$, while the sheaf $j_{w,!}(\theta)$ is monodromic for the generalized monodromy $(\theta, \Ad(w)(\theta))$. Since the left monodromies are different, we get the desired vanishing.
\end{proof}

\section{Jordan decomposition}\label{sectionJordan}

\subsection{Recollections on the categorical trace of the Hecke category}\label{sec:TraceOnCategoricalHeckeCat}

\begin{thm}[\protect{\cite[Theorem 6.1.1]{EteveFreeMonodromic}}]\label{thm:Trace}
There is a canonical equivalence 
$$\Tr(\Frob^*,\HH_{\tot}^{\Gbf}) \cong \DD(\pt/\Gbf^{\Frob}, \Lambda)$$ 
and the map $\tr_{\Frob} : \HH_{\tot}^{\Gbf} \to \Tr(\Frob^*,\HH_{\tot}^{\Gbf})$ is identified with the functor $\ch \pfrak_!$ (see Section \ref{sec:Definition-horocycle} for the notation).
\end{thm}

The following is an immediate consequence of \Cref{prop:ComparisonTraceAndDeligneLusztig}.

\begin{corol} \label{cor:TraceVsDeligneLusztig}
There are canonical isomorphisms
$$\tr_{\Frob}(\Delta_{w, s}) \cong e_{s}\RGamma_c(Y(\dot{w}), \Lambda)[\ell(w)]$$
and 
$$\tr_{\Frob}(\nabla_{w, s}) \cong e_{s}\RGamma(Y(\dot{w}), \Lambda)[\ell(w)]$$
where $s \in \Ch(\Tbf)$ is fixed by $w\Frob$ and $e_{s}$ denotes the idempotent in $\Lambda[\Tbf^{w\Frob}]$ onto the block containing $s$. 
\end{corol}

\subsubsection{}
Let us also recall the more general parabolic case. This case is based on Lusztig's theory of parabolic character sheaves \cite{LusztigParabolicCS1}, which is reviewed in Section \ref{sec:parabolic-character-sheaves}.

In terms of traces, we can consider as before the Hecke category $\HH^{\Lbf}$ for the Levi $\Lbf$ and the category $\HH^{\Gbf}_{\Frob}$ as an $\HH^{\Lbf}$-bimodule (where the right action is again twisted by $\Frob$). There is a natural morphism of pairs 
$$(\HH^{\Lbf}_{\tot},\HH^{\Gbf}_{\tot, \Frob}) \to (\HH^{\Gbf}_{\tot}, \HH^{\Gbf}_{\tot, \Frob}).$$
This morphism induces a morphism between categorical traces 
$$\Tr(\HH^{\Gbf}_{\tot,\Frob},\HH^{\Lbf}_{\tot}) \to \Tr(\HH^{\Gbf}_{\tot,\Frob},\HH^{\Gbf}_{\tot}).$$
\begin{thm}\label{thm:TracesAndParabolicCharacterSheaves}
There is a canonical isomorphism 
$$\Tr(\HH^{\Gbf}_{\tot,\Frob},\HH^{\Lbf}_{\tot}) \cong \DD(\frac{\UP \backslash \Gbf/\Frob(\UP)}{\Ad_{\Frob}(\Lbf)}, \Lambda)$$ 
and the induced map on traces makes the following diagram commutative
\[\begin{tikzcd}
	{\Tr(\HH^{\Gbf}_{\tot,\Frob},\HH^{\Lbf}_{\tot})} & {\Tr(\HH^{\Gbf}_{\tot,\Frob},\HH^{\Gbf}_{\tot})} \\
	{\DD(\frac{\UP \backslash \Gbf/\Frob(\UP)}{\Ad_{\Frob}(\Lbf)}, \Lambda)} & {\DD(\frac{\Gbf}{\Ad_{\Frob}\Gbf}, \Lambda).}
	\arrow[from=1-1, to=1-2]
	\arrow["\wr"', from=1-1, to=2-1]
	\arrow["\wr", from=1-2, to=2-2]
	\arrow["{{\ch_{\Pbf \subset \Gbf}}}"', from=2-1, to=2-2]
\end{tikzcd}\]
\end{thm}
The proof of this theorem is sketched in \Cref{thm:calculation-traces} and the argument is essentially the $\Frob$-twisted version of \cite{LiNadlerYun}. 

Let $\Pbf = \Lbf\Ubf_{\Pbf}$ be a standard parabolic as above and let $\dot{v}$ be a lift of an element $v \in {^\Pbf}\Wbf^{\Frob(\Pbf)}$ such that $\Ad(\dot{v})\circ \Frob(\Lbf) = \Lbf$. We denote by $\HH^{\Pbf, v}_{\tot}$ the full subcategory of $\HH^{\Gbf}_{\tot}$ of sheaves supported on $\Pbf v \Frob(\Pbf)$. 
\begin{corol}
There is a canonical equivalence $\Tr(\HH^{\Pbf, v}_{\tot,\Frob},\HH^{\Lbf}_{\tot}) \cong \DD(\frac{\UP \backslash \Pbf v \Frob(\Pbf)/\Frob(\UP)}{\Ad_{\Frob}\Lbf})$ making the following diagram commutative 
\[\begin{tikzcd}
	{\Tr(\HH^{\Pbf, v}_{\tot,\Frob},\HH^{\Lbf}_{\tot})} & {\Tr(\HH^{\Gbf}_{\tot,\Frob},\HH^{\Gbf}_{\tot})} \\
	{\DD(\frac{\UP \backslash \Pbf v \Frob(\Pbf)/\Frob(\UP)}{\Ad_{\Frob}\Lbf}, \Lambda)} & {\DD(\frac{\Gbf}{\Ad_{\Frob}\Gbf}, \Lambda).}
	\arrow[from=1-1, to=1-2]
	\arrow["\wr"', from=1-1, to=2-1]
	\arrow["\wr", from=1-2, to=2-2]
	\arrow["{{\ch_{\Pbf \subset \Gbf}}j_{v,!}}"', from=2-1, to=2-2]
\end{tikzcd}\]
\end{corol} 

\begin{rque}
By \Cref{prop:ComparisonTraceAndDeligneLusztig}, the functor $\ch_{\Pbf \subset \Gbf}j_{v,!}$ is isomorphic to the usual parabolic Deligne--Lusztig induction. 
\end{rque}

\subsection{$\mathcal{O}((\Tbf^{\vee}\sslash \Wbf)^{\Frob^{\vee}})$-linearity}

Recall that the Hecke category is naturally an object of $\Alg(\Pr_{\Ccal(\Tbf)\sslash \Wbf})$ by Lemma \ref{lem:LinearityAlgebra}. We consider the map 
of schemes 
$$\Ccal(\Tbf) \to \Tbf^{\vee}$$ 
introduced in Section \ref{subsect:Linearity}. This induces a map 
$$\Ccal(\Tbf) \sslash \Wbf \to \Tbf^{\vee} \sslash \Wbf.$$
We can forget part of the structure and view $\HH_{\tot}^{\Gbf}$ as an object of $\Alg(\Pr_{\Tbf^{\vee}\sslash \Wbf})$. Recall that we take the trace of Frobenius in $\Pr_{\Lambda}$, rather than in $\Pr_{\Tbf^{\vee}\sslash \Wbf}$. Since $\DD_{\qcoh}(\Tbf^{\vee}\sslash \Wbf)$ is symmetric monoidal, we deduce that
$$\Tr(\Frob^*,\HH_{\tot}^{\Gbf})$$ 
is naturally a module over
$$\Tr(\Frob^{\vee,*},\DD_{\qcoh}(\Tbf^{\vee}\sslash \Wbf)).$$

\begin{lem}[\protect{\cite{BenZviNadlerFrancis}, \cite[Example 3.7.3]{ToyModel}}]
There is a canonical isomorphism of symmetric monoidal categories 
$$\Tr(\Frob^{\vee,*},\DD_{\qcoh}(\Tbf^{\vee}\sslash \Wbf)) \cong \DD_{\qcoh}((\Tbf^{\vee}\sslash \Wbf)^{\Frob^{\vee}}).$$
\end{lem}
Here $(\Tbf^{\vee}\sslash \Wbf)^{\Frob^{\vee}}$ denotes the \emph{derived} scheme of $\Frob^{\vee}$-fixed points on $\Tbf^{\vee}\sslash\Wbf$. Thus the following diagram is derived Cartesian:
\[\begin{tikzcd}
	{(\Tbf^{\vee} \sslash \Wbf)^{\Frob^{\vee}}} & {\Tbf^{\vee} \sslash \Wbf} \\
	{\Tbf^{\vee} \sslash \Wbf} & {\Tbf^{\vee} \sslash \Wbf \times \Tbf^{\vee} \sslash \Wbf}
	\arrow[from=1-1, to=1-2]
	\arrow[from=1-1, to=2-1]
	\arrow["{\id \times \Frob^{\vee}}", from=1-2, to=2-2]
	\arrow["\Delta"', from=2-1, to=2-2]
\end{tikzcd}\]
where $\Delta$ is the diagonal map. Furthermore, the canonical trace map
$$\DD_{\qcoh}(\Tbf^{\vee}\sslash \Wbf) \to \Tr(\Frob^{\vee,*},\DD_{\qcoh}(\Tbf^{\vee}\sslash \Wbf)) \cong \DD_{\qcoh}((\Tbf^{\vee}\sslash \Wbf)^{\Frob^{\vee}})$$
is simply the pullback map along $(\Tbf^{\vee}\sslash \Wbf)^{\Frob^{\vee}} \to \Tbf^{\vee}\sslash \Wbf$.

\begin{rque}
Let us comment on the importance of symmetric monoidality here. A priori, the trace of a monoidal category $C$ does not have a monoidal structure. However, if $C$ is symmetric monoidal, then so is $\Tr(\Frob,C)$. For any algebra $(A,\Frob_A)$ over $(C,\Frob)$---that is, when the endomorphisms $\Frob_A$ and $\Frob$ are compatible---the trace $\Tr(\Frob_A,A)$ is linear over $\Tr(\Frob,C)$. We refer to \cite[Section 3.3]{ToyModel} and \cite[Section 4.2]{HoyoisScherotzkeSibilla} for a discussion.
\end{rque}

\begin{rque}
If $\Gbf$ has connected center then it follows that $(\Tbf^{\vee}\sslash \Wbf)^{\Frob^{\vee}}$ is reduced by \cite[Theorem 3.9]{Li}.
\end{rque}

\begin{thm}\label{thm:LinearityForReps}
Let $\Pbf = \Lbf\Ubf_{\Pbf}$ be a standard parabolic subgroup of $\Gbf$. The category $\DD(\frac{\Ubf_{\Pbf} \backslash \Gbf/\Frob(\Ubf_{\Pbf})}{\Ad_{\Frob}\Lbf}, \Lambda)$ is linear over $(\Tbf^{\vee}\sslash \Wbf)^{\Frob^{\vee}}$ and all functors $\ch_{\Pbf \subset \Gbf}$, $\hc^!_{\Pbf \subset \Gbf}$ and $\hc^*_{\Pbf \subset \Gbf}$ are $(\Tbf^{\vee}\sslash \Wbf)^{\Frob^{\vee}}$-linear.
\end{thm}

\begin{proof}
Consider the map of algebras and bimodules $(\HH^{\Lbf}_{\tot}, \HH_{\tot,\Frob}^{\Gbf}) \to (\HH_{\tot}^{\Gbf}, \HH_{\tot,\Frob}^{\Gbf})$, where $\HH_{\tot,\Frob}^{\Gbf}$ again denotes the twist of the right module structure by $\Frob^*$. This map induces the following map on traces:
$$\Tr(\HH_{\tot,\Frob}^{\Gbf},\HH^{\Lbf}_{\tot}) \to \Tr(\Frob^*,\HH_{\tot}^{\Gbf}).$$
By \Cref{thm:TracesAndParabolicCharacterSheaves}, there is an equivalence 
$$\Tr(\HH_{\tot,\Frob}^{\Gbf},\HH^{\Lbf}_{\tot}) \cong \DD(\frac{\Ubf_{\Pbf} \backslash \Gbf/\Frob(\Ubf_{\Pbf})}{\Ad_{\Frob}\Lbf}, \Lambda)$$ 
and the map induced on traces is $\ch_{\Pbf \subset \Gbf}$ .

The algebra $\HH^{\Lbf}_{\tot}$ is an algebra in $\Pr_{\Tbf^{\vee}\sslash \Wbf}$ and the map of bimodules $(\HH^{\Lbf}_{\tot}, \HH_{\tot,\Frob}^{\Gbf}) \to (\HH_{\tot}^{\Gbf}, \HH_{\tot,\Frob}^{\Gbf})$ is $\Tbf^{\vee}\sslash \Wbf$-linear. It follows that the map of traces is then linear over $(\Tbf^{\vee}\sslash \Wbf)^{\Frob^{\vee}}$. This shows the linearity statement for $\ch_{\Pbf \subset \Gbf}$. Since $\hc^*_{\Pbf \subset \Gbf}$ is right adjoint to $\ch_{\Pbf \subset \Gbf}$, it is lax-linear over $(\Tbf^{\vee}\sslash \Wbf)^{\Frob^{\vee}}$ and since $\DD_{\qcoh}((\Tbf^{\vee}\sslash \Wbf)^{\Frob^{\vee}})$ is rigid it follows that $\hc^*_{\Pbf \subset \Gbf}$ is $(\Tbf^{\vee}\sslash \Wbf)^{\Frob^{\vee}}$ linear. The statement for $\hc^!_{\Pbf \subset \Gbf}$ follows from duality.
\end{proof}

\subsection{Geometric series}\label{sec:GeometricSeries}

Recall that we have identified the following sets in Section \ref{sec:SemisimplePoints}
$$(\Tbf^{\vee}\sslash \Wbf)^{\Frob^{\vee}}_{\tors, \Lambda}(\Lambda) \cong \Gbf^*_{\mathrm{ss}}/\sim_{\Lambda, \geo}.$$
More precisely, the points $(\Tbf^{\vee}\sslash \Wbf)^{\Frob^{\vee}}_{\tors, \Lambda}(\Lambda)$ enumerate the connected components of the classical truncation $t_0((\Tbf^{\vee}\sslash \Wbf)^{\Frob^{\vee}})$. In particular, we shall identify these connected components with $\Wbf$-orbits in $\Ch(\Tbf)$.
The linearity over $(\Tbf^{\vee}\sslash \Wbf)^{\Frob^{\vee}}$ of $\DD(\pt/\Gbf^{\Frob}, \Lambda) \cong \Tr(\Frob^*,\HH^{\Gbf}_{\tot})$ corresponds to a map 
$$\Ocal((\Tbf^{\vee}\sslash \Wbf)^{\Frob^{\vee}}) \to \End(\id_{\DD(\pt/\Gbf^{\Frob}, \Lambda)}) \cong Z(\Lambda[\Gbf^{\Frob}]).$$
Given a connected component $\ofrak$ of $(\Tbf^{\vee}\sslash \Wbf)^{\Frob^{\vee}}$, we obtain an idempotent of $\Ocal((\Tbf^{\vee}\sslash \Wbf)^{\Frob^{\vee}})$, namely the characteristic function of this component. Its image in $Z(\Lambda[\Gbf^{\Frob}])$ yields a central idempotent, which we denote by $e_{\ofrak, \Lambda}$. The collection $(e_{\ofrak, \Lambda})$ is a complete collection of orthogonal central idempotents of $\Lambda[\Gbf^{\Frob}]$.

\begin{lem}[Decomposition in geometric series]\label{lem:GeometricSeries}
There is a direct sum decomposition 
$$\DD(\pt/\Gbf^{\Frob}, \Lambda) = \bigoplus_{\ofrak \in \Gbf^*_{\mathrm{ss}}/\sim_{\Lambda, \geo}} \DD^{\ofrak}(\pt/\Gbf^{\Frob}, \Lambda)$$
where $\DD^{\ofrak}(\pt/\Gbf^{\Frob}, \Lambda)$ is the direct summand cut out by the idempotent $e_{\ofrak, \Lambda}$. They satisfy the following properties.
\begin{enumerate}
\item Let $(w, s)$ be a pair such that $w\Frob(s) = s$ and $s \in \ofrak$, the object $\tr_{\Frob}(\Delta_{w,s})$ lies in the summand $\DD^{\ofrak}(\pt/\Gbf^{\Frob}, \Lambda)$ of $\DD(\pt/\Gbf^{\Frob}, \Lambda)$.
\item There is an equality in $\Qlb[\Gbf^{\Frob}]$, for all conjugacy classes $\ofrak \in \Gbf^*/{\sim_{\Zlb, \geo}}$,
$$e_{\ofrak, \Zlb} = \sum_{\ofrak' \sim_{\Qlb, \geo} \ofrak} e_{\ofrak', \Qlb}.$$ 
\end{enumerate}
Furthermore Property $(i)$ characterizes the summands $\DD^{\ofrak}(\pt/\Gbf^{\Frob}, \Lambda)$ uniquely. 
\end{lem}

\begin{proof}
We get a direct sum decomposition since the idempotents $(e_{\ofrak, \Lambda})$ satisfy 
$$\sum_{\ofrak \in \Gbf^*/\sim_{\Lambda, \geo}} e_{\ofrak,\Lambda} = 1.$$
Point $(i)$ follows from the $(\Tbf^{\vee}\sslash \Wbf)^{\Frob^{\vee}}$-linearity of the functor $\ch$. The fact that they characterize the idempotents follows from the fact that the cohomology complexes of the Deligne--Lusztig varieties generate $\DD(\Rep_{\Lambda}(\Gbf^{\Frob}))$. Point $(ii)$ follows from the comparison of the connected components of $(\Tbf^{\vee}\sslash \Wbf)^{\Frob^{\vee}}$ over $\Qlb$ and $\Zlb$. 
\end{proof}

\begin{rque}
Lemma \ref{lem:GeometricSeries} is not new. Over $\Qlb$, the decomposition comes from the partition of $\Irr_{\Qlb}(\Gbf^{\Frob})$ into geometric Lusztig series constructed in \cite{DeligneLusztig}. Over $\Zlb$ and $\Flb$, this statement was proved in \cite{BroueMichel} by considering the sums of idempotents in $\Qlb[\Gbf^{\Frob}]$ appearing in $(ii)$ and showing that they are in fact defined over $\Zlb$.
\end{rque}

\begin{lem}[Trace of a permuted direct sum]\label{lem:trace-permuted-summands}
Let $\Acal=\bigoplus_{i\in I}\Acal_i$ be a finite direct sum of pairwise orthogonal monoidal summands, and let $F$ be a monoidal endofunctor permuting them. Then
$$\Tr(F,\Acal)\simeq\bigoplus_{i\,:\,F(i)=i}\Tr(F|_{\Acal_i},\Acal_i).$$
\end{lem}

\begin{rque}\label{rk:Series}
We have constructed this decomposition using linearity over $(\Tbf^{\vee}\sslash \Wbf)^{\Frob^{\vee}}$. Another way to construct it is to observe that there is a direct sum decomposition 
$$\HH^{\Gbf}_{\tot} \colondef= \bigoplus_{\mathfrak{o}} \HH_{\mathfrak{o}_{\tot}}^{\Gbf}$$ 
of monoidal categories where $\mathfrak{o}$ is a $\Wbf$-orbit in $\Ch(\Tbf)$ and $\HH_{\mathfrak{o},\tot}^{\Gbf} = \oplus_{s, s' \in \mathfrak{o}} \HH_{[s,s']}^{\Gbf}$. Let us further group them according to the orbits of $\Frob$ on the set of $\Wbf$-orbits $(\mathfrak{o})$, namely,
$$\HH^{\Gbf}_{\tot} \colondef= \bigoplus_{\Omega} \HH_{\Omega}^{\Gbf}$$
where $\Omega$ is an $\Frob$-orbit of $\Wbf$-orbits and $\HH_{\Omega}^{\Gbf} \colondef= \bigoplus_{\mathfrak{o} \in \Omega} \HH_{\mathfrak{o}}^{\Gbf}$. Each monoidal category $\HH_{\Omega}^{\Gbf}$ is stable under $\Frob^*$. The first assertion follows from \Cref{lem:trace-permuted-summands}, and the second follows by applying the same decomposition to all $\Frob$-stable summands; see also \cite[Section 3.5.6]{EteveThesis}:
\begin{enumerate}
\item $\Tr(\Frob^*,\HH_{\Omega}^{\Gbf}) = 0$ if $\Omega$ is not reduced to a single element, 
\item Assume that $\mathfrak{o}$ is $\Frob$-stable so that $\HH_{\mathfrak{o},\tot}^{\Gbf}$ is stable under $\Frob^*$, and let $\ofrak \in \Gbf^*_{\mathrm{ss}}/\sim_{\Lambda, \geo}$ be the corresponding $\Gbf^*$-orbit of semisimple elements, then 
$$\Tr(\Frob^*,\HH_{\mathfrak{o},\tot}^{\Gbf}) \cong \DD^{\ofrak}(\pt/\Gbf^{\Frob}, \Lambda).$$
\end{enumerate}
\end{rque}

\subsection{Rational series}\label{sec:RationalSeries}

If the center of $\Gbf$ is disconnected, there is a discrepancy between the rational and geometric conjugacy classes of semisimple elements in $(\Gbf^*)^{\Frob^*}$. This gives rise to the theory of rational series, which refines that of geometric series. We refer to \cite{BroueMichel} and \cite{BonnafeRouquier} for the classical approach. We fix an $\Frob$-stable $\Wbf$-orbit $\mathfrak{o}$ in $\Ch(\Tbf)$. Recall that $\Xibf_{\mathfrak{o}}$ denotes the connected component of $\Xibf$ corresponding to $\mathfrak{o}$.

Recall that the Hecke category $\HH_{\ofrak,\tot}^{\Gbf}$ is the totalization of the functor $\HH^{\Gbf}_{\ofrak} : \Xibf_{\ofrak}^{\opp} \to \Morita$. Recall also that $\Xibf_{\ofrak}$ is viewed as a finite algebraic stack over $\Fqb$. As $\ofrak$ is stable under $\Frob$, so is $\Xibf_{\ofrak}$ and we can form its stack of Frobenius fixed points defined as the following pullback 
\[\begin{tikzcd}
	{\Lcal_{\Frob}\Xibf_{\ofrak}} & {\Xibf_{\ofrak}} \\
	{\Xibf_{\ofrak}} & {\Xibf_{\ofrak} \times \Xibf_{\ofrak}.}
	\arrow[from=1-1, to=1-2]
	\arrow[from=1-1, to=2-1]
	\arrow["{\id \times \Frob}", from=1-2, to=2-2]
	\arrow["\Delta"', from=2-1, to=2-2]
\end{tikzcd}\]
It is sometimes also called the $\Frob$-loop stack (see for instance \cite{BenZviNadlerCharacterTheory}) hence our notation. It is still a finite stack over $\Fqb$ described as follows 
\begin{enumerate}
\item its objects are pairs $(\beta, s)$, where $\beta : s \to \Frob(s)$ is a map in $\Xibf_{\ofrak}$, 
\item its morphisms from $(\beta_1,s_1)$ to $(\beta_2,s_2)$ are diagrams
\[\begin{tikzcd}
	{s_1} & {\Frob(s_1)} \\
	{s_2} & {\Frob(s_2)}
	\arrow["{\beta_1}", from=1-1, to=1-2]
	\arrow["\alpha"', from=1-1, to=2-1]
	\arrow["{\Frob(\alpha)}", from=1-2, to=2-2]
	\arrow["{\beta_2}"', from=2-1, to=2-2]
\end{tikzcd}\]\\
where $\alpha : s_1 \to s_2$ is a morphism in $\Xibf_{\ofrak}$. 
\end{enumerate}

\subsubsection{}
We will apply the construction of \Cref{sec:local-traces}. 
The Frobenius equivariance of the functor $\HH_{\ofrak}^{\Gbf} : \Xibf_{\ofrak}^{\opp} \to \Morita$ yields a functor 
$$\Tr^{\loc}(\Frob^*,\HH_{\ofrak}^{\Gbf}) : \Lcal_{\Frob}\Xibf_{\ofrak}^{\opp} \to \mathrm{Pr}_{\Lambda}$$
defined as follows 
\begin{enumerate}
\item an object $(\beta, s)$ is sent to the category
$$\Tr^{\loc}(\Frob^*,\HH_{\ofrak}^{\Gbf})(s, \beta) :=\Tr(\HH^{\Gbf, \beta}_{[s,\Frob(s)]},\HH^{\Gbf, \circ}_{[s,s]}),$$ 
\item a morphism from $(s_1,\beta_1)$ to $(s_2,\beta_2)$ in $\Lcal_{\Frob}\Xibf_{\ofrak}^{\opp}$ is represented by a morphism $\alpha:(s_2,\beta_2)\to(s_1,\beta_1)$ in $\Lcal_{\Frob}\Xibf_{\ofrak}$; thus $\alpha:s_2\to s_1$ and $\beta_1\alpha=\Frob(\alpha)\beta_2$. The block $\HH^{\Gbf,\alpha}_{[s_2,s_1]}$ is an $(\HH^{\Gbf,\circ}_{[s_2,s_2]},\HH^{\Gbf,\circ}_{[s_1,s_1]})$-bimodule, and convolution gives the structure isomorphism
$$\HH^{\Gbf,\alpha}_{[s_2,s_1]}\otimes_{\HH^{\Gbf,\circ}_{[s_1,s_1]}}\HH^{\Gbf,\beta_1}_{[s_1,\Frob(s_1)],\Frob}
\xrightarrow{\ \sim\ }
\HH^{\Gbf,\beta_2}_{[s_2,\Frob(s_2)],\Frob}\otimes_{\HH^{\Gbf,\circ}_{[s_2,s_2]}}\HH^{\Gbf,\alpha}_{[s_2,s_1]}.$$
Here the subscript $\Frob$ means that the right diagonal action is transported by $\Frob^*$; on underlying blocks the right-hand tensor product is convolution with $\HH^{\Gbf,\Frob(\alpha)}_{[\Frob(s_2),\Frob(s_1)]}$. This is the specialization of \Cref{prop:functoriality-local-traces} and yields a map
$$\Tr^{\loc}(\Frob^*,\HH_{\ofrak}^{\Gbf})(s_1, \beta_1) \xrightarrow{\Tr^{\loc}(\alpha)} \Tr^{\loc}(\Frob^*,\HH_{\ofrak}^{\Gbf})(s_2, \beta_2).$$
\end{enumerate} 

We now relate the functor $\Tr^{\loc}$ to the category $\Tr(\HH_{\ofrak, \Frob}^{\Gbf},\HH_{\ofrak}^{\Gbf})$. For any $(s,\beta) \in \Lcal_{\Frob}(\Xibf_{\ofrak})$, there is a morphism of pairs
$$(\HH^{\Gbf,\circ}_{[s,s]}, \HH^{\Gbf,\beta}_{[s,\Frob(s)]}) \rightarrow (\HH_{\ofrak, \tot}^{\Gbf}, \HH_{\ofrak, \tot, \Frob}^{\Gbf})$$
induced by the obvious inclusions. Corresponding to this morphism, we get a map 
$$\incl_{(s, \beta)} : \Tr(\HH^{\Gbf,\beta}_{[s,\Frob(s)]},\HH^{\Gbf,\circ}_{[s,s]}) \to \Tr(\HH_{\ofrak, \tot, \Frob}^{\Gbf},\HH_{\ofrak, \tot}^{\Gbf}).$$
For any map $\alpha : (s_1, \beta_1) \to (s_2, \beta_2)$, the natural triangle 
\[\begin{tikzcd}
	{(\HH^{\Gbf,\circ}_{[s_1,s_1]}, \HH^{\Gbf,\beta_1}_{[s_1,\Frob(s_1)]})} & {(\HH_{\ofrak, \tot}^{\Gbf}, \HH_{\ofrak, \tot, \Frob}^{\Gbf})} \\
	{(\HH^{\Gbf,\circ}_{[s_2,s_2]}, \HH^{\Gbf,\beta_2}_{[s_2,\Frob(s_2)]})}
	\arrow[from=1-1, to=1-2]
	\arrow["{\HH^{\alpha}_{[s_1,s_2]}}"', from=1-1, to=2-1]
	\arrow[from=2-1, to=1-2]
\end{tikzcd}\]
is naturally commutative. 
It follows that there is a natural map 
\begin{equation}\label{eq:natural-transfo-colimit-xibf}
\varinjlim_{\Lcal_{\Frob}\Xibf_{\ofrak}^{\opp}} \Tr^{\loc} \to  \Tr(\HH_{\ofrak, \tot, \Frob}^{\Gbf},\HH_{\ofrak, \tot}^{\Gbf}).
\end{equation}

\begin{thm}[\Cref{thm:local-to-global-trace}]\label{thm:local-global-rational-series}
The map \Cref{eq:natural-transfo-colimit-xibf} is an equivalence. 
\end{thm}

\subsubsection{}
It follows from \Cref{thm:local-global-rational-series} that the category $\Tr(\HH_{\ofrak, \tot, \Frob},\HH_{\ofrak, \tot})$ naturally carries an action of $\DD(\pt)_{\Lcal_{\Frob}\Xibf_{\ofrak}} \cong \DD(\Lcal_{\Frob}\Xibf_{\ofrak})$. Recall from \Cref{lem:RationalConjClasses} that the connected components of $\Lcal_{\Frob}\Xibf_{\ofrak}$ are in canonical bijection with the $(\Gbf^*)^{\Frob^*}$-orbits in the $\Gbf^*$-orbit containing $\ofrak$. Hence, for any such $(\Gbf^*)^{\Frob^*}$-orbit, denoted $\ofrak_{\rat}'$, there is a canonical idempotent $e_{\ofrak_{\rat}', \Lambda} \in Z(\Lambda[\Gbf^{\Frob}])$ and a corresponding direct summand
$$\DD^{\ofrak_{\rat}'}(\pt/\Gbf^{\Frob}) \subset \DD^{\ofrak}(\pt/\Gbf^{\Frob}).$$

\begin{lem}\label{lem:rational-series-structure}
The idempotents $e_{\ofrak_{\rat}', \Lambda}$ satisfy the following properties : 
\begin{enumerate}
\item $\sum_{\ofrak_{\rat}'} e_{\ofrak_{\rat}', \Lambda} = e_{\ofrak, \Lambda}$,
\item Let $\ofrak_{\Zlb} \subset \Ch_{\Zlb}(\Tbf)$ be a $\Frob$-stable $\Wbf$-orbit and let $\ofrak_{\rat, \Zlb}' \subset \ofrak_{\Zlb}$ be a rational orbit, then we have 
$$\sum_{\ofrak_{\rat, \Qlb} \sim_{\rat, \Zlb} \ofrak_{\rat,\Zlb}'} e_{\ofrak_{\rat, \Qlb}, \Qlb} = e_{\ofrak_{\rat, \Zlb}', \Zlb}.$$
\item The category $\DD^{\ofrak_{\rat}'}(\pt/\Gbf^{\Frob})$ is generated by the objects $\tr_{\Frob}(\Delta_{w,\Frob(s)})$ where $w\Frob(s) = s$ and $w$ lies in a block $\beta \in {_s}\underline{\Wbf}_{\Frob(s)}$ that corresponds to $\ofrak_{\rat}'$ under \Cref{lem:RationalConjClasses}.
\end{enumerate}
\end{lem}

\begin{proof}
Follows from unraveling the definitions.
\end{proof}

\begin{rque}
This construction recovers the theory of rational series originally constructed in \cite{DeligneLusztig} for characteristic $0$ coefficients and \cite[Theorem 2.2]{BroueMichel} for modular coefficients. 
\end{rque}

\subsection{Series and Deligne-Lusztig induction}\label{sec:series-and-induction}

We also construct geometric and rational series for parabolic character sheaves and discuss their compatibility with the horocycle and character functor. For this entire subsection, we fix $\Pbf = \Lbf\Ubf_{\Pbf}$ a standard parabolic,  $v \in {^\Pbf}\Wbf^{\Frob(\Pbf)}$ and a lift $\dot{v}$ satisfying $\Ad(\dot{v}) \circ \Frob(\Lbf) = \Lbf$. In particular, we are going to reprove that rational series commutes with Deligne--Lusztig induction, again see \cite{DeligneLusztig} (for the $\Qlb$-version) and \cite{BroueMichel} for the modular version. 

Recall that the pair $(\Lbf, \dot{v})$ determines a well-defined Deligne--Lusztig induction functor
$$\RR_{\Lbf, \dot{v}} : \DD(\pt/\Lbf^{\dot{v}\Frob}) \to \DD(\pt/\Gbf^{\Frob}).$$

\begin{thm}[Compatibility of Deligne--Lusztig induction with series]\label{thm:compat-DL-and-series}
\begin{enumerate}
\item The functor $\RR_{\Lbf, \dot{v}}$ is linear over $(\Tbf^{\vee}\sslash \Wbf_{\Lbf})^{v\Frob^{\vee}} \to (\Tbf^{\vee}\sslash \Wbf)^{\Frob^{\vee}}$. 
\item Let $\ofrak \subset \Ch(\Tbf)$ be a $\Frob$-stable $\Wbf$-orbit and let $\ofrak_{\Lbf} \subset \ofrak$ be a $\Wbf_{\Lbf}$-orbit that is stable under $\Ad(\dot{v})\Frob$, then the induced functor 
$$\DD^{\ofrak_{\Lbf}}(\pt/\Lbf^{\dot{v}\Frob}) \to \DD^{\ofrak}(\pt/\Gbf^{\Frob})$$
is linear over $\Lcal_{v\Frob}\Xibf_{\ofrak_{\Lbf}}^{\Lbf} \xrightarrow{\cdot v} \Lcal_{\Frob}\Xibf_{\ofrak}$ (see \Cref{constr:translation-by-v}).
\end{enumerate}
\end{thm}

\begin{lem}\label{lem:coset-group-vs-Levi}
The subset $\Wbf^{\Lbf}v\Frob(\Wbf^{\Lbf}) \cap {_s}\Wbf_{\Frob(s)}$ of $\Wbf$ is a single $\Wbf^{\Lbf}_s$-left coset and, equivalently, a single right coset for $\Frob(\Wbf_s^{\Lbf})$.
\end{lem}

\begin{proof}
This is a direct calculation. We have 
\begin{align*}
\Wbf^{\Lbf}v\Frob(\Wbf^{\Lbf}) \cap {_s}\Wbf_{\Frob(s)} &= \Wbf^{\Lbf}v\Frob(\Wbf^{\Lbf})v^{-1}v \cap {_s}\Wbf_{\Frob(s)} \\
&= \Wbf^{\Lbf}\cdot v \cap {_s}\Wbf_{\Frob(s)} \\
&= (\Wbf^{\Lbf} \cap {_s}\Wbf_{v\Frob(s)})\cdot v \\
&= {_s}\Wbf_{v\Frob(s)}^{\Lbf}\cdot v.
\end{align*}
\end{proof}

\begin{construction}\label{constr:translation-by-v}
It follows from Lemma \ref{lem:coset-group-vs-Levi} that there is a well-defined map $\Lcal_{v\Frob}\Xibf_{\ofrak^{\Lbf}}^{\Lbf} \to \Lcal_{\Frob}\Xibf_{\ofrak}$ given by
\begin{align*}
(s, \beta) &\mapsto (s, \beta\cdot v) \\
\alpha &\mapsto \alpha.
\end{align*}
\end{construction} 

The map of \Cref{constr:translation-by-v} can be obtained in the following way, consider the diagram of groupoids 
\begin{equation}\label{diag:morphism-of-xi}
\begin{tikzcd}
	{\Xibf_{\ofrak^{\Lbf}}^{\Lbf}} & {\Xibf_{\ofrak}} \\
	{\Xibf_{\ofrak^{\Lbf}}^{\Lbf}} & {\Xibf_{\ofrak}}
	\arrow["{\mathrm{incl}}", from=1-1, to=1-2]
	\arrow["{\Ad(v) \Frob}"', from=1-1, to=2-1]
	\arrow["{\Ad(v)}"{description}, Rightarrow, from=1-2, to=2-1]
	\arrow["\Frob", from=1-2, to=2-2]
	\arrow["{\mathrm{incl}}", from=2-1, to=2-2]
\end{tikzcd}
\end{equation}
where $\Ad(v)$ is viewed in the middle as a natural transformation making the diagram $2$-commutative. Passing to loop space along the vertical maps yields the map of \Cref{constr:translation-by-v}. 

\begin{lem}\label{lem:identification-twist-by-delta-v}
There is a natural isomorphism of functors $\HH^{\Lbf}_{\ofrak^{\Lbf}, \tot} \to \HH^{\Ad(v)\Lbf}_{\Ad(v)\ofrak^{\Lbf}, \tot}$
$$\Ad(v)^* \cong \Ad(\Delta_v)$$
where $\Delta_v = \oplus_{s \in \ofrak} \Delta_{\dot{v}, s}$ and the functor $\Ad(\Delta_v)$ is computed in $\HH^{\Gbf}_{\tot}$. 
\end{lem}

\begin{proof}
Since both functors $\Ad(v)^*$ and $\Ad(\Delta_v)$ send tilting sheaves to tilting sheaves, it is enough to prove it for tilting sheaves by formality (see \Cref{thm:StructureCategoryviaTilting}). This can then be verified on Soergel bimodules. As $\V(\Delta_v) \cong \oplus_{s \in \ofrak} \RR_{\Tbf, v}$ is the $\RR_{\Tbf}$-bimodule whose right action is twisted by $\Ad(v)$, the lemma follows from the fact that on $\RR_{\Tbf}$-bimodules, there is an isomorphism of functors:  
$$\RR_{\Tbf, v} \otimes_{\RR_{\Tbf}} - \otimes_{\RR_{\Tbf}} \RR_{\Tbf, v^{-1}} \cong \Ad(v).$$
\end{proof}

\begin{rque}
Applying \Cref{lem:identification-twist-by-delta-v} to the case $\Lbf = \Gbf$ also yields that $\Ad(\dot{v})^* \cong \Ad(\Delta_v)$ as functors $\HH^{\Gbf}_{\ofrak, \tot} \to \HH^{\Gbf}_{\ofrak, \tot}$ for any element $w \in \Wbf$. 
\end{rque}

\begin{lem}\label{lem:comparison-bimodules}
The map 
$$\HH^{\Lbf}_{\ofrak^{\Lbf}, \tot, v\Frob} \xrightarrow{\cdot \Delta_v} \HH^{\Pbf,v}_{\ofrak^{\Lbf}, \tot, \Frob}$$
is an isomorphism of $\HH^{\Lbf}_{\ofrak^{\Lbf}, \tot}$-bimodules, where $\HH^{\Pbf,v}_{\ofrak^{\Lbf}, \tot} = \bigoplus_{s,s' \in \ofrak^{\Lbf}} \HH^{\Pbf,v}_{[s,v\Frob(s')]}$ and $\cdot \Delta_v$ is the multiplication in $\HH^{\Gbf}_{\tot}$ on the right by the element $\Delta_{\dot{v}} = \oplus_{s \in \ofrak} \Delta_{\dot{v},s}$.
\end{lem}

\begin{proof}
Since $\Delta_v$ is invertible in $\HH^{\Gbf}_{\ofrak, \tot}$ this functor is invertible. The fact that it is compatible with the left and right actions of $\HH^{\Lbf}_{\ofrak^{\Lbf},\tot}$ is immediate in view of Lemma \ref{lem:identification-twist-by-delta-v}.
\end{proof}

\begin{lem}
There is a natural commutative diagram of functors from Diagram \ref{diag:morphism-of-xi} to $\Morita$ given by 
\begin{equation}\label{diag:mapping-into-morita-diag}
\begin{tikzcd}
	{\HH^{\Lbf}_{\ofrak^{\Lbf}}} & {\HH^{\Gbf}_{\ofrak^{\Gbf}}} \\
	{\HH^{\Lbf}_{\ofrak^{\Lbf}}} & {\HH^{\Gbf}_{\ofrak^{\Gbf}}}
	\arrow["{\mathrm{\incl}}", from=1-1, to=1-2]
	\arrow["{\Ad(\dot{v}\Frob)^*}"', from=1-1, to=2-1]
	\arrow["{\Frob^*}", from=1-2, to=2-2]
	\arrow["{\cdot \Delta_{\dot{v}}}"{description}, Rightarrow, from=2-1, to=1-2]
	\arrow["{\mathrm{\incl}}"', from=2-1, to=2-2]
\end{tikzcd}
\end{equation}
where $\incl$ is induced by the natural inclusion $\HH^{\Lbf}_{\ofrak^{\Lbf}, \tot} \to \HH^{\Gbf}_{\ofrak, \tot}$, and $\cdot \Delta_v$ is induced by multiplication on the right by $\Delta_v = \oplus_{s \in \ofrak} \Delta_{\dot{v}, s}$. 
\end{lem}

\begin{proof}
In view of Lemma \ref{lem:comparison-bimodules}, the diagram \ref{diag:mapping-into-morita-diag} is isomorphic to the diagram 
\[\begin{tikzcd}
	{\HH^{\Lbf}_{\ofrak^{\Lbf}}} & {\HH^{\Gbf}_{\ofrak^{\Gbf}}} \\
	{\HH^{\Lbf}_{\ofrak^{\Lbf}}} & {\HH^{\Gbf}_{\ofrak^{\Gbf}}}
	\arrow["{\mathrm{\incl}}", from=1-1, to=1-2]
	\arrow["{\HH^{\Pbf,v}_{\Frob}}"', from=1-1, to=2-1]
	\arrow["{\Frob^*}", from=1-2, to=2-2]
	\arrow["\incl"{description}, Rightarrow, from=2-1, to=1-2]
	\arrow["{\mathrm{\incl}}"', from=2-1, to=2-2]
\end{tikzcd}\]
where $\HH^{\Pbf,v}_{\Frob}$ is viewed as natural transformation 
\[\begin{tikzcd}
	{\Xibf_{\ofrak^{\Lbf}}^{\Lbf}} && \Morita \\
	\\
	{\Xibf_{\ofrak^{\Lbf}}^{\Lbf}} && {\Morita.}
	\arrow["{\HH_{\ofrak^{\Lbf}}^{\Lbf}}", from=1-1, to=1-3]
	\arrow["{\Ad(v)\Frob}"', from=1-1, to=3-1]
	\arrow["{\HH^{\Pbf,v}_{\ofrak^{\Lbf}}}"{description},	 Rightarrow, from=1-3, to=3-1]
	\arrow[equals, from=1-3, to=3-3]
	\arrow["{\HH_{\ofrak^{\Lbf}}^{\Lbf}}"', from=3-1, to=3-3]
\end{tikzcd}\]

\end{proof}

\begin{proof}[Proof of \Cref{thm:compat-DL-and-series}]
Passing the diagram \ref{diag:mapping-into-morita-diag} to totalizations yields a morphism of pairs 
$$(\HH^{\Lbf}_{\ofrak^{\Lbf}, \tot}, \HH^{\Lbf}_{\ofrak^{\Lbf}, \tot, \Ad(\dot{v})\Frob}) \xrightarrow{(\mathrm{\incl}, - \cdot \Delta_v)} (\HH^{\Gbf}_{\ofrak, \tot}, \HH^{\Gbf}_{\ofrak, \tot, \Frob}).$$
By \Cref{lem:LinearityAlgebra}, this map of pairs is a map of pairs in $\mathrm{Pr}_{(\Tbf^{\vee}\sslash \Wbf)}$. The same argument as in \Cref{lem:GeometricSeries} establishes point $(i)$ of the theorem.

For the second point of the theorem, consider again diagram \ref{diag:mapping-into-morita-diag}. Passing to local traces yields the following commutative diagram:
\[\begin{tikzcd}
	{\Lcal_{v\Frob}\Xibf_{\ofrak^{\Lbf}}^{\Lbf,\opp}} && {\mathrm{\Pr}} \\
	\\
	{\Lcal_{\Frob}\Xibf_{\ofrak^{\Gbf}}^{\Gbf, \opp}} && {\mathrm{\Pr}}
	\arrow["{{\Tr^{\loc}(\Ad(\dot{v})\Frob,\HH^{\Lbf}_{\ofrak^{\Lbf}})}}", from=1-1, to=1-3]
	\arrow["{{\cdot v}}", from=1-1, to=3-1]
	\arrow["{{\cdot \Delta_v}}"{description}, Rightarrow, from=1-3, to=3-1]
	\arrow[equals, from=1-3, to=3-3]
	\arrow["{{\Tr^{\loc}(\Frob,\HH^{\Gbf}_{\ofrak^{\Gbf}})}}"', from=3-1, to=3-3]
\end{tikzcd}\]
As the $\Lcal_{\Frob}\Xibf_{\ofrak}$ (resp. $\Lcal_{v\Frob}\Xibf_{\ofrak^{\Lbf}}^{\Lbf}$)-actions come from totalization of the horizontal maps in the above diagram, we get the desired functoriality. 
\end{proof}

\subsection{Series and parabolic character sheaves}

Let $\Pbf = \Lbf\UP$, $v, \dot{v}$ and $\ofrak$ be as in \Cref{sec:series-and-induction}. Recall from Theorem \ref{thm:TracesAndParabolicCharacterSheaves} that we have an equivalence 
$$\Tr(\HH^{\Gbf}_{\ofrak, \tot, \Frob},\HH^{\Lbf}_{\ofrak, \tot}) \cong \DD^{\ofrak}(\frac{\UP\backslash \Gbf /\Frob(\UP)}{\Ad_{\Frob}(\Lbf)}).$$

\begin{defi}\label{def:groupoid-relative-loop-space}
Let $\Lcal^{\Lbf}_{\Frob}\Xibf_{\ofrak}^{\Gbf}$ be the fiber product 
$$\Lcal^{\Lbf}_{\Frob}\Xibf_{\ofrak}^{\Gbf} = \Xibf_{\ofrak}^{\Lbf} \times_{(\incl,\,\incl\circ\Frob),\, \Xibf_{\ofrak}^{\Gbf} \times \Xibf_{\ofrak}^{\Gbf},\, \Delta} \Xibf_{\ofrak}^{\Gbf}.$$
Here $\Delta : \Xibf_{\ofrak}^{\Gbf} \to \Xibf_{\ofrak}^{\Gbf} \times \Xibf_{\ofrak}^{\Gbf}$ is the diagonal, and the other map is the pair consisting of the natural inclusion $\Xibf_{\ofrak}^{\Lbf} \subset \Xibf_{\ofrak}^{\Gbf}$ and the composite
$$\Xibf_{\ofrak}^{\Lbf} \xrightarrow{\Frob} \Xibf_{\ofrak}^{\Frob(\Lbf)}\xrightarrow{\incl} \Xibf_{\ofrak}^{\Gbf}.$$
As a groupoid, it is the groupoid of pairs $(s, \beta) \in \Lcal_{\Frob}\Xibf_{\ofrak}^{\Gbf}$ and morphisms are commutative diagrams 
\[\begin{tikzcd}
	s & {\Frob(s)} \\
	{s'} & {\Frob(s')}
	\arrow["\beta", from=1-1, to=1-2]
	\arrow["\alpha"', from=1-1, to=2-1]
	\arrow["{\Frob(\alpha)}", from=1-2, to=2-2]
	\arrow["{\beta'}"', from=2-1, to=2-2]
\end{tikzcd}\]
where $\alpha$ is a morphism in $\Xibf_{\ofrak}^{\Lbf}$. 
\end{defi} 

The local trace construction of \Cref{sec:local-traces} yields a functor 
\begin{align*}
\Lcal^{\Lbf}_{\Frob}\Xibf_{\ofrak}^{\Gbf, \opp} &\to \mathrm{Pr}_{\Lambda} \\
(s, \beta) &\mapsto \Tr(\HH^{\Gbf, \beta}_{[s,\Frob(s)]},\HH^{\Lbf, \circ}_s) ,
\end{align*}
whose totalization is $\DD^{\ofrak}(\frac{\UP\backslash \Gbf /\Frob(\UP)}{\Ad_{\Frob}(\Lbf)})$. 

\begin{lem}
Let $\ofrak^{\Lbf} \subset \ofrak$ be a $\Ad(v)\Frob$ stable $\Wbf_{\Lbf}$-orbit of $\ofrak$. The map of \Cref{constr:translation-by-v} factors as 
\begin{equation}\label{eq:lem-parabolic-char-sheaves-and-series}
\Lcal_{v\Frob}\Xibf_{\ofrak^{\Lbf}}^{\Lbf} \xrightarrow{\cdot v} \Lcal^{\Lbf}_{\Frob}\Xibf_{\ofrak}^{\Gbf} \to \Lcal_{\Frob}\Xibf_{\ofrak}^{\Gbf}.
\end{equation}
Furthermore, the functors 
$$\DD^{\ofrak^{\Lbf}}(\pt/\Lbf^{\dot{v}\Frob}) \xrightarrow{j_{v,!}k_{\dot{v}}^*} \DD^{\ofrak}(\frac{\UP\backslash \Gbf /\Frob(\UP)}{\Ad_{\Frob}(\Lbf)}) \xrightarrow{\ch_{\Pbf \subset \Gbf}} \DD^{\ofrak}(\pt/\Gbf^{\Frob})$$
are linear over the maps \eqref{eq:lem-parabolic-char-sheaves-and-series}. 
\end{lem}

\begin{proof}
The existence of the factorization (\ref{eq:lem-parabolic-char-sheaves-and-series}) is immediate. The linearity of the map on traces follows from the functoriality statement \Cref{prop:functoriality-local-traces}, and its compatibility with totalization follows from \Cref{thm:local-to-global-trace}.
\end{proof}

\begin{lem}
\begin{enumerate}
\item There is a canonical disjoint union decomposition
$$\Lcal^{\Lbf}_{\Frob}\Xibf_{\ofrak}^{\Gbf} = \bigsqcup_{\ofrak^{\Lbf}\subset \ofrak} \Lcal^{\Lbf}_{\Frob}\Xibf_{\ofrak^{\Lbf}}^{\Gbf}$$
where the disjoint union ranges over all $\Wbf_{\Lbf}$-orbits in $\ofrak$, and $\Lcal^{\Lbf}_{\Frob}\Xibf_{\ofrak^{\Lbf}}^{\Gbf}$ is the subgroupoid of pairs $(s,\beta) \in \Lcal^{\Lbf}_{\Frob}\Xibf_{\ofrak}^{\Gbf}$ with $s \in \ofrak^{\Lbf}$.
\item For $s \in \ofrak^{\Lbf} \subset \ofrak$, there is a canonical isomorphism of groupoids
$$\Lcal^{\Lbf}_{\Frob}\Xibf_{\ofrak^{\Lbf}}^{\Gbf} = \frac{{_s}\underline{\Wbf}_{\Frob(s)}}{\Ad_{\Frob}(\Gammabf_s^{\Lbf})}.$$
\end{enumerate}
\end{lem}

\begin{proof}
For the first point, two elements  $(s, \beta)$ and $(s', \beta')$ can be conjugate in $\Lcal^{\Lbf}_{\Frob}\Xibf_{\ofrak}^{\Gbf}$ only if $s$ and $s'$ are $\Wbf_{\Lbf}$-conjugate. This establishes the disjoint union decomposition. 

The second point is immediate in view of \Cref{def:groupoid-relative-loop-space}.
\end{proof}

\subsection{Unipotent representations of $\Hbffrak_{\ofrak}$}\label{sec:unipotent-representations-endoscopic-side}

Let $\ofrak \subset \Ch(\Tbf)$ be a $\Frob$-stable orbit. Recall that we have a groupoid $\Hbffrak_{\ofrak}$ with a map 
$$\Hbffrak_{\ofrak} \to \Xibf_{\ofrak}$$
that is an $\Hbf_{s}^{\circ}$-gerbe (for some choice of $s \in \ofrak$). Since $\ofrak$ is stable under $\Frob$, there is a well-defined Frobenius morphism
$$\Frob : \Hbffrak_{\ofrak} \to \Hbffrak_{\ofrak}.$$

Consider the groupoid $\Lcal_{\Frob}\Hbffrak_{\ofrak}$, the corresponding $\Frob$-loop-space (or $\Frob$-inertia stack). By functoriality it comes equipped with a map 
$$\pi_{\ofrak} : \Lcal_{\Frob}\Hbffrak_{\ofrak} \to \Lcal_{\Frob}\Xibf_{\ofrak}.$$
We now describe the stack $\Lcal_{\Frob}\Hbffrak_{\ofrak}$ and the fibers of $\pi_{\ofrak}$. 
\begin{lem}
Let $s \xrightarrow{\beta} \Frob(s)$ be a point of $\Lcal_{\Frob}\Xibf_{\ofrak}$; see \Cref{sec:RationalSeries}. The fiber $\pi^{-1}(s,\beta)$ is canonically isomorphic to the stack $\frac{\Hbf_{[s,\Frob(s)]}^{\beta}}{\Ad_{\Frob}(\Hbf_s^{\circ})}$, so there are pullback diagrams
\[\begin{tikzcd}
	{\frac{\Hbf_{[s,\Frob(s)]}^{\beta}}{\Ad_{\Frob}(\Hbf_s^{\circ})}} & {\Lcal_{\Frob}(\Hbffrak_{\ofrak})} & {\frac{\Hbf_{[s,\Frob(s)]}}{\Ad_{\Frob}(\Hbf_s)}} & {\Lcal_{\Frob}(\Hbffrak_{\ofrak})} \\
	\pt & {\Lcal_{\Frob}\Xibf_{\ofrak}} & {\frac{{_s}\underline{\Wbf}_{\Frob(s)}}{\Ad_{\Frob}(\Gammabf_s)}} & {\Lcal_{\Frob}\Xibf_{\ofrak}}
	\arrow[from=1-1, to=1-2]
	\arrow[from=1-1, to=2-1]
	\arrow["{\pi_{\ofrak}}", from=1-2, to=2-2]
	\arrow[from=1-3, to=1-4]
	\arrow[from=1-3, to=2-3]
	\arrow["{\pi_{\ofrak}}", from=1-4, to=2-4]
	\arrow["{(s,\beta)}"', from=2-1, to=2-2]
	\arrow[from=2-3, to=2-4]
\end{tikzcd}\]
\end{lem}

\begin{proof}
This is clear from the definitions.
\end{proof}

\begin{rque}
Given a base point $\dot{\beta} \in \Hbf_{[s,\Frob(s)]}^{\beta}$, we get an isomorphism 
$$\frac{\Hbf_{[s,\Frob(s)]}^{\beta}}{\Ad_{\Frob}(\Hbf_{s}^{\circ})} \cong \frac{\Hbf_{s}^{\circ}}{\Ad_{\Frob_{\dot{\beta}}}(\Hbf_s^{\circ})}$$
where $\Frob_{\dot{\beta}} := \Ad(\dot{\beta}) \circ \Frob : \Hbf_s^{\circ} \to \Hbf_s^{\circ}$ is the corresponding twisted Frobenius endomorphism.
\end{rque}

\begin{defilem}\label{def:unipotent-endoscopic-side}
For $\dot{\beta} \in \Hbf_{[s,\Frob(s)]}^{\beta}$, we denote by
$$\DD^{\unip}(\frac{\Hbf_{[s,\Frob(s)]}^{\beta}}{\Hbf_{s}^{\circ}}) \subset \DD(\frac{\Hbf_{[s,\Frob(s)]}^{\beta}}{\Hbf_{s}^{\circ}})$$
the full subcategory of objects that belong to $\DD^{\unip}(\pt/\Hbf_s^{\Frob_{\dot{\beta}}})$. This category does not depend on the choice of $\dot{\beta}$. 
\end{defilem}

\begin{proof}
If $\dot{\beta}'$ is another lift of $\beta$, then $\dot{\beta} = x\dot{\beta}'$ for a unique $x \in \Hbf_s^{\circ}$. The well-definedness of $\DD^{\unip}(\frac{\Hbf_{[s,\Frob(s)]}^{\beta}}{\Hbf_{s}^{\circ}})$ follows from the fact that $\DD^{\unip}(\pt/\Hbf_s^{\Frob_{\dot{\beta}}})$ is stable under conjugation. 
\end{proof}

\begin{defi}
We denote by $\DD^{\unip}(\Lcal_{\Frob}(\Hbffrak_{\ofrak}))$ the full subcategory of sheaves $A \in \DD(\Lcal_{\Frob}(\Hbffrak_{\ofrak}))$ such that for all $(s,\beta) \in \Lcal_{\Frob}\Xibf_{\ofrak}$, the restriction to $\frac{\Hbf_{[s,\Frob(s)]}^{\beta}}{\Ad_{\Frob}(\Hbf_s^{\circ})}$ lies in $\DD^{\unip}(\frac{\Hbf_{[s,\Frob(s)]}^{\beta}}{\Ad_{\Frob}(\Hbf_s^{\circ})})$. 
\end{defi} 

\begin{rque}
Note that there is an action of $\DD(\Lcal_{\Frob}\Xibf_{\ofrak})$ on $\DD(\Lcal_{\Frob}(\Hbffrak_{\ofrak}))$ and that $\DD^{\unip}(\Lcal_{\Frob}(\Hbffrak_{\ofrak}))$ is stable under this action. 
\end{rque}

The next theorem is a variant of Theorem \ref{thm:Trace}, we record it separately.

\begin{thm}\label{thm:equiv-trace-endoscopic-side}
There is a natural equivalence 
$$\Tr(\Frob^*, \HH^{\Hbf}_{\ofrak, \tot}) \cong \DD^{\unip}(\Lcal_{\Frob}(\Hbffrak_{\ofrak})).$$
\end{thm}

\begin{rque}\label{rque:local-trace-endoscopic-side}
Note that the $\DD(\Lcal_{\Frob}\Xibf_{\ofrak})$-action on $\DD^{\unip}(\Lcal_{\Frob}(\Hbffrak_{\ofrak}))$ coming from the map 
$$\Lcal_{\Frob}(\Hbffrak_{\ofrak}) \to \Lcal_{\Frob}(\Xibf_{\ofrak})$$
is canonically equivalent to the one coming from the local trace formalism, see \Cref{thm:local-to-global-trace}.
\end{rque}

\subsection{Parabolic character sheaves for $\Hbffrak_{\ofrak}$}

Let $\Pbf = \Lbf\UP \subset \Gbf$ be a parabolic subgroup, let $v \in {^\Pbf}\Wbf^{\Frob(\Pbf)}$, and let $\dot{v}$ be a lift of $v$ such that $\Lbf$ is $\Ad(\dot{v}) \Frob$-stable. Let $\ofrak \subset \Ch(\Tbf)$ be a $\Frob$-stable $\Wbf$-orbit. The goal of this section is to define parabolic character sheaves for $\Hbffrak_{\ofrak}$. Recall from \Cref{sec:endoscopic-parabolic-groups} that there is a diagram of groupoids
\[\begin{tikzcd}
	{\Hbffrak^{\Lbf}_{\ofrak^{\Lbf}}} & {\Hbffrak^{\Pbf}_{\ofrak^{\Lbf}}} & {\Hbffrak^{\Gbf}_{\ofrak^{\Gbf}}} \\
	{\Xibf_{\ofrak^{\Lbf}}^{\Lbf}} & {\Xibf_{\ofrak^{\Lbf}}^{\Lbf}} & {\Xibf_{\ofrak^{\Gbf}}^{\Gbf}}
	\arrow[from=1-1, to=2-1]
	\arrow["{\alpha_{\Pbf}^{\Hbf}}"', from=1-2, to=1-1]
	\arrow["{\beta_{\Pbf}^{\Hbf}}", from=1-2, to=1-3]
	\arrow[from=1-2, to=2-2]
	\arrow[from=1-3, to=2-3]
	\arrow[equals, from=2-2, to=2-1]
	\arrow[from=2-2, to=2-3]
\end{tikzcd}\]
Recall as well that the fiber of $\Hbffrak_{\ofrak^?}^{?}$ over a block $\beta : s \to s'$ is $\Hbf_{[s,s']}^{?, \beta}$ where $? \in \{\Lbf, \Pbf, \Gbf\}$.  

\begin{defi}
We define $\Lcal_{\Frob}^{?}\Hbffrak_{\ofrak}^{\Gbf}$ to be the fiber product of groupoids 
\[\begin{tikzcd}
	{\Lcal_{\Frob}^{?}\Hbffrak_{\ofrak}^{\Gbf}} & {\Hbffrak_{\ofrak^{\Gbf}}^{\Gbf}} \\
	{\Xibf_{\ofrak}^{?}} & {\Hbffrak_{\ofrak^{\Gbf}}^{\Gbf} \times \Hbffrak_{\ofrak^{\Gbf}}^{\Gbf}}
	\arrow[from=1-1, to=1-2]
	\arrow[from=1-1, to=2-1]
	\arrow["\Delta", from=1-2, to=2-2]
	\arrow["{\incl \times \Frob}"', from=2-1, to=2-2]
\end{tikzcd}\]
where $? \in \{\Lbf, \Pbf\}$.
\end{defi} 

There is a natural map of groupoids
\begin{equation}\label{eq:map-to-xi-endoscopic-side}
\Lcal_{\Frob}^{?}\Hbffrak_{\ofrak}^{\Gbf} \to \Lcal_{\Frob}^{\Lbf}\Xibf_{\ofrak}^{\Gbf}.
\end{equation}
\begin{lem}
For $(s, \beta) \in \Lcal_{\Frob}^{\Lbf}\Xibf_{\ofrak}^{\Gbf}$, there is a canonical Cartesian diagram
\[\begin{tikzcd}
	{\frac{\Hbf^{\Gbf}_{[s,\Frob(s)]}}{\Ad_{\Frob}(\Hbf^{?}_s)}} & {\Lcal_{\Frob}^{?}\Hbffrak_{\ofrak}^{\Gbf}} \\
	{\frac{{_s}\Wbf_{\Frob(s)}}{\Ad_{\Frob}(\Gammabf_s^{\Lbf})}} & {\Lcal_{\Frob}^{?}\Xibf_{\ofrak}^{\Gbf}}
	\arrow[from=1-1, to=1-2]
	\arrow[from=1-1, to=2-1]
	\arrow[from=1-2, to=2-2]
	\arrow[from=2-1, to=2-2]
\end{tikzcd}\]
\end{lem}
\begin{proof}
Immediate from the definition. 
\end{proof}

\begin{defi}
For $(s, \beta) \in \Lcal_{\Frob}^{\Lbf}\Xibf_{\ofrak}^{\Gbf}$, we define
$$\DD_{(\UP,\UP)}\!\left(\frac{\Hbf^{\Gbf}_{[s,\Frob(s)]}}{\Ad_{\Frob}(\Hbf^{\Lbf}_s)}\right)
:=\DD\!\left(\frac{\Ubf_s^{\Hbf,\Pbf}\backslash\Hbf^{\Gbf}_{[s,\Frob(s)]}/\Ubf_{\Frob(s)}^{\Hbf,\Frob(\Pbf)}}{\Ad_{\Frob}(\Hbf^{\Lbf}_s)}\right).$$
Thus the notation $(\UP,\UP)$ in this subsection is only a shorthand for the two fiberwise endoscopic unipotent radicals; it never denotes the subgroup $\UP\subset\Gbf$.
\end{defi} 

\begin{lem}\label{lem:lusztig-stratif-endoscopic-side}
For $(s, \beta) \in \Lcal_{\Frob}^{\Lbf}\Xibf_{\ofrak}^{\Gbf}$, the category $\DD_{(\UP,\UP)}(\frac{\Hbf^{\Gbf}_{[s,\Frob(s)]}}{\Ad_{\Frob}(\Hbf^{\Lbf}_s)})$ has a natural semiorthogonal decomposition coming from Lusztig stratification. A stratum $t$ in Lusztig's stratification has stable datum $(I_{t,\infty},w_{t,\infty})$; let $\Hbf^{\Lbf_t}$ denote the corresponding stable endoscopic Levi. The graded piece attached to $t$ is of the form
$$\DD\!\left(\pt/\bigl((\Hbf^{\Lbf_t})^{w_{t,\infty}\Frob}\rtimes\Vbf_t\bigr)\right),$$
where $\Vbf_t$ is a connected unipotent group.
\end{lem}

\begin{proof}
Apply \Cref{corol:StructureOfStrata} fiberwise to the endoscopic parabolic just defined. Its stable pair $(I_{t,\infty},w_{t,\infty})$ gives $\Hbf^{\Lbf_t}$ and the displayed stabilizer. The semiorthogonal decomposition comes from Lusztig's  stratification.
\end{proof}

\begin{defi}\label{defi:unipotent-parabolic-cs-endoscopic-side}
We define $\DD_{(\UP,\UP)}^{\unip}(\frac{\Hbf^{\Gbf}_{[s,\Frob(s)]}}{\Ad_{\Frob}(\Hbf^{\Lbf}_s)})$ to be the full subcategory of $\DD_{(\UP,\UP)}(\frac{\Hbf^{\Gbf}_{[s,\Frob(s)]}}{\Ad_{\Frob}(\Hbf^{\Lbf}_s)})$ for which the $!$- and $*$-restrictions to the strata of Lemma \ref{lem:lusztig-stratif-endoscopic-side} are unipotent in the sense of Definition \ref{def:unipotent-endoscopic-side}.

We denote by $\DD^{\unip}_{(\UP, \UP)}(\Lcal^{\Lbf}_{\Frob}\Hbffrak_{\ofrak}^{\Gbf})$ the full subcategory of $\DD(\Lcal^{\Lbf}_{\Frob}\Hbffrak_{\ofrak}^{\Gbf})$ of objects whose restriction to $\frac{\Hbf^{\Gbf}_{[s,\Frob(s)]}}{\Ad_{\Frob}(\Hbf^{\Lbf}_s)}$ lie in $\DD_{(\UP,\UP)}^{\unip}(\frac{\Hbf^{\Gbf}_{[s,\Frob(s)]}}{\Ad_{\Frob}(\Hbf^{\Lbf}_s)})$. 
\end{defi} 

\begin{defi}\label{defi:DL-induction-endoscopic-side}
The maps $\alpha_{\Pbf}^{\Hbf}$ and $\beta_{\Pbf}^{\Hbf}$ induce natural maps 
$$\DD^{\unip}_{(\UP, \UP)}(\Lcal^{\Lbf}_{\Frob}\Hbffrak_{\ofrak}^{\Gbf}) \xrightarrow{\alpha_{\Pbf}^{\Hbf,*}} \DD(\Lcal^{\Pbf}_{\Frob}\Hbffrak_{\ofrak}^{\Gbf}) \xrightarrow{\beta_{\Pbf,!}^{\Hbf}} \DD^{\unip}(\Lcal_{\Frob}\Hbffrak_{\ofrak}^{\Gbf}).$$
We denote by $\ch_{\Pbf \subset \Gbf}^{\Hbf} = \beta_{\Pbf,!}^{\Hbf}\alpha_{\Pbf}^{\Hbf,*}$ the composition.
\end{defi}

\begin{thm}\label{thm:functoriality-DL-induction-endoscopic-side}
There is a canonical equivalence $\Tr(\HH^{\Hbf^{\Gbf}}_{\ofrak, \Frob},\HH^{\Hbf^{\Lbf}}_{\ofrak}) \cong \DD^{\unip}_{(\UP, \UP)}(\Lcal^{\Lbf}_{\Frob}\Hbffrak_{\ofrak}^{\Gbf})$ and a commutative diagram 
\[\begin{tikzcd}
	{\Tr(\HH^{\Hbf^{\Gbf}}_{\ofrak, \Frob},\HH^{\Hbf^{\Lbf}}_{\ofrak}) } & {\Tr(\HH^{\Hbf^{\Gbf}}_{\ofrak, \Frob},\HH^{\Hbf^{\Gbf}}_{\ofrak}) } \\
	{\DD^{\unip}_{(\UP, \UP)}(\Lcal^{\Lbf}_{\Frob}\Hbffrak_{\ofrak}^{\Gbf})} & {\DD^{\unip}(\Lcal_{\Frob}\Hbffrak_{\ofrak}^{\Gbf})}
	\arrow[from=1-1, to=1-2]
	\arrow["\wr"', from=1-1, to=2-1]
	\arrow["\wr", from=1-2, to=2-2]
	\arrow["{\ch_{\Pbf \subset \Gbf}^{\Hbf}}"', from=2-1, to=2-2]
\end{tikzcd}\]
where the top map is induced by the natural map of pairs $(\HH^{\Hbf^{\Lbf}}_{\ofrak}, \HH^{\Hbf^{\Gbf}}_{\ofrak, \Frob}) \to (\HH^{\Hbf^{\Gbf}}_{\ofrak}, \HH^{\Hbf^{\Gbf}}_{\ofrak, \Frob})$. 
\end{thm}

\begin{proof}
Using the formalism of local traces \Cref{sec:local-traces}, both $\Tr(\HH^{\Hbf^{\Gbf}}_{\ofrak, \Frob},\HH^{\Hbf^{\Lbf}}_{\ofrak})$ and $\DD^{\unip}_{(\UP, \UP)}(\Lcal^{\Lbf}_{\Frob}\Hbffrak_{\ofrak}^{\Gbf})$ are fibered over $\DD(\Lcal^{\Lbf}_{\Frob}\Xibf^{\Gbf}_{\ofrak})$ hence it is enough to prove the statement fiberwise. This is then an application of \Cref{thm:TracesAndParabolicCharacterSheaves}.
\end{proof}

Recall the definition of the scheme $\Hbf^{\Gbf}_{[s,\Frob(s)]} = \Hbf_{[s,s]}^{\Gbf} \times^{\Nbf_{[s,s]}} \Nbf_{[s,\Frob(s)]}$, where $\Nbf_{[s,\Frob(s)]}$ is the union of all connected components of $\Nbf$ above elements $w \in {_s}\Wbf_{\Frob(s)}$. Associated with $\Pbf$, there are (disconnected) parabolic subgroups $\Hbf^{\Pbf}_{[s,s]}$ of $\Hbf^{\Gbf}_{[s,s]}$ and $\Hbf^{\Pbf}_{[\Frob(s), \Frob(s)]}$ of $\Hbf^{\Gbf}_{[\Frob(s), \Frob(s)]}$. We consider the stratification of $\Hbf^{\Gbf}_{[s,\Frob(s)]}$ into $(\Hbf^{\Pbf}_{[s,s]}, \Hbf^{\Pbf}_{[\Frob(s), \Frob(s)]})$-double cosets. These double cosets are indexed by
$$\Wbf^{\Lbf}_s \backslash {_s}\Wbf_{\Frob(s)} /\Wbf_{\Frob(s)}^{\Lbf}$$
and the above choice of $v$ determines a unique such double coset. 

Let us denote by $\Hbf^{\Pbf,v}_{[s, \Frob(s)]} \subset \Hbf^{\Gbf}_{[s,\Frob(s)]}$ the corresponding double coset. By our assumption on $v$, it contains a unique $\Hbf^{\Pbf}_{[s,s]}$-left coset. Denote by $\Nbf^{\Lbf}_s \subset \Hbf^{\Lbf}_{[s,s]}$ the normalizer of $\Tbf$ and $\Nbf^{\Pbf, v}_{[s, \Frob(s)]}$ the union of connected components of $\Nbf(\Tbf)$ above $\Wbf^{\Lbf}v\Frob(\Wbf^{\Lbf}) \cap {_s}\Wbf_{\Frob(s)}$.
\begin{lem}
There is a canonical isomorphism 
$$\Hbf^{\Pbf,v}_{[s, \Frob(s)]} = \Hbf^{\Pbf}_{[s,s]} \times^{\Nbf^{\Lbf}_s} \Nbf_{[s,\Frob(s)]}^{\Pbf,v}.$$
\end{lem}

\begin{proof}
Multiplication defines a map from the RHS to the LHS and it is clear that this map is an isomorphism by counting Bruhat strata.
\end{proof}

Note that multiplication yields a canonical isomorphism 
$$\Nbf_{[s,\Frob(s)]}^{\Pbf,v} = \Nbf^{\Lbf}_{[s, v\Frob(s)]} \times^{\Tbf} \Nbf_v$$ 
where $\Nbf_v$ is the fiber of $\Nbf$ above $v$. It follows that the choice of $\dot{v} \in \Nbf_v$, yields a map 
\begin{align*}
\Hbf^{\Pbf,v}_{[s, \Frob(s)]} &= \Hbf^{\Pbf}_{[s,s]} \times^{\Nbf^{\Lbf}_s} \Nbf_{[s,\Frob(s)]}^{\Pbf,v} \\
&= \Hbf^{\Pbf}_{[s,s]} \times^{\Nbf^{\Lbf}_s} \Nbf^{\Lbf}_{[s, v\Frob(s)]} \times^{\Tbf} \Nbf_v \\
&\to \Hbf^{\Lbf}_{[s,s]} \times^{\Nbf^{\Lbf}_s} \Nbf^{\Lbf}_{[s, v\Frob(s)]} \times^{\Tbf} \Nbf_v \\
&\xrightarrow{k_{\dot{v}}} \Hbf^{\Lbf}_{[s,s]} \times^{\Nbf^{\Lbf}_s} \Nbf^{\Lbf}_{[s, v\Frob(s)]} \\
&= \Hbf^{\Lbf}_{[s,v\Frob(s)]}.
\end{align*}
In the above morphism, the third map is induced by the projection $\Hbf^{\Pbf}_{[s,s]} \to \Hbf^{\Lbf}_{[s,s]}$. It is clear that this composition is equivariant for both the left (resp. right) action of $\Hbf^{\Lbf}_{[s,s]}$ (resp. $\Hbf^{\Lbf}_{[\Frob(s), \Frob(s)]}$). 

The correspondence 
$$\Hbf^{\Lbf}_{[s,v\Frob(s)]} \xleftarrow{k_{\dot{v}}} \Hbf^{\Pbf,v}_{[s,\Frob(s)]} \xrightarrow{j_v} \Hbf^{\Gbf}_{[s,\Frob(s)]}$$
then yields a map 
$$\DD^{\unip}(\Lcal_{\dot{v}\Frob}\Hbffrak^{\Lbf}_{\ofrak^{\Lbf}}) \xrightarrow{j_{v,!}k_{\dot{v}}^*} \DD_{(\UP, \UP)}^{\unip}(\Lcal^{\Lbf}_{\Frob}\Hbffrak_{\ofrak}).$$

\begin{rque}
The composition 
$$\DD^{\unip}(\Lcal_{\dot{v}\Frob}\Hbffrak^{\Lbf}_{\ofrak^{\Lbf}}) \xrightarrow{j_{v,!}k_{\dot{v}}^*} *\DD_{(\UP, \UP)}^{\unip}(\Lcal^{\Lbf}_{\Frob}\Hbffrak_{\ofrak}) \xrightarrow{\ch^{\Hbf}_{\Pbf \subset \Gbf}} \DD^{\unip}(\Lcal_{\Frob}\Hbffrak_{\ofrak})$$
is a version of Deligne--Lusztig induction for disconnected groups. Restricting to connected components of both the source and target recovers the constructions of \cite{DigneMichelNonConnexe}.
\end{rque}

\subsection{Endoscopy}\label{sec:EndoscopyReps}

We now come to the main result of this section, Theorem \ref{thm:EndoscopyRep}, which establishes the endoscopic equivalence for representations of $\Gbf^{\Frob}$. So far, we have not used any hypothesis on $\ell$; we now restrict to the following setting:
\begin{enumerate}
\item if $\Lambda = \Qlb$, then there is no condition on $\ell$, 
\item if $\Lambda \in \{\Zlb, \Flb\}$, then $\ell$ is not a torsion prime for $\Gbf^*$; see \Cref{def:torsion-prime}.
\end{enumerate}

We will now combine categorical traces with the endoscopic equivalence.
\begin{thm}\label{thm:EndoscopyRep}
There is a canonical $\DD(\Lcal_{\Frob}\Xibf_{\ofrak})$-linear equivalence, depending only on the pinning $\phi$ of $\Gbf$,
$$\EE^{\Frob}_{\ofrak} : \DD^{\ofrak}(\frac{\Gbf}{\Ad_{\Frob}\Gbf}) \cong \DD^{\unip}(\Lcal_{\Frob}\Hbffrak_{\ofrak}).$$
Furthermore, under the equivalences $\Tr(\Frob, \HH_{\ofrak}^{\Gbf}) \cong \DD^{\ofrak}(\frac{\Gbf}{\Ad_{\Frob}\Gbf})$ (see \Cref{thm:Trace}) and $\Tr(\Frob^*, \HH^{\Hbf}_{\ofrak}) \cong \DD^{\unip}(\Lcal_{\Frob}\Hbffrak_{\ofrak})$ (see \Cref{thm:equiv-trace-endoscopic-side}), the following diagram commutes:
\begin{equation}\label{eq:diagram-commutation-into-trace-endoscopy}
\begin{tikzcd}
	{\HH^{\Gbf}_{\ofrak, \tot}} & {\HH^{\Hbf}_{\ofrak, \tot}} \\
	{\DD^{\ofrak}(\frac{\Gbf}{\Ad_{\Frob}\Gbf})} & {\DD^{\unip}(\Lcal_{\Frob}\Hbffrak_{\ofrak})}
	\arrow["{\EE_{\ofrak}}", no head, from=1-1, to=1-2]
	\arrow[from=1-1, to=2-1]
	\arrow[from=1-2, to=2-2]
	\arrow["{\EE_{\ofrak}^{\Frob}}"', from=2-1, to=2-2]
\end{tikzcd}
\end{equation}
\end{thm}

\begin{proof}
The endoscopic equivalence for Hecke categories, \Cref{thm:EndoscopyHeckeCat}, yields an equivalence 
$$\HH^{\Gbf}_{\ofrak,\tot} \cong \HH^{\Hbf}_{\ofrak,\tot}$$
which is compatible with $\Frob$ by \Cref{lem:FrobeniusEquivariance}. This yields an equivalence 
$$\Tr(\Frob^*, \HH^{\Gbf}_{\ofrak,\tot}) \cong \Tr(\Frob^*, \HH^{\Hbf}_{\ofrak, \tot}).$$
This is the desired equivalence $\EE^{\Frob}_{\ofrak}$, furthermore the diagram in \Cref{thm:EndoscopyRep} is commutative. The equivalence $\EE_{\ofrak}^{\Frob} $ depends only on the equivalence $\EE_{\ofrak}$ (and its compatibility with $\Frob$), both of which depend only on the pinning of $\Gbf$ by \Cref{corol:dependency-on-pinning}. Finally, the $\DD(\Lcal_{\Frob}\Xibf_{\ofrak})$-linearity follows from the fact that the endoscopic equivalence  is an equivalence of functors 
$$\Xibf_{\ofrak}^{\opp} \to \Morita$$
compatible with $\Frob$, therefore both categories $\Tr(\Frob^*, \HH^{\Gbf}_{\ofrak,\tot})$ and $\Tr(\Frob^*,\HH^{\Hbf}_{\ofrak, \tot})$ can be realized via local traces. By the construction of \Cref{sec:RationalSeries} and \Cref{rque:local-trace-endoscopic-side}, the equivalence $\EE^{\Frob}_{\ofrak}$ is $\DD(\Lcal_{\Frob}\Xibf_{\ofrak})$-linear. 
\end{proof}

\begin{corol}\label{corol:RationalEndoscopy}
With the notations of Theorem \ref{thm:EndoscopyRep}, let $\ofrak'_{\rat} \subset \ofrak$ be a rational conjugacy class and denote by $\Lcal_{\Frob}\Hbffrak_{\ofrak}^{\ofrak'_{\rat}}$ the corresponding connected component, the equivalence $\EE_{\ofrak}^{\Frob}$ restricts to a $\DD(\Lcal_{\Frob}\Xibf_{\ofrak}^{\ofrak'_{\rat}})$-linear equivalence 
$$\EE_{\ofrak'_{\rat}}^{\Frob} : \DD^{\ofrak'_{\rat}}(\frac{\Gbf}{\Ad_{\Frob}\Gbf}) \cong \DD^{\unip}(\Lcal_{\Frob}\Hbffrak_{\ofrak}^{\ofrak'_{\rat}}).$$
\end{corol}

\begin{proof}
This follows from the $\Lcal_{\Frob}\Xibf_{\ofrak}$-linearity in Theorem \ref{thm:EndoscopyRep} by projecting on the corresponding connected components. 
\end{proof}

\subsubsection{Morita equivalent description}

To make the equivalence $\EE_{\ofrak}^{\Frob}$ more amenable to computation, let us unravel it through the Morita equivalences \Cref{lem:MoritaEquivalence1}.
\begin{equation}\label{eq:morita-equiv-for-traces}
\HH_{\ofrak, \tot}^{\Gbf} \sim \HH^{\Gbf}_{[s,s]}, \ \HH_{\ofrak, \tot}^{\Hbf} \sim \HH_{[s,s]}^{\Hbf}.
\end{equation}

\begin{lem}\label{lem:Morita-equivalence-for-bimodule-for-trace}
Under the Morita equivalences \Cref{eq:morita-equiv-for-traces}, the $\HH_{\ofrak, \tot}^?$ bimodule $\HH_{\ofrak, \tot, \Frob}^?$ corresponds to $\HH_{[s,\Frob(s)]}^?$ for $? \in \{\Hbf, \Gbf\}$.  
\end{lem}

\begin{proof}
We prove the assertion for $\HH^{\Gbf}$. The bimodule corresponding to $\HH^{\Gbf}_{\ofrak, \tot, \Frob}$ is
\begin{align*}
{_s}\mathbb{M} &\otimes_{\HH_{\ofrak,\tot}^{\Gbf}} \HH^{\Gbf}_{\ofrak, \tot, \Frob} \otimes_{\HH_{\ofrak,\tot}^{\Gbf}} \mathbb{M}_s \\
&\cong {_s}\mathbb{M} \otimes_{\HH_{\ofrak,\tot}^{\Gbf}}  \mathbb{M}_{\Frob(s), \Frob} \\
&\cong \HH^{\Gbf}_{[s,\Frob(s)]}
\end{align*}
\end{proof}

It follows from \Cref{lem:Morita-equivalence-for-bimodule-for-trace}, that we have a commutative diagram of traces 
\begin{equation}\label{eq:diagram-morita-equiv-traces}
\begin{tikzcd}
	{\Tr(\Frob^*, \HH_{\mathfrak{o}}^{\Gbf})} & {\Tr(\Frob^*, \HH^{\Hbf}_{\mathfrak{o}})} \\
	{\Tr(\HH^{\Gbf}_{[s,\Frob(s)]}, \HH^{\Gbf}_{[s,s]})} & {\Tr(\HH^{\Hbf}_{[s,\Frob(s)]},\HH^{\Hbf}_{[s,s]}).}
	\arrow["{\EE_{\ofrak}^{\Frob}}", from=1-1, to=1-2]
	\arrow["\wr"', from=1-1, to=2-1]
	\arrow["\wr", from=1-2, to=2-2]
	\arrow["{\EE_{s}^{\Frob}}"', from=2-1, to=2-2]
\end{tikzcd}
\end{equation}

where the vertical maps are induced by the Morita equivalences. 

\begin{corol}\label{cor:endoscopy-from-a-morita-point-of-view}
There is a commutative diagram 
\[\begin{tikzcd}
	{\DD^{\ofrak}(\frac{\Gbf}{\Ad_{\Frob}\Gbf})} & {\DD^{\unip}(\Lcal_{\Frob}\Hbffrak_{\ofrak})} \\
	& {\DD^{\unip}(\frac{\Hbf_{[s,\Frob(s)]}}{\Ad_{\Frob}(\Hbf_s)})}
	\arrow["{\EE_{\ofrak}^{\Frob}}", from=1-1, to=1-2]
	\arrow["{\EE_s^{\Frob}}"', from=1-1, to=2-2]
	\arrow[from=1-2, to=2-2]
\end{tikzcd}\]\\
where the vertical map is restriction along $\frac{\Hbf_{[s,\Frob(s)]}}{\Ad_{\Frob}(\Hbf_s)}) \subset \Lcal_{\Frob}\Hbffrak_{\ofrak}$.
\end{corol}

\begin{proof}
This is clear from the Diagram \ref{eq:diagram-morita-equiv-traces}.
\end{proof}

\begin{corol}\label{cor:endoscopy-rat-series-morita-point-of-view}
Let $\ofrak'_{\rat} \subset \ofrak$, and let $[\beta] \subset {_s}\underline{\Wbf}_{\Frob(s)}$ be the corresponding $\Ad_{\Frob}(\Gammabf_s)$-orbit (see \Cref{lem:RationalConjClasses}). The equivalence $\EE_{s}^{\Frob}$ restricts to an equivalence
$$\EE_{s, \ofrak'_{\rat}}^{\Frob} : \DD^{\ofrak'_{\rat}}(\frac{\Gbf}{\Ad_{\Frob}(\Gbf)}) \to \DD^{\unip}(\frac{\Hbf_{[s,\Frob(s)]}^{[\beta]}}{\Ad_{\Frob}(\Hbf_{s})}),$$
where $\Hbf_{[s,\Frob(s)]}^{[\beta]} = \sqcup_{\gamma \in [\beta]} \Hbf_{[s,\Frob(s)]}^{\gamma}$.
\end{corol}

\begin{proof}
This is a direct combination of \Cref{corol:RationalEndoscopy} and \Cref{cor:endoscopy-from-a-morita-point-of-view}. 
\end{proof}

Let $s \in \ofrak$ and $\ofrak'_{\rat} \subset \ofrak$ be as in \Cref{cor:endoscopy-rat-series-morita-point-of-view}. The stack $\frac{\Hbf_{[s,\Frob(s)]}^{[\beta]}}{\Ad_{\Frob}(\Hbf_{s})}$ is connected. For $\dot{\beta} \in \Hbf_{[s,\Frob(s)]}^{[\beta]}$, there is an isomorphism
$$\frac{\Hbf_{[s,\Frob(s)]}^{[\beta]}}{\Ad_{\Frob}(\Hbf_{s})} \cong \pt/\Hbf_{s}^{\Frob_{\dot{\beta}}}$$
where $\Frob_{\dot{\beta}} = \Ad(\dot{\beta}) \circ \Frob$. The group $\Hbf_{s}^{\Frob_{\dot{\beta}}}$ is the group of rational points of the disconnected group $\Hbf_{s}$. 

\subsubsection{Endoscopy for parabolic character sheaves}

We assume that $\ell$ satisfies the assumptions as in Section \ref{sec:EndoscopyReps}. We fix $\Pbf = \Lbf\UP$, $v, \dot{v}$ and $\ofrak^{\Lbf} \subset \ofrak$ as in \Cref{sec:series-and-induction}. 

\begin{thm}\label{thm:compat-DL-induction}
There is an equivalence of categories 
$$\EE^{\Lbf, \Frob}_{\ofrak} : \DD^{\ofrak}(\frac{\UP \backslash \Gbf /\Frob(\UP)}{\Ad_{\Frob}(\Lbf)}) \cong \DD^{\unip}_{(\UP, \UP)}(\Lcal^{\Lbf}_{\Frob}\Hbffrak_{\ofrak}^{\Gbf})$$
making the following diagram commutative 

\[\begin{tikzcd}
	{\DD^{\ofrak^{\Lbf}}(\frac{\Lbf}{\Ad_{\dot v\Frob}\Lbf})} & {\DD^{\unip}(\Lcal_{v\Frob}\Hbffrak_{\ofrak^{\Lbf}}^{\Lbf})} \\
	{\DD^{\ofrak}(\frac{\UP\backslash \Gbf /\Frob(\UP)}{\Ad_{\Frob}\Lbf})} & {\DD^{\unip}_{(\UP,\UP)}(\Lcal_{\Frob}^{\Lbf}\Hbffrak_{\ofrak})} \\
	{\DD^{\ofrak}(\frac{\Gbf}{\Ad_{\Frob}\Gbf})} & {\DD^{\unip}(\Lcal_{\Frob}\Hbffrak_{\ofrak})}
	\arrow["{\EE^{v\Frob}_{\ofrak^{\Lbf}}}", from=1-1, to=1-2]
	\arrow["{j_{v,!}k_{\dot{v}}^*}"', from=1-1, to=2-1]
	\arrow["{j_{v,!}k_{\dot{v}}^*}", from=1-2, to=2-2]
	\arrow["{\EE^{\Lbf, \Frob}_{\ofrak}}", from=2-1, to=2-2]
	\arrow["{\ch_{\Pbf\subset\Gbf}}"', from=2-1, to=3-1]
	\arrow["{\ch_{\Pbf\subset\Gbf}^{\Hbf}}", from=2-2, to=3-2]
	\arrow["{\EE^{\Frob}_{\ofrak}}"', from=3-1, to=3-2]
\end{tikzcd}\]

\end{thm}

\begin{proof}
This is a combination of Theorems \ref{thm:Trace}, \ref{thm:TracesAndParabolicCharacterSheaves}, \ref{thm:equiv-trace-endoscopic-side}, \ref{thm:functoriality-DL-induction-endoscopic-side}, \ref{thm:EndoscopyRep} and \ref{thm:EndoscopyHeckeCat}.
\end{proof}

\subsubsection{Comparison with \cite{BonnafeRouquier} and \cite{BonnafeDatRouquier}}

We keep the same notations as before. We fix $s \in \ofrak^{\Lbf}$.

Assume that $s \in \Gbf^*$ is central, then $\ofrak = \{s\}$ and $s$ defines a character $\Lcal_s$ of $\Gbf^{\Frob}$.
\begin{lem}\label{lem:endocsopy-rep-trivial}
Assume that $s \in \Gbf^*$ is central, then the functor $\EE^{\Frob}_{\ofrak}$ is isomorphic to 
$$\DD^{\ofrak}(\frac{\Gbf}{\Ad_{\Frob}\Gbf}) \xrightarrow{\Lcal_s} \DD^{\unip}(\frac{\Gbf}{\Ad_{\Frob}\Gbf}).$$
\end{lem} 

\begin{proof}
This is clear by \Cref{lem:endoscopy-trivial-case}.
\end{proof}

Assume that $\Lbf^* = \Zbf(s) \subset \Gbf^*$ is a (connected) Levi subgroup of $\Gbf$, then \cite{BonnafeRouquier} and \cite{BonnafeDatRouquier} have defined an equivalence 
$$\EE^{\mathrm{BDR}} : \DD^{\ofrak}(\frac{\Gbf}{\Ad_{\Frob}\Gbf}) \cong \DD^{\unip}(\frac{\Hbf_s}{\Ad_{\dot{v}\Frob}\Hbf_s})$$ 
given by 
$$\DD^{\unip}(\frac{\Lbf}{\Ad_{\dot{v}\Frob}\Lbf}) \xrightarrow{\otimes \Lcal_s^{\Lbf}} \DD^{\ofrak^{\Lbf}}(\frac{\Lbf}{\Ad_{\dot{v}\Frob}\Lbf}) \xrightarrow{\mathrm{R}_{\Pbf, \dot{v}}} \DD^{\ofrak}(\frac{\Gbf}{\Ad_{\Frob}\Gbf}).$$

\begin{thm}\label{thm:comparison-with-BDR}
Assume that $\Lbf^* = \Zbf(s) \subset \Gbf^*$ is a (connected) Levi subgroup of $\Gbf$ then 
$$\EE^{\mathrm{BDR}} \cong \EE^{\Frob}_{\ofrak}.$$
\end{thm}

\begin{proof}
This is direct application of \Cref{lem:endocsopy-rep-trivial} and \Cref{thm:compat-DL-induction}. 
\end{proof}

\subsection{A variant of Theorem \ref{thm:EndoscopyRep}}\label{sec:Variant}

We keep the hypothesis of the previous section. Let $s \in (\Gbf^*)^{\Frob^*}$ be a semisimple element and choose an $\Frob^*$-maximal torus $\Sbf^*$ such that $s \in \Sbf^*$ and a pinning $(\Sbf^*, \tilde{\Bbf}^*, \Gbf^*, X'^*)$ of $\Gbf^*$. Let $\Sbf \subset \Gbf$ be an $\Frob$-maximal torus such that $\Sbf^*$ is identified with the dual of $\Sbf$ and fix a pinning $(\Sbf,\tilde{\Bbf}, \Gbf, X')$ of $\Gbf$ dual to the pinning of $\Gbf^*$. Let $\Gbf_s$ be the endoscopic group of $s$ obtained as in Section \ref{sec:EndoscopicGroup}. We have 
\begin{enumerate}
\item $\pi_0(\Gbf_s) = \Gammabf_s$, 
\item $\Gbf_s^{\circ, *} = \Zbf^{\circ}(s)$. 
\end{enumerate}
This group is defined canonically through its pinning relative to $\Gbf$. Since $s$ is fixed by $\Frob^*$, $\Zbf(s)$ is stable under $\Frob^*$ and $\Gbf_s$ acquires an isogeny 
$$\Frob : \Gbf_s \to \Gbf_s$$
dual to $\Frob^*$. 

\begin{thm}\label{thm:EndoscopyRepCanonical}
There is a $\DD(\frac{\Gammabf_s}{\Ad_{\Frob}(\Gammabf_s)})$-linear equivalence 
$$\DD^{\ofrak}(\pt/\Gbf^{\Frob}, \Lambda) \cong \DD^{\unip}(\frac{\Gbf_s}{\Ad_{\Frob}\Gbf_s}, \Lambda).$$
\end{thm}

\begin{proof}
We identify both sides of \Cref{thm:EndoscopyRepCanonical} with those of \Cref{thm:EndoscopyRep}; the left-hand side requires no modification. There is a unique element $\bar g^*\in\Gbf^*_{\ad}$ carrying the pinning $(\Sbf^*,\tilde{\Bbf}^*,\Gbf^*,X'^*)$ to $(\Tbf^*,\Bbf^*,\Gbf^*,X^*)$. Choose a lift $g^*\in\Gbf^*$ and put $t=g^*sg^{*,-1}$. Likewise, there is a unique $\bar g\in\Gbf_{\ad}$ carrying $(\Sbf,\tilde{\Bbf},\Gbf,X')$ to $(\Tbf,\Bbf,\Gbf,X)$, and we choose a lift $g\in\Gbf$. Two lifts differ by a central element, so they induce the same conjugation maps and all equivalences below are independent of these choices. The resulting conjugation map yields an isomorphism of relatively pinned reductive groups
$$\Hbf_t = \Gbf_s$$
and an isomorphism of $\Hbf_t = \Gbf_s$-bimodules
$${_t}\Hbf_{\Frob(t)} = \Gbf_{s,\Frob}$$
where $\Gbf_{s,\Frob}$ is the $\Gbf_s$-bimodule where the right multiplication is twisted by $\Frob$. 
There is therefore an equivalence 
$$\DD(\frac{{_t}\Hbf_{\Frob(t)}}{\Ad_{\Frob}\Hbf_t}) \cong \DD(\frac{\Gbf_s}{\Ad_{\Frob}\Gbf_s}).$$
Combining this with Theorem \ref{thm:EndoscopyRep} yields the desired equivalence. 
The linearity over $\DD(\frac{\Gammabf_s}{\Ad_{\Frob}(\Gammabf_s)})$ is immediate. 
\end{proof}

\begin{rque}
Using the $\DD(\frac{\Gammabf_s}{\Ad_{\Frob}(\Gammabf_s)})$-linearity, we obtain, as in Corollary \ref{corol:RationalEndoscopy}, an equivalence between each rational series and the category of unipotent representations of a (possibly disconnected) reductive group.
\end{rque}

\subsection{Properties of the endoscopic equivalence}\label{sec:PropertiesEndoscopic}

In this section we assume that $\ell$ is not a torsion prime for $\Gbf^*$ if $\Lambda = \Zlb, \Flb$. 

\subsubsection{$t$-exactness}

\begin{prop}[$t$-exactness]\label{lem:tExactness}
The equivalence $\EE_{\ofrak}^{\Frob}$ of Theorem \ref{thm:EndoscopyRep} is $t$-exact.
\end{prop}

\begin{proof}
The Deligne--Lusztig varieties are affine: this follows from \cite{DeligneLusztig} when $q>h$, where $h$ is the Coxeter number of $\Gbf$, and from \cite{DeligneLusztigAreAffine} in general. The objects $\tr_{\Frob}(\Delta_w)$ lie in $\DD^{\geq 0}(\pt/\Gbf^{\Frob})$ and generate this category under colimits and extensions. Similarly, $\DD^{\leq 0}(\pt/\Gbf^{\Frob})$ is generated by the objects $\tr_{\Frob}(\nabla_w)$ under limits and extensions. The same statements hold for $\Hbf$. By \Cref{thm:EndoscopyHeckeCat}, the equivalence $\EE_{\ofrak}$ sends (co)standard objects to (co)standard objects. Therefore, by the commutativity of diagram \ref{eq:diagram-commutation-into-trace-endoscopy}, the equivalence $\EE_{\ofrak}^{\Frob}$ sends $\DD^{? 0, \ofrak}(\pt/\Gbf^{\Frob})$ into $\DD^{? 0, \unip}(\Lcal_{\Frob}\Hbffrak_{\ofrak})$, where $? \in \{\leq, \geq\}$. It follows that $\EE_{\ofrak}^{\Frob}$ is $t$-exact.
\end{proof}

\subsubsection{Inversion of $\ell$}\label{sec:compatibility-with-inversion-of-ell}

In this subsection we assume that $\Gbf$ has connected center. This hypothesis is most likely superfluous but it makes the statement much clearer, we leave it to the interested reader to generalize Proposition \ref{prop:InversionOfEll} to the disconnected case. We fix a point $\dot{\beta} \in \Nbf(\Tbf)$ lifting the minimal element of ${_s}\Wbf_{\Frob(s)}$. 

Let $s \in \Ch_{\Zlb}(\Tbf)$ and consider the endoscopic equivalence of Theorem \ref{thm:EndoscopyRep}
$$\DD^{\ofrak}(\pt/\Gbf^{\Frob}, \Zlb) \cong \DD^{\unip}(\pt/\Hbf^{\Frob_{\dot{\beta}}}, \Zlb).$$
After inverting $\ell$, we get an equivalence 
\begin{equation}\label{eq:DecompoInvertell}
\DD^{\ofrak}(\pt/\Gbf^{\Frob}, \Zlb)[\frac{1}{\ell}] \cong  \DD^{\unip}(\pt/\Hbf^{\Frob_{\dot{\beta}}}, \Zlb)[\frac{1}{\ell}].
\end{equation}
Both sides decompose as
\begin{enumerate}
\item $\DD^{\ofrak}(\pt/\Gbf^{\Frob}, \Zlb)[\frac{1}{\ell}] = \bigoplus_{s' \in (\Gbf^*_{\mathrm{ss}})^{\Frob^*}/\sim_{\Qlb, \geo}, s' \sim_{\Zlb,\geo} s} \DD^{\ofrak'}(\pt/\Gbf^{\Frob}, \Qlb)$, 
\item $\DD^{\unip}(\pt/\Hbf^{\Frob_{\dot{\beta}}}, \Zlb)[\frac{1}{\ell}] = \bigoplus_{s' \in (\Hbf^*_{\mathrm{ss}})^{\Frob_{\dot{\beta}}^*}/\sim_{\Qlb, \geo}, s' \sim_{\Zlb,\geo} 1} \DD^{\ofrak'}(\pt/\Hbf^{\Frob_{\dot{\beta}}}, \Qlb).$
\end{enumerate}

Since $\Hbf^* = \Zbf(s) \subset \Gbf^*$, multiplication by $s$ yields a map 
\begin{equation}\label{eq:CompareHvsG}
\{s' \in (\Hbf^*_{\mathrm{ss}})^{\Frob_{\dot{\beta}}^*}/\sim_{\Qlb, \geo}, s' \sim_{\Zlb,\geo} 1\} \to \{s' \in (\Gbf^*_{\mathrm{ss}})^{\Frob^*}/\sim_{\Qlb, \geo}, s' \sim_{\Zlb,\geo} s\}.
\end{equation}

\begin{lem}
The map \ref{eq:CompareHvsG} is an isomorphism. 
\end{lem}

\begin{proof}
Let $s' \in \Hbf_s^*$ be a point in the LHS, then $s'$ is a torsion point of order a power of $\ell$. Conversely let $s' \in \Gbf^*$ be a point in the RHS, then $s' = ss_{\ell}$ where $s_{\ell}$ is a semisimple point of order a power of $\ell$. Since $s$ has order prime to $\ell$, the lemma follows from the observation that $\Zbf_{\Zbf_{\Gbf^*}(s)}(s_{\ell}) = \Zbf_{\Gbf^*}(ss_{\ell})$.
\end{proof}

\begin{lem}
Under the bijection \ref{eq:CompareHvsG}, the equivalence \ref{eq:DecompoInvertell} preserves each summand.
\end{lem}

\begin{proof}
Both categories appearing in the equivalence \ref{eq:DecompoInvertell} are linear over $((\Tbf^{\vee}_{\Zlb}\sslash \Wbf)^{\Frob^{\vee}})[\frac{1}{\ell}]$ and the decomposition into direct summands comes from the various connected components of this scheme.
\end{proof}

Let $s' \in (\Hbf^*_{s,\mathrm{ss}})^{\Frob_{\dot{\beta}}^*}$ be such that $s' \sim_{\Qlb, \geo} 1$ in $\Hbf^*_s$ and let $\Hbf_{s'}$ be the endoscopic group of $s'$ relative to $\Hbf$. This group is relatively pinned to $\Hbf_s$ and thus relatively pinned to $\Gbf$. Correspondingly, $\Hbf_{s'}$ is also the endoscopic group for the point $ss' \in \Gbf^*$ pinned relatively to $\Gbf$. 

\begin{prop}[Inversion of $\ell$]\label{prop:InversionOfEll}
There is a commutative diagram 
\[\begin{tikzcd}
	{\DD^{\ofrak}(\pt/\Gbf^{\Frob}, \Zlb)[\frac{1}{\ell}]} & {\DD^{\unip}(\pt/\Hbf_s^{\Frob_{\dot{\beta}}}, \Zlb)[\frac{1}{\ell}]} \\
	{\DD^{s\cdot\ofrak'}(\pt/\Gbf^{\Frob}, \Qlb)} & {\DD^{\ofrak'}(\pt/\Hbf_s^{\Frob_{\dot{\beta}}}, \Qlb)} \\
	& {\DD^{\unip}(\pt/\Hbf_{s'}^{\Frob_{\dot{\beta}'}}, \Qlb)}
	\arrow["{\EE_s^{\Frob}}", from=1-1, to=1-2]
	\arrow[from=1-1, to=2-1]
	\arrow[from=1-2, to=2-2]
	\arrow["{e_{\ofrak',\Qlb}\EE_s^{\Frob}}", from=2-1, to=2-2]
	\arrow["{\EE_{s\cdot\ofrak'}^{\Frob}}"', from=2-1, to=3-2]
	\arrow["{\EE_{\ofrak'}^{\Frob}}", from=2-2, to=3-2]
\end{tikzcd}\]
where the top vertical maps are the projections on the corresponding direct summands and the two vertical bottom maps are the endoscopic equivalence for $\Lambda = \Qlb$ and $ss'$ in $\Gbf$ and $s'$ in $\Hbf$.
\end{prop}

\begin{proof}
This follows from the corresponding statement about Hecke categories and applying traces to Lemma \ref{lem:CompatInvertEll}. 
\end{proof}

\subsection{Compatibility with automorphisms of $\Gbf^{\Frob}$}\label{sec:compat-with-automorphisms}

We now study how automorphisms of $\Gbf^{\Frob}$ act on the Jordan decomposition. 

\subsubsection{Outer automorphisms of $\Gbf^{\Frob}$}
\begin{defi}
We let $\Isog(\Gbf, \Frob)$ be the monoid of isogenies $\tau : \Gbf \to \Gbf$ that commute with $\Frob$, preserve the pair $(\Bbf,\Tbf)$, and have purely inseparable kernel.
\end{defi}

There are three specific collections of isogenies in $\Isog(\Gbf, \Frob)$ that we can distinguish : 
\begin{enumerate}
\item graph automorphisms $\Gbf \to \Gbf$,
\item inner automorphisms $\Ad(t) : \Gbf \to \Gbf$ for $t \in \Tbf_{\ad}^{\Frob}$, 
\item isogenies coming from roots of $\Frob$.  
\end{enumerate}
Any element in $\Isog(\Gbf, \Frob)$ restricts to an automorphism $\Gbf^{\Frob} \to \Gbf^{\Frob}$. 

\begin{thm}[\protect{\cite[2.5.12]{AutomFiniteGroupsLieType}}]
Assume that $\Gbf$ is almost simple and simply connected, the map 
$$\Isog(\Gbf, \Frob) \to \Aut(\Gbf^{\Frob}) \to \Out(\Gbf^{\Frob})$$
is surjective. 
\end{thm} 

There is a natural map of sets 
\begin{align*}
\Isog(\Gbf, \Frob) &\to \Tbf_{\ad}^{\Frob} \\
\tau &\mapsto t_{\tau},
\end{align*}
that sends $\tau$ to the element $t_{\tau}$ giving the relative position of the two pinnings $(\phi, \tau(\phi))$, i.e. $\tau(\phi) = \phi_{t_{\tau}}$. 

\subsubsection{Statement}
Let us fix a $\Wbf$-orbit $\ofrak \subset \Ch(\Tbf)$ and $s \in \ofrak$. Recall that there is a canonical embedding, see Lemma \ref{lem:ControlCentralizer}, 
$$\Gammabf_{s} \subset \pi_1(\Gbf^*_{\der}).$$

\begin{lem}\label{lem:invariance-of-embedding}
For any $s,s' \in \ofrak$, we have $\Gammabf_s = \Gammabf_{s'}$ as subgroups of $\pi_1(\Gbf^*_{\der})$ (but not necessarily as subquotients of $\Wbf$). 
\end{lem}

\begin{proof}
Let $\gamma : s \to s'$ be a morphism in $\Xibf$. Let $\dot{\gamma} \in \Nbf(\Tbf)$ be any lift of $\gamma$. Conjugation by $\dot\gamma$ identifies the two component groups, so inside $\Zbf(\Gbf^*_{\mathrm{sc}})$ we have
$$\Ad(\dot{\gamma})(\Gamma_s)=\Gamma_{s'}.$$
But $\Ad(\dot{\gamma})$ acts trivially on the center $\Zbf(\Gbf^*_{\mathrm{sc}})$. Hence $\Gamma_s=\Gamma_{s'}$ as subgroups of $\pi_1(\Gbf^*_{\der})$.
\end{proof}

\begin{corol}
As a subgroup $\pi_1(\Gbf^*_{\der})$, the group $\Gammabf_s$ is $\Frob$-stable. 
\end{corol}

\begin{proof}
Immediate from \Cref{lem:invariance-of-embedding} and the fact that $\ofrak$ is $\Frob$-stable. 
\end{proof}

Let $\beta : s \to \Frob(s)$ be a morphism in $\Xibf_{\ofrak}$. As subgroups of $\Gammabf_{s}$, we have
$$\Gammabf_s^{\beta\Frob} = \Gammabf_{s}^{\Frob}$$
where the RHS is obtained by restricting the Frobenius of $\pi_{1}(\Gbf^*_{\der})$. 

Recall that we have a canonical perfect pairing $\pi_0(\Zbf(\Gbf)) \times \pi_1(\Gbf_{\der}^*) \to \Q/\Z$. This induces a canonical perfect pairing
$$H^1(\Frob, \pi_0(\Zbf(\Gbf))) \times \pi_1(\Gbf^*_{\der})^{\Frob} \to \Q/\Z.$$
We thus obtain a natural map 
\begin{equation}\label{eq:map-defining-twisting}
\Tbf_{\ad}^{\Frob} \to H^1(\Frob,  \pi_0(\Zbf(\Gbf))) \to \Hom(\pi_1(\Gbf^*_{\der})^{\Frob}, \Q/\Z) \to \Hom(\Gammabf_s^{\Frob}, \Q/\Z),
\end{equation}
where the last map is the restriction along $\Gammabf_s^{\Frob} \subset \pi_1(\Gbf^*_{\der})^{\Frob}$. 

\begin{thm}\label{thm:compat-automorphisms}
\begin{enumerate}
\item There is a canonical rank one local system $\Lcal_{\tau}$ on $\Lcal_{\Frob}\Xibf_{\ofrak}^{\opp}$ whose restriction to $\pt/\Gammabf_s^{\Frob} \xrightarrow{(s, \beta)} \Lcal_{\Frob}\Xibf_{\ofrak}^{\opp}$, is the character of $\Gammabf_s^{\Frob}$ given by the image of $t_{\tau}$ under \Cref{eq:map-defining-twisting}. 
\item There is a commutative diagram 
\[\begin{tikzcd}
	{\DD^{\ofrak}(\frac{\Gbf}{\Ad_{\Frob}\Gbf})} & {\DD^{\unip}(\Lcal_{\Frob}\Hbffrak_{\ofrak})} \\
	{\DD^{\tau(\ofrak)}(\frac{\Gbf}{\Ad_{\Frob}\Gbf})} & {\DD^{\unip}(\Lcal_{\Frob}\Hbffrak_{\tau(\ofrak)})}
	\arrow["{\EE^{\Frob}_{\ofrak}}", equals, from=1-1, to=1-2]
	\arrow["{\tau^*}", from=1-1, to=2-1]
	\arrow["{\tau^*}", from=1-2, to=2-2]
	\arrow["{\Lcal_{\tau} \otimes \EE^{\Frob}_{\tau(\ofrak)}}"', equals, from=2-1, to=2-2]
\end{tikzcd}\]
\end{enumerate}
\end{thm} 

\begin{rque}
It follows from \Cref{thm:compat-automorphisms} that the induced bijection 
$$\Irr_{\Lambda}^{\ofrak}(\frac{\Gbf}{\Ad_{\Frob}\Gbf}) \cong \Irr_{\Lambda}^{\unip}(\Lcal_{\Frob}\Hbffrak_{\ofrak})$$
is equivariant for the image of $\Isog(\Gbf,\Frob)$ in $\Out(\Gbf^{\Frob})$. This recovers most of \cite[Theorem B]{SpaethTypeD}.
\end{rque}

\subsubsection{Proof of \Cref{thm:compat-automorphisms}} 

In this section, we abbreviate $t=t_{\tau}$. Recall that the endoscopic equivalence depends on the pinning $\phi$ of $\Gbf$, and that we write $\EE_{\phi}$ when we want to emphasize this dependence (see \Cref{sec:changing-the-pinning}).

\begin{lem}\label{lem:diagram-computation-character}
There is a natural commutative diagram of functors $\Xibf_{\ofrak}^{\opp} \to \Morita$
\[\begin{tikzcd}
	{\HH^{\Gbf}_{\ofrak}} && {\HH^{\Hbf}_{\ofrak}} \\
	{\HH^{\Gbf}_{\tau(\ofrak)}} && {\HH^{\Hbf}_{\tau(\ofrak)}} \\
	{\HH^{\Gbf}_{\tau(\ofrak)}} && {\HH^{\Hbf}_{\tau(\ofrak)}}
	\arrow["{{\EE_{\ofrak, \phi}}}", from=1-1, to=1-3]
	\arrow["{{\tau^*}}"', from=1-1, to=2-1]
	\arrow["{{\tau^*}}"{description}, from=1-3, to=2-3]
	\arrow["{{\EE_{\tau(\ofrak), \tau(\phi)}}}"', from=2-1, to=2-3]
	\arrow[equals, from=2-1, to=3-1]
	\arrow[equals, from=2-3, to=3-3]
	\arrow["{{\Lcal_t \otimes \EE_{\tau(\ofrak), \phi}}}"', from=3-1, to=3-3]
\end{tikzcd}\]
\end{lem}

\begin{proof}
The first part is the functoriality of the construction and the second part is \Cref{lem:changing-pinning}. 
\end{proof}

In the diagram of \Cref{lem:diagram-computation-character}, all the terms are compatible with $\Frob$ by assumption on $\tau$. It follows that there is a commutative diagram of functors $\Lcal_{\Frob}\Xibf_{\ofrak}^{\opp} \to \Pr_{\Lambda}$ : 
\[\begin{tikzcd}
	{\Tr^{\loc}(\Frob^*, \HH^{\Gbf}_{\ofrak})} &&& {\Tr^{\loc}(\Frob^*, \HH^{\Hbf}_{\ofrak})} \\
	{\Tr^{\loc}(\Frob^*, \HH^{\Gbf}_{\tau(\ofrak)})} &&& {\Tr^{\loc}(\Frob^*, \HH^{\Hbf}_{\tau(\ofrak)})}
	\arrow["{{\EE_{\ofrak}^{\Frob}}}", from=1-1, to=1-4]
	\arrow["{{ \tau^*}}"', from=1-1, to=2-1]
	\arrow["{{\tau^*}}", from=1-4, to=2-4]
	\arrow["{{\EE_{\tau(\ofrak)}^{\Frob} \otimes \Tr^{\loc}(\Frob, \Lcal_t)}}"', from=2-1, to=2-4]
\end{tikzcd}\]

Using \Cref{thm:local-to-global-trace}, it is enough to compute $\Tr^{\loc}(\Frob, \Lcal_t)$. We can view the natural transformation $\otimes \Lcal_t$ as a natural transformation $\id_C \Rightarrow \id_C$ where $C : \Xibf_{\ofrak}^{\opp} \to \Morita$ is the constant diagram with values $\DD(\Lambda)$. The Frobenius equivariance of $\Lcal_t$ amounts to the datum of a collection of compatible isomorphisms 
$$\Frob_{\beta} : \Lcal^{\beta} \xrightarrow{\sim} \Lcal^{\Frob(\beta)}.$$

Passing to the trace, we get that $\Tr^{\loc}(\Frob, \Lcal_t)$ is a natural transformation $\id_{\Tr^{\loc}(\Frob, C)} \to \id_{\Tr^{\loc}(\Frob, C)}$. Unwinding, $\Tr^{\loc}(\Frob, C) : \Lcal_{\Frob}\Xibf_{\ofrak}^{\mathrm{op}} \to \Pr_{\Lambda}$ is the constant diagram. It follows that $$\Tr^{\loc}(\Frob, \Lcal_t)(s, \beta): \Tr^{\loc}(\Frob, C)(s,\beta) \cong \DD(\Lambda) \to \DD(\Lambda) \cong \Tr^{\loc}(\Frob, C)(s,\beta) $$
is the datum of an object in $\DD(\Lambda)$, which is nothing else than $\Lcal_t^{\beta}$ itself by $\DD(\Lambda)$-linearity of the traces we consider. The naturality of $\Tr^{\loc}(\Frob, \Lcal_t)$ amounts to a collection of compatible isomorphisms 
$$\Lcal^{\beta}_t \xrightarrow{t_{\alpha}} \Lcal^{\beta'}_t$$
for all maps $\alpha : (s, \beta) \to (s', \beta')$ in $\Lcal_{\Frob}\Xibf_{\ofrak}^{\opp}$. The collection of all the maps $t_{\alpha}$ is equivalent to the datum of a rank one local system on $\Lcal_{\Frob}\Xibf_{\ofrak}$ which yields the desired local system. It remains to identify this local system. Let $\alpha : (s, \beta) \to (s',\beta')$ and consider the diagram 
\[\begin{tikzcd}
	s & {\Frob(s)} \\
	{s'} & {\Frob(s')}
	\arrow["\beta", from=1-1, to=1-2]
	\arrow["\alpha"', from=1-1, to=2-1]
	\arrow["{\Frob(\alpha)}", from=1-2, to=2-2]
	\arrow["{\beta'}"', from=2-1, to=2-2]
\end{tikzcd}\]

 The map $t_{\alpha}$ is given by the following composition 
\begin{align*}
\Lcal_t^{\beta} &\cong \Lcal_t^{\beta}\otimes \Lcal_t^{\alpha} \otimes \Lcal_t^{\alpha^{-1}} \\
&\xrightarrow{\Frob_{\alpha}} \Lcal_t^{\beta}\otimes \Lcal_t^{\alpha} \otimes \Lcal_t^{\Frob(\alpha^{-1})} \\
&\cong \Lcal_t^{\beta'}.
\end{align*}

Assume that $s = s'$ and $\beta = \beta'$, so that $\alpha \in \Gammabf_s^{\Frob}$. It follows that $t_{\alpha} = \Tr(\Frob_{\alpha}, \Lcal^{\alpha}_t)$. Since $\Lcal_t^{\alpha}$ is the fiber at $t$ of a local system $\Lcal^{\alpha}$ on $\Tbf_{\ad}$ induced by 
$$X_*(\Tbf_{\ad}) \xrightarrow{\ev(-,\alpha)} \Lambda^{\times}$$ 
by compatibility with Grothendieck function sheaf dictionary, $t_{\alpha}$ is the image of $t$ under the map \ref{eq:map-defining-twisting}.

\subsection{Compatibility with isogenies}\label{sec:compatibility-with-isogenies}
We let $f : \Gbf \to \Gbf'$ be a finite central isogeny that commutes with $\Frob$. We shall use the notations of Section \ref{sec:PrelimEndoscopicGroupsAndIsogeny}. For this entire section, we fix $s' \in \Ch(\Tbf')$, with image $s \in \Ch(\Tbf)$ and we denote by $\ofrak' \subset \Ch(\Tbf')$ and $\ofrak \subset \Ch(\Tbf)$ their corresponding $\Wbf$-orbits. 

Since $f$ commutes with $\Frob$, the monoidal functor $f_*$ commutes with $\Frob_*$ and we have an induced map on traces 
$$\Tr(\Frob^*, \HH^{\Gbf}) \xrightarrow{\Tr(f_*)} \Tr(\Frob^*, \HH^{\Gbf'}).$$ We shall also denote by $f$ the map 
$$\pt/\Gbf^{\Frob} \to \pt/(\Gbf')^{\Frob}.$$

\begin{lem}\label{lem:calculation-pushfoward-trace}
There is a commutative diagram 
\[\begin{tikzcd}
	{\Tr(\Frob^*, \HH^{\Gbf})} & {\DD(\pt/\Gbf^{\Frob})} \\
	{\Tr(\Frob^*, \HH^{\Gbf'})} & {\DD(\pt/(\Gbf')^{\Frob})}
	\arrow["\sim", from=1-1, to=1-2]
	\arrow["{{\Tr(f_*)}}", from=1-1, to=2-1]
	\arrow["{{f_*}}", from=1-2, to=2-2]
	\arrow["\sim"', from=2-1, to=2-2]
\end{tikzcd}\]
\end{lem}

\begin{proof}
The map $f_*$ is monoidal hence we view $\HH^{\Gbf'}$ as a $(\HH^{\Gbf'},\HH^{\Gbf})$-bimodule, its dual is, by rigidity, again $\HH^{\Gbf'}$. The induced map on traces is given by 
$$\HH^{\Gbf'} \otimes_{\HH^{\Gbf}} \HH^{\Gbf}_{\Frob} \cong \HH^{\Gbf'}_{\Frob} \xrightarrow{\id} \HH^{\Gbf'}_{\Frob} \cong \HH^{\Gbf'}_{\Frob} \otimes_{\HH^{\Gbf'}} \HH^{\Gbf'}.$$
The unit and counit are then given by 
$$\unit = f_* : \HH^{\Gbf} \to \HH^{\Gbf'} \otimes_{\HH^{\Gbf'}} \HH^{\Gbf'} \cong \HH^{\Gbf'}$$
and 
$$\counit = \mult : \HH^{\Gbf'} \otimes_{\HH^{\Gbf}} \HH^{\Gbf'} \rightarrow \HH^{\Gbf'}$$
The induced functoriality is then given by 
\begin{align*}
\Tr(\HH^{\Gbf}_{\Frob},\HH^{\Gbf}) &\xrightarrow{\unit} \Tr(\HH^{\Gbf}_{\Frob} \otimes_{\HH^{\Gbf}} \HH^{\Gbf'} \otimes_{\HH^{\Gbf'}} \HH^{\Gbf'},\HH^{\Gbf}) \\
&\xrightarrow{\cycl} \Tr(\HH^{\Gbf'} \otimes_{\HH^{\Gbf}} \HH^{\Gbf}_{\Frob} \otimes_{\HH^{\Gbf}} \HH^{\Gbf'},\HH^{\Gbf'}) \\
&\xrightarrow{\mathrm{commut}} \Tr(\HH^{\Gbf'}_{\Frob} \otimes_{\HH^{\Gbf'}} \HH^{\Gbf'} \otimes_{\HH^{\Gbf}} \HH^{\Gbf'},\HH^{\Gbf'}) \\
&\xrightarrow{\counit} \Tr(\HH^{\Gbf'}_{\Frob},\HH^{\Gbf'}).
\end{align*}
Using \Cref{thm:calculation-traces} and \Cref{lem:essential-surjectivity-of-trace-map}, this series of functors is equivalent to the following composition in the top row of the following diagram 
\[\begin{tikzcd}
	{\DD(\frac{\Gbf}{\Ad_{\Frob}\Gbf})} & {\DD(\frac{\Gbf \times^{\Gbf}\Gbf' \times^{\Gbf'}\Gbf'}{\Ad_{\Frob}\Gbf})} & {\DD(\frac{\Gbf \times^{\Gbf}\Gbf' \times^{\Gbf}\Gbf'}{\Ad_{\Frob}\Gbf'})} & {\DD(\frac{\Gbf' \times^{\Gbf'}\Gbf' \times^{\Gbf}\Gbf'}{\Ad_{\Frob}\Gbf})} & {\DD(\frac{\Gbf'}{\Ad_{\Frob}\Gbf'})} \\
	{\DD(\frac{\Gbf}{\Ad_{\Frob}\Gbf})} & {\DD(\frac{\Gbf'}{\Ad_{\Frob}\Gbf})} & {\DD(\frac{\Gbf'}{\Ad_{\Frob}\Gbf})} & {\DD(\frac{\Gbf'}{\Ad_{\Frob}\Gbf})} & {\DD(\frac{\Gbf'}{\Ad_{\Frob}\Gbf'})}
	\arrow["{\unit}", from=1-1, to=1-2]
	\arrow[equals, from=1-1, to=2-1]
	\arrow["\cycl", from=1-2, to=1-3]
	\arrow["\wr"{description}, from=1-2, to=2-2]
	\arrow["{\mathrm{commut}}", from=1-3, to=1-4]
	\arrow["\wr"{description}, from=1-3, to=2-3]
	\arrow["\counit", from=1-4, to=1-5]
	\arrow["\wr"{description}, from=1-4, to=2-4]
	\arrow[equals, from=1-5, to=2-5]
	\arrow["{f_*}"', from=2-1, to=2-2]
	\arrow[equals, from=2-2, to=2-3]
	\arrow[equals, from=2-3, to=2-4]
	\arrow["{f_*}"', from=2-4, to=2-5]
\end{tikzcd}\]

\end{proof}

Let us now restrict to the orbit $\ofrak$ and analyze the corresponding map under the Morita equivalence of \Cref{lem:Morita-equivalence-for-bimodule-for-trace} and the bimodule $\HH^{\Gbf'}_{[s',s']}$ constructed in \Cref{sec:compat-with-isog-hecke-cat} (beware that we have switched $s$ and $s'$). Consider the map 
$$\beta^{\Gbf} : \HH^{\Gbf'}_{[s',\Frob(s')]} \otimes_{\HH^{\Gbf'}_{[s',s']}} \HH^{\Gbf'}_{[s',s]} \xrightarrow{\sim}  \HH^{\Gbf'}_{[s',s]} \otimes_{\HH^{\Gbf}_{[s,s]}} \HH^{\Gbf}_{[s,\Frob(s)]}$$ 
defined so that the following diagram is commutative 
\[\begin{tikzcd}
	{ \HH^{\Gbf'}_{[s',\Frob(s')]} \otimes_{\HH^{\Gbf'}_{[s',s']}} \HH^{\Gbf'}_{[s',s]}} & { \HH^{\Gbf'}_{[s',\Frob(s')]} \otimes_{\HH^{\Gbf'}_{[\Frob(s'),\Frob(s')]}} \HH^{\Gbf'}_{[\Frob(s'),\Frob(s)]}} \\
	{\HH^{\Gbf'}_{[s',s]} \otimes_{\HH^{\Gbf}_{[s,s]}} \HH^{\Gbf}_{[s,\Frob(s)]}} & {\HH^{\Gbf'}_{[s',\Frob(s)]}}
	\arrow["{\id \otimes\Frob^*}", from=1-1, to=1-2]
	\arrow["\beta"', from=1-1, to=2-1]
	\arrow["{\mathrm{act}}", from=1-2, to=2-2]
	\arrow["{\mathrm{act}}"', from=2-1, to=2-2]
\end{tikzcd}\]
where $\mathrm{act}$ is the action map coming from the bimodule structure on $\HH^{\Gbf'}_{[s',s]}$. Note that all three maps different from $\beta$ are isomorphisms.

On the endoscopic side, we have a map
$$\Hbf'_{[s',\Frob(s')]} \times^{\Hbf_{[s',s']}'} \Hbf'_{[s',s]} \xrightarrow{\id \times \Frob} \Hbf'_{[s',\Frob(s')]} \times^{\Hbf_{[\Frob(s'),\Frob(s')]}'} \Hbf'_{[\Frob(s'),\Frob(s)]} \cong \Hbf'_{[s',\Frob(s)]} \cong \Hbf'_{[s',s]} \times^{\Hbf_{[s,s]}} \Hbf_{[s,\Frob(s)]}.$$
Pushing forward along this map yields a map of bimodules
$$\beta^{\Hbf} : \HH^{\Hbf'}_{[s',\Frob(s')]} \otimes_{\HH^{\Hbf'}_{[s',s']}} \HH^{\Hbf'}_{[s',s]} \xrightarrow{\sim} \HH^{\Hbf'}_{[s',s]} \otimes_{\HH^{\Hbf}_{[s,s]}} \HH^{\Hbf}_{[s,\Frob(s)]}.$$

We then obtain a commutative diagram of pairs 
\[\begin{tikzcd}
	{(\HH^{\Gbf}_{[s,s]}, \HH^{\Gbf}_{[s,\Frob(s)]})} & {(\HH^{\Hbf}_{[s,s]}, \HH^{\Hbf}_{[s,\Frob(s)]})} \\
	{(\HH^{\Gbf'}_{[s',s']}, \HH^{\Gbf'}_{[s',\Frob(s')]})} & {(\HH^{\Hbf'}_{[s',s']}, \HH^{\Hbf'}_{[s',\Frob(s')]})}
	\arrow["\sim"{description}, from=1-1, to=1-2]
	\arrow["{(\HH^{\Gbf'}_{[s',s]}, \beta^{\Gbf})}"', from=1-1, to=2-1]
	\arrow["{(\HH^{\Hbf'}_{[s',s]}, \beta^{\Hbf})}", from=1-2, to=2-2]
	\arrow["\sim"{description}, from=2-1, to=2-2]
\end{tikzcd}\]\\
Passing to traces, we obtain the following commutative diagram:
\[\begin{tikzcd}
	{\DD^{\ofrak}(\pt/\Gbf^{\Frob})} & {\DD^{\unip}(\frac{\Hbf_{[s,\Frob(s)]}}{\Ad_{\Frob}\Hbf_{[s,s]}})} \\
	{\DD^{\ofrak'}(\pt/\Gbf^{'\Frob})} & {\DD^{\unip}(\frac{\Hbf'_{[s',\Frob(s')]}}{\Ad_{\Frob}\Hbf'_{[s',s']}}).}
	\arrow["\sim"{description}, from=1-1, to=1-2]
	\arrow["{f_*}", from=1-1, to=2-1]
	\arrow["{\Tr(\HH^{\Hbf'}_{[s',s]}, \beta^{\Hbf})}", from=1-2, to=2-2]
	\arrow["\sim"{description}, from=2-1, to=2-2]
\end{tikzcd}\]
Finally, consider the scheme $\Hbf_{[s,\Frob(s)]} \times \Hbf'_{[s',\Frob(s')]} \times \Hbf_{[s',s]}$ seen as a $(\Hbf_{[s,s]}, \Hbf'_{[s',s']})$-bimodule for the action given by 
$$h'\cdot(x,y,z)\cdot h=\bigl(\Frob(h)xh^{-1},\,\Frob(h')y(h')^{-1},\,h'zh^{-1}\bigr).$$
Let $\Xbf_{[s',s]}$ be the substack of $(\Hbf_{[s,\Frob(s)]} \times \Hbf'_{[s',\Frob(s')]} \times \Hbf_{[s',s]})/(\Hbf_{[s,s]} \times \Hbf'_{[s',s']})$ cut out by the equation 
$$\Frob(z)x = yz$$
where $(x,y,z) \in \Hbf_{[s,\Frob(s)]} \times \Hbf'_{[s',\Frob(s')]} \times \Hbf_{[s',s]}$.
This equation is invariant under the displayed action: writing $(x',y',z')$ for the image triple, one has
$$\Frob(z')x'=\Frob(h')\Frob(z)xh^{-1}=\Frob(h')yzh^{-1}=y'z'.$$
The stack $\Xbf_{[s',s]}$ yields a correspondence 
$$\frac{\Hbf_{[s,\Frob(s)]}}{\Ad_{\Frob}\Hbf_{[s,s]}} \xleftarrow{\alpha_{[s',s]}} \Xbf_{[s',s]} \xrightarrow{\beta_{[s',s]}} \frac{\Hbf'_{[s',\Frob(s')]}}{\Ad_{\Frob}\Hbf'_{[s',s']}}.$$
\begin{corol} 
There is a commutative diagram 
\[\begin{tikzcd}
	{\DD^{\ofrak}(\pt/\Gbf^{\Frob})} & {\DD^{\unip}(\frac{\Hbf_{[s,\Frob(s)]}}{\Ad_{\Frob}\Hbf_{[s,s]}})} \\
	{\DD^{\ofrak'}(\pt/\Gbf^{'\Frob})} & {\DD^{\unip}(\frac{\Hbf'_{[s',\Frob(s')]}}{\Ad_{\Frob}\Hbf'_{[s',s']}}).}
	\arrow["{\EE_s^{\Frob}}"{description}, from=1-1, to=1-2]
	\arrow["{f_*}", from=1-1, to=2-1]
	\arrow["{\beta_{[s',s],!}\alpha^*_{[s',s]}}", from=1-2, to=2-2]
	\arrow["{\EE_{s'}^{\Frob}}"{description}, from=2-1, to=2-2]
\end{tikzcd}\]
\end{corol} 

\begin{proof}
The only part not yet proved is the identification of the map $\Tr(\HH^{\Hbf'}_{[s',s]}, \beta^{\Hbf})$; this follows from the same calculation as in \Cref{lem:calculation-pushfoward-trace}.
\end{proof}

\subsection{Application to the classification of irreducible representations of $\Gbf^{\Frob}$}

As a consequence of Theorem \ref{thm:EndoscopyRep}, we can obtain a classification of the irreducible representations of $\Gbf^{\Frob}$. 

\begin{defi}\label{defi:JordanIsoClasses}
Let $\mathbb{J}(\Gbf^{\Frob}, \Lambda)$ be the set 
$$\{(s,\rho), s \in (\Gbf^{*}_{\sss,\Lambda})^{\Frob^*}, \rho \in \Irr_{\Lambda}^{\unip}(\Gbf^{\Frob}_s)\} / \sim$$
where $\Lambda \in \{\Flb, \Qlb\}$, and the equivalence relation is given by 
$$(s, \rho) \sim (s', \rho')$$
if and only if the two pairs lie in the same orbit of the rational-conjugacy action groupoid: a morphism carrying $s$ to $s'$ transports $\Gbf_s^{\Frob}$ to $\Gbf_{s'}^{\Frob}$ and carries the isomorphism class of $\rho$ to that of $\rho'$. Equivalently, $\sim$ is the equivalence relation generated by all such transports. 
\end{defi}

\begin{thm}\label{thm:classificationOfIrrep}
Assume that $\ell$ satisfies the condition of \Cref{sec:EndoscopyReps}.
\begin{enumerate}
\item There is a bijection, which depends only on the pinning of $\Gbf$
$$\Irr_{\Lambda}(\Gbf^{\Frob}) \cong \mathbb{J}(\Gbf^{\Frob}, \Lambda).$$
\item Assume further that $\Gbf$ has connected center and $p$ and $\ell$ are good for $\Gbf$, then there is a canonical bijection
$$\mathbb{J}(\Gbf^{\Frob}, \Flb) \cong \{(s,\rho), s \in (\Gbf^{*}_{\sss,\Flb})^{\Frob^*}, \rho \in \Irr^{\unip}_{\Qlb}(\Gbf^{\Frob}_s)\}/\sim$$
where the equivalence relation is as in \Cref{defi:JordanIsoClasses}. 
\item Assume that $p$ and $\ell$ are good for $\Gbf$. With the basic set and the ordering specified in \cite[Theorem B]{DudasBrunatTaylor} and transported through \cite[Corollary 3.7]{FengSpaeth}, the decomposition matrix of $\Gbf^{\Frob}$ is unitriangular.
\end{enumerate}
\end{thm}

\begin{rque}
In the statement of Theorem \ref{thm:classificationOfIrrep}, the assumption $\ell$ is good for $\Gbf$ is almost always redundant given the assumption on $\ell$ in Section \ref{sec:EndoscopyReps}. 
\end{rque}

\begin{proof}
Point $(i)$ is a direct consequence of Theorems \ref{thm:EndoscopyRep}, \ref{thm:EndoscopyRepCanonical} and the compatibility with reduction mod $\ell$ which is Proposition \ref{prop:InversionOfEll}. For point $(ii)$, since $\Gbf$ has connected center, the group $\Gbf_s$ is connected and by $(i)$, the statement is \cite[Theorem B]{DudasBrunatTaylor}. Point $(iii)$ is a direct combination of \cite[Theorem B]{DudasBrunatTaylor} and \cite[Corollary 3.7]{FengSpaeth}.
\end{proof}

\begin{rque}
Point $(iii)$ of \Cref{thm:classificationOfIrrep} was originally conjectured by Geck \cite{GeckHiss}. 
\end{rque}

\begin{conj}
\Cref{thm:classificationOfIrrep} holds without assumptions on $p$, $\ell$ or $\Gbf$ (besides $p \neq \ell$). 
\end{conj}

\subsection{Generalizations}

To conclude this section, we state a few possible generalizations and strengthenings of the equivalence of Theorem \ref{thm:EndoscopyRep}, all of which are currently conjectural.

\begin{conj}
There exists an endoscopic equivalence as in Theorem \ref{thm:EndoscopyRep} without any hypothesis on $\ell$. 
\end{conj}
Ultimately, the hypothesis on $\ell$ was used twice: in the proof of the Endomorphismensatz \ref{thm:Endomorphismensatz} and in the trivialization of the Whittaker model of $\Hbf_{\mathfrak{o}}$ in Section \ref{sec:TrivializationWhittModelOfH}. The machinery of categorical traces introduces no new restrictions on $\ell$. One could expect a Soergel-type description of the category $\HH^{\Gbf}_{[s,s]}$ even for bad $\ell$. In the Betti context, results in this direction include \cite{TaylorSheaves}.

\begin{conj}\label{conj:SplendidEquiv}
The endoscopic equivalence of Theorem \ref{thm:EndoscopyRep} induces a Morita equivalence between the corresponding representation categories, and the associated derived equivalence is splendid in the sense of Rickard \cite{Rickard}.
\end{conj}
This conjecture is known in the setting of \cite{BonnafeDatRouquier}. Here the Morita assertion concerns the equivalence of the corresponding abelian representation categories, whereas ``splendid'' qualifies its derived enhancement. The latter statement requires one to realize the derived equivalence by a bounded complex of bimodules whose terms are direct summands of finite direct sums of permutation bimodules. The argument of \cite{BonnafeDatRouquier} starts from Rickard's theorem \cite{Rickard}, which gives such a permutation-module model for the compactly supported cohomology complex of a variety with an action of $\Gbf^F$. Our equivalence is constructed by abstract methods that do not make it easy to identify a complex realizing the derived equivalence. We expect Conjecture \ref{conj:SplendidEquiv} to follow from a motivic refinement of all the endoscopic equivalences. Indeed, \cite{Voevodsky} shows that the category of mixed Artin--Tate motives on $\Spec(F)$, for a field $F$, is equivalent to the category of permutation modules of the absolute Galois group of $F$. Replacing $\Spec(F)$, which we loosely interpret as $\pt/\Gal(\overline{F}/F)$, by the stack $\pt/\Gbf^{\Frob}$ suggests constructing a similar category of motives on $\pt/\Gbf^{\Frob}$, which should be equivalent to the category of permutation modules of $\Gbf^{\Frob}$.

\section{Endomorphisms of the Gelfand--Graev representation}\label{sectionGelfandGraev}

In this section, the following hypotheses are in force: $\Gbf$ has connected center and $\ell$ is good for $\Gbf$, see \Cref{def:torsion-prime}. 

\subsection{Statement of the theorems}

Recall that in Section \ref{sec:WhittakerModelForG} we fixed a nondegenerate character $\phi : \Ubf^- \to \A^1$, which we further assume to be stable under $\Frob$. We also chose a primitive $p$-th root of unity $\zeta$ in $\Lambda$. These choices give rise to a character
$$\psi : (\Ubf^-)^{\Frob} \xrightarrow{\phi} \Fq \xrightarrow{\tr_{\Fq/\Fp}} \Fp \xrightarrow{\zeta} \Lambda^{\times}$$
of $(\Ubf^-)^{\Frob}$. The Gelfand--Graev representation attached to $\psi$ is the representation
\begin{equation}
\Gamma_{\psi} \colondef = \ind_{(\Ubf^-)^{\Frob}}^{\Gbf^{\Frob}}(\psi).
\end{equation}

\begin{rque}
This representation is called ``the'' Gelfand--Graev representation of $\Gbf^{\Frob}$ because, under our assumption that the center of $\Gbf$ is connected, all pairs $(\Ubf^-,\psi)$ consisting of the unipotent radical of an $\Frob$-stable Borel subgroup and a nondegenerate character of $(\Ubf^-)^{\Frob}$ are conjugate under $\Gbf^{\Frob}$.
\end{rque}

The goal of this section is to propose new proofs of Theorems \ref{thmDudas} and \ref{thmShottonLi}. Let $w \in \Wbf$ and denote by ${^*}R_w : \DD(\Rep(\Gbf^{\Frob})) \to \DD(\Rep(\Tbf^{w\Frob}))$ the usual Deligne--Lusztig restriction functors, recalled in Section \ref{sec:TraceOnCategoricalHeckeCat}, and by ${^!}R_w$ the functor $\D \circ {^*}R_w \circ \D$ where $\D$ is Verdier duality. 

\begin{thm}[\cite{Dudas}]\label{thmDudas}
There is a $\Tbf^{w\Frob}$-linear isomorphism 
$${^*}R_w(\Gamma_{\psi}) \cong \Lambda[\Tbf^{w\Frob}][-\ell(w)].$$
\end{thm}

The isomorphism in this theorem is noncanonical because it involves several choices. However, one can deduce from it a canonical morphism
$$\Cur_w : \End(\Gamma_{\psi}) \to \Lambda[\Tbf^{w\Frob}]$$
which is called the $w$-Curtis morphism. This morphism was first introduced by Curtis \cite{Curtis} when $\Lambda = \Qlb$. Integrality properties of this morphism were studied in \cite{BonnafeKessar} and fully established by the above theorem. Putting together all these morphisms, we define the full Curtis morphism 
$$\Cur : \End(\Gamma_{\psi}) \to \bigoplus_{w \in \Wbf} \Lambda[\Tbf^{w\Frob}].$$
It follows from the conservativity of the Deligne--Lusztig restriction functors that this morphism is injective; see also \cite{BonnafeKessar}. This statement holds without any additional restriction on $\ell$, as does our proof.

Recall that we have chosen a topological generator of $\pi_1^t(\Gm)$. This choice yields a trivialization of the roots of unity in $\Fqb$ and provides an isomorphism 
$$\Lambda[\Tbf^{w\Frob}] \cong \Ocal((\Tbf^{\vee})^{w\Frob^{\vee}})$$
where $(\Tbf^{\vee})^{w\Frob^{\vee}}$ is the scheme of $w\Frob^{\vee}$-fixed points on $\Tbf^{\vee}$. 
Let $(\Tbf^{\vee}\sslash \Wbf)^{\Frob^{\vee}}$ denote the scheme of $\Frob^{\vee}$-fixed points on $\Tbf^{\vee}\sslash \Wbf$. The natural map $\Tbf^{\vee} \to \Tbf^{\vee} \sslash \Wbf$ yields a map on fixed points
$$(\Tbf^{\vee})^{w\Frob^{\vee}} \to (\Tbf^{\vee}\sslash \Wbf)^{\Frob^{\vee}}.$$
We denote by 
$$\Cur_w^{\spec} : \Ocal((\Tbf^{\vee}\sslash \Wbf)^{\Frob^{\vee}}) \to \Ocal((\Tbf^{\vee})^{w\Frob^{\vee}})$$
the corresponding map on algebras, and we call this morphism the `$w$-spectral Curtis morphism'. We also denote by 
$$\Cur^{\spec} : \Ocal((\Tbf^{\vee}\sslash \Wbf)^{\Frob^{\vee}}) \to \bigoplus_{w \in \Wbf} \Ocal((\Tbf^{\vee})^{w\Frob^{\vee}})$$ 
the direct sum of all these maps. 

\begin{thm}[\cite{Li}, \cite{ShottonLi}]\label{thmShottonLi}
There is a unique isomorphism 
$$\End(\Gamma_{\psi}) \cong \Ocal((\Tbf^{\vee}\sslash \Wbf)^{\Frob^{\vee}})$$
making the following diagram commutative
\[\begin{tikzcd}
	{\End(\Gamma_{\psi})} & {\Ocal((\Tbf^{\vee}\sslash \Wbf)^{\Frob^{\vee}})} \\
	{\bigoplus_{w\in W}\Lambda[\Tbf^{w\Frob}]} & {\bigoplus_{w\in W}\Ocal((\Tbf^{\vee})^{w\Frob^{\vee}}).}
	\arrow[from=1-1, to=1-2]
	\arrow["\Cur"', from=1-1, to=2-1]
	\arrow["{{\Cur^{\spec}}}", from=1-2, to=2-2]
	\arrow[from=2-1, to=2-2]
\end{tikzcd}\]
\end{thm} 

\begin{rque}
The uniqueness in Theorem \ref{thmShottonLi} follows from the fact that both maps $\Cur$ and $\Cur^{\spec}$ are injective.
\end{rque} 

In Section \ref{subsect:Linearity}, we have produced a canonical $\Ocal((\Tbf^{\vee}\sslash \Wbf)^{\Frob^{\vee}})$-linear structure on $\DD(\Rep(\Gbf^{\Frob}))$ which is compatible with the Deligne--Lusztig restriction functors. Hence we already have a canonical morphism 
\begin{equation}\label{eq:MonodromyMapForGelfandGraev}
\Ocal((\Tbf^{\vee}\sslash \Wbf)^{\Frob^{\vee}}) \to \End(\Gamma_{\psi}).
\end{equation}
Theorem \ref{thmShottonLi} is equivalent to the assertion that map (\ref{eq:MonodromyMapForGelfandGraev}) is an isomorphism. To give a new proof of both theorems, we study the object $\hc(\Gamma_{\psi})$ carefully. Theorem \ref{thmDudas} will follow from its fibers, whereas Theorem \ref{thmShottonLi} will follow from its endomorphism algebra.

\begin{rque} 
Unlike the proofs of \cite{Li} and \cite{ShottonLi}, our proof is entirely geometric and does not rely on character computations. The strategy used here to control the endomorphism algebra of Whittaker objects has already appeared in geometric representation theory: in \cite{BezrukavnikovTolmachov} for (unipotent) character sheaves and in \cite{LiNadlerYun} in the Betti setting.
\end{rque} 

\begin{rque}
The versions of Theorems \ref{thmDudas} and \ref{thmShottonLi} for general $\Lambda$ follow from the statement for $\Lambda = \Zlb$. From now on, we assume that $\Lambda = \Zlb$. 
\end{rque}

\subsection{Horocycle transform of $\Gamma_{\psi}$ and proof of Theorem \ref{thmDudas}}\label{sec:ProofDudas}

We consider the object $\hc(\Gamma_{\psi})$. We recall that $\pfrak : \Ubf \backslash \Gbf/\Ubf \to \frac{\Ubf \backslash \Gbf/\Ubf}{\Ad_{\Frob}\Tbf}$ is the quotient map. 

\begin{lem}\label{lem:HorocycleOfGelfandGraev}
There is a natural isomorphism 
$$\hc(\Gamma_{\psi}) \cong \pfrak_!\mathbb{T}[\dim \Tbf]$$
where $\mathbb{T}$ is the big tilting sheaf of Section \ref{subsec:BigTilting}. 
\end{lem}

\begin{proof}
The argument closely follows that of \cite[Theorem 5.5.1]{BezrukavnikovTolmachov}. Namely the same argument shows that 
$$\hc(\Gamma_{\psi}) \cong \pfrak_!(\Av_{\Ubf}\Av_{\psi}\delta_1)[\dim \Tbf]$$
where 
$$\Av_{\Ubf} : \DD((\Ubf^-, \Lcal_{\psi})\backslash \Gbf/\Ubf) \to \DD(\Ubf \backslash \Gbf/\Ubf)$$ 
and 
$$ \Av_{\psi} : \DD(\Ubf \backslash \Gbf/\Ubf) \to \DD((\Ubf^-, \Lcal_{\psi})\backslash \Gbf/\Ubf)$$ 
are the averaging functors introduced in Definition \ref{def:AveragingFunctors} and $\delta_1$ is the skyscraper sheaf at $1$. 
Recall from Lemma \ref{lem:MultiplicityOne} that $\mathbb{T} = \Av_{\Ubf}\Av_{\psi}(\Delta_1)$. The object $\Delta_1$ lies in $\DD(\Tbf, \RR_{\Tbf})$. As the map $\pfrak$ is $\Tbf$-equivariant for the action of $\Tbf$ by right translations and the averagings are done with respect to actions by left translations, we get 
$$\pfrak_!(\Av_{\Ubf}\Av_{\psi}\delta_1) * \Delta_1 \cong \pfrak_!(\Av_{\Ubf}\Av_{\psi}(\delta_1 * \Delta_1)) \cong \pfrak_!(\Av_{\Ubf}\Av_{\psi}\Delta_1) = \pfrak_!\mathbb{T}.$$
\end{proof}

\begin{proof}[Proof of Theorem \ref{thmDudas}]
The normalization of Deligne--Lusztig restriction used here gives an isomorphism of functors
$$i_w^*\hc\cong{}^!R_w.$$
By \Cref{lem:HorocycleOfGelfandGraev},
$${}^!R_w(\Gamma_{\psi})
\cong \pfrak_!i_w^*\mathbb{T}[\dim\Tbf].$$
As $\mathbb{T}$ has a $\Delta$-filtration with each $\Delta_w$ appearing exactly once, we have
$$i_w^*\mathbb{T}\cong
\nu_w^*\!\left(\bigoplus_{\chi\in\Ch(\Tbf)}
(L_{\Tbf}\otimes\Lcal_{\chi})\right)[\dim\Tbf+\ell(w)].$$
Therefore \cite[Lemma 2.8.3]{EteveFreeMonodromic} gives
$${}^!R_w(\Gamma_{\psi})\cong
\Lambda[\Tbf^{w\Frob}][\ell(w)].$$
The Theorem follow immediately from duality. 
\end{proof}

\subsection{Endomorphisms of $\hc(\Gamma_{\psi})$}

Recall from Theorem \ref{thm:Endomorphismensatz} that we have 
$$\End(\mathbb{T}) \cong \Ocal(\Ccal(\Tbf) \times_{\Ccal(\Tbf)\sslash \Wbf} \Ccal(\Tbf)).$$
We denote by $\mathbf{\Xcal}$ the intersection of $\Ccal(\Tbf) \times_{\Ccal(\Tbf)\sslash \Wbf} \Ccal(\Tbf)$ with the graph of Frobenius in $\Ccal(\Tbf) \times \Ccal(\Tbf)$, namely there is a Cartesian square 
\[\begin{tikzcd}
	\mathbf{\Xcal} & {\Ccal(\Tbf)} \\
	{\Ccal(\Tbf) \times_{\Ccal(\Tbf)\sslash \Wbf} \Ccal(\Tbf)} & {\Ccal(\Tbf) \times \Ccal(\Tbf).}
	\arrow[from=1-1, to=1-2]
	\arrow[from=1-1, to=2-1]
	\arrow["{\Frob^{\vee}\times \id }", from=1-2, to=2-2]
	\arrow[from=2-1, to=2-2]
\end{tikzcd}\]
In this subsection we present the fixed locus using the graph $\Frob^{\vee}\times\id$, so a graph point is $(\Frob^{\vee}(y),y)$. Swapping the two factors identifies this presentation canonically with the $\id\times\Frob^{\vee}$ convention used for loop stacks in \Cref{sec:local-traces,sec:GeometricSeries}.
Consider the map $\Tbf^{\vee} \times \Tbf^{\vee} \to \Tbf^{\vee}$ given by $(x,y) \mapsto x\Frob^{\vee}(y^{-1})$ and denote by $\phi$ the map induced on $\Ccal(\Tbf) \times \Ccal(\Tbf) \to \Ccal(\Tbf)$. 
\begin{lem}\label{lem:CartesianDiag}
The following diagram is Cartesian
\[\begin{tikzcd}
	\mathbf{\Xcal} & {\Ccal(\Tbf) \times_{\Ccal(\Tbf)\sslash \Wbf}\Ccal(\Tbf)} & {\Ccal(\Tbf) \times \Ccal(\Tbf)} \\
	1 & {\Ccal(\Tbf).}
	\arrow[from=1-1, to=1-2]
	\arrow[from=1-1, to=2-1]
	\arrow[from=1-2, to=1-3]
	\arrow[from=1-2, to=2-2]
	\arrow["\phi", from=1-3, to=2-2]
	\arrow[from=2-1, to=2-2]
\end{tikzcd}\]
\end{lem}
\begin{proof}
Denote by $\mathbf{\Ycal}$ the pullback fitting into the previous diagram, i.e. 
\[\begin{tikzcd}
	\mathbf{\Ycal} & {\Ccal(\Tbf) \times_{\Ccal(\Tbf)\sslash \Wbf}\Ccal(\Tbf)} & {\Ccal(\Tbf) \times \Ccal(\Tbf)} \\
	1 & {\Ccal(\Tbf).}
	\arrow[from=1-1, to=1-2]
	\arrow[from=1-1, to=2-1]
	\arrow[from=1-2, to=1-3]
	\arrow[from=1-2, to=2-2]
	\arrow["\phi", from=1-3, to=2-2]
	\arrow[from=2-1, to=2-2]
\end{tikzcd}\]
We must show that $\mathbf{\Xcal} = \mathbf{\Ycal}$ (as subschemes of $\Ccal(\Tbf) \times \Ccal(\Tbf)$). 
Consider the composition $\mathbf{\Xcal} \to \Ccal(\Tbf) \times \Ccal(\Tbf) \xrightarrow{\phi} \Ccal(\Tbf)$. This composition factors through $1 \in \Ccal(\Tbf)$ and we have a commutative diagram 
\[\begin{tikzcd}
	\mathbf{\Xcal} & {\Ccal(\Tbf) \times_{\Ccal(\Tbf) \sslash \Wbf} \Ccal(\Tbf)} \\
	1 & {\Ccal(\Tbf)}
	\arrow[from=1-1, to=1-2]
	\arrow[from=1-1, to=2-1]
	\arrow[from=1-2, to=2-2]
	\arrow[from=2-1, to=2-2]
\end{tikzcd}\]
which yields a map $\mathbf{\Xcal} \to \mathbf{\Ycal}$. Conversely, consider the composition $\mathbf{\Ycal} \to \Ccal(\Tbf) \times \Ccal(\Tbf) \xrightarrow{\phi} \Ccal(\Tbf)$ as this composition factors through $1$, we get that $\mathbf{\Ycal} \to 1$ factors through $\Ccal(\Tbf) \to \Ccal(\Tbf) \times \Ccal(\Tbf)$. We then have a commutative diagram 
\[\begin{tikzcd}
	\mathbf{\Ycal} & {\Ccal(\Tbf) } \\
	{\Ccal(\Tbf) \times_{\Ccal(\Tbf) \sslash \Wbf} \Ccal(\Tbf)} & {\Ccal(\Tbf) \times \Ccal(\Tbf)}
	\arrow[from=1-1, to=1-2]
	\arrow[from=1-1, to=2-1]
	\arrow[from=1-2, to=2-2]
	\arrow[from=2-1, to=2-2]
\end{tikzcd}\]
which yields a map $\mathbf{\Ycal} \to \mathbf{\Xcal}$. It is easily seen that these maps are inverse of each other. 
\end{proof}

\begin{lem}\label{lem:EndomorphimspTilting}
There is a canonical isomorphism 
$$\End(\pfrak_!\mathbb{T}) \cong \Ocal(\mathbf{\Xcal}).$$
\end{lem}

\begin{proof}
We have the following sequence of isomorphisms
\begin{align*}
\End(\pfrak_!\mathbb{T}) &\cong \Hom(\mathbb{T}, \pfrak^!\pfrak_!\mathbb{T}) \\
&\cong \Hom(\mathbb{T}, \Lambda \otimes_{\Ocal(\Ccal(\Tbf))} \mathbb{T}) \\
&\cong \V(\Lambda \otimes_{\Ocal(\Ccal(\Tbf))} \mathbb{T}) \\
&= \Lambda \otimes_{\Ocal(\Ccal(\Tbf))}  \V(\mathbb{T}) \\
&\cong \Ocal(\mathbf{\Xcal}).
\end{align*}
Let us detail this calculation. 
\begin{enumerate}
\item The first line is the $(\pfrak_!, \pfrak^!)$ adjunction.
\item The second line follows from the fact that for a monodromic sheaf $A$ on a space $\mathbf{X}$, the sheaf $\qfrak^!\qfrak_!A = \Lambda \otimes_{\Ocal(\Ccal(\Tbf))} A$ where $\qfrak : \mathbf{X} \to \mathbf{X}/\Tbf$ is the natural quotient and the $\Ocal(\Ccal(\Tbf))$-structure on $A$ is given by its monodromy, see \cite[Section 2.5]{EteveFreeMonodromic} and $\Lambda$ corresponds to the point $1 \in \Ccal(\Tbf)$. In the case of the above calculation, the action is given by $\Ad_{\Frob}\Tbf$. 
\item The third line follows from the definition of $\V$.
\item The fourth line follows from the linearity of $\V$ over $\Ocal(\Ccal(\Tbf) \times \Ccal(\Tbf))$. 
\item The last line follows from Lemma \ref{lem:CartesianDiag} and the Endomorphismensatz \ref{thm:Endomorphismensatz}. 
\end{enumerate}
\end{proof}

Consider the natural map $\Ccal(\Tbf) \to \Tbf^{\vee}$ and denote by $\mathbf{X}$ the following pullback 
\[\begin{tikzcd}
	\mathbf{X} & {\Tbf^{\vee}} \\
	{\Tbf^{\vee} \times_{\Tbf^{\vee}\sslash \Wbf}\Tbf^{\vee}} & {\Tbf^{\vee} \times \Tbf^{\vee}}
	\arrow[from=1-1, to=1-2]
	\arrow[from=1-1, to=2-1]
	\arrow["{\Frob^{\vee}\times \id }", from=1-2, to=2-2]
	\arrow[from=2-1, to=2-2]
\end{tikzcd}\]
There is a natural map $\mathbf{\Xcal} \to \mathbf{X}$ induced by $\Ccal(\Tbf) \to \Tbf^{\vee}$. 
\begin{lem}\label{lem:XvsXcal}
The map $\mathbf{\Xcal} \to \mathbf{X}$ is an isomorphism. 
\end{lem}

\begin{proof}
For each $w\in\Wbf$, let
$$\mathbf X_w:=\ker(w\Frob^{\vee}-\id:\Tbf^{\vee}\to\Tbf^{\vee}).$$
Choose $d>0$ such that $w\Frob^{\vee}$ has $d$-th power $(\Frob^{\vee})^d$. Then $\mathbf X_w$ is a closed subgroup of the finite group scheme $\ker((\Frob^{\vee})^d-id)$; hence $\mathbf X_w$ is finite over $\Lambda$. For every geometric point $x$ of $\mathbf X$, equality of the two images in $\Tbf^{\vee}\sslash\Wbf$ gives $w\Frob^{\vee}(x)=x$ for some $w$. Thus the support of $\mathbf X$ is contained in the finite union $\bigcup_w\mathbf X_w$. It follows that $\Xbf$ is finite over $\Lambda$. 

As the component $\Ccal(\Tbf)_s$ is the formal completion of $\Tbf^{\vee}$ at $s$ (after the finite coefficient extension over which $s$ is defined), $\mathbf{\Xcal}$ is the base change of $\mathbf X$ to the union of the completed local schemes at its torsion support. At every supporting pair $(x,y)$,
$$\widehat{\Ocal}_{\mathbf X,(x,y)}
\cong \Ocal_{\mathbf X,(x,y)}\widehat{\otimes}_{\Ocal(\Tbf^{\vee}\times\Tbf^{\vee})}
\bigl(\widehat{\Ocal}_{\Tbf^{\vee},x}\widehat{\otimes}_{\Lambda}\widehat{\Ocal}_{\Tbf^{\vee},y}\bigr),$$
and these completed local rings are exactly those defining $\mathbf{\Xcal}$, hence the map $\Xcal \to \Xbf$ is an isomorphism.
\end{proof}

Consider the natural $\Wbf \times \Wbf$ action on $\Tbf^{\vee} \times \Tbf^{\vee}$. The subscheme $\Tbf^{\vee} \times_{\Tbf^{\vee} \sslash \Wbf} \Tbf^{\vee}$ is stable under this action. Consider the inclusion $\Wbf \xrightarrow{\Frob \times \id} \Wbf \times \Wbf$, then $\Tbf^{\vee} \xrightarrow{\Frob^{\vee} \times \id} \Tbf^{\vee} \times \Tbf^{\vee}$ is stable under the action of $\Wbf$. It then follows that $\mathbf{X}$ is also stable under the action of $\Wbf$. The diagram defining $\mathbf{X}$ lies above the diagram 
\[\begin{tikzcd}
	{(\Tbf^{\vee}\sslash \Wbf)^{\Frob^{\vee}}} & {\Tbf^{\vee}\sslash \Wbf} \\
	{\Tbf^{\vee}\sslash \Wbf} & {\Tbf^{\vee}\sslash \Wbf \times \Tbf^{\vee}\sslash \Wbf}
	\arrow[from=1-1, to=1-2]
	\arrow[from=1-1, to=2-1]
	\arrow["{\Frob^{\vee}\times \id }", from=1-2, to=2-2]
	\arrow["\Delta"', from=2-1, to=2-2]
\end{tikzcd}\]
where $\Delta$ is the diagonal. Hence, there is a canonical map $\mathbf{X} \to (\Tbf^{\vee}\sslash \Wbf)^{\Frob^{\vee}}$. This map is $\Wbf$-equivariant, so it induces a map $\mathbf{X}\sslash \Wbf \to (\Tbf^{\vee}\sslash \Wbf)^{\Frob^{\vee}}$.

\begin{rque}
Note that this map is essentially controlling the commutation of taking $\Wbf$-invariant (which is a limit in the category of algebras) with taking schematic $\Frob^{\vee}$-invariants (which is a colimit in the category of algebras). 
\end{rque} 

\begin{lem}\label{lem:InvariantvsCoinvariant}
The natural map $\mathbf{X}\sslash \Wbf \to  (\Tbf^{\vee}\sslash \Wbf)^{\Frob^{\vee}}$ is an isomorphism. 
\end{lem}

\begin{proof}
We first consider the diagram 
\[\begin{tikzcd}
	\mathbf{X} & {\Tbf^{\vee}} \\
	& {\Tbf^{\vee} \times_{\Tbf^{\vee}\sslash \Wbf}\Tbf^{\vee}} & {\Tbf^{\vee} \times \Tbf^{\vee}} \\
	{(\Tbf^{\vee}\sslash \Wbf)^{\Frob^{\vee}}} & {\Tbf^{\vee}\sslash \Wbf} \\
	& {\Tbf^{\vee}\sslash \Wbf} & {\Tbf^{\vee}\sslash \Wbf \times \Tbf^{\vee}\sslash \Wbf.}
	\arrow[from=1-1, to=1-2]
	\arrow[from=1-1, to=2-2]
	\arrow[from=1-1, to=3-1]
	\arrow["{\Frob^{\vee} \times \id}"{description}, from=1-2, to=2-3]
	\arrow[shift right, from=1-2, to=3-2]
	\arrow[from=2-2, to=2-3]
	\arrow[shift left, from=2-2, to=4-2]
	\arrow[from=2-3, to=4-3]
	\arrow[from=3-1, to=3-2]
	\arrow[from=3-1, to=4-2]
	\arrow["{\Frob^{\vee} \times \id}"{description}, from=3-2, to=4-3]
	\arrow["\Delta", from=4-2, to=4-3]
\end{tikzcd}\]
Note that the top, bottom and front faces are Cartesian. We now prove that the back face is Cartesian. Consider a commutative diagram 
\[\begin{tikzcd}
	\mathbf{Z} \\
	& \mathbf{X} & {\Tbf^{\vee}} \\
	& {(\Tbf^{\vee}\sslash \Wbf)^{\Frob^{\vee}}} & {\Tbf^{\vee}\sslash \Wbf;}
	\arrow[from=1-1, to=2-3]
	\arrow[from=1-1, to=3-2]
	\arrow[from=2-2, to=2-3]
	\arrow[from=2-2, to=3-2]
	\arrow[from=2-3, to=3-3]
	\arrow[from=3-2, to=3-3]
\end{tikzcd}\]
Consider the compositions 
$$\mathbf{Z} \to \Tbf^{\vee} \xrightarrow{\Frob^{\vee} \times \id} \Tbf^{\vee} \times \Tbf^{\vee} \to \Tbf^{\vee}\sslash \Wbf \times \Tbf^{\vee}\sslash \Wbf$$ 
and 
$$\mathbf{Z} \to (\Tbf^{\vee}\sslash \Wbf)^{\Frob^{\vee}} \to \Tbf^{\vee}\sslash \Wbf \xrightarrow{\Delta} \Tbf^{\vee}\sslash \Wbf \times \Tbf^{\vee}\sslash \Wbf.$$ 
It follows from the commutativity of the cube above and the Cartesianity of the front face that we have a commutative diagram
\[\begin{tikzcd}
	\mathbf{Z} && {\Tbf^{\vee}} \\
	& {\Tbf^{\vee} \times_{\Tbf^{\vee}\sslash \Wbf} \Tbf^{\vee}} & {\Tbf^{\vee}\times \Tbf^{\vee}} \\
	& {\Tbf^{\vee} \sslash \Wbf} & {\Tbf^{\vee} \sslash \Wbf \times \Tbf^{\vee}\sslash \Wbf}
	\arrow[from=1-1, to=1-3]
	\arrow["{\exists!}"{description}, from=1-1, to=2-2]
	\arrow[from=1-1, to=2-3]
	\arrow[from=1-1, to=3-2]
	\arrow["{\Frob^{\vee}\times \id }"{description}, from=1-3, to=2-3]
	\arrow[from=2-2, to=2-3]
	\arrow[from=2-2, to=3-2]
	\arrow[from=2-3, to=3-3]
	\arrow["\Delta"', from=3-2, to=3-3]
\end{tikzcd}\]
Hence we get a commutative diagram 
\[\begin{tikzcd}
	\mathbf{Z} & {\Tbf^{\vee}} \\
	{\Tbf^{\vee} \times_{\Tbf^{\vee}\sslash \Wbf} \Tbf^{\vee}} & {\Tbf^{\vee} \times \Tbf^{\vee}}
	\arrow[from=1-1, to=1-2]
	\arrow[from=1-1, to=2-1]
	\arrow[from=1-2, to=2-2]
	\arrow[from=2-1, to=2-2]
\end{tikzcd}\]
and therefore a map $\mathbf{Z} \to \mathbf{X}$. The commutativity of the square
\[\begin{tikzcd}
	\mathbf{Z} & {\Tbf^{\vee}} \\
	{(\Tbf^{\vee}\sslash \Wbf)^{\Frob^{\vee}}} & {\Tbf^{\vee}\sslash \Wbf}
	\arrow[from=1-1, to=1-2]
	\arrow[from=1-1, to=2-1]
	\arrow[from=1-2, to=2-2]
	\arrow[from=2-1, to=2-2]
\end{tikzcd}\]
follows from the fact that the maps from the back face of the above cube to the front face are closed immersions. The uniqueness of the map $\mathbf{Z} \to \mathbf{X}$ also follows from the fact that the map $\mathbf{X} \to \Tbf^{\vee} \times_{\Tbf^{\vee} \sslash \Wbf} \Tbf^{\vee}$ is a closed immersion. This establishes that the back face of the above cube is Cartesian.
\medskip

Since we have assumed that $\Gbf$ has connected center, it follows from the Pittie-Steinberg Theorem \cite{PittieSteinberg} that the map 
$$\Tbf^{\vee} \to \Tbf^{\vee} \sslash \Wbf$$
is finite and faithfully flat, more precisely, the algebra $\Ocal(\Tbf^{\vee})$ is free of finite rank $|\Wbf|$ over $\Ocal(\Tbf^{\vee}\sslash \Wbf)$. Since the back face of the cube is Cartesian, we then have that the map $\Ocal((\Tbf^{\vee}\sslash \Wbf)^{\Frob^{\vee}}) \to \Ocal(\mathbf{X})$ is finite and faithfully flat. 
By \cite[Theorem 3.9]{Li}, $(\Tbf^{\vee}\sslash \Wbf)^{\Frob^{\vee}}$ is finite and torsion free over $\Zlb$ hence it is a finite free $\Zlb$-module and thus so is $\Ocal(\mathbf{X})$. As $\Ocal(\mathbf{X}\sslash \Wbf) = \Ocal(\mathbf{X})^{\Wbf}$, it follows that $\Ocal(\mathbf{X}\sslash \Wbf)$ is $\ell$-torsion free and thus $\Ocal(\mathbf{X}\sslash \Wbf)$ is a finite flat $\Zlb$-algebra. 
Moreover, since inverting $\ell$ is exact and over $\Qlb$ the map $\mathbf{X}\sslash \Wbf \to (\Tbf^{\vee}\sslash \Wbf)^{\Frob^{\vee}}$ is an isomorphism, we deduce that both $\Ocal(\mathbf{X}\sslash \Wbf)$ and $\Ocal((\Tbf^{\vee}\sslash \Wbf)^{\Frob^{\vee}})$ are finite flat $\Zlb$-algebras of the same rank. To show that they are isomorphic, it is enough to show that after reducing mod $\ell$ the map is an isomorphism. 
Since $\mathbf{X}$ is flat over $\Zlb$, there is a natural injective map $\Ocal(\mathbf{X}_{\Flb}\sslash \Wbf) \leftarrow \Ocal((\mathbf{X}\sslash \Wbf)_{\Flb})$. We then have a sequence of maps
$$\Ocal((\Tbf^{\vee}\sslash \Wbf)^{\Frob^{\vee}}_{\Flb}) \to \Ocal((\mathbf{X}\sslash \Wbf)_{\Flb}) \to \Ocal(\mathbf{X}_{\Flb}\sslash \Wbf) \to \Ocal(\mathbf{X}_{\Flb}).$$
The composition is the reduction mod $\ell$ of the map $\mathbf{X} \to (\Tbf^{\vee}\sslash \Wbf)^{\Frob^{\vee}}$ hence it is faithfully flat and, in particular, injective. The middle and right maps are injective. This implies that the first map is an injection; since both the source and target of this map are finite $\Flb$-algebras of the same dimension, it is therefore an isomorphism. 
\end{proof}

\subsection{Proof of Theorem \ref{thmShottonLi}}\label{sec:ProofShottonLi}

The rest of this section is dedicated to the proof of Theorem \ref{thmShottonLi}. By functoriality and Lemmas \ref{lem:HorocycleOfGelfandGraev}, \ref{lem:EndomorphimspTilting} and \ref{lem:XvsXcal}, there is a map 
\begin{equation}\label{eq:CurtisGeometric}
\End(\Gamma_{\psi}) \to \End(\pfrak_!\mathbb{T}) = \Ocal(\mathbf{X}).
\end{equation}

By the discussion in the previous section, the target of this map has a $\Wbf$-action. 
\begin{lem}\label{lem:InvariantOfW}
The map (\ref{eq:CurtisGeometric}) factors through $\Ocal(\mathbf{X})^{\Wbf}$.
\end{lem}
Let us assume Lemma \ref{lem:InvariantOfW} and deduce the theorem. It follows from Lemma \ref{lem:InvariantvsCoinvariant} that the map (\ref{eq:CurtisGeometric}) is then 
\begin{equation}\label{eq:CurtisGeometricRefined}
\End(\Gamma_{\psi}) \to \Ocal((\Tbf^{\vee}\sslash \Wbf)^{\Frob^{\vee}}).
\end{equation}
Moreover, both $\DD(\pt/\Gbf^{\Frob})$ and $\DD(\frac{\Ubf \backslash \Gbf/\Ubf}{\Ad_{\Frob}\Tbf})$ are linear over $\Ocal((\Tbf^{\vee}\sslash \Wbf)^{\Frob^{\vee}})$, and $\hc$ is linear over the same algebra. It follows that map (\ref{eq:CurtisGeometricRefined}) is also linear. In particular, it is a section of map (\ref{eq:MonodromyMapForGelfandGraev}), and is therefore surjective.
Finally, consider the commutative diagram 
\[\begin{tikzcd}
	{\End(\Gamma_{\psi})} & {\Ocal((\Tbf^{\vee}\sslash \Wbf)^{\Frob^{\vee}})} & {\End(\pfrak_!\mathbb{T})} \\
	& {\bigoplus_w \End(i_w^*\pfrak_!\mathbb{T})} \\
	& {\bigoplus_{w} \Lambda[\Tbf^{w\Frob}]}
	\arrow[from=1-1, to=1-2]
	\arrow["\Cur"{description}, from=1-1, to=2-2]
	\arrow[from=1-2, to=1-3]
	\arrow["{\Cur^{\spec}}"{description}, from=1-2, to=2-2]
	\arrow[from=1-3, to=2-2]
	\arrow[Rightarrow, no head, from=2-2, to=3-2]
\end{tikzcd}\]
where the vertical maps are deduced by functoriality. Since the collection of functors $(i_w^*)$ is conservative on $\DD(\frac{\Ubf \backslash \Gbf/\Ubf}{\Ad_{\Frob}\Tbf})$, all these maps are injections. This implies that the map (\ref{eq:CurtisGeometricRefined}) is injective and is thus an isomorphism. 

\begin{proof}[Proof of Lemma \ref{lem:InvariantOfW}]
The morphism $\pfrak$ is a smooth surjective quotient map with connected fibers, hence the natural map
$$\End(\pfrak_!\mathbb T)\hookrightarrow\End(\pfrak^!\pfrak_!\mathbb T)$$
is an isomorphism. We construct a $\Wbf$-action on the target and show below that its restriction to the image agrees with the action on $\End(\pfrak_!\mathbb T)$ given by \Cref{lem:EndomorphimspTilting}.

By \Cref{thm:Trace}, the functor $\pfrak^!\hc : \DD(\pt/\Gbf^{\Frob}) \to \HH$ is $\Frob$-central. Hence the object $\pfrak^!\hc(\Gamma_{\psi}) = \pfrak^!\pfrak_!\mathbb{T}[\dim \Tbf]$ comes equipped with a canonical $\Frob$-central structure. Thus for all $w \in \Wbf$ there is an isomorphism
\begin{equation}\label{eq:Centrality}
\Delta_w * \pfrak^!\pfrak_!\mathbb{T} * \Frob^*(\Delta_w^{-1}) \cong \pfrak^!\pfrak_!\mathbb{T}.
\end{equation}
This defines an automorphism of $\End(\pfrak^!\pfrak_!\mathbb{T})$ as follows 
$$\End(\pfrak^!\pfrak_!\mathbb{T}) \to \End(\Delta_w * \pfrak^!\pfrak_!\mathbb{T} * \Frob^*(\Delta_w^{-1})) \to \End(\pfrak^!\pfrak_!\mathbb{T})$$ 
where the first map is induced by the functor $(-) \to (\Delta_w *- * \Frob^*(\Delta_w^{-1}))$ and the second is the isomorphism (\ref{eq:Centrality}). The same argument as in \cite[Lemma 9.6]{BezrukavnikovRicheSoergelTheory} shows that the collection of these morphisms yields an action of $\Wbf$. If $f\in\End(\Gamma_\psi)$ and $c_w$ denotes the centrality isomorphism in \eqref{eq:Centrality}, naturality of the central structure is the commutative square expressed by
$$c_w\circ\bigl(\Delta_w*\pfrak^!\hc(f)*\Frob^*(\Delta_w^{-1})\bigr)
=\pfrak^!\hc(f)\circ c_w.$$
Thus the image of every $f$ is fixed by the above $\Wbf$-action. It follows that the morphism \eqref{eq:CurtisGeometric} factors through $\Ocal(\mathbf X)^{\Wbf}$.

It remains to check that the action we have defined is the same as the one constructed from the identification of Lemma \ref{lem:EndomorphimspTilting}. Tracking the isomorphisms in the lemma, the main input is the Endomorphismensatz \ref{thm:Endomorphismensatz} and the $\Wbf \times \Wbf$-action constructed on $\End(\mathbb{T})$ \ref{rque:WWAction}. That action is canonical, see \Cref{rque:WWAction}.
\begin{enumerate}
\item There is an isomorphism of functors $\HH \to \HH$ 
$$\Delta_v * \pfrak^!\pfrak_!(-) * \Frob^*(\Delta_v^{-1}) \cong  \pfrak^!\pfrak_! (\Delta_v *(-) * \Frob^*(\Delta_v^{-1})).$$
\item There is an isomorphism $\Delta_v * \mathbb{T} * \Frob^*(\Delta_v^{-1}) \cong \mathbb{T}$ making the following diagram commutative 
\[\begin{tikzcd}
	{\pfrak^!\pfrak_!\Delta_v * \mathbb{T} * \Frob^*(\Delta_v^{-1})} & {\pfrak^!\pfrak_!\mathbb{T}} \\
	{\Delta_v * \pfrak^!\pfrak_! \mathbb{T} * \Frob^*(\Delta_v^{-1})} & {\pfrak^!\pfrak_! \mathbb{T} }
	\arrow[from=1-1, to=1-2]
	\arrow[from=1-1, to=2-1]
	\arrow[Rightarrow, no head, from=1-2, to=2-2]
	\arrow[from=2-1, to=2-2]
\end{tikzcd}\]
where the left vertical map is the isomorphism of functors applied to $\mathbb{T}$ and the bottom map is the centrality isomorphism of $\pfrak^!\mathbb{T}$. 
\end{enumerate}
We sketch the construction of the first isomorphism. It is enough to construct it for tilting sheaves. For such sheaves, we use the monoidality of $\V$, together with $\V(\Delta_v)=\Ocal(\Ccal(\Tbf))$ (with bimodule structure given by the graph of $v$) and $\V(\pfrak^!\pfrak_!(-))=\Lambda\otimes_{\Ocal(\Ccal(\Tbf))}-$. The second isomorphism can again be constructed using $\V$ and the Endomorphismensatz \ref{thm:Endomorphismensatz}.
Finally there is a commutative diagram 
\[\begin{tikzcd}
	{\End(\pfrak^!\pfrak_!\mathbb{T})} & {\End(\Delta_v * \pfrak^!\pfrak_!\mathbb{T} * \Frob^*(\Delta_v^{-1}))} & {\End(\pfrak^!\pfrak_!\mathbb{T})} \\
	& {\End(\pfrak^!\pfrak_!(\Delta_v * \mathbb{T} * \Frob^*(\Delta_v^{-1}))} \\
	{\End(\mathbb{T})} & {\End(\Delta_v *\mathbb{T} * \Frob^*(\Delta_v^{-1}))} & {\End(\mathbb{T})}
	\arrow[from=1-1, to=1-2]
	\arrow[from=1-2, to=1-3]
	\arrow[from=2-2, to=1-2]
	\arrow[from=3-1, to=1-1]
	\arrow[from=3-1, to=3-2]
	\arrow[from=3-2, to=2-2]
	\arrow[from=3-2, to=3-3]
	\arrow[from=3-3, to=1-3]
\end{tikzcd}\]
where the vertical maps are the ones coming from the functoriality of $\pfrak^!\pfrak_!$ and the horizontal maps are the ones defining the action of $\Wbf$ and $\Wbf \times \Wbf$ discussed previously. 
\end{proof}

\appendix

\section{Some trace calculations}\label{appendix:trace-calculation}

\subsection{Recollections on categorical traces}\label{sec:recollection-traces}

\subsubsection{}
We start by recalling the notion of categorical traces that we will use. 
\begin{defi}
Let $\Ccal \in \Pr_{\Lambda}$ be a monoidal presentable category and let $\Frob : \Ccal \to \Ccal$ be a continuous $\otimes$-endomorphism of $\Ccal$. 
\begin{enumerate}
\item the categorical trace of $\Frob$ on $\Ccal$ is the object of $\Pr_{\Lambda}$ given by 
$$\Tr(\Frob, \Ccal) := \Ccal \otimes_{\Ccal \otimes \Ccal^{\rev}} \Ccal_{\Frob}$$ 
where $\Ccal^{\rev}$ is the category $\Ccal$ with $\otimes$-structure given by $x \otimes^{\rev} y = y \otimes x$ and $\Ccal_{\Frob}$ is the category $\Ccal$ with the right module structure twisted by $\Frob$,
\item the $\Frob$-categorical center of $\Ccal$ is the category 
$$\Zcal_{\Frob}(\Ccal) := \Fun^L_{\Ccal \otimes \Ccal^{\rev}}(\Ccal, \Ccal_{\Frob})$$
with the same notations as above. 
\end{enumerate}
More generally if $\Mcal$ is a $\Ccal$-bimodule, define its Hochschild homology (which we will also denote by $\Tr$) to be 
$$\Tr(\Mcal, \Ccal) := \Ccal \otimes_{\Ccal \otimes \Ccal^{\rev}} \Mcal \in \mathrm{Pr}_{\Lambda}.$$
\end{defi} 
Thus, for a general bimodule, our convention is always ``bimodule first, algebra second.'' For a monoidal endomorphism $F$ of $\Ccal$, we use the shorthand
$$\Tr(F,\Ccal):=\Tr(\Ccal_F,\Ccal).$$
These notions categorify the usual Hochschild homology and cohomology; we refer to \cite{BenZviNadlerCharacterTheory}, \cite{HoyoisScherotzkeSibilla}, and \cite{ToyModel} for a longer discussion. We always assume that the monoidal categories under consideration are unital. Therefore, for every $\Ccal$-bimodule $\Mcal$, there is a natural map
$$\tr_{\Mcal} : \Mcal \to \Tr(\Mcal, \Ccal)$$
given by $m \mapsto 1_{\Ccal} \otimes m$. Note that we have not defined Hochschild cohomology for general bimodules as we will not need it thanks to the following theorems. 

\begin{thm}[\protect{\cite[Proposition 3.13]{BenZviNadlerCharacterTheory}}]
Let $\Ccal \in \Pr_{\Lambda}$ be a quasi-rigid dualizable monoidal category equipped with a pivotal structure (see \emph{loc. cit.} for a definition), then for all $\Frob : \Ccal \to \Ccal$ continuous monoidal endomorphism of $\Ccal$ there is a canonical isomorphism 
$$\Tr(\Frob, \Ccal) \cong \Zcal_{\Frob} (\Ccal).$$
\end{thm}

\begin{thm}[\protect{\cite[Theorem 4.4.2, Theorem 4.4.8]{EteveFreeMonodromic}}]
The category $\HH_{\tot}$ is quasi-rigid dualizable and equipped with a pivotal structure. 
\end{thm}

\subsubsection{Functoriality} 

Let $(\Ccal, \Mcal)$ and $(\Dcal, \Ncal)$ be two pairs consisting of quasi-rigid compactly generated monoidal categories $\Ccal$ and $\Dcal$ and let $\Mcal$ and $\Ncal$ be dualizable bimodules.

\begin{defi}\label{def:morphism-of-pairs}
A morphism of pairs $(\Ccal, \Mcal) \to (\Dcal, \Ncal)$ is the datum of a pair $(\Acal, \phi)$, where $\Acal$ is a $(\Dcal,\Ccal)$-bimodule that is right-adjointable as a bimodule (see \cite[3.6.2]{ToyModel}) and $\phi$ is a morphism of $(\Dcal, \Ccal)$-bimodules
$$\phi : \Acal \otimes_{\Ccal} \Mcal \to \Ncal \otimes_{\Dcal} \Acal.$$
\end{defi} 

\begin{construction}[see \protect{\cite[Section 3.6]{ToyModel}}]
Given a morphism of pairs $(\Ccal, \Mcal) \xrightarrow{(\Acal, \phi)} (\Dcal, \Ncal)$ as in \Cref{def:morphism-of-pairs}, there is a well-defined map
$$\Tr(\Acal, \phi) : \Tr(\Mcal, \Ccal) \to \Tr(\Ncal, \Dcal).$$
The map is constructed as follows. Let $\Acal^{R}$ denote the $(\Ccal,\Dcal)$-bimodule right adjoint to $\Acal$. Then $\Tr(\Acal, \phi)$ is defined as
\begin{align*}
\Tr(\Mcal, \Ccal) &\xrightarrow{\mathrm{unit}} \Tr(\Acal^R \otimes_{\Dcal} \Acal \otimes_{\Ccal} \Mcal, \Ccal) \\
&\xrightarrow{\phi} \Tr(\Acal^R \otimes_{\Dcal} \Ncal \otimes_{\Dcal} \Acal, \Ccal) \\
&\xrightarrow{\mathrm{cycl}} \Tr(\Ncal \otimes_{\Dcal} \Acal \otimes_{\Ccal} \Acal^{R}, \Dcal) \\
&\xrightarrow{\mathrm{counit}} \Tr(\Ncal, \Dcal).
\end{align*}
Furthermore, the formation of $\Tr(\Acal, \phi)$ is compatible with composition. 
\end{construction}

We will often consider the following situation. Let $f : \Ccal \to \Dcal$ be a monoidal functor between compactly generated rigid categories. Using $f$, we view $\Dcal$ as a $(\Dcal, \Ccal)$-bimodule. The rigidity of $\Ccal$ and $\Dcal$ implies that $\Dcal$ is right adjoinable as a bimodule (see \cite[3.6.10]{ToyModel}), and
\begin{equation}\label{eq:defi-right-adjoint-module}
\Dcal^R = \Dcal
\end{equation}
where $\Dcal^R$ is the $(\Ccal, \Dcal)$-bimodule that is right adjoint to $\Dcal$ and the RHS of (\ref{eq:defi-right-adjoint-module}) is $\Dcal$ viewed as a $(\Ccal, \Dcal)$-bimodule via the morphism $f$. Moreover, we have 
$$\Dcal^{RR} = \Dcal.$$

Assume now that we also have bimodules $(\Ccal, \Mcal)$ and $(\Dcal, \Ncal)$ as above. Let us denote $\Acal = \Dcal$ viewed as $(\Dcal, \Ccal)$-bimodule. Extending $\Acal$ to a morphism of pairs amounts to providing a map 
$$\phi : \Dcal \otimes_{\Ccal} \Mcal \to \Ncal$$
of $(\Dcal, \Ccal)$-bimodules. 

\begin{lem}[\protect{\cite[3.10.4]{ToyModel}}]\label{lem:computation-right-adjoint}
Let $f : \Ccal \to \Dcal$ be a morphism as above, and let $(\Dcal, \phi) : (\Ccal, \Mcal) \to (\Dcal, \Ncal)$ be a morphism of pairs. If $\phi$ admits a continuous right adjoint, then $\Tr(\Dcal, \phi)$ admits a continuous right adjoint and
$$\Tr(\Dcal, \phi)^R = \Tr(\Dcal, \phi^R),$$
where on the RHS $\Dcal$ is viewed as a $(\Ccal, \Dcal)$-bimodule. 
\end{lem}

\subsection{Local traces}\label{sec:local-traces}

Let $\Xibf$ be a finite discrete groupoid and consider a functor 
$$\Ccal : \Xibf^{\mathrm{op}} \to \Morita$$
such that
\begin{enumerate}
\item for all $\xi \in \Xibf$, the monoidal category $\Ccal_{\xi}$ is rigid and compactly generated, 
\item for all $\beta : \xi \to \xi'$, the corresponding $(\Ccal_{\xi}, \Ccal_{\xi'})$-bimodule $\Ccal^{\beta}$ is compactly generated (and thus dualizable as a bimodule by rigidity) and free of rank one over $\Ccal_{\xi}$. 
\end{enumerate}
Let $\Frob : \Xibf \to \Xibf$ be an automorphism of $\Xibf$ and let $\Frob_{\Ccal}$ be an isomorphism of functors $\Ccal \Rightarrow \Ccal\circ \Frob^{\mathrm{op}}$ making the following diagram commute 
\[\begin{tikzcd}
	{\Xibf^{\mathrm{op}}} && \Morita \\
	\\
	{\Xibf^{\mathrm{op}}.} && \Morita
	\arrow["\Ccal", from=1-1, to=1-3]
	\arrow["\Frob"', from=1-1, to=3-1]
	\arrow["{{{\Frob_{\Ccal}}}}"{description}, Rightarrow, from=1-1, to=3-3]
	\arrow[equals, from=1-3, to=3-3]
	\arrow["\Ccal"', from=3-1, to=3-3]
\end{tikzcd}\]
It follows that $\Frob_{\Ccal}$ is the datum of 
\begin{enumerate}
\item for all $\xi \in \Xibf$ of a monoidal equivalence $\Frob_{\xi} : \Ccal_{\xi} \to \Ccal_{\Frob(\xi)}$, 
\item for all $\beta : \xi \to \xi'$ of a map of $(\Ccal_{\xi}, \Ccal_{\xi'})$-bimodule $\Frob^{\beta} : \Ccal^{\beta} \to \Ccal^{\Frob(\beta)}$, where the actions on the target are transported along $\Frob_{\xi}$ and $\Frob_{\xi'}$. 
\end{enumerate}

\begin{defi}
We let $\Lcal_{\Frob}\Xibf$ be the $\Frob$-loop stack of $\Xibf$, defined by the following pullback:
\[\begin{tikzcd}
	{\Lcal_{\Frob}\Xibf} & \Xibf \\
	\Xibf & {\Xibf \times \Xibf.}
	\arrow[from=1-1, to=1-2]
	\arrow[from=1-1, to=2-1]
	\arrow["{(\id \times \Frob)}", from=1-2, to=2-2]
	\arrow["\Delta"', from=2-1, to=2-2]
\end{tikzcd}\]
It is the groupoid whose objects are pairs $(\xi, \beta)$, where $\beta :\xi \to \Frob(\xi)$ is a map, and whose morphisms are commutative diagrams
\[\begin{tikzcd}
	\xi & {\Frob(\xi)} \\
	{\xi'} & {\Frob(\xi').}
	\arrow["\beta", from=1-1, to=1-2]
	\arrow["\alpha"', from=1-1, to=2-1]
	\arrow["{\Frob(\alpha)}", from=1-2, to=2-2]
	\arrow["{\beta'}"', from=2-1, to=2-2]
\end{tikzcd}\]
\end{defi} 

\begin{construction}
From the datum of $(\Ccal, \Frob_{\Ccal})$, by taking trace, we obtain a functor  
$$\Tr^{\loc}(\Frob_{\Ccal}, \Ccal) : \Lcal_{\Frob}\Xibf^{\mathrm{op}} \to \mathrm{Pr}_{\Lambda}.$$
It is defined as follows 
\begin{enumerate}
\item for a pair $(\xi, \beta)$, 
$$\Tr^{\loc}(\xi, \beta) := \Tr(\Ccal^{\beta}_{\Frob},\Ccal_{\xi}),$$
where $\Ccal^{\beta}_{\Frob}$ has the underlying category of $\Ccal^{\beta}$, its original left $\Ccal_{\xi}$-action, and the right $\Ccal_{\xi}$-action obtained by composing $\Frob_{\xi}:\Ccal_{\xi}\to\Ccal_{\Frob(\xi)}$ with the original right $\Ccal_{\Frob(\xi)}$-action;
\item a morphism $(\xi,\beta)\to(\xi',\beta')$ in $\Lcal_{\Frob}\Xibf^{\mathrm{op}}$ is represented by a morphism $\alpha:(\xi',\beta')\to(\xi,\beta)$ in $\Lcal_{\Frob}\Xibf$. It gives the morphism of pairs
$$(\Ccal_{\xi},\Ccal^{\beta}_{\Frob})\longrightarrow
(\Ccal_{\xi'},\Ccal^{\beta'}_{\Frob})$$
whose right-adjointable bimodule is $\Ccal^{\alpha}$ and whose structure isomorphism is
$$\phi_{\alpha}:\Ccal^{\alpha}\otimes_{\Ccal_{\xi}}\Ccal^{\beta}_{\Frob}
\xrightarrow{\ \sim\ }
\Ccal^{\beta'}_{\Frob}\otimes_{\Ccal_{\xi'}}\Ccal^{\alpha}.$$
This isomorphism is induced by the equality $\beta\alpha=\Frob(\alpha)\beta'$, the composition constraints of $\Ccal$, and the equivariance $\Frob_{\Ccal}$. We denote the map on local traces by $\Tr(\Ccal^{\alpha},\phi_{\alpha})$.
\end{enumerate}
\end{construction}

Attached to the functor $\Ccal$, there is a monoidal category $\Ccal_{\tot} = \bigoplus_{\beta \in \Arr(\Xibf)} \Ccal^\beta$ where $\Arr(\Xibf)$ is the set of all arrows in $\Xibf$ and the monoidal structure is induced by the convolution map of bimodules. The category $\Ccal_{\tot}$ is still compactly generated and rigid. The datum of the natural transformation $\Frob$ yields a monoidal endofunctor 
$$\Frob_{\tot} : \Ccal_{\tot} \to \Ccal_{\tot}.$$
We can therefore consider the category $\Tr(\Frob_{\tot}, \Ccal_{\tot})$. 

\begin{construction}
Let $\xi \in \Xibf$, consider the inclusion $\Ccal_{\xi} \to \Ccal_{\tot}$. There is a natural map of pairs $(\Ccal_{\xi}, \Ccal^{\beta}) \to (\Ccal_{\tot}, \Ccal_{\tot, \Frob})$ given by $(\Ccal_{\tot}, \Frob)$ where $\Ccal_{\tot}$ is viewed as a $(\Ccal_{\tot}, \Ccal_{\xi})$-bimodule via the inclusion $\Ccal_{\xi} \to \Ccal_{\tot}$. By the rigidity of $\Ccal$ and $\Ccal_{\xi}$, it is right adjoinable as a bimodule. The map $\Frob$ is the right adjoint to the map 
$$\Ccal_{\tot} \otimes_{\Ccal_{\xi}} \Ccal^{\beta}_{\Frob} \cong \Ccal_{\tot, \Frob} \xrightarrow{\Frob} \Ccal_{\tot,\Frob} \cong \Ccal_{\tot, \Frob} \otimes_{\Ccal_{\tot}} \Ccal_{\tot}$$
induced by multiplication in $\Ccal$. By the rigidity of $\Ccal$, this map is again continuous and a morphism of bimodules. It follows that we have a well defined map 
$$\incl_{\xi, \beta} : \Tr(\Ccal^{\beta}_{\Frob},\Ccal_{\xi}) \to \Tr(\Frob, \Ccal_{\tot}).$$

Let $\alpha : (\xi, \beta) \to (\xi', \beta')$ be a morphism in $\Lcal_{\Frob}\Xibf^{\mathrm{op}}$. We have a commutative diagram of pairs 
\[\begin{tikzcd}
	{(\Ccal_{\xi},\Ccal^{\beta})} & {(\Ccal_{\tot}, \Ccal_{\tot,\Frob})} \\
	{(\Ccal_{\xi'},\Ccal^{\beta'})}
	\arrow[from=1-1, to=1-2]
	\arrow["{(\Ccal^{\alpha}, \phi_{\alpha})}"', from=1-1, to=2-1]
	\arrow[from=2-1, to=1-2]
\end{tikzcd}\]

In turn, this defines a map 
\begin{equation}\label{eq:local-to-global-trace}
\varinjlim_{\Lcal_{\Frob}\Xibf^{\mathrm{op}}} \Tr^{\loc} \to \Tr(\Frob, \Ccal_{\tot}).
\end{equation}
\end{construction}

\begin{thm}\label{thm:local-to-global-trace}
The map (\ref{eq:local-to-global-trace}) is an equivalence. 
\end{thm}

\begin{proof}
Consider the category $\Ccal_{\tot}^{\circ} = \bigoplus_{(\xi, \beta)} \Ccal_{\xi}$. It is equipped with a natural bimodule, namely $\Mcal_{\tot} = \bigoplus_{(\xi,\beta)} \Ccal^{\beta}$ and we have 
$$\Tr(\Mcal_{\tot}, \Ccal_{\tot}^{\circ}) = \bigoplus_{(\xi, \beta)} \Tr^{\loc}(\xi,\beta).$$
Furthermore, the natural map 
$$\Ccal_{\tot}^{\circ} \to \Ccal_{\tot}$$ 
is monoidal and compatible with $\Frob$ and thus defines a map on traces 
\begin{equation}\label{eq:proof-part-1}
\incl : \Tr(\Mcal_{\tot}, \Ccal_{\tot}^{\circ}) \to \Tr(\Ccal_{\tot, \Frob}, \Ccal)
\end{equation}
whose restriction to the summand $\Tr^{\loc}(\xi,\beta)$ is $\incl_{\xi, \beta}$. 

The collection of isomorphisms $\Tr(\Ccal^{\alpha}, \phi_{\alpha})$ define an action of $\Xibf^{\opp}$ on $\Tr(\Mcal_{\tot}, \Ccal_{\tot}^{\circ})$ and we tautologically have 
$$\varinjlim_{\Lcal_{\Frob}\Xibf^{\opp}} \Tr^{\loc} \cong (\Tr(\Mcal_{\tot}, \Ccal_{\tot}^{\circ}))_{\Xibf^{\opp}}.$$
Moreover, the map (\ref{eq:proof-part-1}) factors as 
$$\Tr(\Mcal_{\tot}, \Ccal_{\tot}^{\circ}) \xrightarrow{\pr} \varinjlim_{\Lcal_{\Frob}\Xibf^{\opp}} \Tr^{\loc}  \to \Tr(\Ccal_{\tot, \Frob}, \Ccal_{\tot})$$
where the second map is (\ref{eq:local-to-global-trace}). To prove that (\ref{eq:local-to-global-trace}) is an equivalence, by Barr-Beck, it is enough to construct an isomorphism of comonads $\incl^R\incl \cong \pr^R\pr$. By Lemma \ref{lem:computation-right-adjoint}, the functor $\incl^R\incl$ is given by the $\Ccal_{\tot}^{\circ}$-bimodule $\Ccal_{\tot} \otimes_{\Ccal^{\circ}_{\tot}} \Ccal_{\tot} \cong \Ccal_{\tot}$ (and the monad structure is induced by the multiplication in $\Ccal_{\tot}$). On the other hand, the functor $\pr^R\pr$ is given by the composition 
$$\Tr(\Mcal_{\tot}, \Ccal_{\tot}^{\circ}) \xrightarrow{\id \otimes \Lambda_{\Xibf}} \Tr(\Mcal_{\tot}, \Ccal_{\tot}^{\circ}) \otimes \DD(\Xibf^{\opp}) \xrightarrow{\mathrm{act}} \Tr(\Mcal_{\tot}, \Ccal_{\tot}^{\circ}).$$
This functor is given by the $\Ccal^{\circ}_{\tot}$-bimodule $\oplus_{\alpha} \Ccal^{\alpha} = \Ccal_{\tot}$. It follows that the two bimodules agree and the map (\ref{eq:local-to-global-trace}) is an equivalence. 
\end{proof}

Let us now examine the functoriality of the local trace construction. Let $(\Xibf,\Frob)$ and $(\Xibf',\Frob')$ be finite groupoids with automorphisms, let $u:\Xibf\to\Xibf'$ satisfy $u\Frob=\Frob'u$, and let $\Ccal$ and $\Dcal$ be equivariant Morita diagrams as above. Suppose that $\Acal:\Ccal\Rightarrow u^*\Dcal$ is an equivariant right-adjointable morphism of Morita diagrams. Explicitly, at $\xi$ it gives a $(\Dcal_{u(\xi)},\Ccal_{\xi})$-bimodule $\Acal_{\xi}$ and, at a loop $\beta:\xi\to\Frob(\xi)$, an isomorphism
$$\phi_{\beta}:\Acal_{\xi}\otimes_{\Ccal_{\xi}}\Ccal^{\beta}_{\Frob}
\xrightarrow{\ \sim\ }
\Dcal^{u(\beta)}_{\Frob'}\otimes_{\Dcal_{u(\xi)}}\Acal_{\xi},$$
compatible with morphisms of loops. Then the trace maps $\Tr(\Acal_{\xi},\phi_{\beta})$ assemble into a natural transformation
$$\Tr^{\loc}(\Frob_{\Ccal},\Ccal)\longrightarrow
(\Lcal_{\Frob}u)^*\Tr^{\loc}(\Frob'_{\Dcal},\Dcal).$$
This natural transformation in particular induces a map 
\begin{equation}\label{eq:morphism-between-local-traces}
\varinjlim_{\Lcal_{\Frob}\Xibf^{\mathrm{op}}}\Tr^{\loc}(\Frob_{\Ccal},\Ccal) \to \varinjlim_{\Lcal_{\Frob'}\Xibf^{'\mathrm{op}}} \Tr^{\loc}(\Frob_{\Dcal}',\Dcal).
\end{equation}

Correspondingly, passing to totalizations, taking $\Acal_{\tot} = \bigoplus_{\xi \in \Xibf} \Acal_{\xi}$, this is a $((u^*\Dcal)_{\tot}, \Ccal_{\tot})$-bimodule which is right adjoinable and compatible with $\Frob$. Moreover the map $\Xibf \to \Xibf'$ yields a natural monoidal map $(u^*\Dcal)_{\tot} \to \Dcal_{\tot}$ thus the bimodule $\Dcal_{\tot} \otimes_{(u^*\Dcal)_{\tot}} \Acal_{\tot}$ yields a map on traces 
$$\Tr(\Frob, \Ccal_{\tot}) \to \Tr(\Frob', \Dcal_{\tot}).$$

\begin{prop}[Functoriality of local traces]\label{prop:functoriality-local-traces}
There is a commutative diagram 
\[\begin{tikzcd}
	{\varinjlim_{\Lcal_{\Frob}\Xibf^{\mathrm{op}}}\Tr^{\loc}(\Frob_{\Ccal},\Ccal)} & {\varinjlim_{\Lcal_{\Frob'}\Xibf^{'\mathrm{op}}} \Tr^{\loc}(\Frob_{\Dcal}',\Dcal)} \\
	{\Tr(\Frob, \Ccal_{\tot})} & {\Tr(\Frob', \Dcal_{\tot})}
	\arrow[from=1-1, to=1-2]
	\arrow[from=1-1, to=2-1]
	\arrow[from=1-2, to=2-2]
	\arrow[from=2-1, to=2-2]
\end{tikzcd}\]
\end{prop}

\begin{proof}
It is enough to prove that for all pair $(\xi, \beta) \in \Lcal_{\Frob}\Xibf^{\mathrm{op}}$, the following diagram commutes (functorialy in the pair $(\xi, \beta)$) 
\[\begin{tikzcd}
	{\Tr^{\loc}(\Frob_{\Ccal}, \Ccal)(\xi,\beta)} & {\varinjlim_{\Lcal_{\Frob'}\Xibf^{'\mathrm{op}}} \Tr^{\loc}(\Frob_{\Dcal}',\Dcal)} \\
	{\Tr(\Frob, \Ccal_{\tot})} & {\Tr(\Frob', \Dcal_{\tot})}
	\arrow[from=1-1, to=1-2]
	\arrow[from=1-1, to=2-1]
	\arrow[from=1-2, to=2-2]
	\arrow[from=2-1, to=2-2]
\end{tikzcd}\]
where the top horizontal map is the composition
$$\Tr^{\loc}(\Frob_{\Ccal}, \Ccal)(\xi,\beta) \to \varinjlim_{\Lcal_{\Frob}\Xibf^{\mathrm{op}}}\Tr^{\loc}(\Frob_{\Ccal},\Ccal) \to \varinjlim_{\Lcal_{\Frob'}\Xibf^{'\mathrm{op}}} \Tr^{\loc}(\Frob_{\Dcal}',\Dcal) \to \Tr(\Frob', \Dcal_{\tot}).$$
That map factors as 
$$\Tr^{\loc}(\Frob_{\Ccal}, \Ccal)(\xi,\beta) \to \Tr^{\loc}(\Frob_{\Dcal}', \Dcal)(u(\xi,\beta)) \to \varinjlim_{\Lcal_{\Frob'}\Xibf^{'\mathrm{op}}} \Tr(\Frob', \Dcal_{\tot}).$$
That composition is then induced by the bimodule $\Dcal_{\tot} \otimes_{\Dcal_{u(\xi)}} \Acal_{\xi}$ (and the morphism $\phi$ is clear). 

On the other hand the composition along the left vertical map and bottom horizontal one is induced by the bimodule $\Dcal_{\tot} \otimes_{(u^*\Dcal)_{\tot}} \Acal_{\tot} \otimes_{\Ccal_{\tot}} \Ccal_{\tot}$. Furthermore, there is a canonical map of $(\Dcal_{\tot}, \Ccal_{\xi})$-bimodules  
$$\Dcal_{\tot} \otimes_{\Dcal_{u(\xi)}} \Acal_{\xi} \to \Dcal_{\tot} \otimes_{(u^*\Dcal)_{\tot}} \Acal_{\tot} \otimes_{\Ccal_{\tot}} \Ccal_{\tot}$$ 
given by $d \otimes a \mapsto d \otimes a \otimes 1_{\Ccal_{\tot}}$. Since $\Ccal_{\xi}$ acts by $0$ on $\Acal_{xi'}$ for $\xi' \neq \xi$, it follows that both bimodule induce the same functor on after passing to traces. 
\end{proof}

\subsection{Finite orbit stacks}

\begin{defi}
\begin{enumerate}
\item Let $\Xbf$ be an algebraic stack, we say that $\Xbf$ is a finite orbit stack if the underlying topological space $|\Xbf|$ is finite.
\item Let $\Xbf$ be an algebraic stack with an action of $\Gbf$, $\Xbf$ is a $\Gbf$-finite orbit stack if $\Xbf/\Gbf$ is a finite orbit stack. 
\end{enumerate}
\end{defi}

\begin{nota} 
Let $\Xbf$ be a stack with an action of $\Gbf \times \Gbf$ (we think of $\Xbf$ as a $\Gbf$-bimodule), we define 
$$\HH^{\Xbf} = \bigoplus_{s,s' \in \Ch(\Tbf)} \DD(\Ubf \backslash \Xbf/\Ubf, \RR_{\Tbf \times\Tbf})_{[s,s']}$$
the `Hecke category' of $\Xbf$. It is naturally a $\HH^{\Gbf}$-bimodule. 
\end{nota} 

\begin{lem}[\protect{\cite[Lemma 6.2.4]{EteveFreeMonodromic}}]\label{lem:CategoricalKunneth}
Let $\Xbf$ and $\Ybf$ be two $\Tbf$-stacks, and assume that $\Xbf$ is a finite $\Tbf$-orbit stack. Then there is a canonical equivalence
$$\HH^{\Xbf} \otimes_{\HH^{\Tbf}} \HH^{\Ybf} \cong \HH^{\Xbf \times^{\Tbf} \Ybf}.$$
\end{lem}

\begin{lem}\label{lem:RestrictednessOfTermsInTrace}
Let $\Xbf$ be a $\Gbf \times \Gbf$-stack such that for all $w \in \Wbf$, the stack $\frac{\Ubf \backslash \Xbf/\Ubf}{\Ad_w \Tbf}$ is a finite orbit stack, then the stack $\frac{\Ubf \backslash \Gbf /\Ubf \times^{\Tbf} \Ubf \backslash \Xbf /\Ubf}{\Ad \Tbf}$ is a finite orbit stack. 
\end{lem}

\begin{proof}
It is enough to prove that, under the assumption on $\Xbf$, the stack $\frac{\Ubf \backslash \Bbf w \Bbf /\Ubf \times^{\Tbf} \Ubf \backslash \Xbf /\Ubf}{\Ad \Tbf}$ is a finite orbit stack for all $w \in \Wbf$. Let $\dot{w} \in \Nbf(\Tbf)$ be a lift of $w$. This choice yields a splitting $\Ubf \backslash \Bbf w \Bbf/\Ubf \cong \Tbf.w \times \pt/(\Ubf \cap {^w}\Ubf)$. It follows that the map 
$$\frac{\Ubf \backslash \Bbf w \Bbf /\Ubf \times^{\Tbf} \Ubf \backslash \Xbf /\Ubf}{\Ad \Tbf} \to \frac{\Tbf.w \times^{\Tbf} \Ubf \backslash \Xbf /\Ubf}{\Ad \Tbf} \cong \frac{\Ubf \backslash \Xbf /\Ubf}{\Ad_w\Tbf}$$
is a $\pt/(\Ubf \cap {^{w}}\Ubf)$-fibration. By assumption the target has finitely many points hence so does the source. 
\end{proof}

\begin{corol}\label{cor:Restrictedness}
Let $\Xbf$ be a $\Gbf \times \Gbf$-stack such that for all $w \in \Wbf$, the stack $\frac{\Ubf \backslash \Xbf/\Ubf}{\Ad_w \Tbf}$ is a finite orbit stack, then for all $n \geq 1$, the stack  $\frac{\Ubf \backslash \Gbf /\Ubf \times^{\Tbf} \dots \times^{\Tbf} \Ubf \backslash \Gbf /\Ubf \times^{\Tbf} \Ubf \backslash \Xbf /\Ubf}{\Ad \Tbf}$ with $n$ copies of $\Gbf$ is a finite orbit stack.
\end{corol}

\begin{proof}
This is an immediate induction on Lemma \ref{lem:RestrictednessOfTermsInTrace}.
\end{proof}

\subsection{Trace calculations}\label{sec:traces-calculations}

Let $\Xbf$ and $\Ybf$ be two stacks with actions $\Gbf \times \Gbf$, consider the correspondence 
\[\begin{tikzcd}
	& {\Ubf \backslash \Xbf  \times^{\Ubf}  \Ybf /\Ubf } \\
	{\Ubf \backslash \Xbf /\Ubf \times \Ubf \backslash \Ybf /\Ubf } && {\Ubf \backslash \Xbf \times^{\Gbf} \Ybf /\Ubf }
	\arrow["p"{description}, from=1-2, to=2-1]
	\arrow["m"{description}, from=1-2, to=2-3]
\end{tikzcd}\]
and the `convolution' functor 
$$\HH^{\Xbf} \otimes \HH^{\Ybf} \to \HH^{\Xbf \times^{\Gbf} \Ybf}$$
given by $(A \otimes B) \mapsto m_!p^*(A \widehat{\boxtimes} B)[\dim \Tbf]$. Compare with the actual convolution in $\Gbf$, see \Cref{sec:Hecke-cat-generalities}. 

\begin{prop}\label{prop:RelativeKunneth}
The convolution map induces a fully faithful functor
\begin{equation}\label{eq:ConvolutionForBimodules}
\HH^{\Xbf} \otimes_{\HH^{\Gbf}} \HH^{\Ybf} \to \HH^{\Xbf \times^{\Gbf} \Ybf}.
\end{equation}
\end{prop}
\begin{proof}
By definition of $\HH^{\Xbf} \otimes_{\HH^{\Gbf}} \HH^{\Ybf}$, when we compute it in the category of $\HH^{\Tbf}$-modules, we have
$$\HH^{\Xbf} \otimes_{\HH^{\Gbf}} \HH^{\Ybf} = \varinjlim_{[n] \in \Delta^{\opp}} \HH^{\Xbf} \otimes_{\HH^{\Tbf}} (\HH^{\Gbf})^{\otimes_{\HH^{\Tbf}} n} \otimes_{\HH^{\Tbf}} \HH^{\Ybf},$$
where the face maps are the left action, partial convolutions, and the right action, while the degeneracy maps insert the unit of $\HH^{\Gbf}$. For $n=0$, the middle factor is $\HH^{\Tbf}$, so the term is $\HH^{\Xbf}\otimes_{\HH^{\Tbf}}\HH^{\Ybf}$.
Using the categorical Künneth formula of Lemma \ref{lem:CategoricalKunneth}, and setting $\Xbf_n = \Xbf /\Ubf \times^{\Tbf}  \Ubf \backslash \Gbf /\Ubf \times^{\Tbf} \dots \times^{\Tbf}  \Ubf \backslash \Ybf$, we have 
$$\HH^{\Xbf} \otimes_{\HH^{\Tbf}} (\HH^{\Gbf})^{\otimes_{\HH^{\Tbf}} n} \otimes_{\HH^{\Tbf}} \HH^{\Ybf} = \HH^{\Xbf_n}.$$
Let $\Ybf_n =  \Xbf \times^{\Bbf}  \Gbf \times^{\Bbf} \dots \times^{\Bbf} \Ybf$ and $q_n : \Ybf_n \to \Xbf_n$. Note that this is a unipotent map. Arguing as in \cite[Lemma 6.3.2]{EteveFreeMonodromic}, the collection of functors $q_n^*$ define a fully faithful morphism of simplicial objects $(\HH^{\Xbf_n} \to \HH^{\Ybf_n}).$
The simplicial stack $(\Ybf_n)$ appears as the \v{C}ech nerve of the map $\Xbf \times^{\Bbf} \Ybf \to \Xbf \times^{\Gbf} \Ybf$ and by descent for categories of free monodromic sheaves, the colimit $\varinjlim_n \HH^{\Ybf_n}$ is identified with $\HH^{\Xbf \times^{\Gbf} \Ybf}$. Finally the convolution map $\HH^{\Xbf} \otimes \HH^{\Ybf} \to \HH^{\Xbf \times^{\Gbf} \Ybf}$ is identified with the inclusion functor in the colimit. 
\end{proof}

We prove a slight generalization of \cite[Theorem 1.3.1]{EteveFreeMonodromic}, the structure of the proof is very similar to the one in \emph{loc. cit.}. Let $\Xbf$ be a scheme with an action of $\Gbf \times \Gbf$ such that for all $w \in \Wbf$, the stack $\frac{\Ubf \backslash \Xbf/\Ubf}{\Ad_w \Tbf}$ is a finite orbit stack.

\begin{thm}\label{thm:calculation-traces}
Under the current assumptions on $\Xbf$, there is a natural fully faithful map  
\begin{equation}\label{eq:map-traces}
\Tr(\HH^{\Xbf}, \HH^{\Gbf}) \hookrightarrow \DD(\frac{\Xbf}{\Ad(\Gbf)})
\end{equation}
and the following canonical map $\HH^{\Xbf} \to \DD(\frac{\Xbf}{\Ad(\Gbf)})$ is given by $\alpha_{\Xbf,!}\beta^*_{\Xbf}$ where $\alpha_{\Xbf}$ and $\beta_{\Xbf}$ are the following maps
$$\Ubf \backslash \Xbf / \Ubf \xleftarrow{\beta_{\Xbf}} \frac{\Xbf}{\Ad(\Ubf)} \xrightarrow{\alpha_{\Xbf}} \frac{\Xbf}{\Ad(\Gbf)}.$$
\end{thm}

\begin{proof}
Recall that $\Tr(\HH^{\Xbf}, \HH^{\Gbf}) = \HH^{\Gbf} \otimes_{\HH^{\Gbf} \otimes \HH^{\Gbf, \rev}} \HH^{\Xbf}$. As in \cite[Section 6.3]{EteveFreeMonodromic}, let us resolve $\HH^{\Gbf}$ as a $\HH^{\Gbf} \otimes \HH^{\Gbf}$-bimodule, namely 
$$\HH^{\Gbf} = \varinjlim_{\Delta^{\opp}} (\HH^{\Gbf})^{\otimes_{\HH^{\Tbf}} (n+2)}$$
where the transition maps are given by partial multiplications. We thus have
\begin{align*}
\Tr(\HH^{\Xbf}, \HH^{\Gbf}) &= \HH^{\Gbf} \otimes_{\HH^{\Gbf} \otimes \HH^{\Gbf, \rev}} \HH^{\Xbf} \\
&= \varinjlim_{n} (\HH^{\Gbf})^{\otimes_{\HH^{\Tbf}} (n+2)} \otimes_{\HH^{\Gbf} \otimes \HH^{\Gbf, \rev}} \HH^{\Xbf} \\
&= \varinjlim_{n} (\HH^{\Gbf})^{\otimes_{\HH^{\Tbf}} n} \otimes_{\HH^{\Tbf} \otimes \HH^{\Tbf}} \HH^{\Xbf}.
\end{align*}

Let $\Xbf_n = \frac{\Ubf \backslash \Gbf /\Ubf \times^{\Tbf} \dots \times^{\Tbf} \Ubf \backslash \Gbf /\Ubf \times^{\Tbf} \Ubf\backslash\Xbf/\Ubf}{\Ad(\Tbf)}$, where there are $n$ copies of $\Gbf$ (note that this stack is not the same as in the proof of Proposition \ref{prop:RelativeKunneth}). Using the categorical Künneth formula of Lemma \ref{lem:CategoricalKunneth}, we have 
$$(\HH^{\Gbf})^{\otimes_{\HH^{\Tbf}} n} \otimes_{\HH^{\Tbf} \otimes \HH^{\Tbf}} \HH^{\Xbf} = \bigoplus_{s,s'\in \Ch(\Tbf)} \DD(\Xbf_n, \RR_{\Tbf \times \Tbf})_{[s,s']}.$$
By assumption on $\Xbf$ and using Corollary \ref{cor:Restrictedness}, the stack $\Xbf_n$ is a finite orbit stack and using \cite[Lemma 2.8.3]{EteveFreeMonodromic}, the forgetful map 
$$ \bigoplus_{s,s'\in \Ch(\Tbf)} \DD(\Xbf_n, \RR_{\Tbf \times \Tbf})_{[s,s']} \to \DD(\Xbf_n, \Lambda)$$
is an equivalence. 

Let $\Ybf_n = \frac{\Gbf \times^{\Bbf} \dots \times^{\Bbf} \Gbf}{\Ad(\Bbf)}$, where there are $n$-copies of $\Gbf$ (again this stack differs from the proof of Proposition \ref{prop:RelativeKunneth}) and let $q_n : \Ybf_n \to \Xbf_n$. Arguing as in \cite[Lemma 6.3.2]{EteveFreeMonodromic}, the collection of functors $(q_n^*)$ yields a fully faithful map of simplicial objects $(\DD(\Xbf_n, \Lambda))_n \to (\DD(\Ybf_n, \Lambda))_n$. Again since the simplicial scheme $(\Ybf_n)$ is the \v{C}ech nerve of the map $\frac{\Xbf}{\Ad(\Bbf)} \to \frac{\Xbf}{\Ad(\Gbf)}$, we get a fully faithful functor 
$$\Tr(\HH^{\Xbf}, \HH^{\Gbf}) \to \DD(\frac{\Xbf}{\Ad(\Gbf)}, \Lambda).$$

It follows from the construction of the isomorphism, that the natural functor 
$$\HH^{\Xbf} \to \Tr(\HH^{\Xbf}, \HH^{\Gbf}) \to \DD(\frac{\Xbf}{\Ad(\Gbf)}, \Lambda)$$
is given by $\alpha_{\Xbf,!}\beta_{\Xbf}^*$. 
\end{proof}

\begin{lem}\label{lem:essential-surjectivity-of-trace-map}
The map (\ref{eq:map-traces})
is essentially surjective if $\ell$ does not divide $|\Wbf|$. In the cases where the pair $(\Gbf, \Xbf)$ is either $(\Gbf, \Gbf_{\Frob})$ or $(\Lbf, \Ubf_{\Pbf}\backslash \Gbf/\Frob(\Ubf_{\Pbf}))$ where $\Pbf = \Lbf\Ubf_{\Pbf}$ is a parabolic of $\Gbf$, the functor (\ref{eq:map-traces}) is essentially surjective.
\end{lem} 

\begin{proof}
This is essentially the same argument as in \cite{EteveFreeMonodromic}. We prove that the essential image of the functor $\alpha_{\Xbf,!}\beta_{\Xbf}^*$ generates the category. In the two cases $(\Gbf, \Gbf_{\Frob})$ and $(\Lbf, \Ubf_{\Pbf}\backslash \Gbf/\Frob(\Ubf_{\Pbf}))$, this is a direct calculation that reduces to the fact that the cohomology of Deligne--Lusztig varieties generate $\DD(\Rep \Gbf^{\Frob})$, see \Cref{prop:Generation}. 

Assume now that $\ell$ does not divide $|\Wbf|$. By \cite{MirkovicVilonen}, we have
$$ \alpha_{\Xbf,!}\beta_{\Xbf}^*\beta_{\Xbf,!}\alpha_{\Xbf}^* \cong \Spr_{\Gbf} *_!^{\Gbf} -$$ 
where $\Spr_{\Gbf} \in \Perv(\frac{\Gbf}{\Ad(\Gbf)})$ is the Springer sheaf and $*_!^{\Gbf}$ is the natural action of $\DD(\frac{\Gbf}{\Ad(\Gbf)})$ on $\DD(\frac{\Xbf}{\Ad(\Gbf)})$. Under the assumption on $\ell$, skyscraper sheaf $\delta_1$ is a direct summand of $\Spr_{\Gbf}$ and therefore the functor $\beta_{\Xbf,!}\alpha_{\Xbf}^*$ is conservative. It then follows from the Barr--Beck--Lurie theorem, that $\DD(\frac{\Xbf}{\Ad(\Gbf)})$ is identified with the category of comodules over $\beta_{\Xbf,!}\alpha_{\Xbf}^*\alpha_{\Xbf,!}\beta_{\Xbf}^*$ and the statement follows.  
\end{proof}

\begin{conj}
\Cref{lem:essential-surjectivity-of-trace-map} holds without assumption on $\ell$. 
\end{conj}

\bibliographystyle{alpha}
\bibliography{biblio.bib}

\end{document}